\documentclass{elsart}

\usepackage{graphicx}
\usepackage{psfrag}
\usepackage{color}

\usepackage{amssymb}
\usepackage{amsmath}
\usepackage{dsfont}
\usepackage{rotating}

\newcommand{\Path}{\boldsymbol{\Psi}}

\newcommand{\tens}[1]{#1}
\newcommand{\K}{\mathbf{K}}

\newcommand{\D}{\mathbf{D}}

\newcommand{\U}{\mathbf{U}}

\renewcommand{\S}{\mathbf{S}} 
\newcommand{\Q}{\mathbf{Q}} 
\newcommand{\q}{\mathbf{q}} 
\newcommand{\f}{\mathbf{f}} 
\newcommand{\g}{\mathbf{g}}
\newcommand{\h}{\mathbf{h}} 
\newcommand{\A}{\mathbf{A}} 
\newcommand{\B}{\mathbf{B}} 
\newcommand{\F}{\mathbf{F}} 
\newcommand{\G}{\mathbf{G}} 
\newcommand{\x}{\mathbf{x}} 
\newcommand{\xxi}{\boldsymbol{\xi}} 
 
\renewcommand{\P}{\mathbf{P}} 
\newcommand{\w}{\mathbf{w}}

\newcommand{\halb}{\frac{1}{2}}

\newcommand{\ar}{\phi_1\rho_1}

\newcommand{\ub}{\textbf{u}_\textbf{1}}

\renewcommand{\v}{\mathbf{v}}

\begin{document}

\begin{frontmatter}



\title{Lagrangian ADER-WENO Finite Volume Schemes on Unstructured Tetrahedral Meshes for Conservative and Nonconservative Hyperbolic Systems in 3D}


\author[UNITN]{Walter Boscheri}
\ead{walter.boscheri@unitn.it}
\author[UNITN]{Michael Dumbser}
\ead{michael.dumbser@unitn.it}
\address[UNITN]{Laboratory of Applied Mathematics \\ Department of Civil, Environmental and Mechanical Engineering \\ University of Trento, Via Mesiano 77, I-38123 Trento, Italy}

\begin{abstract}
In this paper we present a new family of high order accurate Arbitrary-Lagrangian-Eulerian (ALE) one-step ADER-WENO finite volume schemes for the solution of 
nonlinear systems of conservative and non-conservative hyperbolic partial differential equations with stiff source terms on moving tetrahedral meshes in three 
space dimensions. 
A WENO reconstruction technique is used to achieve high order of accuracy in space, while an element-local space-time Discontinuous Galerkin finite element predictor 
on moving meshes is used to obtain a high order accurate one-step time discretization. Within the space-time predictor the physical element is mapped 
onto a reference element using an isoparametric approach, where the space-time basis and test functions are given by the Lagrange interpolation polynomials 
passing through a predefined set of space-time nodes. 
Since our algorithm is cell-centered, the final mesh motion is computed by using a suitable node solver algorithm. A rezoning step as well as a flattener 
strategy are used in some of the test problems to avoid mesh tangling or excessive element deformations that may occur when the computation involves strong 
shocks or shear waves. 
The ALE algorithm presented in this article belongs to the so-called \textit{direct} ALE methods because the final Lagrangian finite volume scheme is 
based directly on a space-time conservation formulation of the governing PDE system, with the rezoned geometry taken already into account during the 
computation of the fluxes. 

We apply our new high order unstructured ALE schemes to the 3D Euler equations of compressible gas dynamics, for which a set of classical numerical 
test problems has been solved and for which convergence rates up to sixth order of accuracy in space and time have been obtained. 
We furthermore consider the equations of classical ideal magnetohydrodynamics (MHD) as well as the non-conservative seven-equation Baer-Nunziato model 
of compressible multi-phase flows with stiff relaxation source terms. 

\end{abstract}
\begin{keyword}
Arbitrary-Lagrangian-Eulerian (ALE) finite volume schemes \sep 
WENO reconstruction on moving unstructured tetrahedral meshes \sep 
high order of accuracy in space and time \sep 
stiff source terms \sep 
local rezoning \sep 
conservation laws and nonconservative hyperbolic PDE \sep Euler equations \sep MHD equations \sep compressible multi-phase flows \sep Baer-Nunziato model 
\end{keyword}
\end{frontmatter}


\section{Introduction}
\label{sec.introduction}

Any Lagrangian method aims at following the fluid motion as closely as possible, with a computational mesh that is moving with the local fluid velocity. Therefore the Lagrangian approach allows 
material interfaces and contact waves to be precisely located and tracked during the computation, achieving a much more accurate resolution of these waves compared to classical Eulerian methods 
on fixed grids. 
For this reason a lot of research has been carried out in the last decades in order to develop Lagrangian methods. Already John von Neumann and Richtmyer were working on Lagrangian schemes in the 
1950ies \cite{Neumann1950}, using a formulation of the governing equations in primitive variables, which was also used later in \cite{Benson1992,Caramana1998}. However, most of the modern 
Lagrangian finite volume schemes use the conservation form of the equations based on the physically conserved quantities like mass, momentum and total energy in order to compute shock waves
properly, see e.g. \cite{munz94,Smith1999,Maire2007,Despres2009}. Lagrangian schemes can be also classified according to the location of the physical variables on the mesh: when all variables 
are defined on a collocated grid the so-called \textit{cell-centered} approach is adopted \cite{Depres2012,ShashkovCellCentered,Maire2007,Maire2008,Maire2010}, while in the 
\textit{staggered mesh} approach \cite{StagLag,LoubereSedov3D} the velocity is defined at the cell interfaces and the other variables at the cell center. 

Cell-centered Lagrangian Godunov-type schemes of the Roe-type and of the HLL-type for the Euler equations of compressible gas dynamics have first been considered by Munz in \cite{munz94}. 
A cell-centered Godunov-type scheme has also been introduced by Carr\'e et al. in \cite{Despres2009}, who developed a Lagrangian finite volume algorithm on general multi-dimensional 
unstructured meshes. The resulting finite volume scheme is node based and compatible with the mesh displacement. In the work of Despr\'es et al. \cite{DepresMazeran2003,Despres2005} the 
physical part of the system of equations is coupled and evolved together with the geometrical part, hence obtaining a weakly hyperbolic system of conservation laws that is solved using a 
node-based finite volume scheme. Furthermore they presented a cell-centered Lagrangian method \cite{Depres2012} that is translation invariant and suitable for curved meshes. 
In \cite{Maire2009,Maire2010,Maire2011} Maire proposed first and second order accurate cell-centered Lagrangian schemes in two- and three- space dimensions on general polygonal grids, 
where the time derivatives of the fluxes are obtained using a node-centered solver that may be considered as a multi-dimensional extension of the Generalized Riemann problem methodology 
introduced by Ben-Artzi and Falcovitz \cite{Artzi}, Le Floch et al. \cite{Raviart.GRP.1,Raviart.GRP.2} and Titarev and Toro \cite{Toro:2006a,toro3,titarevtoro}. The node solver algorithm 
developed for hydrodynamics by Maire in \cite{Maire2009} is used also in this paper and applied to both Euler and MHD equations on moving tetrahedral meshes. Since Lagrangian schemes may
lead to severe mesh deformation after a finite time, it is necessary to remesh (or at least to rezone) the computational grid from time to time. A very popular approach consists therefore 
in Lagrangian remesh and remap schemes, such as the family of cell-centered ALE remap algorithms introduced by Shashkov et al. and Maire et al. in 
\cite{ShashkovCellCentered,ShashkovRemap1,ShashkovRemap3,ShashkovRemap4,ShashkovRemap5,MaireMM2}. 
In \cite{ShashkovMultiMat1,ShashkovMultiMat2,ShashkovMultiMat3,ShashkovMultiMat4} purely Lagrangian and Arbitrary-Lagrangian-Eulerian (ALE) numerical schemes with remapping for multi-phase 
and multi-material flows are discussed. All the Lagrangian schemes listed so far are at most second order accurate in space and time. 

Higher order of accuracy in space was first achieved in \cite{chengshu1,chengshu2,chengshu3,chengshu4} by Cheng and Shu, who introduced a third order accurate  essentially non-oscillatory (ENO)  reconstruction operator into Godunov-type Lagrangian finite volume schemes. High order of accuracy in time was guaranteed either by the use of a Runge-Kutta or by a Lax-Wendroff-type time stepping. 
The mesh velocity is simply computed as the arithmetic average of the corner-extrapolated values in the cells adjacent to a mesh vertex. Such a node solver algorithm is very simple and general and 
can be  easily applied to different complicated nonlinear systems of hyperbolic PDE in multiple space dimensions. Cheng and Toro \cite{chengtoro} also investigated Lagrangian ADER-WENO schemes in 
one space dimension. In the finite element framework high order Lagrangian schemes have been developed for example by Scovazzi et al. \cite{scovazzi1,scovazzi2}. In \cite{Lagrange1D} Dumbser et 
al. presented high order ADER-WENO Lagrangian finite volume schemes for hyperbolic balance laws with stiff source terms. In this case the high order of accuracy in time was achieved by using the 
local space-time Galerkin predictor method proposed in \cite{DumbserEnauxToro,HidalgoDumbser} for the Eulerian case, whereas a high order WENO reconstruction algorithm was used to obtain high 
order of accuracy in space. In \cite{Lagrange2D,LagrangeNC} Boscheri and Dumbser extended this algorithm to unstructured triangular meshes 
for conservative and non-conservative hyperbolic systems with stiff source terms. In \cite{LagrangeMHD} three different node solver algorithms have been applied to the Euler equations of compressible 
gas dynamics as well as to the equations for magnetohydrodynamics and have been compared with each other. The multidimensional HLL Riemann solver presented in \cite{BalsaraMultiDRS} for the Eulerian  framework on fixed grids has been used as a node solver for the computation of the mesh velocity in \cite{LagrangeMHD} and for the computation of the space-time  fluxes of a high order Lagrangian 
finite volume scheme in \cite{LagrangeMDRS}. In the latter reference it has been shown that the use of a multi-dimensional Riemann solver allows the use of larger time steps in multiple space 
dimensions and therefore leads to a computationally more efficient scheme compared to a method based on classical one-dimensional Riemann solvers. 

In literature there are also other methods using a Lagrangian approach and these schemes are at least briefly mentioned in the following. For example, also meshless particle schemes, such as 
the smooth particle hydrodynamics (SPH) method, belong to the category of fully Lagrangian schemes, see e.g. \cite{Monaghan1994,SPHLagrange,SPHWeirFlow,SPH3D,Dambreak3D}. SPH is generally used 
to follow the fluid motion in very complex deforming domains. Since it is a particle method, no rezoning or remeshing has to be applied. Furthermore, also semi-Lagrangian methods should be 
mentioned. They are typically adopted to solve transport equations \cite{ALE2000Belgium,ALE1996FV}. Although these schemes use a fixed mesh, as in the classical Eulerian approach, 
the Lagrangian trajectories of the fluid are followed backward in time in order to compute the numerical solution at the the new time level, see for example 
\cite{Casulli1990,CasulliCheng1992,LentineEtAl2011,HuangQiul2011,QuiShu2011,BoscheriDumbser}. There is also the class of 
Arbitrary-Lagrangian-Eulerian (ALE) methods \cite{Hirt1974,Peery2000,Smith1999,Feistauer1,Feistauer2,Feistauer3,Feistauer4}, where the mesh moves with a velocity that does not necessarily 
have to coincide with the local fluid velocity. This method is often used for fluid-structure interaction (FSI) problems, but it is also used together with Lagrangian remap schemes. For the 
sake of generality, the scheme presented in this article uses an ALE approach so that the local mesh velocity can in principle be chosen independently from the local fluid velocity. 

In this paper we extend the algorithm presented in \cite{Lagrange2D,LagrangeNC} to moving unstructured tetrahedral meshes in three space dimensions. To the knowledge of the authors, this is the first
better than second order accurate Lagrangian finite volume method on three-dimensional tetrahedral meshes ever presented. 
We consider the Euler equations of compressible gas  dynamics as well as the ideal classical MHD equations and the non-conservative seven-equation Baer-Nunziato model of compressible multi-phase 
flows with stiff source terms. The node solver proposed by Maire in \cite{Maire2009} is applied, as well as the node solver of Cheng and Shu \cite{chengshu1}. 

The rest of this article is structured as follows: in Section \ref{sec.numethod} we describe the proposed numerical scheme in detail, while in Section \ref{sec.validation} we show numerical 
convergence studies up to sixth order of accuracy in space and time as well as numerical results for several classical test problems for all of the above-mentioned hyperbolic systems. 
Finally, in Section \ref{sec.concl} we give some concluding remarks and an outlook to future research and developments.

\section{Numerical method}
\label{sec.numethod} 

In this paper we consider nonlinear systems of hyperbolic balance laws which may also contain non-conservative products and stiff source terms. A general formulation that is suitable to write  
the above mentioned systems reads
\begin{equation}
\label{eqn.pde.nc} 
  \frac{\partial \Q}{\partial t} + \nabla \cdot \tens{\F}(\Q) + \tens{\B}(\Q) \cdot \nabla \Q = \S(\Q), \qquad \x \in \Omega \subset \mathds{R}^3, t \in \mathds{R}_0^+, 
\end{equation} 
where $\Q=(q_1,q_2,...,q_\nu)$ denotes the vector of conserved variables, $\tens{\F} = (\mathbf{f}, \mathbf{g}, \mathbf{h})$ is the conservative nonlinear flux tensor, $\tens{\B}=(\B_1,\B_2,\B_3)$ contains the purely non-conservative part of the system written in block-matrix notation and $\S(\Q)$ represents a nonlinear algebraic source term that is allowed to be stiff. We furthermore introduce the abbreviation $\P = \P(\Q,\nabla \Q) = \B(\Q) \cdot \nabla \Q$ to ease notation in some parts of the manuscript. 

In a Lagrangian framework the computational domain $\Omega(t) \subset \mathds{R}^3$ is time-dependent and is discretized at the current time $t^n$ by a set of tetrahedral elements $T^n_i$. $N_E$ denotes the total number of elements contained in the domain and the union of all elements is called the \textit{current tetrahedrization} $\mathcal{T}^n_{\Omega}$ of the domain 
\begin{equation}
\mathcal{T}^n_{\Omega} = \bigcup \limits_{i=1}^{N_E}{T^n_i}. 
\label{trian}
\end{equation}

\begin{figure}[!htbp]
\begin{center}
\includegraphics[width=0.80\textwidth]{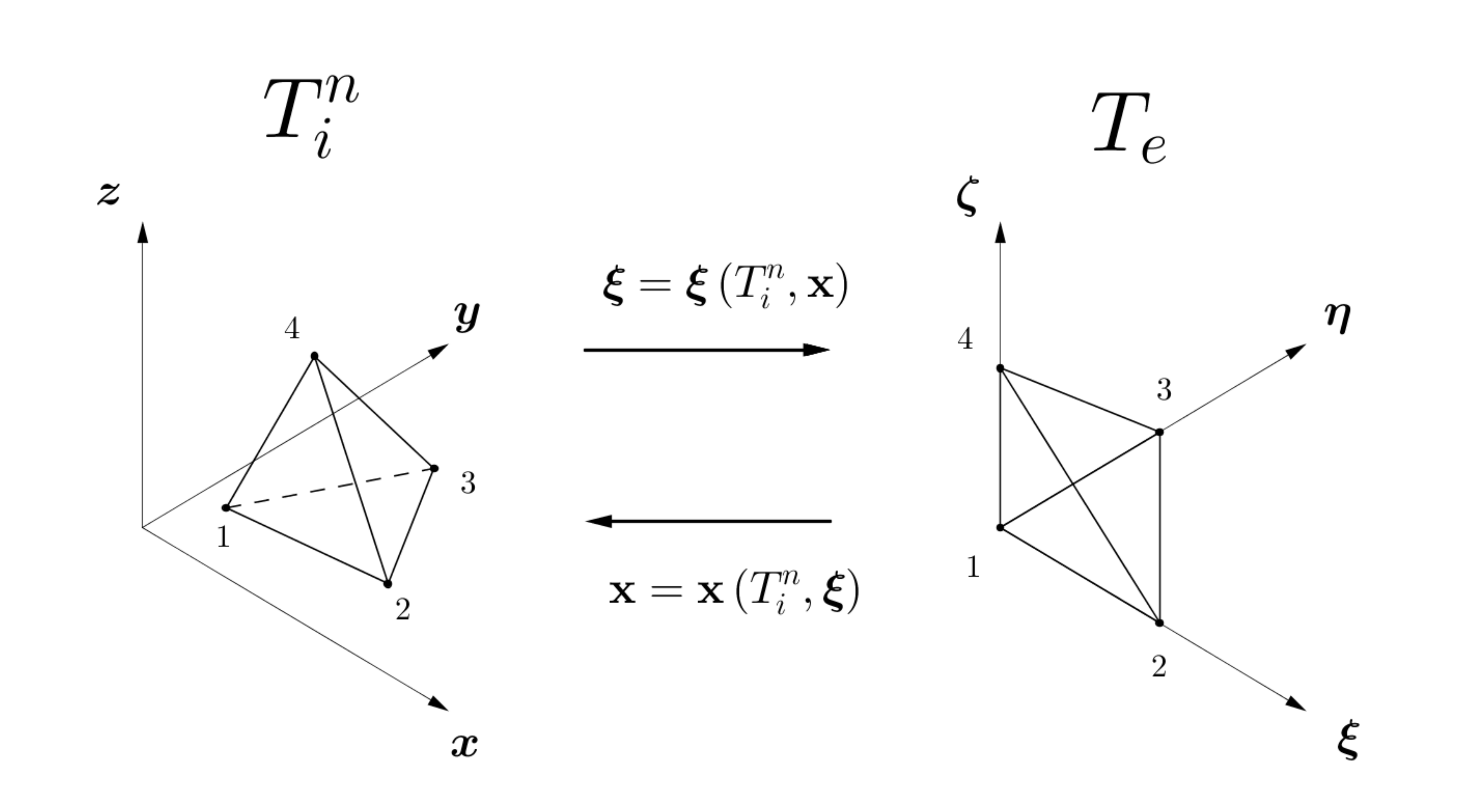}  
\caption{Spatial mapping from the physical element $T^n_i$ defined with $\mathbf{x}=(x,y,z)$ to the unit reference tetrahedron $T_e$ in $\boldsymbol{\xi} = (\xi, \eta, \zeta)$.}
\label{fig.refSystem}
\end{center}
\end{figure}

Since we are dealing with a moving computational domain where the mesh configuration continuously changes in time, it is more convenient to map the physical element $T^n_i$ to a reference element $T_e$ via a \textit{local} reference coordinate system $\xi-\eta-\zeta$. The spatial reference element $T_e$ is the unit tetrahedron shown in Figure \ref{fig.refSystem} and is defined by the nodes $\boldsymbol{\xi}_{e,1}=(\xi_{e,1},\eta_{e,1},\zeta_{e,1})=(0,0,0)$, $\boldsymbol{\xi}_{e,2}=(\xi_{e,2},\eta_{e,2},\zeta_{e,2})=(1,0,0)$, $\boldsymbol{\xi}_{e,3}=(\xi_{e,3},\eta_{e,3},\zeta_{e,3})=(0,1,0)$ and $\boldsymbol{\xi}_{e,4}=(\xi_{e,4},\eta_{e,4},\zeta_{e,4})=(0,0,1)$, where $\boldsymbol{\xi} = (\xi, \eta, \zeta)$ is the vector of the spatial coordinates in the reference system, while the position vector $\mathbf{x}=(x,y,z)$ is defined in the physical system. Let furthermore $\mathbf{X}^n_{k,i} = (X^n_{k,i},Y^n_{k,i},Z^n_{k,i})$ be the vector of physical spatial coordinates of the $k$-th vertex of tetrahedron $T^n_i$. Then the linear mapping from $T^n_i$ to $T_e$ is given by
\begin{equation} 
 \mathbf{x} = \mathbf{X}^n_{1,i} + 
\left( \mathbf{X}^n_{2,i} - \mathbf{X}^n_{1,i} \right) \xi + 
\left( \mathbf{X}^n_{3,i} - \mathbf{X}^n_{1,i} \right) \eta + 
\left( \mathbf{X}^n_{4,i} - \mathbf{X}^n_{1,i} \right) \zeta.
 \label{xietaTransf} 
\end{equation} 

As usual for finite volume schemes, data are stored and evolved in time as piecewise constant cell averages. They are defined at each time level $t^n$ within the control volume $T^n_i$ as
\begin{equation}
  \Q_i^n = \frac{1}{|T_i^n|} \int_{T^n_i} \Q(\mathbf{x},t^n) d\x,     
 \label{eqn.cellaverage}
\end{equation}  
with $|T_i^n|$ denoting the volume of tetrahedron $T_i^n$. In the next Section \ref{sec.weno} a WENO reconstruction technique is described and used to obtain piecewise higher order polynomials 
$\mathbf{w}_h(\x,t^n)$ from the known cell averages $\Q_i^n$. High order of accuracy in time is achieved later in Section \ref{sec.lst} by applying a local space-time Galerkin predictor method
to the reconstruction polynomials $\mathbf{w}_h(\x,t^n)$.

\subsection{Polynomial WENO reconstruction} 
\label{sec.weno} 
The WENO reconstruction operator produces piecewise polynomials $\mathbf{w}_h(\x,t^n)$ of degree $M$. The  $\mathbf{w}_h(\x,t^n)$ are computed for each control volume $T^n_i$ from the known 
cell averages within a so-called \textit{reconstruction stencil} $\mathcal{S}_i^s$, which is composed of an appropriate neighborhood of element $T^n_i$ and contains a prescribed total number 
$n_e$ of tetrahedra. We do not use the original \textit{pointwise} WENO method first introduced by Shu et al. \cite{shu_efficient_weno,HuShuTri,ZhangShu3D}, but we adopt the 
\textit{polynomial} formulation proposed in \cite{friedrich,kaeserjcp,Dumbser2007204,DumbserKaeser06b} and also used in \cite{MixedWENO2D,MixedWENO3D}, which is relatively simple to code and 
which allows the scheme to reach very high order of accuracy even on unstructured tetrahedral meshes in three space dimensions.    

\begin{figure}[!htbp]
\begin{center}
\begin{tabular}{ccc} 
\includegraphics[width=0.33\textwidth]{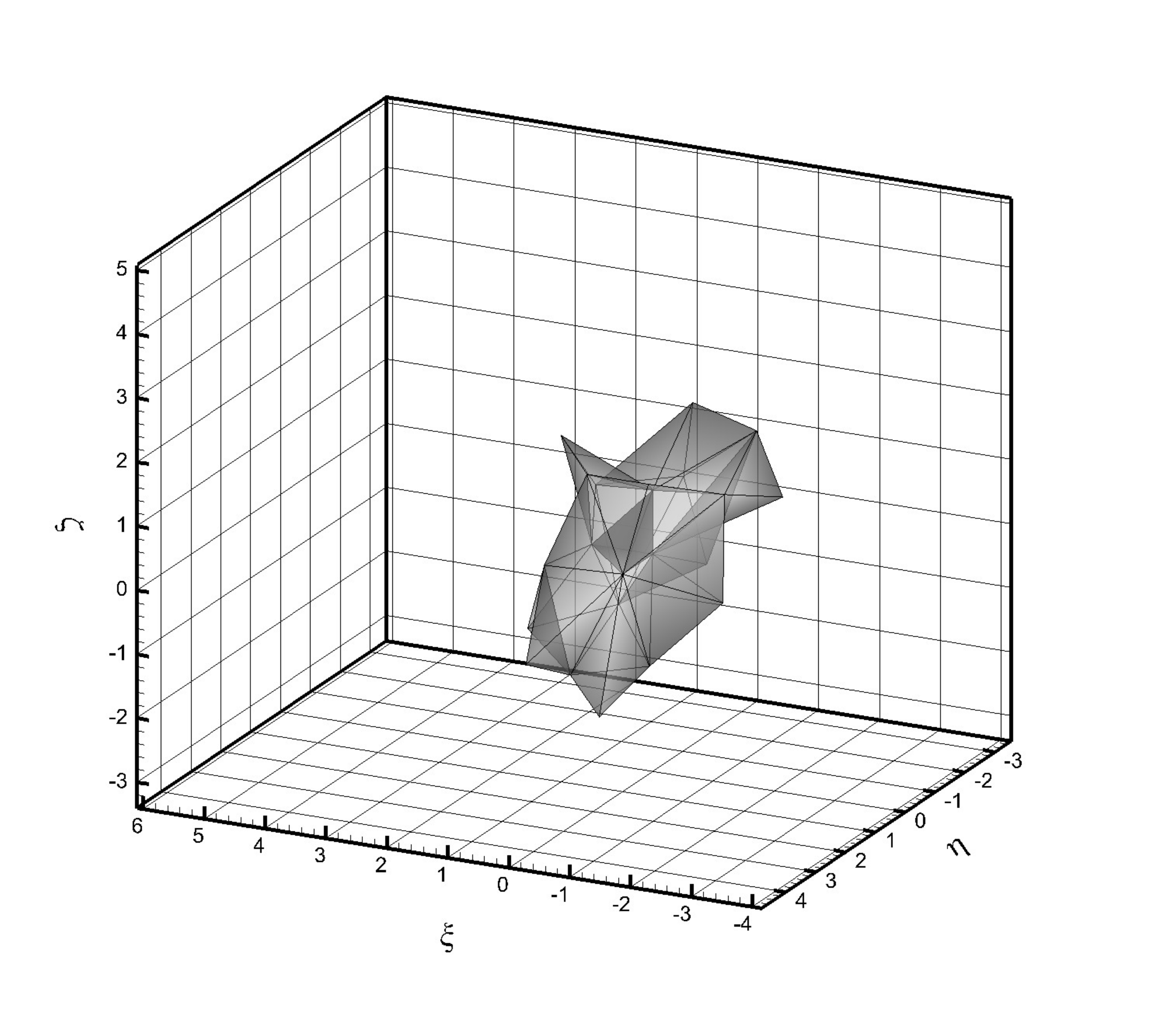}  &           
\includegraphics[width=0.33\textwidth]{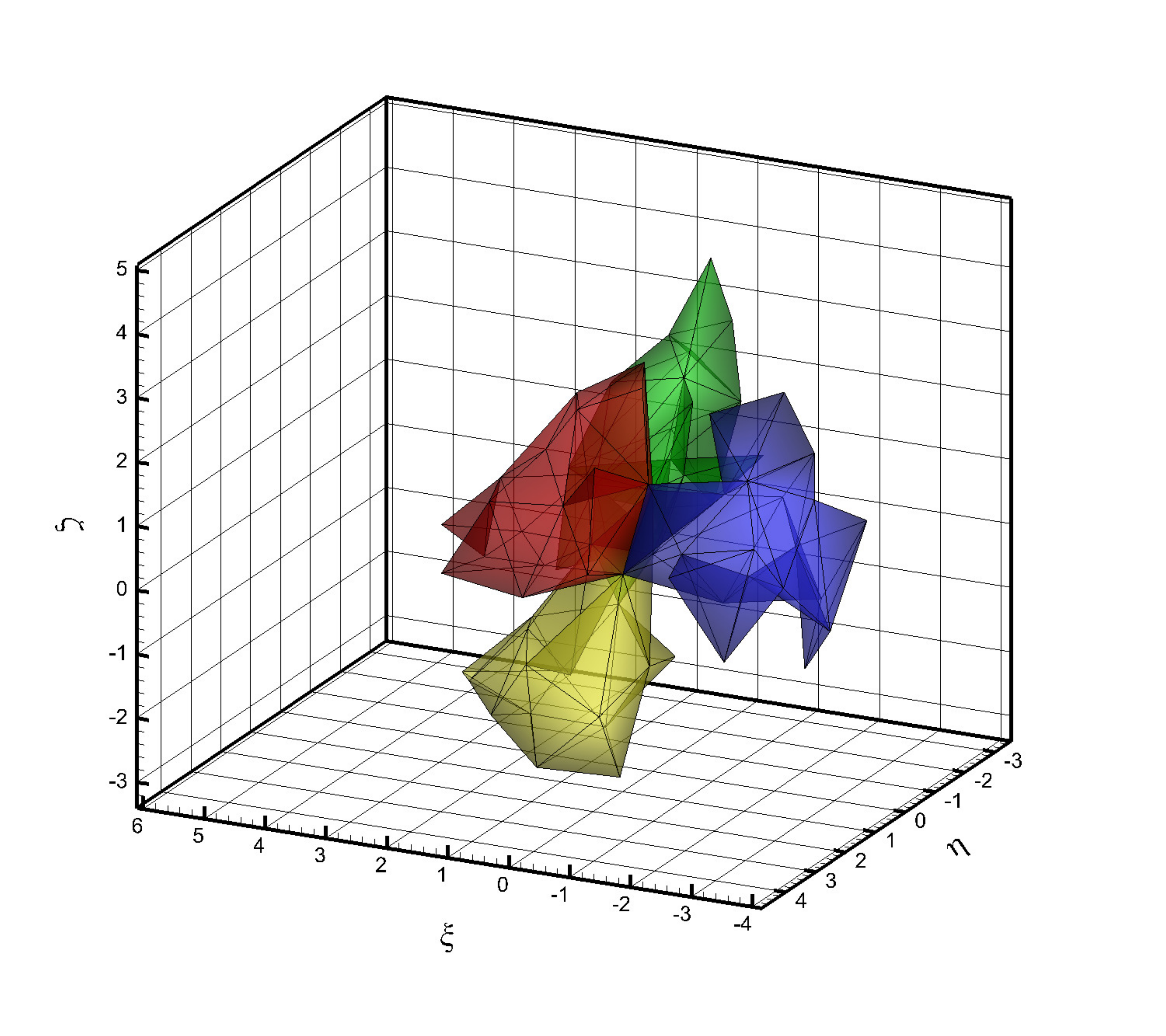}  &
\includegraphics[width=0.33\textwidth]{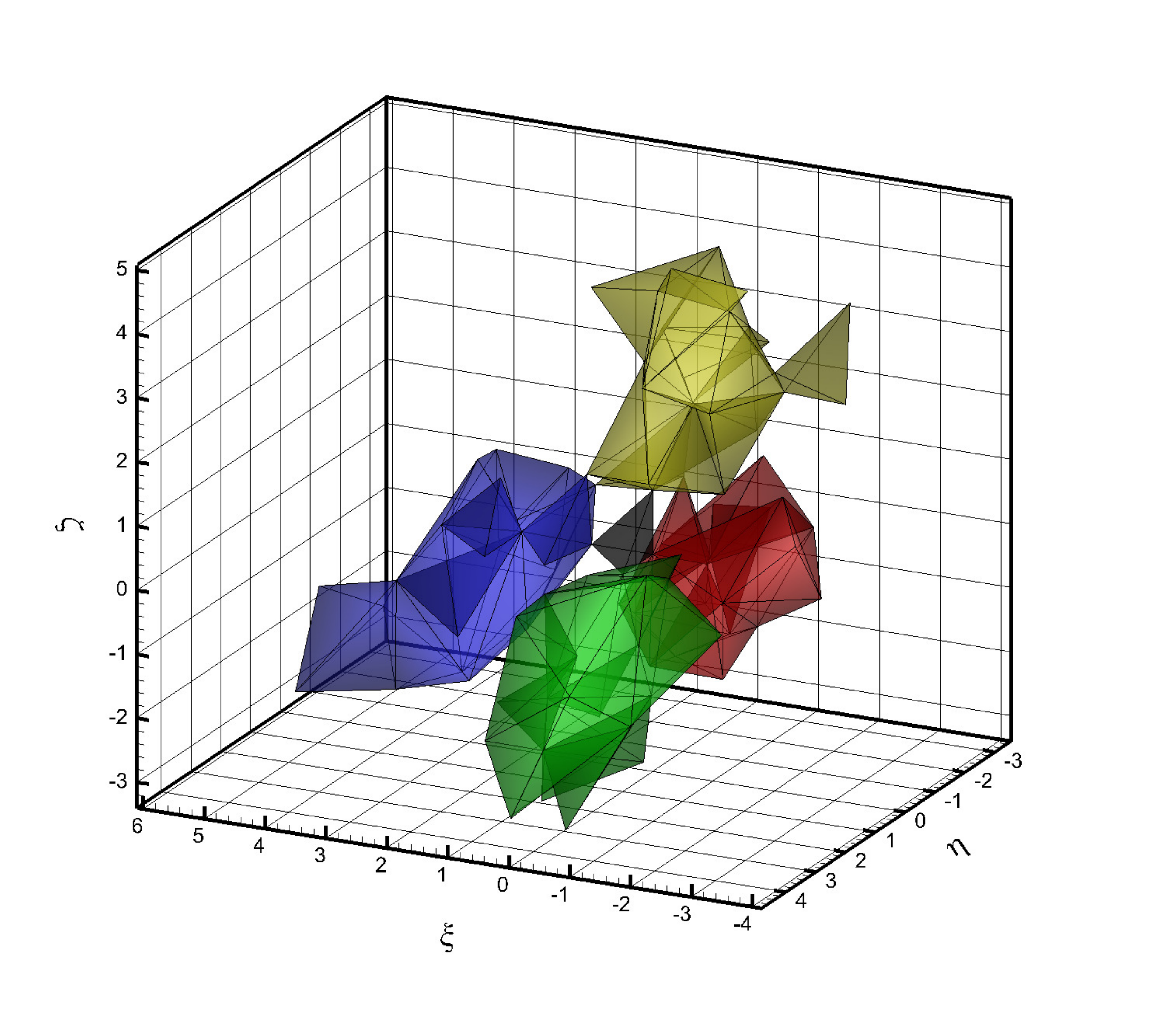}  \\           
\end{tabular} 
\caption{Three-dimensional WENO reconstruction stencils in the reference coordinate system with $M=2$ and $n_e=30$. One central stencil (left), four forward stencils (center) and four backward stencils (right). Tetrahedron $T_i^n$ is highlighted in black.} 
\label{fig.Stencil}
\end{center}
\end{figure}

According to \cite{DumbserKaeser06b}, in three space dimensions we always use nine reconstruction stencils, hence $1 \leq s \leq 9$. Specifically, we consider one central stencil given 
by $s=1$, four forward stencils with $s \in \left\{2,3,4,5 \right\}$ and four backward stencils with $s \in \left\{ 6,7,8,9 \right\}$, as depicted in Figure \ref{fig.Stencil}. Each forward and 
backward stencil sector is spanned by one point and three vectors: the four forward stencils are defined by a vertex $k$ of the tetrahedron $T_i^n$ and the triplet of vectors connecting $k$ to the three 
vertices of the opposite face, while the backward sectors are defined by the negative vectors of the forward stencils and the opposite face barycenter. Each type of stencil is filled by 
recursively adding neighbor elements until the prescribed total number $n_e$ is reached. An element belongs to the stencil if its barycenter is located in the corresponding sector. For the 
central stencil we use a simple Neumann-type neighbor search algorithm that recursively adds direct face neighbors to the stencil, until the desired number $n_e$ is reached. For the 
remaining eight one-sided stencils we use a Voronoi-type search algorithm, which fills the stencil starting from the vertex neighborhood of the tetrahedron and then using recursively vertex 
neighbors of stencil elements. Each stencil contains a total number of elements $n_e$ that depends on the reconstruction degree $M$, hence 
\begin{equation}
\mathcal{S}_i^s = \bigcup \limits_{j=1}^{n_e} T^n_{m(j)}, 
\label{stencil}
\end{equation}
where $1 \leq j \leq n_e$ is a local index which progressively counts the elements in the stencil number $s$ and $m(j)$ represents a mapping from the local index $j$ to the global index of 
the element in $\mathcal{T}^n_{\Omega}$. As explained in \cite{StencilRec1990,Olliver2002,kaeserjcp}, the total number of elements $n_e$ must be greater than the smallest number 
$\mathcal{M} = (M+1)(M+2)(M+3)/6$ needed to reach the formal order of accuracy $M+1$. As suggested in \cite{DumbserKaeser06b,Dumbser2007204} we typically take $n_e = 3 \mathcal{M}$ in 
three space dimensions. 

The high order reconstruction polynomial for each candidate stencil $\mathcal{S}_i^s$ for tetrahedron $T_i^n$ is written in terms of the \textit{orthogonal} Dubiner-type basis functions $\psi_l(\xi,\eta,\zeta)$ \cite{Dubiner,orth-basis,CBS-book} on the reference tetrahedron $T_e$, i.e.
\begin{equation}
\label{eqn.recpolydef} 
\w^s_h(\x,t^n) = \sum \limits_{l=1}^\mathcal{M} \psi_l(\boldsymbol{\xi}) \hat \w^{n,s}_{l,i} := \psi_l(\boldsymbol{\xi}) \hat \w^{n,s}_{l,i},   
\end{equation}
where the mapping to the reference coordinate system is given by \eqref{xietaTransf} and $\hat \w^{n,s}_{l,i}$ denote the \textit{unknown} degrees of freedom (expansion coefficients) of the reconstruction polynomial on stencil $\mathcal{S}_i^s$ for element $T_i^n$ at time $t^n$. In the rest of the paper we will use classical tensor index notation based on the Einstein summation convention, which implies summation over two equal indices.

Integral conservation is required for the reconstruction on each element $T_j^n$ of the stencil $\mathcal{S}_i^s$, hence 
\begin{equation}
\label{intConsRec}
\frac{1}{|T^n_j|} \int \limits_{T^n_j} \psi_l(\boldsymbol{\xi}) \hat \w^{n,s}_{l,i} d\x = \Q^n_j, \qquad \forall T^n_j \in \mathcal{S}_i^s.     
\end{equation}
Inserting the transformation \eqref{xietaTransf} into the above expression \eqref{intConsRec}, an analytical integration formula can be obtained that is a function of the four physical 
vertex coordinates $\mathbf{X}^n_{k,j}$ of the tetrahedron. The resulting algebraic expressions of the integrals appearing in \eqref{intConsRec} can be obtained for example at the aid 
of a symbolic computer algebra system like MAPLE. Up to $M=3$ we use the aforementioned analytical integration, while for higher reconstruction degrees the integrals in \eqref{intConsRec} 
are simply evaluated using Gaussian quadrature formulae of suitable order, see \cite{stroud} for details, since the analytical expressions become too cumbersome. The reconstruction matrix, 
which is given by the integrals of the linear system \eqref{intConsRec}, depends on the geometry of the tetrahedral elements in stencil $\mathcal{S}_i^s$. Therefore, since in the Lagrangian 
framework the mesh is moving in time the reconstruction matrix can \textit{not} be inverted and stored once and for all during a preprocessing stage, like in the Eulerian case. As a 
consequence, we assemble and solve the small reconstruction system \eqref{intConsRec} for each element $T^n_i$ directly at the beginning of each time step $t^n$ using optimized LAPACK 
subroutines. This makes the Lagrangian WENO reconstruction \textit{computationally more expensive} but at the same time also \textit{much less memory consuming} compared to the original 
Eulerian WENO algorithm presented in \cite{DumbserKaeser06b,Dumbser2007204}, since no reconstruction matrices are stored. 

While the mesh is moving in time, we always assume that the connectivity of the mesh and therefore also the topology of each reconstruction stencil remains constant in time. Therefore, 
the definition of the stencils $\mathcal{S}_i^s$ does \textit{not} need to be updated during the simulation. This is a very important simplification, since the stencil search may be quite 
time consuming in three space dimensions. 

Since each stencil $\mathcal{S}_i^s$ is filled with a total number of $n_e = 3 \mathcal{M}$ elements, system \eqref{intConsRec} results in an overdetermined linear system that has to be 
solved properly by either using a constrained least-squares technique (LSQ), see \cite{DumbserKaeser06b}, or a more sophisticated singular value decomposition (SVD) algorithm. In order 
to avoid ill-conditioned reconstruction matrices, see \cite{abgrall_eno}, each element in a stencil $\mathcal{S}_i^s$ is first mapped to the reference coordinate system $\xi-\eta-\zeta$ 
associated with element $T_i^n$ by using the transformation \eqref{xietaTransf} before solving \eqref{intConsRec}. 

As stated by the Godunov theorem \cite{godunov}, linear monotone schemes are at most of order one and if the scheme is required to be higher order accurate and non-oscillatory, 
it must be \textit{nonlinear}. Therefore a nonlinear formulation has to be used for the final WENO reconstruction polynomial. We first measure the \textit{smoothness} of each 
reconstruction polynomial obtained on stencil $\mathcal{S}_i^s$ by a so-called oscillation indicator $\boldsymbol{\sigma}_s$ \cite{shu_efficient_weno},
\begin{equation}
\boldsymbol{\sigma}_s = \Sigma_{lm} \hat \w^{n,s}_{l,i} \hat w^{n,s}_{m,i},
\end{equation}
which is computed on the reference element using the (universal) oscillation indicator matrix $\Sigma_{lm}$, which, according to \cite{DumbserKaeser06b}, is given by 
\begin{equation}
\Sigma_{lm} = \sum \limits_{ 1 \leq \alpha + \beta + \gamma \leq M}  \, \, \int \limits_{T_e} \frac{\partial^{\alpha+\beta+\gamma} \psi_l(\xi,\eta,\zeta)}{\partial \xi^\alpha \partial \eta^\beta \partial \zeta^\gamma} \cdot 
                                                                         \frac{\partial^{\alpha+\beta+\gamma} \psi_m(\xi,\eta,\zeta)}{\partial \xi^\alpha \partial \eta^\beta \partial \zeta^\gamma} d\xi d\eta d\zeta.   
\end{equation}  
The nonlinearity is then introduced into the scheme by the WENO weights $\boldsymbol{\omega}_s$, which read
\begin{equation}
\label{eqn.weno.weights}
\tilde{\boldsymbol{\omega}}_s = \frac{\lambda_s}{\left(\boldsymbol{\sigma}_s + \epsilon \right)^r}, \qquad 
\boldsymbol{\omega}_s = \frac{\tilde{\boldsymbol{\omega}}_s}{\sum_k \tilde{\boldsymbol{\omega}}_k},  
\end{equation} 
with the parameters $r=8$ and $\epsilon=10^{-14}$. According to \cite{DumbserKaeser06b} the linear weights are chosen as $\lambda_1=10^5$ for the central stencil and $\lambda_s=1$ for the one-sided 
stencils ($2 \leq s \leq 9$). Formula \eqref{eqn.weno.weights} is intended to be read componentwise. For a WENO reconstruction based on characteristic variables see \cite{Dumbser2007204}. 
A weighted nonlinear combination of the reconstruction polynomials obtained on each candidate stencil $\mathcal{S}_i^s$ yields the final WENO reconstruction polynomial 
and its coefficients: 
\begin{equation}
\label{eqn.weno} 
 \w_h(\x,t^n) = \sum \limits_{l=1}^{\mathcal{M}} \psi_l(\xxi) \hat \w^{n}_{l,i}, \qquad \textnormal{ with } \qquad  
 \hat \w^{n}_{l,i} = \sum_s \boldsymbol{\omega}_s \hat \w^{n,s}_{l,i}.   
\end{equation}

\paragraph*{Positivity preserving technique.} \label{sec.flattener}
Several phenomena in physics and engineering as well as many classical benchmark test cases in the Lagrangian framework are involving strong shock waves, which may lead to a loss of positivity for density and pressure in the numerical scheme. Such a problem typically occurs after carrying out the high-order reconstruction algorithm presented in Section \ref{sec.weno}, which is designed to be \textit{essentially} but not \textit{absolutely} non-oscillatory. For this reason we rely on the positivity preserving technique of Balsara \cite{BalsaraFlattener}, where a flattener variable is computed in 
order to smear out the oscillations and to bring back density and pressure values to their physically admissible range if the positivity constraint has been violated. In \cite{BalsaraFlattener} 
the equations for both hydrodynamics and magnetohydrodynamics have been considered on two- and three-dimensional Cartesian grids and in this paper we extend the method to moving unstructured 
tetrahedral meshes. 

First we have to detect those regions of the computational domain $\Omega(t)$ which are characterized by strong shocks. Let us consider a tetrahedron $T_i^n$ and its Neumann neighborhood 
$\mathcal{N}_i$, i.e. all the elements $T_j^n$ that are attached to a face of $T_i^n$. Let furthermore $\Q_i^n$ and $\Q_j^n$ be the vectors of conserved variables of element $T_i^n$ and 
its direct neighbor $T_j^n$, respectively, and let $\rho^n$ denote the density and $p^n$ the pressure. A shock can be identified by comparing 
the divergence of the velocity field $\nabla \cdot \mathbf{v}^n$ with the minimum of the sound speed $c_{i,\min}^n$ obtained by considering the element $T_i^n$ itself as well as its neighborhood 
$\mathcal{N}_i$. Hence, 
\begin{equation}
\nabla \cdot \mathbf{v}^n = \frac{1}{|T_i^n|}\sum \limits_{T_j^n \in \mathcal{N}_i}{ S_j^n\left(\mathbf{v}_j^n - \mathbf{v}_i^n \right)\cdot \mathbf{n}_{ij}^n }, \qquad c_{i,\min}^n = \min \limits_{T_j^n \in \mathcal{N}_i}{\left(c_i^n,c_j^n\right)},
\label{eq:divV}
\end{equation}
where $|T_i^n|$ represents as usual the volume of the tetrahedron $T_i^n$, $S_j^n$ denotes the surface shared between element $T_i^n$ and the neighbor $T_j^n$, $\mathbf{n}_{ij}^n$ is the associated unit  normal vector w.r.t. the surface $S_j^n$ and $c_{i,j}^n=\sqrt{\frac{\gamma p_{i,j}^n}{\rho_{i,j}^n}}$ are the sound speeds of $T_i^n$ and the neighbor element $T_j^n$, respectively, 
with $\gamma$ representing the ratio of specific heats. The divergence of the velocity field is estimated from the cell-averaged states $\Q_{i,j}^n$ and \textit{not} from the reconstructed 
states $\w_h(\x,t^n)$ obtained from \eqref{eqn.weno}.

The flattener variable $f_i^n$ is then computed according to \cite{BalsaraFlattener} as
\begin{equation}
f_i^n = \min {\left[ 1, \max {\left(0, -\frac{\nabla \cdot \mathbf{v}^n + k_1 c_{i,\min}^n}{k_1 c_{i,\min}^n}\right)}\right]},
\label{eq:flattener}
\end{equation}
with the coefficient $k_1$ that is set to the value of $k_1=0.1$ for all our computations. For rarefaction waves the divergence of the velocity field is positive, i.e. $\nabla \cdot \mathbf{v}^n\geq 0$, hence obtaining $f_i^n=0$ and leaving the reconstruction polynomial $\w_h(\x,t^n)$ as it is. Even when shocks of modest strength occur, i.e. $-k_1 c_{i,\min}^n \leq \nabla \cdot \mathbf{v}^n \leq 0$, the reconstruction remains untouched. 

In the work of Balsara \cite{BalsaraFlattener} the flattener variable is propagated even to those elements that are about to be crossed by a shock, but have still to enter the wave, i.e. the neighbors of an element which has already experienced the shock. Due to the more complex computational domain on unstructured meshes, we propose to define a \textit{node based} flattener $\tilde f_k^n$: the Voronoi neighborhood $\mathcal{V}_k$ of each vertex $k$ of element $T_i^n$ is also considered in order to propagate the flattener, hence taking into account all those elements that share vertex $k$ of tetrahedron $T_i^n$. The node based flattener results in the maximum value among the flattener values $f_j^n$ of the attached tetrahedra that has been previously computed according to \eqref{eq:flattener}:
\begin{equation}
\tilde f_k^n = \max \limits_{j \in \mathcal{V}_k} {f_j^n}.
\label{eq:flattenerNode}
\end{equation}
Each element $T_i^n$ is then assigned again with the maximum value of the node based flattener among the set $\mathcal{K}_i$ of the four vertices that define the tetrahedron $T_i^n$, i.e. 
\begin{equation}
f_i^n = \max \limits_{k \in \mathcal{K}_i} {\tilde f_k^n}. 
\label{eq:flattenerFinal}
\end{equation}
Once the flattener variable has been computed for each element of the computational domain $\Omega(t)$, the WENO reconstruction polynomials are corrected with the following expression:
\begin{equation}
\w_h(\x,t^n) := (1 - f_i^n) \psi_l(\xxi) \hat \w^{n}_{l,i} + f_i^n \cdot \Q_i^n.
\label{eq:WENOcorr}
\end{equation}

If positivity is still violated even after using \eqref{eq:WENOcorr}, then $f_i^n:=1$ is set, thus recovering a (positivity preserving) first order finite volume scheme. 
This strategy resembles to some extent the recently developed MOOD algorithm of Diot et al. \cite{MOODorg,MOODhighorder,ADERMOOD}, however, it is still used as an \textit{a priori} limiter here, 
while the MOOD approach uses an innovative \textit{a posteriori} limiting philosophy. The development of high order Lagrangian MOOD schemes will be the topic of future research.  
The presented flattener technique is by default switched off and has been used only for those test problems where it was absolutely necessary in order to run the simulation to the final time. 
We therefore explicitly state in Section \ref{sec.validation} if the flattener has been used. 

\subsection{Local space-time Discontinuous Galerkin predictor on moving curved tetrahedra} 
\label{sec.lst} 
The reconstructed polynomials $\w_h(\x,t^n)$ computed at the current time $t^n$ are then \textit{evolved} during one time step \textit{locally} within each element $T_i(t)$ without requiring any 
neighbor information. As a result, one obtains piecewise space-time polynomials of degree $M$, denoted by $\q_h(\x,t)$. This allows the scheme to achieve also high 
order of accuracy in time. For this purpose an \textit{element-local} weak space-time formulation of the governing PDE \eqref{eqn.pde.nc} is used. This approach has first been developed in 
the Eulerian framework on fixed grids by Dumbser et al. in \cite{DumbserEnauxToro,Dumbser20088209,USFORCE2,HidalgoDumbser}. Later, it has been extended to the Lagrangian framework on moving grids 
in 1D and 2D in \cite{Lagrange1D,Lagrange2D,LagrangeNC}. Here, we extend this approach for the first time to moving tetrahedral meshes in 3D. As already done in the past 
\cite{DumbserEnauxToro,HidalgoDumbser,DumbserZanotti} we use the local space-time Discontinuous Galerkin predictor method, since it is able to handle also \textit{stiff source terms}.  

Let $\mathbf{x}=(x,y,z)$ and $\boldsymbol{\xi}=(\xi,\eta,\zeta)$ be the spatial coordinate vectors defined in the physical and in the reference system, respectively, and let $\mathbf{\tilde{x}}=(x,y,z,t)$ and $\boldsymbol{\tilde{\xi}}=(\xi,\eta,\zeta,\tau)$ be the corresponding space-time coordinate vectors. Let furthermore $\theta_l=\theta_l(\boldsymbol{\tilde{\xi}})=\theta_l(\xi,\eta,\zeta,\tau)$ be a space-time basis function defined by the Lagrange interpolation polynomials passing through the space-time nodes $\boldsymbol{\tilde{\xi}}_m=(\xi_m,\eta_m,\zeta_m,\tau_m)$, which are defined by the tensor product of the spatial nodes of classical conforming high order finite elements and the Gauss-Legendre quadrature points in time.

Since the Lagrange interpolation polynomials define a \textit{nodal} basis, the functions $\theta_l$ satisfy the following interpolation property:
\begin{equation}
 \theta_l(\boldsymbol{\tilde{\xi}}_m) = \delta_{lm}, 
 \label{nodprop}
\end{equation} 
where $\delta_{lm}$ denotes the usual Kronecker symbol.
According to \cite{Dumbser20088209} the local solution $\q_h$, the fluxes $\F_h = (\f_{h}, \g_h, \h_h)$, the source term $\S_h$ and the non-conservative product $\P_h = \B(\q_h) \cdot \nabla \q_h$ are approximated within the space-time element $T_i(t) \times [t^n;t^{n+1}]$ with  
\begin{equation*}
\q_h=\q_h(\boldsymbol{\tilde{\xi}}) = \theta_{l}(\boldsymbol{\tilde{\xi}}) \, \widehat{\q}_{l,i}, \qquad \qquad  
\F_h=\F_h(\boldsymbol{\tilde{\xi}}) = \theta_{l}(\boldsymbol{\tilde{\xi}}) \, \widehat{\F}_{l,i},  
\end{equation*}
\begin{equation} 
\S_h=\S_h(\boldsymbol{\tilde{\xi}}) = \theta_{l}(\boldsymbol{\tilde{\xi}}) \, \widehat{\S}_{l,i}, \qquad \qquad  
\P_h=\P_h(\boldsymbol{\tilde{\xi}}) = \theta_{l}(\boldsymbol{\tilde{\xi}}) \, \widehat{\P}_{l,i}.
\label{thetaSol}
\end{equation}

Because of the interpolation property \eqref{nodprop} we evaluate the degrees of freedom for $\F_h$, $\S_h$ and $\P_h$ in a \textit{pointwise} manner from $\q_h$ as 
\begin{equation}
  \widehat{\F}_{l,i} = \F(\widehat{\q}_{l,i}), \quad 
  \widehat{\S}_{l,i} = \S(\widehat{\q}_{l,i}), \quad 
  \widehat{\P}_{l,i} = \P(\widehat{\q}_{l,i},\nabla \widehat{\q}_{l,i}), \quad 
  \nabla \widehat{\q}_{l,i} = \nabla \theta_{m}(\boldsymbol{\tilde{\xi}}_l) \widehat{\q}_{m,i}.   
\end{equation}  
The degrees of freedom $\nabla \widehat{\q}_{l,i}$ represent the gradient of $\q_h$ in node $\boldsymbol{\tilde{\xi}}_l$.

An \textit{isoparametric} approach is used, where the mapping between the physical space-time coordinate vector $\mathbf{\tilde{x}}$ and the reference space-time coordinate vector $\boldsymbol{\tilde{\xi}}$ is represented by the \textit{same} basis functions $\theta_l$ used for the discrete solution $\q_h$ itself. Therefore 
\begin{equation}
 \x(\boldsymbol{\tilde{\xi}}) = \theta_l(\boldsymbol{\tilde{\xi}}) \, \widehat{\x}_{l,i}, \qquad 
 t(\boldsymbol{\tilde{\xi}}) = \theta_l(\boldsymbol{\tilde{\xi}}) \, \widehat{t}_l,  
 \label{eqn.isoparametric} 
\end{equation} 
where $\widehat{\mathbf{x}}_{l,i} = (\widehat{x}_{l,i},\widehat{y}_{l,i},\widehat{z}_{l,i})$ are the degrees of freedom of the spatial physical coordinates of the moving space-time control volume, which are unknown, while $\widehat{t}_l$ denote the \textit{known} degrees of freedom of the physical time at each space-time node $\tilde{\x}_{l,i} = (\widehat{x}_{l,i}, \widehat{y}_{l,i}, \widehat{z}_{l,i}, \widehat{t}_l)$. The mapping in time is linear and simply reads
\begin{equation}
t = t_n + \tau \, \Delta t, \qquad  \tau = \frac{t - t^n}{\Delta t}, \qquad \Rightarrow \qquad \widehat{t}_l = t_n + \tau_l \, \Delta t, 
\label{timeTransf}
\end{equation} 
where $t^n$ represents the current time and $\Delta t$ is the current time step, which is computed under a classical Courant-Friedrichs-Levy number (CFL) stability condition, i.e.
\begin{equation}
\Delta t = \textnormal{CFL} \, \min \limits_{T_i^n} \frac{d_i}{|\lambda_{\max,i}|}, \qquad \forall T_i^n \in \Omega^n, 
\label{eq:timestep}
\end{equation}
with $d_i$ denoting the insphere diameter of tetrahedron $T_i^n$ and $|\lambda_{\max,i}|$ corresponding to the 
maximum absolute value of the eigenvalues computed from the solution $\Q_i^n$ in $T_i^n$. On unstructured three-dimensional meshes the CFL stability condition must satisfy the inequality $\textnormal{CFL} \leq \frac{1}{3}$.

The Jacobian of the transformation from the physical space-time element to the reference space-time element reads
\begin{equation}
J_{st} = \frac{\partial \mathbf{\tilde{x}}}{\partial \boldsymbol{\tilde{\xi}}} = \left( \begin{array}{cccc} x_{\xi} & x_{\eta} & x_{\zeta} & x_{\tau} \\ y_{\xi} & y_{\eta} & y_{\zeta} & y_{\tau} \\ z_{\xi} & z_{\eta} & z_{\zeta} & z_{\tau} \\ 0 & 0 & 0 & \Delta t \\ \end{array} \right) 
\label{Jac}
\end{equation}
and its inverse is given by 
\begin{equation}
J_{st}^{-1} = \frac{\partial \boldsymbol{\tilde{\xi}}}{\partial \mathbf{\tilde{x}}} = \left( \begin{array}{cccc} \xi_{x} & \xi_{y} & \xi_{z} & \xi_{t} \\ \eta_{x} & \eta_{y} & \eta_{z} & \eta_{t} \\ \zeta_{x} & \zeta_{y} & \zeta_{z} & \zeta_{t} \\ 0 & 0 & 0 & \frac{1}{\Delta t} \\ \end{array} \right).
\label{iJac}
\end{equation}
We point out that in the Jacobian matrix $t_\xi = t_\eta = t_\zeta = 0$ and $t_\tau = \Delta t$, as can be easily derived from the time mapping \eqref{timeTransf}.

In the following we introduce the notation adopted for the nabla operator $\nabla$ in the reference space $\boldsymbol{\xi}=(\xi,\eta,\zeta)$ and in the physical space $\mathbf{x}=(x,y,z)$:
\begin{equation}
 \nabla_{\xxi} = \left( \begin{array}{c} \frac{\partial}{\partial \xi} \\ \frac{\partial}{\partial \eta} \\ \frac{\partial}{\partial \zeta}  \end{array} \right), \qquad 
 \nabla        = \left( \begin{array}{c} \frac{\partial}{\partial x  } \\ \frac{\partial}{\partial y   } \\ \frac{\partial}{\partial z   }  \end{array} \right) = 
  \left( \begin{array}{ccc} \xi_x & \eta _x & \zeta _x \\ \xi_y & \eta_y & \zeta _y \\ \xi_z & \eta_z & \zeta _z \end{array} \right) 
  \left( \begin{array}{c} \frac{\partial}{\partial \xi} \\ \frac{\partial}{\partial \eta} \\ \frac{\partial}{\partial \zeta}  \end{array} \right)  = 
  \left( \frac{\partial \xxi}{\partial \x} \right)^T \nabla_{\xxi},
\label{not.ref.pde}
\end{equation}
and let us furthermore introduce the two integral operators
\begin{eqnarray}
\left[f,g\right]^{\tau} &=& \int \limits_{T_e} f(\xi,\eta,\zeta,\tau) g(\xi,\eta,\zeta,\tau) d\xi d\eta d\zeta, \nonumber \\
\left\langle f,g \right\rangle &=& \int \limits_{0}^{1} \int \limits_{T_e} f(\xi,\eta,\zeta,\tau)g(\xi,\eta,\zeta,\tau) d\xi d\eta d\zeta d\tau,  
\label{intOperators}
\end{eqnarray}
that denote the scalar products of two functions $f$ and $g$ over the spatial reference element $T_e$ at time $\tau$ and over the space-time reference element $T_e\times \left[0,1\right]$, respectively.

The governing PDE \eqref{eqn.pde.nc} is then reformulated in the reference coordinate system $(\xi,\eta,\zeta)$ using the inverse of the associated Jacobian matrix \eqref{iJac} with $\tau_x = \tau_y = 0$ and $\tau_t = \frac{1}{\Delta t}$ according to (\ref{timeTransf}) and adopting the gradient notation illustrated in \eqref{not.ref.pde} above:
\begin{equation}
\frac{\partial \Q}{\partial \tau} + \Delta t \left[ \frac{\partial \Q}{\partial \xxi} \cdot \frac{\partial \xxi}{\partial t} + \left( \frac{\partial \xxi}{\partial \x} \right)^T \nabla_{\xxi} \cdot \F  + \B(\Q) \cdot \left( \frac{\partial \xxi}{\partial \x} \right)^T \nabla_{\xxi} \Q \right] = \Delta t \mathbf{S}(\Q).
\label{PDECG}
\end{equation}
By introducing the following abbreviation
\begin{equation}
 \mathbf{H} = \frac{\partial \Q}{\partial \xxi} \cdot \frac{\partial \xxi}{\partial t} + \left( \frac{\partial \xxi}{\partial \x} \right)^T \nabla_{\xxi} \cdot \F  + \B(\Q) \cdot \left( \frac{\partial \xxi}{\partial \x} \right)^T \nabla_{\xxi} \Q ,
\end{equation}
Eqn. \eqref{PDECG} simplifies to
\begin{equation}
\frac{\partial \Q}{\partial \tau} + \Delta t \mathbf{H} = \Delta t \mathbf{S}(\Q).
\label{PDECGsimple}
\end{equation}
The numerical approximation of $\mathbf{H}$ is computed by the same isoparametric approach used in \eqref{thetaSol} for the solution and the flux representation, i.e.
\begin{equation}
\mathbf{H}_h = \theta_{l}(\boldsymbol{\tilde{\xi}}) \, \widehat{\mathbf{H}}_{l,i}. 
\label{eqn.stdg.rhs} 
\end{equation} 
Inserting \eqref{thetaSol} and \eqref{eqn.stdg.rhs} into \eqref{PDECG}, then multiplying Eqn. \eqref{PDECG} with the space-time test functions $\theta_k(\xxi)$ and integrating the resulting equation over the space-time reference element $T_e \times [0,1]$, one obtains a \textit{weak formulation} of the governing PDE \eqref{eqn.pde.nc}:
\begin{equation}
\left\langle \theta_k,\frac{\partial \theta_l}{\partial \tau} \right\rangle \widehat{\q}_{l,i}  
=  \left\langle \theta_k,\theta_l \right\rangle \Delta t \left( \widehat{\mathbf{S}}_{l,i} - \widehat{\mathbf{H}}_{l,i} \right). \nonumber\\ 
\label{eqn.weak.lag} 
\end{equation}
The term on the left hand side can be integrated by parts in time, yielding
\begin{equation}
 \left[ \theta_k(\xxi,1), \theta_l(\xxi,1)\right]^1 \widehat{\q}_{l,i} - \left\langle \frac{\partial \theta_k}{\partial \tau}, \theta_l \right\rangle \widehat{\q}_{l,i}  
= \left[ \theta_k(\xxi, 0), \psi_l(\xxi) \right]^0 \hat \w^n_{l,i} + \left\langle \theta_k,\theta_l \right\rangle \Delta t \left( \widehat{\mathbf{S}}_{l,i}  - \widehat{\mathbf{H}}_{l,i}\right), 
\label{LagrSTPDECG}
\end{equation}
where the initial condition of the local Cauchy problem has been introduced in a weak form.

Adopting the following more compact matrix-vector notation
\begin{equation}
\mathbf{K}_{1} = \left[ \theta_k(\xxi,1), \theta_l(\xxi,1)\right]^1 - \left\langle \frac{\partial \theta_k}{\partial \tau}, \theta_l \right\rangle, \quad 
\F_0   = \left[ \theta_k(\xxi, 0), \psi_l(\xxi) \right], \quad 
\mathbf{M} = \left\langle \theta_k,\theta_l \right\rangle, 
\end{equation}
the system \eqref{LagrSTPDECG} is reformulated as
\begin{equation}
 \mathbf{K}_1 \widehat{\q}_{l,i} = \F_0 \hat \w^n_{l,i} + \Delta t \mathbf{M}  \left( \widehat{\mathbf{S}}_{l,i} - \widehat{\mathbf{H}}_{l,i} \right).
 \label{DGsystem}
\end{equation}
Eqn. \eqref{DGsystem} constitutes an element-local nonlinear algebraic equation system for the unknown space-time expansion coefficients $\widehat{\q}_{l,i}$ which can be solved using the 
following iterative scheme 
\begin{equation}
 \widehat{\q}_{l,i}^{r+1} - \Delta t \mathbf{K}_1^{-1} \mathbf{M} \, \widehat{\mathbf{S}}^{r+1}_{l,i}  = \mathbf{K}_1^{-1} \left( \F_0 \hat \w^n_{l,i} - \Delta t \mathbf{M} \widehat{\mathbf{H}}_{l,i}^r \right), 
 \label{DGfinal} 
\end{equation}
where $r$ denotes the iteration number. In case of stiff algebraic source terms, the discretization of $\mathbf{S}$ must be \textit{implicit}, see \cite{DumbserEnauxToro,DumbserZanotti,HidalgoDumbser,Lagrange1D}. For an efficient initial guess of this iterative procedure in the case of stiff source terms see \cite{HidalgoDumbser}.

Together with the solution, we also have to evolve in time the geometry of the space-time control volume, i.e. the vertex coordinates of element $T^n_i$, whose motion is described by the ODE system 
\begin{equation}
\frac{d \mathbf{x}}{dt} = \mathbf{V}(\Q,\x,t),
\label{ODEmesh}
\end{equation}
with $\mathbf{V}=\mathbf{V}(\Q,\x,t)$ denoting the local mesh velocity. In this paper we are developing an \textit{Arbitrary-Lagrangian-Eulerian} (ALE) method, which allows the mesh velocity to be chosen independently from the fluid velocity, so that the scheme may reduce either to a pure Eulerian approach in the case where $\mathbf{V}=0$ or to a fully Lagrangian-type algorithm if $\mathbf{V}$ coincides with the local fluid velocity $\mathbf{v}$. Any other choice for the mesh velocity is possible. The velocity inside element $T_i(t)$ is also expressed in terms of the space-time basis functions $\theta_{l}$ as 
\begin{equation}
\mathbf{V}_h= \theta_{l}(\xxi,\tau) \widehat{\mathbf{V}}_{l,i}, 
\label{Vdof}
\end{equation}
with $ \widehat{\mathbf{V}}_{l,i} = \mathbf{V}(\mathbf{\hat \q}_{l,i}, \hat{\x}_{l,i}, \hat t_l)$.

The local space-time DG method is used again to solve Eqn.\eqref{ODEmesh} for the unknown coordinate vector $\widehat{\mathbf{x}}_l=(x_l,y_l,z_l)$, according to \cite{Lagrange1D,Lagrange2D,LagrangeNC}, hence 
\begin{equation}
\K_1 \widehat{\mathbf{x}}_{l,i} = \left[ \theta_k(\xxi,0), \x(\xxi,t^n) \right]^0 + \Delta t \mathbf{M} \, \widehat{\mathbf{V}}_{l,i},
\label{VCG}
\end{equation}
where $\x(\xxi,t^n)$ is given by the mapping \eqref{xietaTransf} based on the known vertex coordinates of tetrahedron $T_i^n$ at time $t^n$. The above expression is then solved by an iterative procedure \textit{together} with Eqn. \eqref{DGfinal} until the residuals of the predicted solution given by \eqref{DGfinal} and the new vertex position $\widehat{\mathbf{x}}^{r+1}_{l,i}$ at iteration $r$ are less than a prescribed tolerance, typically set to $10^{-12}$.

Once we have carried out the above procedure for all the elements of the computational domain, we end up with an \textit{element-local predictor} for the numerical solution $\q_h$, for the fluxes $\mathbf{F}_h=(\f_h,\g_h,\h_h)$, 
for the source term $\S_h$ and also for the mesh velocity $\mathbf{V}_h$. 

Then we have to update the mesh \textit{globally}, by assigning a \textit{unique} velocity vector to each node, since we do not admit discontinuities in the geometry. In the next Section \ref{sec.meshMot} a \textit{local} node solver algorithm for the velocity together with a \textit{rezoning} algorithm will be presented in detail, in order to obtain a uniquely defined vertex location at the new time level $t^{n+1}$.

\subsection{Mesh motion}
\label{sec.meshMot}
Lagrangian schemes have been designed and developed in order to compute the flow variables by moving together with the fluid. As a consequence, the computational mesh continuously changes its configuration in time, following as closely as possible the flow motion. The mesh velocity plays indeed an important role and should be evaluated very accurately using a \textit{node solver} algorithm, which assigns a velocity vector to each vertex of the mesh. A comparison between different node solver techniques can be found in \cite{LagrangeMHD}. Moreover, the flow motion may become very complex, hence highly deforming the computational elements, that are compressed, twisted or even tangled. Therefore, the challenge of any Lagrangian scheme is to preserve at the same time the excellent 
resolution properties of contact waves and material interfaces together with a good mesh quality without invalid elements. A suitable \textit{rezoning algorithm} \cite{KnuppRezoning} is typically 
used to improve the mesh quality together with a so-called \textit{relaxation algorithm} \cite{MaireRezoning} to partially recover the optimal Lagrangian accuracy where the computational elements 
are not distorted too much. In the following we present in detail the three main steps adopted in our ADER-WENO ALE finite volume schemes to move the mesh vertices to the final mesh configuration 
at the new time level $t^{n+1}$: the Lagrangian step, the rezoning step and the relaxation step. 

\subsubsection{The Lagrangian step.}
At the end of the local predictor procedure illustrated in Section \ref{sec.lst}, each vertex $k$ is assigned with several velocity vectors $\mathbf{V}_{k,j}$, each of them coming from the Voronoi neighborhood which is composed by the neighbor elements that share the common node $k$. Moving the same vertex $k$ to the next time level $t^{n+1}$ with different velocities would lead to a discontinuity in the geometry, that is not admissible in our Lagrangian algorithm. Therefore a node solver technique is adopted in order to fix a \textit{unique} velocity for each node of the computational grid. In \cite{LagrangeMHD} Boscheri et al. compare three different node solvers for unstructured triangular meshes with each other and here we extend two of them to the three-dimensional case, in particular 
the node solver $\mathcal{NS}_{cs}$ of Cheng and Shu and the node solver $\mathcal{NS}_{m}$ of Maire. 

Let $\mathcal{V}_k$ be the Voronoi neighborhood of vertex $k$, that is composed by a total number of $N_k$ neighbor elements denoted by $T_j^n$, and let furthermore $m(k)$ represent a mapping from the global node number $k$ defined in $\mathcal{T}^n_{\Omega}$ to the local vertex number in element $T_j^n$. The local velocity $\mathbf{V}_{k,j}$ computed within element $T_j^n$ is evaluated as the time integral of the high order vertex-extrapolated velocity at node $k$, i.e. 
\begin{equation}
\mathbf{V}_{k,j} = \left( \int \limits_{0}^{1} \theta_l(\xi^e_{m(k)}, \eta^e_{m(k)}, \zeta^e_{m(k)}, \tau) d \tau \right) \widehat{\mathbf{V}}_{l,j}. 
\label{NodesVel}
\end{equation}

The node solver $\mathcal{NS}_{cs}$ computes the velocity $\overline{\mathbf{V}}_k$ of vertex $k$ as a \textit{mass weighted} average velocity among its neighborhood and it reads
\begin{equation}
\overline{\mathbf{V}}_k = \frac{1}{\mu_k}\sum \limits_{T_j^n \in \mathcal{V}_k}{\mu_{k,j}\mathbf{V}_{k,j}},
\label{eqnNScs}
\end{equation}
where the local weights $\mu_{k,j}$ are defined as the product between the cell averaged value of density $\rho^n_j$ and the cell volume $|T_j^n|$, hence
\begin{equation}
\mu_{k,j}=\rho^n_j |T_j^n|, \qquad \mu_k = \sum \limits_{T_j^n \in \mathcal{V}_k}{\mu_{k,j}}.
\label{eqn.NScs.weights}
\end{equation}

\begin{figure}[htbp]
	\centering
		\includegraphics[width=0.45\textwidth]{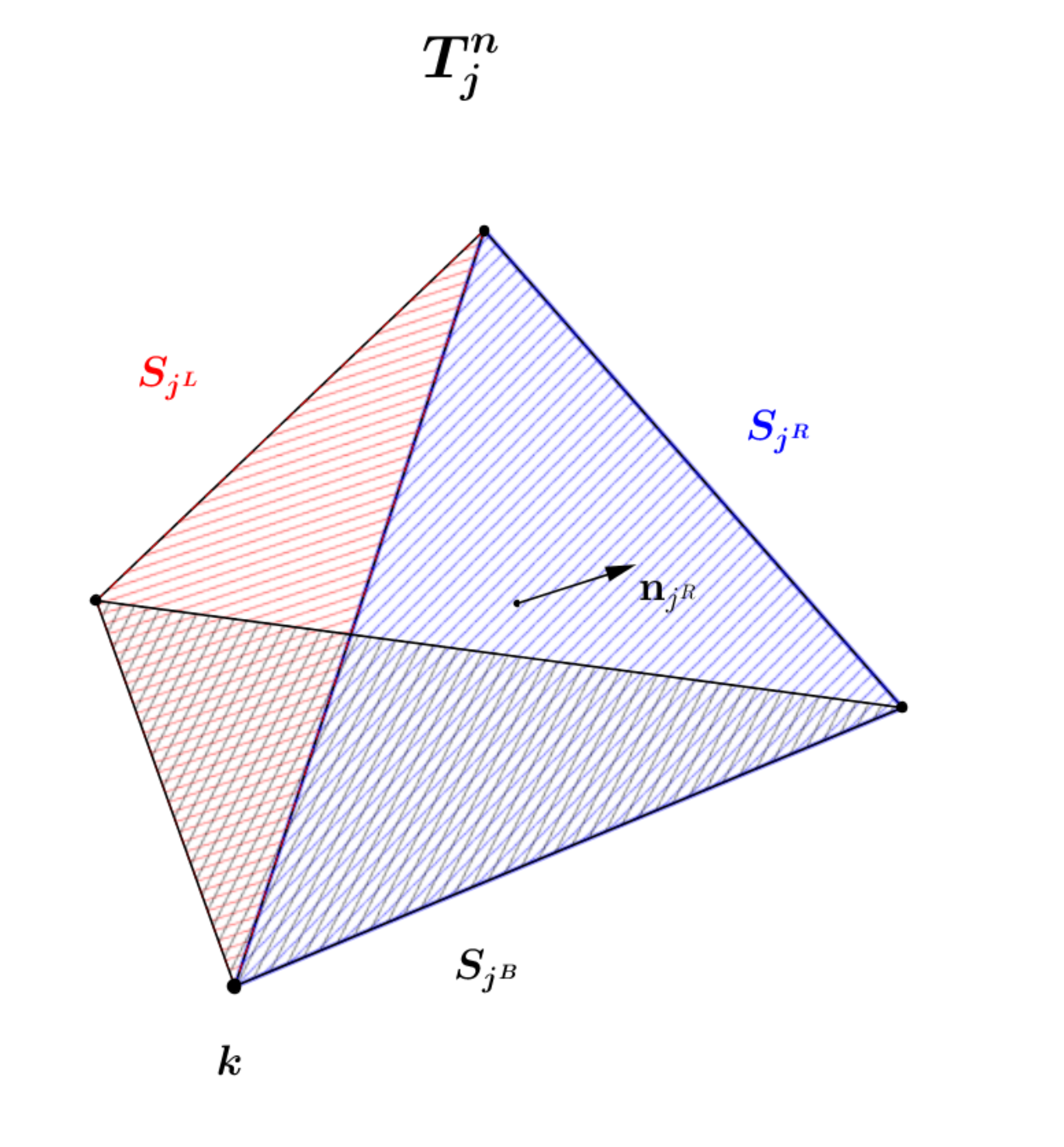}
	\caption{Geometrical notation for the node solver $\mathcal{NS}_{m}$, where only one neighbor element $T_j^n$ of node $k$ is depicted. $S_{j^R},S_{j^L},S_{j^R}$ denote one third of the total area of the faces $R,L,B$ of $T_j^n$ that share vertex $k$, while $\mathbf{n}_{j^R},\mathbf{n}_{j^L},\mathbf{n}_{j^B}$ are the corresponding outward pointing unit normal vectors.}
	\label{NS_Maire}
\end{figure}

In \cite{Maire2009,Maire2010,Maire2011} Maire et al. developed the node solver $\mathcal{NS}_{m}$ for hydrodynamics, while in \cite{Despres2009} Despr\'es presented a similar approach. 
All the details can be found in the above-mentioned references, hence we limit us here only to a brief overview of this node solver algorithm, which is based on the conservation of total 
energy in the equations for compressible hydrodynamics. According to Figure \ref{NS_Maire}, $k$ is the node index, $T_j^n$ denotes the neighbor element $j$ of vertex $k$ and the subscripts 
$({j^R},{j^L},{j^B})$ represent the three faces of tetrahedron $T_j^n$ which share node $k$, ordered adopting a counterclockwise convention. Furthermore $(S_{j^R},S_{j^L},S_{j^R})$ are 
assumed to be one third of the corresponding face areas and $(\mathbf{n}_{j^R},\mathbf{n}_{j^L},\mathbf{n}_{j^B})$ denote the associated outward pointing unit normal vectors. Finally 
$p_j$ is the fluid pressure and $c_j$ is the speed of sound for hydrodynamics. 

The total energy at the generic node $k$ is conserved only if the sum of the forces acting on node $k$ is zero, i.e.
\begin{equation}
\sum \limits_{T_j^n \in \mathcal{V}_k}{\mathbf{F}_{k,j}} = 0.
\label{eqnEnergycons}
\end{equation} 
In \eqref{eqnEnergycons} the sub-cell force $\mathbf{F}_{k,j}$ exerted by each neighbor element $T_j^n$ onto vertex $k$, is evaluated solving approximately three half Riemann problems on the faces $(j^R,j^L,j^B)$. The acoustic Riemann solver of Dukowicz et al. \cite{DukowiczRS} is used to obtain the final expression for the sub-cell force, which reads
\begin{equation}
\mathbf{F}_{k,j} = S_{k,j}p_{k,j}\mathbf{n}_{k,j} - \mathbf{M}_{k,j}\left(\overline{\mathbf{V}}_k - \mathbf{V}_{k,j}\right),
\label{eqnCellForce}
\end{equation}
with $S_{k,j}\mathbf{n}_{k,j}=S_{j^R}\mathbf{n}_{j^R}+S_{j^L}\mathbf{n}_{j^L}+S_{j^B}\mathbf{n}_{j^B}$ denoting the \textit{corner vector} related to node $k$. $\mathbf{V}_{k,j}$ represents the known vertex velocity of cell $j$ according to \eqref{NodesVel}, while $\overline{\mathbf{V}}_k$ denotes the unknown velocity of node $k$. $\mathbf{M}_{k,j}$ is a $(3 \times 3)$ symmetric positive definite matrix that is evaluated as 
\begin{equation}
\mathbf{M}_{k,j} = z_{j^R}S_{j^R}\left(\mathbf{n}_{j^R}\otimes\mathbf{n}_{j^R}\right) + z_{j^L}S_{j^L}\left(\mathbf{n}_{j^L}\otimes\mathbf{n}_{j^L}\right) + z_{j^L}S_{j^B}\left(\mathbf{n}_{j^B}\otimes\mathbf{n}_{j^B}\right),
\label{eqnMkj}
\end{equation}
where $z_{j}=\rho_j c_j$ is the acoustic impedance. The equation for the total energy conservation \eqref{eqnEnergycons} can be reformulated using the expression for the sub-cell force \eqref{eqnCellForce}, hence obtaining a linear algebraic system for the unknown node velocity $\overline{\mathbf{V}}_k$: 
\begin{equation}
\mathbf{M}_k \overline{\mathbf{V}}_k = \sum \limits_{T_j^n \in \mathcal{V}_k}{\left(S_{k,j}p_{k,j}\mathbf{n}_{k,j} + \mathbf{M}_{k,j}\mathbf{V}_{k,j}\right)}, \quad \mathbf{M}_k = \sum \limits_{T_j^n \in \mathcal{V}_k}{\mathbf{M}_{k,j}}.
\label{eq:NS_mLinSyst}
\end{equation}
Since matrix $\mathbf{M}_k$ is always invertible, this system admits a unique solution and the node velocity can always be evaluated. Instead of taking the above-defined acoustic impedance, one can compute it as originally proposed by Dukowicz in \cite{DukowiczRS}:
\begin{equation}
z_{j^+} = \rho_j\left[c_j+\Gamma_j |\left(\overline{\mathbf{V}}_k-\mathbf{V}_{k,j}\right)\cdot \mathbf{n}_{j^+}|\right],
\label{eqnZnonlin}
\end{equation}
where $\Gamma_j=\frac{\gamma+1}{2}$ is a material dependent parameter which is a function of the ratio of specific heats $\gamma$. In this case the system \eqref{eq:NS_mLinSyst} becomes nonlinear, due to the dependency of the acoustic impedance on the unknown node velocity, and a suitable iterative algorithm has to be used to obtain the solution.

Once a unique velocity $\overline{\mathbf{V}}_k$ has been defined for each node $k$ of the mesh, the new \textit{Lagrangian coordinates} $\mathbf{X}^{Lag}_{k}$ are computed as
\begin{equation} 
	\mathbf{X}^{Lag}_{k}	= \mathbf{X}^{n}_{k}	+ \Delta t \, \overline{\mathbf{V}}_k, 
	\label{eqn.vertex.update}
\end{equation}
with $\mathbf{X}^{n}_{k}$ representing the coordinates of node $k$ at the current time level $t^n$.

\subsubsection{The rezoning step.}
The Lagrangian step allows the nodes to follow the fluid motion as closely as possible. However, this may lead to bad quality elements, where the Jacobians become very small or even negative. 
This either drastically decreases the admissible timestep, according to \eqref{eq:timestep}, or even leads to a failure of the computation. Therefore, also a \textit{rezoned} position should 
be computed for each node $k$ in order to improve the \textit{local} mesh quality without taking into account any physical information. We use a different treatment for internal nodes and 
boundary nodes. Specifically, the rezoning algorithm presented in \cite{KnuppRezoning,MaireRezoning} is adopted for inner nodes, while a variant of the feasible set method proposed by 
Berndt et al. \cite{ShashkovFSM} is used for the boundary nodes. 

The rezoning algorithm aims at improving the mesh quality \textit{locally}, i.e. in the Voronoi neighborhood $\mathcal{V}_k$ of node $k$ considering all the neighbor elements $T_j^{n+1}$, which for sake of simplicity will be addressed by $j$. The starting point is the Lagrangian coordinate vector $\mathbf{X}^{Lag}_{k}$ obtained at the end of the Lagrangian step. The rezoning procedure consists in optimizing a goal function $\mathcal{K}_k$ that has to be defined for each node $k$ as
\begin{equation}
\mathcal{K}_k = \sum\limits_{T_j^{n+1} \in \mathcal{V}_k}{ \kappa_{j} },
\end{equation}
where $\kappa_j$ is the condition number of the Jacobian matrix $\mathbf{J}_{j}$ of the mapping from the reference tetrahedron to the physical element $j$:
\begin{equation}
\mathbf{J}_{j} = \left( \begin{array}{ccc} x_{j,2}-x_k & y_{j,2}-y_k & z_{j,2}-z_k \\ x_{j,3}-x_k & y_{j,3}-y_k & z_{j,3}-z_k \\ x_{j,4}-x_k & y_{j,4}-y_k & z_{j,4}-z_k  \end{array} \right).
\label{eq:locJac}
\end{equation}
In \eqref{eq:locJac} the coordinate vector $\mathbf{x}_{j,l}=(x_{j,l},y_{j,l},z_{j,l})$ represents the four nodes $l=1,2,3,4$ of the neighbor tetrahedron $T_j^{n+1}$, which are counterclockwise ordered in such a way that node $k$ corresponds to $l=1$. Then, the condition number of matrix $\mathbf{J}_{j}$ is given by 
\begin{equation}
\kappa_j = \left\| \mathbf{J}_{j}^{-1} \right\| \left\| \mathbf{J}_{j} \right\|.
\label{eq:CondNumb}
\end{equation} 
The goal function $\mathcal{K}_k$  is computed according to \cite{KnuppRezoning} as the sum of the local condition numbers of the neighbors, i.e.
\begin{equation}
\mathcal{K}_k = \sum\limits_{T_j^{n+1} \in \mathcal{V}_k}{ \kappa_{j} },
\end{equation}
and its minimization leads to a \textit{locally} optimal position of the free node $k$. As proposed in \cite{MaireRezoning}, the optimized \textit{rezoned coordinates} $\mathbf{x}_k^{Rez}$ for vertex $k$ are computed using the first step of a Newton algorithm, hence
\begin{equation}
 \mathbf{x}_k^{Rez} = \mathbf{x}_k^{Lag} - \mathbf{H}_k^{-1}\left(\mathcal{K}_k\right) \cdot \nabla \mathcal{K}_k,
\label{eqn.vertex.rez}
\end{equation}
where $\mathbf{H}_k$ and $\nabla \mathcal{K}_k$ represent the Hessian and the gradient of the goal function $\mathcal{K}_k$, respectively: 
\begin{equation}
\mathbf{H}_k = \sum\limits_{T_j^{n+1} \in \mathcal{V}_k}{\left( \begin{array}{ccc} 
\frac{\partial^2 \kappa_{j}}{\partial x^2} & 
\frac{\partial^2 \kappa_{j}}{\partial x \partial y} & 
\frac{\partial^2 \kappa_{j}}{\partial x \partial z} \\ 
\frac{\partial^2 \kappa_{j}}{\partial y \partial x} & 
\frac{\partial^2 \kappa_{j}}{\partial y^2} & 
\frac{\partial^2 \kappa_{j}}{\partial y \partial z} \\ 
\frac{\partial^2 \kappa_{j}}{\partial z \partial x} & 
\frac{\partial^2 \kappa_{j}}{\partial z \partial y} & 
\frac{\partial^2 \kappa_{j}}{\partial z^2} 
\end{array} \right)}, \quad \nabla \mathcal{K}_k = \sum\limits_{T_j^{n+1} \in \mathcal{V}_k}{\left( \frac{\partial \kappa_{j}}{\partial x}, \frac{\partial \kappa_{j}}{\partial y}, \frac{\partial \kappa_{j}}{\partial z}\right)}.
\label{eqn.HessGrad}
\end{equation}

For the boundary nodes we present a simplified but very efficient version of the feasible set method proposed in \cite{ShashkovFSM} for two-dimensional unstructured meshes. The original feasible set method has been designed in order to find the convex polygon on which a vertex can lie without invalid elements in its neighborhood. In three space dimensions such an algorithm becomes very complex and high demanding in terms of computational efforts. In our simplified procedure the rezoned coordinates $\mathbf{x}_k^{Rez,b}$ of the boundary node $k$ is evaluated as a volume weighted average among the barycenter coordinates $\mathbf{x}_{c,j}^{Lag}$ of each neighbor element $j$, which is
\begin{equation}
 \mathbf{x}_k^{Rez,b} = \frac{1}{\alpha_k} \sum\limits_{T_j^{n+1} \in \mathcal{V}_k}{\mathbf{x}_{c,j}^{Lag} \cdot \alpha_{k,j}},
\label{eq:RezBoundary}
\end{equation}
with the weights 
\begin{equation}
\alpha_{k,j}=|T_j^{n+1}|, \qquad \alpha_k = \sum \limits_{T_j^{n+1} \in \mathcal{V}_k}{\alpha_{k,j}} 
\end{equation}
and the barycenter defined as usual as
\begin{equation}
\mathbf{x}_{c,j}^{n+1} = \frac{1}{4}\sum{\mathbf{x}_k^{Lag}}.
\end{equation}

\subsubsection{The relaxation step.}
Since our ALE scheme is supposed to be as Lagrangian as possible, we do not want to rezone the mesh nodes where it is not strictly necessary in order to carry on the computation. Therefore the 
\textit{final node position} $\mathbf{X}_k^{n+1}$ is obtained applying the relaxation algorithm of Galera et al. \cite{MaireRezoning}, that performs a convex combination between the Lagrangian 
position and the rezoned position of node $k$, hence 
\begin{equation}
\mathbf{X}_k^{n+1} = \mathbf{X}_k^{Lag} + \omega_k \left( \mathbf{X}_k^{Rez} - \mathbf{X}_k^{Lag} \right),
\label{eqn.relaxation}
\end{equation}
where $\omega_k$ is a node-based coefficient associated to the deformation of the Lagrangian grid over the time step $\Delta t$. The values for $\omega_k$ are bounded in the interval $[0,1]$, so that when $\omega_k=0$ a fully Lagrangian mesh motion occurs, while if $\omega_k=1$ the new node location is defined by the pure rezoned coordinates $\mathbf{X}_k^{Rez}$. We point out that the coefficient $\omega_k$ is designed to result in $\omega_k=0$ for rigid body motion, namely rigid translation and rigid rotation, where no element deformation occurs. Further details about the computation of $\omega_k$ can be found in \cite{MaireRezoning}.

\subsection{Finite volume scheme}
\label{sec.SolAlg}
\begin{figure}[!htbp]
	\begin{center} 
	\includegraphics[width=1.2\textwidth]{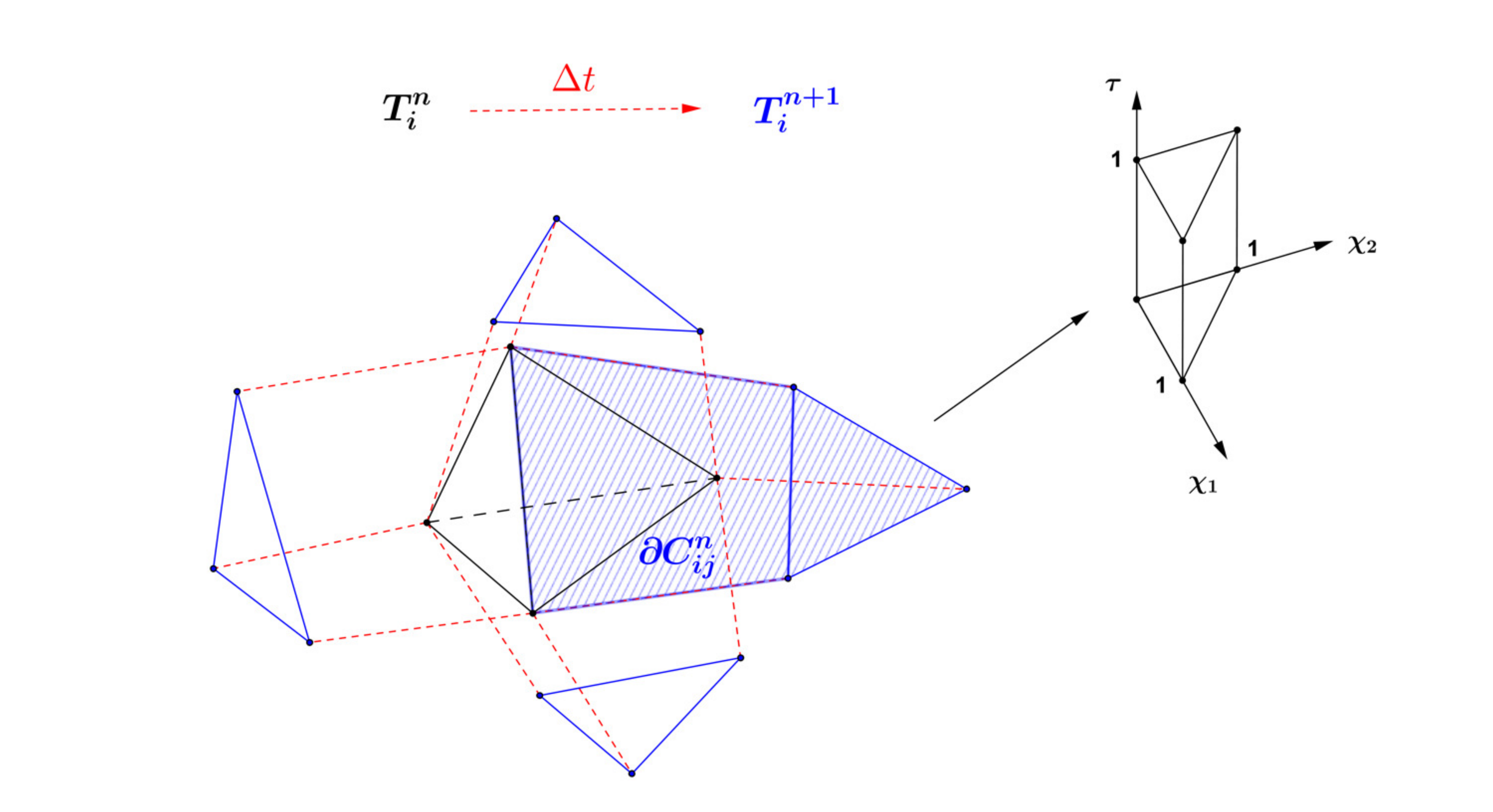}
	\caption{Physical space-time element (a) and parametrization of the lateral space-time sub-face $\partial C_{ij}^n$  (b). The dashed red lines denote the evolution in time of the faces of the tetrahedron, whose configuration at the current time level $t^n$ and at the new time level $t^{n+1}$ is depicted in black and blue, respectively.}
	\label{fig:STelem3D}
	\end{center} 
\end{figure}

In order to develop a Lagrangian finite volume schemes on moving tetrahedra, we adopt the same approach used for our ALE algorithm in two space dimensions presented in \cite{Lagrange2D,LagrangeNC}. There, the governing PDE \eqref{eqn.pde.nc} is reformulated more compactly using a space-time divergence operator $\tilde \nabla$, hence obtaining
\begin{equation}
\tilde \nabla \cdot \tilde{\F} + \tilde \B(\Q) \cdot \tilde \nabla \Q = \mathbf{S}(\Q),  \qquad \tilde \nabla  = \left( \frac{\partial}{\partial x}, \, \frac{\partial}{\partial y}, \, \frac{\partial}{\partial z}, \, \frac{\partial}{\partial t} \right)^T,
\label{eqn.st.pde}
\end{equation}
where the space-time flux tensor $\tilde{\F}$ and the system matrix $\tilde{\B}$ explicitly read
\begin{equation}
\tilde{\F}  = \left( \mathbf{f}, \, \mathbf{g}, \, \mathbf{h}, \, \Q \right), \qquad
\tilde{\B}  = ( \B_1, \B_2, \B_3, 0).  
\label{eqn.st.mat} 
\end{equation}
For the computation of the state vector at the new time level $\Q^{n+1}$, the balance law \eqref{eqn.st.pde} is integrated over a \textit{four-dimensional} space-time control volume $\mathcal{C}^n_i = T_i(t) \times \left[t^{n}; t^{n+1}\right]$, i.e. 
\begin{equation}
 \int \limits_{\mathcal{C}^n_i} \tilde \nabla \cdot \tilde{\F} \, d\mathbf{x} dt + 
\int\limits_{\mathcal{C}^n_i} \tilde{\B}(\Q) \cdot \tilde \nabla \Q \, d\mathbf{x} dt = 
\int\limits_{\mathcal{C}^n_i} \S(\Q) \, d\mathbf{x} dt.   
\label{STPDE}
\end{equation} 
Application of the theorem of Gauss yields
\begin{equation}
\int \limits_{\partial \mathcal{C}^{n}_i} \tilde{\F} \cdot \ \mathbf{\tilde n} \, dS +  
\int\limits_{\mathcal{C}^n_i} \tilde{\B}(\Q) \cdot \tilde \nabla \Q \, d\mathbf{x} dt = 
\int\limits_{\mathcal{C}^n_i} \S(\Q) \, d\mathbf{x} dt,   
\label{STPDEgauss}
\end{equation}
where the space-time volume integral on the left of \eqref{STPDE} has been rewritten as the sum of the fluxes computed over the \textit{three-dimensional} space-time volume $\partial \mathcal{C}^n_i$, given by the evolution of each face of element $T_i(t)$ within the timestep $\Delta t$, as depicted in Figure \ref{fig:STelem3D}. The symbol $\mathbf{\tilde n} = (\tilde n_x,\tilde n_y,\tilde n_z,\tilde n_t)$ denotes the outward pointing space-time unit normal vector on the space-time face $\partial C^n_i$.   

Since the algorithm is required to deal with both, conservative and non-conservative hyperbolic systems, we use a \textit{path-conservative} approach 
to integrate the non-conservative product, see \cite{Toumi1992,Pares2006,Castro2006,Castro2008,Rhebergen2008,ADERNC,USFORCE2,OsherNC,LagrangeNC}. 
One thus obtains 
\begin{equation}
\int \limits_{\partial \mathcal{C}^{n}_i} \left( \tilde{\F} + \tilde{\D} \right) \cdot \ \mathbf{\tilde n} \, dS + \! \! 
\int \limits_{\mathcal{C}^n_i \backslash \partial \mathcal{C}^n_i} \! \! \! \tilde{\B}(\Q) \cdot \tilde \nabla \Q \, d\mathbf{x} dt = 
\int\limits_{\mathcal{C}^n_i} \S(\Q) \, d\mathbf{x} dt,   
\label{I1}
\end{equation}
where a new term $\tilde{\D}$ has been introduced in order to take into account potential jumps of the solution $\Q$ on the space-time element boundaries $\partial \mathcal{C}^n_i$.  This term is computed by the path integral
\begin{equation}
 \tilde{\D} \cdot  \mathbf{\tilde n}  = \halb \int \limits_0^1 \tilde{\B}\left(\Path(\Q^-,\Q^+,s)\right) \cdot \mathbf{\tilde n} \, \frac{\partial \Path}{\partial s} \, ds. 
 \label{eqn.pathint} 
\end{equation}
The integration path $\Path$ in \eqref{eqn.pathint} is chosen to be a simple straight-line segment \cite{Pares2006,Castro2006,USFORCE2,OsherNC}, although other choices are possible. Therefore it reads
\begin{equation}
 \Path = \Path(\Q^-,\Q^+,s) = \Q^- + s (\Q^+ - \Q^-), 
 \label{eqn.segpath} 
\end{equation} 
and the jump term \eqref{eqn.pathint} simply reduces to 
\begin{equation}
 \tilde{\D} \cdot  \mathbf{\tilde n}  = \halb  \left( \int \limits_0^1 \tilde{\B}\left(\Path(\Q^-,\Q^+,s)\right) \cdot \mathbf{\tilde n} \, ds \right) \left( \Q^+ - \Q^- \right), 
 \label{eqn.pathint.seg} 
\end{equation} 
with $\left(\Q^-,\Q^+\right)$ representing the two vectors of conserved variables within element $T_i^n$ and its direct neighbor $T_j^n$, respectively.

Let $\mathcal{N}_i$ denote the \textit{Neumann neighborhood} of tetrahedron $T_i(t)$, which is the set of directly adjacent neighbors $T_j(t)$ that share a common face $\partial T_{ij}(t)$ with tetrahedron $T_i(t)$. The space-time volume $\partial C^n_i$ is composed by four space-time sub-volumes $\partial C^n_{ij}$, each of them defined for each face of tetrahedron $T_i(t)$ as depicted in Figure \ref{fig:STelem3D}, and two more space-time sub-volumes, $T_i^{n}$ and $T_i^{n+1}$, that represent the tetrahedron configuration at times $t^n$ and $t^{n+1}$, respectively. Hence, the space-time volume $\partial C^n_i$ involves overall a total number of six space-time sub-volumes, i.e.
\begin{equation}
\partial C^n_i = \left( \bigcup \limits_{T_j(t) \in \mathcal{N}_i} \partial C^n_{ij} \right) 
\,\, \cup \,\, T_i^{n} \,\, \cup \,\, T_i^{n+1}.  
\label{dCi}
\end{equation} 
Each of the space-time sub-volumes is mapped to a reference element in order to simplify the integral computation. For the configurations at the current and at the new time level, $T_i^{n}$ and $T_i^{n+1}$, we use the mapping \eqref{xietaTransf} with $(\xi,\eta,\zeta) \in \left[0;1\right]$. The space-time unit normal vectors simply read $\mathbf{\tilde n} = (0,0,0,-1)$ for $T_i^{n}$ and $\mathbf{\tilde n} = (0,0,0,1)$ for $T_i^{n+1}$, since these volumes are orthogonal to the time coordinate. For the lateral sub-volumes $\partial C^n_{ij}$ we adopt a linear parametrization to map the physical volume to a four-dimensional space-time reference prism, as shown in Figure \ref{fig:STelem3D}. Starting from the old vertex coordinates $\mathbf{X}_{ik}^n$ and the new ones $\mathbf{X}_{ik}^{n+1}$, that are \textit{known} from the mesh motion algorithm described in Section \ref{sec.meshMot}, the lateral sub-volumes are parametrized using a set of linear basis functions $\beta_k(\chi_1,\chi_2,\tau)$ that are defined on a local reference system $(\chi_1,\chi_2,\tau)$ which is oriented orthogonally w.r.t. the face $\partial T_{ij}(t)$ of tetrahedron $T_i^n$, e.g. the reference time coordinate $\tau$ is orthogonal to the reference space coordinates $(\chi_1,\chi_2)$ that lie on $\partial T_{ij}(t)$. The temporal mapping is simply given by $t = t^n + \tau \, \Delta t$, hence $t_{\chi_1} = t_{\chi_2} = 0$ and $t_\tau = \Delta t$. The lateral space-time volume $\partial C_{ij}^n$ is defined by six vertices of physical coordinates $\mathbf{\tilde{X}}_{ij,k}^n$. The first three vectors $(\mathbf{X}^n_{ij,1},\mathbf{X}^n_{ij,2},\mathbf{X}^n_{ij,3})$ are the nodes defining the common face $\partial T_{ij}(t^n)$ at time $t^n$, while the same procedure applies at the new time level $t^{n+1}$. Therefore the six vectors $\mathbf{\tilde{X}}_{ij,k}^n$ are given by
\begin{equation*}
\mathbf{\tilde{X}}_{ij,1}^n = \left( \mathbf{X}^n_{ij,1}, t^n \right), \qquad 
\mathbf{\tilde{X}}_{ij,2}^n = \left( \mathbf{X}^n_{ij,2}, t^n \right), \qquad
\mathbf{\tilde{X}}_{ij,3}^n = \left( \mathbf{X}^n_{ij,3}, t^n \right),
\end{equation*} 
\begin{equation}
\mathbf{\tilde{X}}_{ij,4}^n = \left( \mathbf{X}^{n+1}_{ij,2}, t^{n+1} \right), \qquad
\mathbf{\tilde{X}}_{ij,5}^n = \left( \mathbf{X}^{n+1}_{ij,1}, t^{n+1} \right), \quad 
\mathbf{\tilde{X}}_{ij,6}^n = \left( \mathbf{X}^{n+1}_{ij,1}, t^{n+1} \right),  
\label{eqn.lateralnodes} 
\end{equation}
and the parametrization for $\partial C_{ij}^n$ reads
\begin{equation}
\partial C_{ij}^n = \mathbf{\tilde{x}} \left( \chi_1,\chi_2,\tau \right) = 
 \sum\limits_{k=1}^{6}{\beta_k(\chi_1,\chi_2,\tau) \, \mathbf{\tilde{X}}_{ij,k}^n },	
\label{SurfPar}
\end{equation}    
with $0 \leq \chi_1 \leq 1$, $0 \leq \chi_2 \leq 1-\chi_1$ and $0 \leq \tau \leq 1$. The basis functions $\beta_k(\chi_1,\chi_2,\tau)$ are given by 
\begin{eqnarray}
\beta_1(\chi_1,\chi_2,\tau) &= (1-\chi_1-\chi_2)(1-\tau),&  \quad \beta_4(\chi_1,\chi_2,\tau) = (1-\chi_1-\chi_2)(\tau) \nonumber \\
\beta_2(\chi_1,\chi_2,\tau) &= \chi_1(1-\tau),&             \quad \beta_5(\chi_1,\chi_2,\tau) = \chi_1\tau, \nonumber \\
\beta_3(\chi_1,\chi_2,\tau) &= \chi_2(1-\tau),&             \quad \beta_6(\chi_1,\chi_2,\tau) = \chi_2\tau.
\label{eq:BetaBaseFunc}
\end{eqnarray}

The coordinate transformation is associated with a matrix $\mathcal{T}$ that reads
\begin{equation}
\mathcal{T} = \left( \hat{\mathbf{e}}, \frac{\partial \mathbf{\tilde{x}}}{\partial \chi_1}, \frac{\partial \mathbf{\tilde{x}}}{\partial \chi_2}, \frac{\partial \mathbf{\tilde{x}}}{\partial \tau} \right)^T,
\label{eqn.JacSTreference}
\end{equation}
with $\hat{\mathbf{e}}=(\hat{\mathbf{e}}_1,\hat{\mathbf{e}}_2,\hat{\mathbf{e}}_3,\hat{\mathbf{e}}_4)$ and where $\hat{\mathbf{e}}_p$ represents the unit vector aligned with the $p$-th axis of the 
physical coordinate system $(x,y,z,t)$. In the following $\tilde{x}_q$ denotes the $q$-th component of vector $\mathbf{\tilde{x}}$.   
The determinant of $\mathcal{T}$ produces at the same time the space-time volume $| \partial C_{ij}^n|$ of the space-time sub-face $\partial C_{ij}^n$ and the space-time normal vectors 
$\mathbf{\tilde n}_{ij}$, as    
\begin{equation}
 \mathbf{\tilde n}_{ij} = \left( \epsilon_{pqrs} \, \hat{\mathbf{e}}_p \, \frac{\partial {\tilde{x}_q}}{\partial \chi_1} \, \frac{\partial {\tilde{x}_r}}{\partial \chi_2} \, \frac{\partial {\tilde{x}_s}}{\partial \tau} \right)/ | \partial C_{ij}^n|,
 \label{n_lateral1}
\end{equation}
where the \textit{Levi-Civita} symbol has been used according to the usual definition
\begin{equation}
\epsilon_{pqrs} = \left\{ \begin{array}{l} +1, \quad \textnormal{if $(p,q,r,s)$ is an \textit{even} permutation of $(1,2,3,4)$}, \\
																					 -1, \quad \textnormal{if $(p,q,r,s)$ is an \textit{odd} permutation of $(1,2,3,4)$}, \\
																					 0,  \quad \textnormal{otherwise,}
                      \end{array}  \right.  
\label{eqn.LeviCivita}
\end{equation}
and with 
\begin{equation}
| \partial C_{ij}^n| = \left\| \epsilon_{pqrs} \, \hat{\mathbf{e}}_p \, \frac{\partial {\tilde{x}_q}}{\partial \chi_1} \, \frac{\partial {\tilde{x}_r}}{\partial \chi_2} \, \frac{\partial {\tilde{x}_s}}{\partial \tau} \right\|. 
\label{eqn.det_lateral}
\end{equation}


The final one-step ALE finite volume scheme takes the following form: 
\begin{equation}
|T_i^{n+1}| \, \Q_i^{n+1} = |T_i^n| \, \Q_i^n - \sum \limits_{T_j \in \mathcal{N}_i} \,\, {\int \limits_0^1 \int \limits_0^1 \int \limits_{0}^{1-\chi_1} 
| \partial C_{ij}^n| \tilde{\G}_{ij} \, d\chi_2 d\chi_1 d\tau}
+ \int \limits_{\mathcal{C}_i^n \backslash \partial \mathcal{C}_i^n}  \left( \S_h - \P_h \right) \, d\mathbf{x} dt, 
\label{PDEfinal}
\end{equation}
where the term $\tilde{\G}_{ij} \cdot \mathbf{\tilde n}_{ij}$ contains the Arbitrary-Lagrangian-Eulerian numerical flux function as well as the path-conservative jump term, hence allowing the discontinuity of the predictor solution $\mathbf{q}_h$ that occurs at the space-time sub-face $\partial C_{ij}^n$ to be properly resolved. The volume and surface integrals appearing in \eqref{PDEfinal} are approximated using multidimensional Gaussian quadrature rules, see \cite{stroud} for details. The term $\tilde{\G}_{ij}$ can be evaluated using a simple ALE Rusanov-type scheme \cite{Lagrange1D} as
\begin{equation}
  \tilde{\G}_{ij} =  
  \frac{1}{2} \left( \tilde{\F}(\q_h^+) + \tilde{\F}(\q_h^-)  \right) \cdot \mathbf{\tilde n}_{ij} +  
  \frac{1}{2} \left( \int \limits_0^1 \tilde{\B}(\Path)\cdot \mathbf{\tilde n} \ ds - |\lambda_{\max}| \mathbf{I} \right) \left( \q_h^+ - \q_h^- \right), 
  \label{eqn.rusanov} 
\end{equation}
where $\q_h^-$ and $\q_h^+$ are the local space-time predictor solution inside element $T_i(t)$ and the neighbor $T_j(t)$, respectively, and $|\lambda_{\max}|$ denotes the maximum absolute value 
of the eigenvalues of the matrix $\tilde{\A} \cdot \mathbf{\tilde n}$ in space-time normal direction. Using the normal mesh velocity $\mathbf{V} \cdot \mathbf{n}$, matrix 
$\tilde{\A}_{\mathbf{\tilde n}}$ reads   
\begin{equation} 
\tilde{\A}_{\mathbf{\tilde n}}  = \tilde{\A} \cdot \mathbf{\tilde n} = \left( \sqrt{\tilde n_x^2 + \tilde n_y^2 + \tilde n_z^2} \, \right) \left[ \left( \frac{\partial \mathbf{F}}{\partial \Q} + \B \right) \cdot \mathbf{n} - 
(\mathbf{V} \cdot \mathbf{n}) \,  \mathbf{I} \right],  
\end{equation}
with $\mathbf{I}$ denoting the $\nu \times \nu$ identity matrix, $\A = \partial \F / \partial \Q + \B$ representing the classical Eulerian system matrix and $\mathbf{n}$ being the spatial 
unit normal vector given by 
\begin{equation}
\mathbf{n} = \frac{(\tilde n_x, \tilde n_y, \tilde n_z)^T}{\sqrt{\tilde n_x^2 + \tilde n_y^2 + \tilde n_z^2}}.  
\end{equation} 


The numerical flux term $\tilde{\G}_{ij}$ can be also computed relying on a more sophisticated Osher-type scheme \cite{osherandsolomon}, introduced by Dumbser et al. for the Eulerian framework in \cite{OsherUniversal,OsherNC} and then extended to moving meshes for conservative \cite{Lagrange1D,Lagrange2D} and non-conservative hyperbolic balance laws \cite{LagrangeNC}. It reads 
\begin{equation}
  \tilde{\G}_{ij} =  
  \frac{1}{2} \left( \tilde{\F}(\q_h^+) + \tilde{\F}(\q_h^-)  \right) \cdot \mathbf{\tilde n}_{ij}  +  
  \frac{1}{2} \left( \int \limits_0^1 \left( \tilde{\B}(\Path) \cdot \mathbf{\tilde n} -  \left| \tilde{\A}_{\mathbf{\tilde n}}(\Path) \right| \right)  \, ds \, \right) \left( \q_h^+ - \q_h^- \right),  
  \label{eqn.osher} 
\end{equation}
where the matrix absolute value operator is computed as usual as
\begin{equation}
 |\mathbf{A}| = \mathbf{R} |\boldsymbol{\Lambda}| \mathbf{R}^{-1},  \qquad |\boldsymbol{\Lambda}| = \textnormal{diag}\left( |\lambda_1|, |\lambda_2|, ..., |\lambda_\nu| \right),  
\end{equation} 
with the right eigenvector matrix $\mathbf{R}$ and its inverse $\mathbf{R}^{-1}$. According to \cite{OsherNC,OsherUniversal} Gaussian quadrature formulae of sufficient accuracy are adopted to evaluate the path integral present in \eqref{eqn.osher}.

Furthermore integration over a closed space-time control volume as done in the scheme presented above automatically respects the geometric conservation law (GCL), since application of Gauss' theorem yields
\begin{equation}
 \int_{\partial \mathcal{C}_i^n} \mathbf{\tilde n} \, dS = 0.
 \label{eqn.gcl} 
\end{equation}

\section{Test problems}
\label{sec.validation} 
\vspace{-2pt}
In order to validate the unstructured three-dimensional ALE ADER-WENO schemes presented in this paper we solve in the following a set of test problems using different hyperbolic systems of 
governing equations that can all be cast into form \eqref{eqn.pde.nc}. We will consider the Euler equations of compressible gas dynamics, the equations of ideal classical magnetohydrodynamics 
(MHD) as well as the Baer-Nunziato model of compressible multi-phase flows with relaxation source terms, hence dealing with both conservative and non-conservative hyperbolic PDE.

For the Euler and ideal classical MHD equations we always use the node solver $\mathcal{NS}_{m}$ to compute the mesh velocity, according to \eqref{eq:NS_mLinSyst}, while the simple algorithm $\mathcal{NS}_{cs}$ given by \eqref{eqnNScs} is adopted for the Baer-Nunziato model. For each of the test cases of the Euler and MHD equations we choose the local mesh velocity as the local fluid velocity, 
hence   
\begin{equation}
 \mathbf{V} = \mathbf{v}. 
\end{equation}
Furthermore we normally do \textit{not} use the flattener technique illustrated in Section \ref{sec.flattener}, but we explicitly write when it has been activated.

\subsection{The Euler equations of compressible gas dynamics} 
Let $\Q=(\rho, \rho u, \rho v, \rho w, \rho E)$ be the vector of conserved variables with $\rho$ denoting the fluid density, $\mathbf{v}=(u,v,w)$ representing the velocity vector and $\rho E$ being the total energy density. Let furthermore $p$ be the fluid pressure and $\gamma$ the ratio of specific heats of the ideal gas, so that the speed of sound is $c=\sqrt{\frac{\gamma p}{\rho}}$. The three-dimensional Euler equations of compressible gas dynamics can be cast into form \eqref{eqn.pde.nc}, with the state vector $\Q$ previously defined and the flux tensor $\tens{\F}=(\f,\g,\h)$ given by   
\begin{equation}
\label{eulerTerms}
\f = \left( \begin{array}{c} \rho u \\ \rho u^2 + p \\ \rho uv \\ \rho uw \\ u(\rho E + p) \end{array} \right), \quad \g = \left( \begin{array}{c} \rho v \\ \rho uv \\ \rho v^2 + p  \\ \rho vw \\ v(\rho E + p) \end{array} \right), \quad \h = \left( \begin{array}{c} \rho w \\ \rho uw \\ \rho vw \\ \rho w^2 + p  \\ w(\rho E + p) \end{array} \right). 
\end{equation}
The term $\tens{\B}$ appearing in \eqref{eqn.pde.nc} is zero for this hyperbolic conservation law, because the system does not involve any non-conservative term. The system is then closed by 
the equation of state for an ideal gas, which reads 
\begin{equation}
\label{eqn.eos} 
p = (\gamma-1)\left(\rho E - \frac{1}{2} \rho \mathbf{v}^2 \right).  
\end{equation}

\subsubsection{Numerical convergence studies}
\label{sec.conv}
The convergence studies of our Lagrangian WENO finite volume schemes are carried out using the Euler equations of compressible gas dynamics \eqref{eulerTerms} for the solution of a smooth convected isentropic vortex first proposed on unstructured meshes by Hu and Shu \cite{HuShuTri} in two space dimensions. The initial computational domain for the three-dimensional case is the box $\Omega(0)=[0;10]\times[0;10]\times[0;5]$ with periodic boundary conditions imposed on each face. The initial condition is the same given in \cite{HuShuTri} where we set to zero the $z-$aligned velocity component $w$ and it is given as a linear superposition of a homogeneous background field and some perturbations $\delta$: 
\begin{equation}
\label{ShuVortIC}
\U = (\rho, u, v, w, p) = (1+\delta \rho, 1+\delta u, 1+\delta v, 1+\delta w, 1+\delta p).
\end{equation}
The perturbation of the velocity vector $\mathbf{v}=(u,v,w)$ as well as the perturbation of temperature $T$ read
\begin{equation}
\label{ShuVortDelta}
\left(\begin{array}{c} \delta u \\ \delta v \\ \delta w \end{array}\right) = \frac{\epsilon}{2\pi}e^{\frac{1-r^2}{2}} \left(\begin{array}{c} -(y-5) \\ \phantom{-}(x-5) \\ 0 \end{array}\right), \qquad \delta T = -\frac{(\gamma-1)\epsilon^2}{8\gamma\pi^2}e^{1-r^2},
\end{equation} 
where the radius of the vortex has been defined on the $x-y$ plane as $r^2=(x-5)^2+(y-5)^2$, the vortex strength is $\epsilon=5$ and the ratio of specific heats is set to $\gamma=1.4$. The entropy perturbation is assumed to be zero, i.e. $S=\frac{p}{\rho^\gamma}=0$, while the perturbations for density and pressure are given by 
\begin{equation}
\label{rhopressDelta}
\delta \rho = (1+\delta T)^{\frac{1}{\gamma-1}}-1, \quad \delta p = (1+\delta T)^{\frac{\gamma}{\gamma-1}}-1. 
\end{equation} 

\begin{figure}[!b]
\begin{center}
\begin{tabular}{ccc} 
\includegraphics[width=0.33\textwidth]{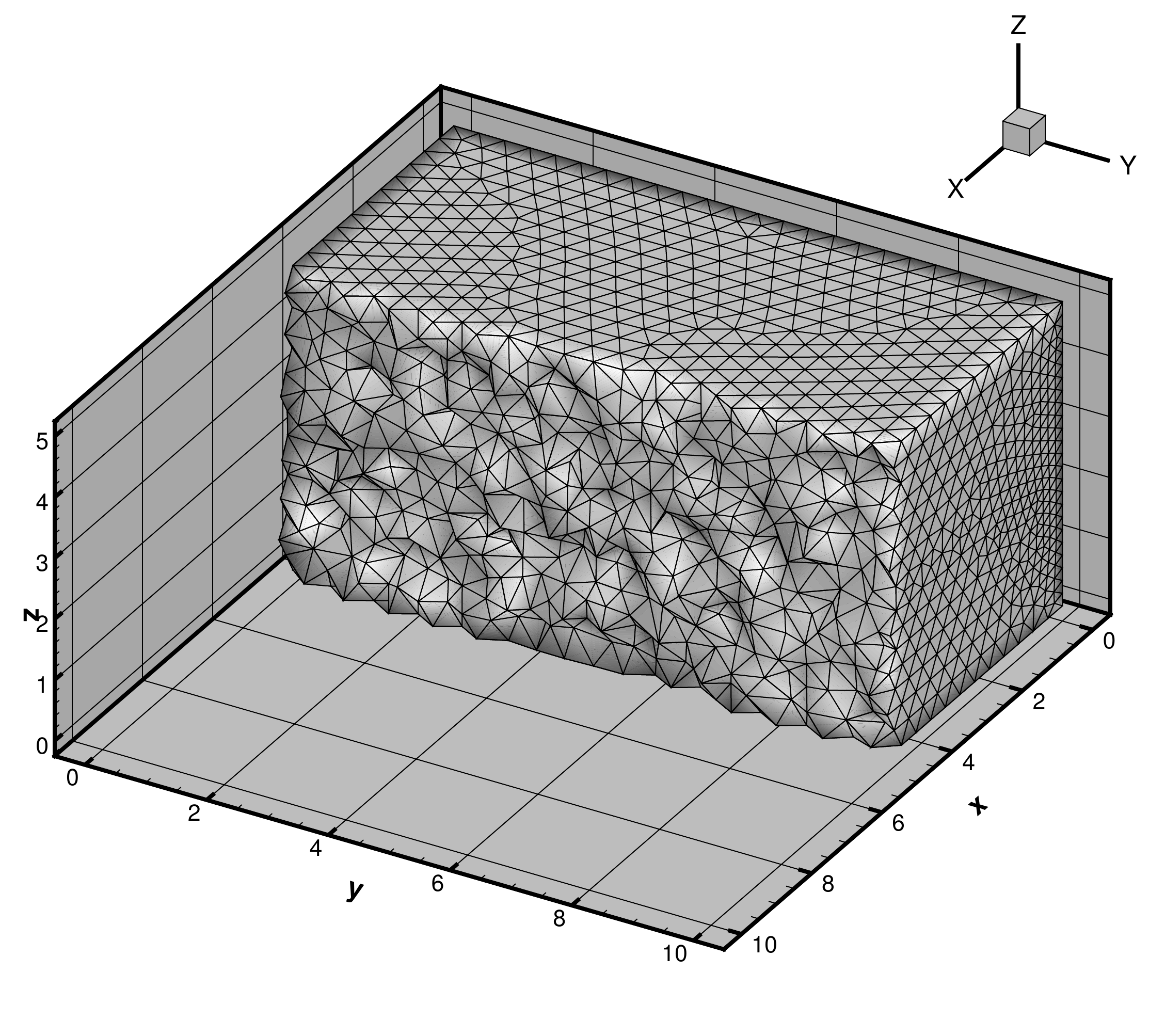}  &           
\includegraphics[width=0.33\textwidth]{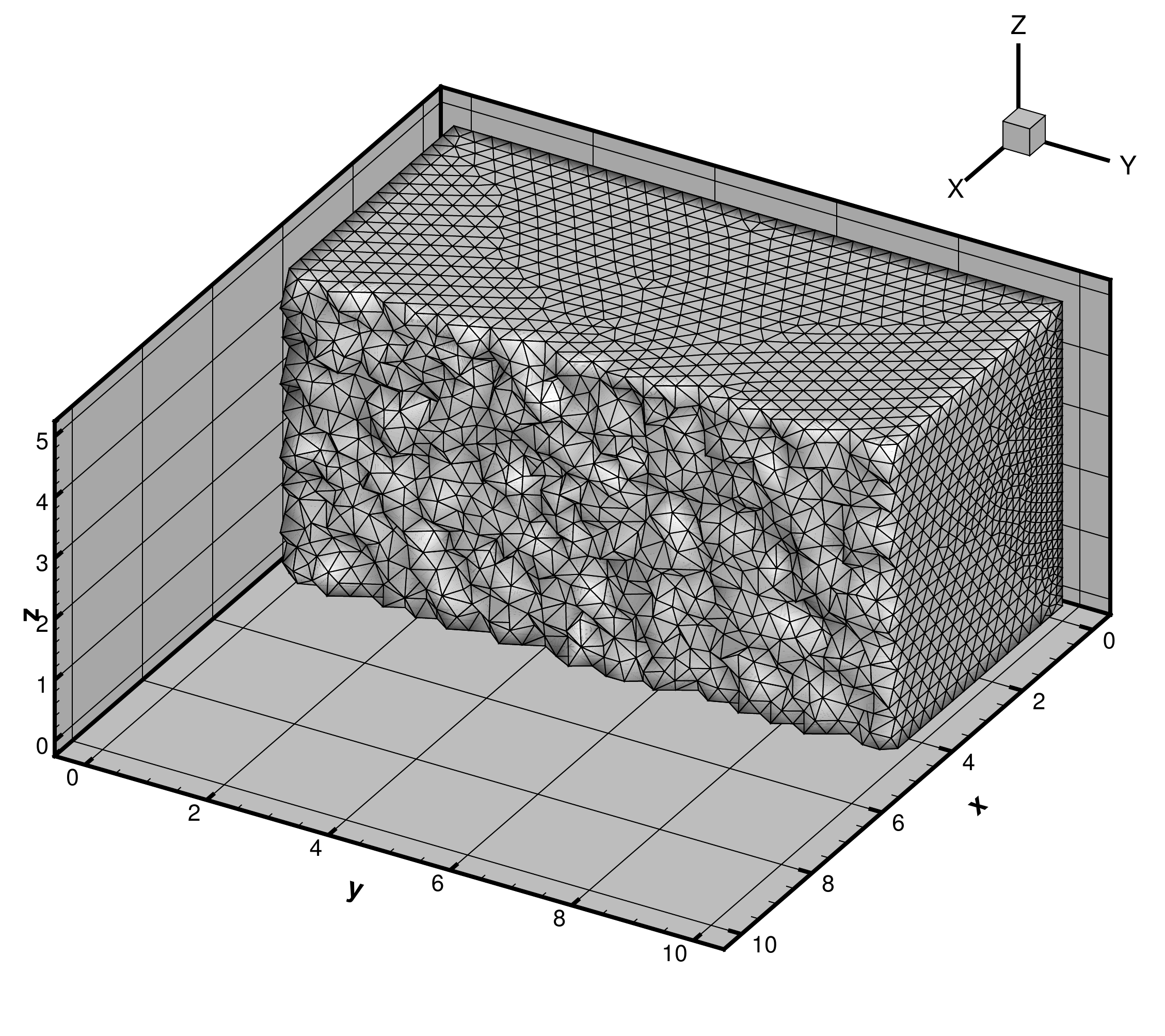}  &
\includegraphics[width=0.33\textwidth]{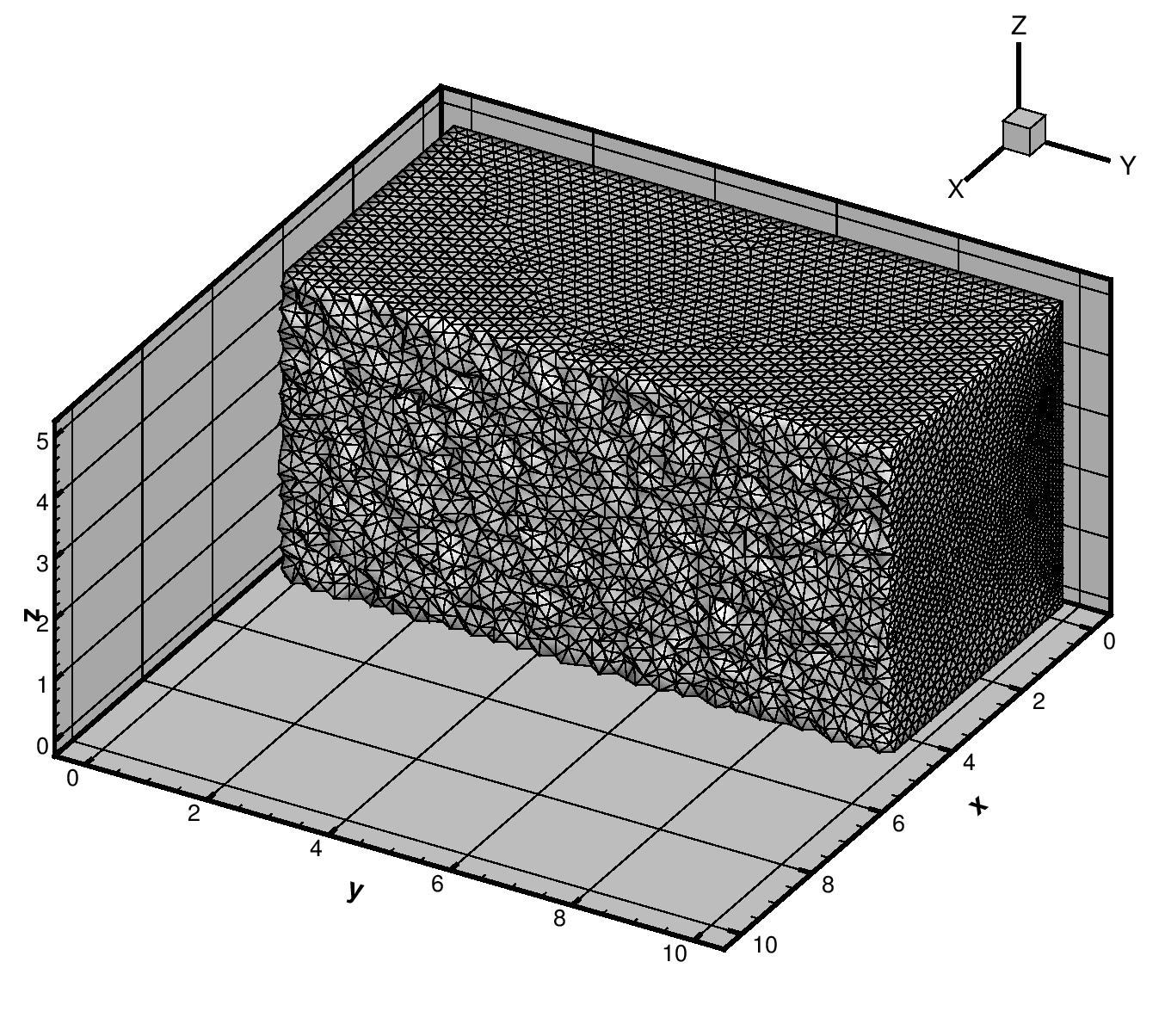}  \\           
\end{tabular} 
\caption{Sequence of tetrahedral meshes used for the numerical convergence studies.} 
\label{fig.SVmesh}
\end{center}
\end{figure}

The vortex is furthermore convected with constant velocity $\v_c=(1,1,1)$. As done in \cite{Lagrange2D}, the final time of the simulation is chosen to be $t_f=1.0$, otherwise the deformations occurring in the mesh due to the Lagrangian motion would stretch and twist the tetrahedral elements so highly that a rezoning stage would be necessary. Here we want the convergence studies to be done with a pure Lagrangian motion, hence no rezoning procedure is admitted and the final time $t_f$ has been set to a sufficiently small value. The exact solution $\Q_e$ can be simply computed as the time-shifted initial condition, e.g. $\Q_e(\x,t_f)=\Q(\x-\v_c t_f,0)$, with the convective mean velocity $\v_c$ previously defined. The error is measured at time $t_f$ using the continuous $L_2$ norm with the high order reconstructed solution $\w_h(\x,t_f)$, hence
\begin{equation}
  \epsilon_{L_2} = \sqrt{ \int \limits_{\Omega(t_f)} \left( \Q_e(\x,t_f) - \w_h(\x,t_f) \right)^2 d\x },  
	\label{eqnL2error}
\end{equation}
where $h(\Omega(t_f))$ represents the mesh size which is taken to be the maximum diameter of the circumspheres of the tetrahedral elements in the final domain configuration $\Omega(t_f)$. Figure \ref{fig.SVmesh} shows some of the successively refined meshes at the initial time $t=0$ used for this test case, while Table \ref{tab.convEul} reports the numerical convergence rates obtained 
with first to sixth order ADER-WENO schemes. The Osher-type flux \eqref{eqn.osher} has been used in all computations. 

\begin{table}  
\caption{Numerical convergence results for the compressible Euler equations using the first up to sixth order version of the three-dimensional Lagrangian one-step WENO finite volume 
schemes presented in this article. The error norms refer to the variable $\rho$ (density) at time $t=1.0$.}  
\begin{center} 
\begin{small}
\renewcommand{\arraystretch}{1.0}
\begin{tabular}{cccccc} 
\hline
  $h(\Omega(t_f))$ & $\epsilon_{L_2}$ & $\mathcal{O}(L_2)$ & $h(\Omega,t_f)$ & $\epsilon_{L_2}$ & $\mathcal{O}(L_2)$ \\ 
\hline
  \multicolumn{3}{c}{$\mathcal{O}1$} & \multicolumn{3}{c}{$\mathcal{O}2$}  \\
\hline
3.43E-01 & 1.081E-01 & -   & 2.89E-01  & 2.214E-02 & -   \\ 
2.85E-01 & 9.159E-02 & 0.9 & 2.16E-01  & 1.202E-02 & 2.1 \\ 
2.09E-01 & 6.875E-02 & 0.9 & 1.52E-01  & 5.865E-03 & 2.0 \\ 
1.47E-01 & 4.899E-02 & 1.0 & 1.13E-01  & 3.254E-03 & 2.0 \\ 
\hline 
  \multicolumn{3}{c}{$\mathcal{O}3$} & \multicolumn{3}{c}{$\mathcal{O}4$}  \\
\hline
2.89E-01 & 1.718E-02 & -    & 2.89E-01 & 4.116E-03 & -  \\ 
2.17E-01 & 7.641E-03 & 2.8  & 2.17E-01 & 1.369E-03 & 3.8 \\ 
1.52E-01 & 2.601E-03 & 3.1  & 1.52E-01 & 3.273E-04 & 4.1 \\ 
1.13E-01 & 1.049E-03 & 3.1  & 1.13E-01 & 9.802E-05 & 4.1\\ 
\hline 
  \multicolumn{3}{c}{$\mathcal{O}5$} & \multicolumn{3}{c}{$\mathcal{O}6$}  \\
\hline
2.89E-01  & 2.272E-03 & -   & 2.89E-01 & 1.015E-03 & -    \\ 
2.17E-01  & 6.605E-04 & 4.3 & 2.17E-01 & 2.312E-04 & 5.1  \\ 
1.52E-01  & 1.234E-04 & 4.8 & 1.52E-01 & 3.090E-05 & 5.7  \\ 
1.13E-01  & 2.932E-05 & 4.9 & 1.13E-01 & 6.576E-06 & 5.2  \\ 
\hline 
\end{tabular}
\end{small}
\end{center}
\label{tab.convEul}
\end{table} 

In order to identify the most expensive part of the algorithm in terms of computational efficiency, we also collect the times used for carrying on the WENO reconstruction, the local space-time predictor and the Lagrangian flux evaluation. We run the simulation in parallel on four Intel Core i7-2600 CPUs with a clock-speed of 3.40GHz. We consider a coarse grid with a characteristic mesh size of $h=0.042$ containing a total number of elements of $N_E=60157$ and we perform the isentropic vortex test case presented in this Section using the Osher-type numerical flux \eqref{eqn.osher} until the final time of $t_f=1.0$. Table \ref{tab.CompCost} reports the computational cost of each part of the algorithm for second, third and fourth order accurate Lagrangian finite volume schemes. The most expensive part of the algorithm is the flux evaluation, since in the Lagrangian framework no quadrature-free approach is possible, due to the continuous evolution of the geometry configuration that does not allow the flux computation to be treated as done for the Eulerian case in \cite{Dumbser2007204}, where the space-time basis used for the flux integrals in \eqref{PDEfinal} can be integrated on the reference space-time element in the pre-processing step and stored only once. As the order of accuracy increases the relative cost of the WENO reconstruction procedure also increases because the reconstruction stencils become larger, while the local space-time predictor step is the least expensive part of the whole algorithm.

\begin{table}  
\caption{Computational cost of the second, third and fourth order version of the ALE WENO finite volume schemes discussed in this paper. The times used for the WENO reconstruction, the local space-time predictor and the flux evaluation are given in percentage w.r.t. the total time of the computation.
}  
\begin{center} 
\begin{small}
\renewcommand{\arraystretch}{1.0}
\begin{tabular}{c|ccc} 
\hline
  \textit{Component of the algorithm} & $\mathcal{O}(2)$ & $\mathcal{O}(3)$ & $\mathcal{O}(4)$ \\ 
\hline
  WENO Reconstruction        & 22 \%              & 30 \%              & 40 \%  \\
  Space-Time Predictor      & 5 \%               & 9  \%              & 3 \%   \\
  Flux Evaluation            & 73 \%              & 61 \%              & 57 \%  \\
\hline
\hline
 \textit{Total time $[s]$}   & 135                & 423                & 2040   \\
\hline
\end{tabular}
\end{small}
\end{center}
\label{tab.CompCost}
\end{table}

\subsubsection{The Sod shock tube problem} 
\label{sec.RP3D}

\begin{figure}[!htbp]
\begin{center}
\begin{tabular}{cc} 
\includegraphics[width=0.47\textwidth]{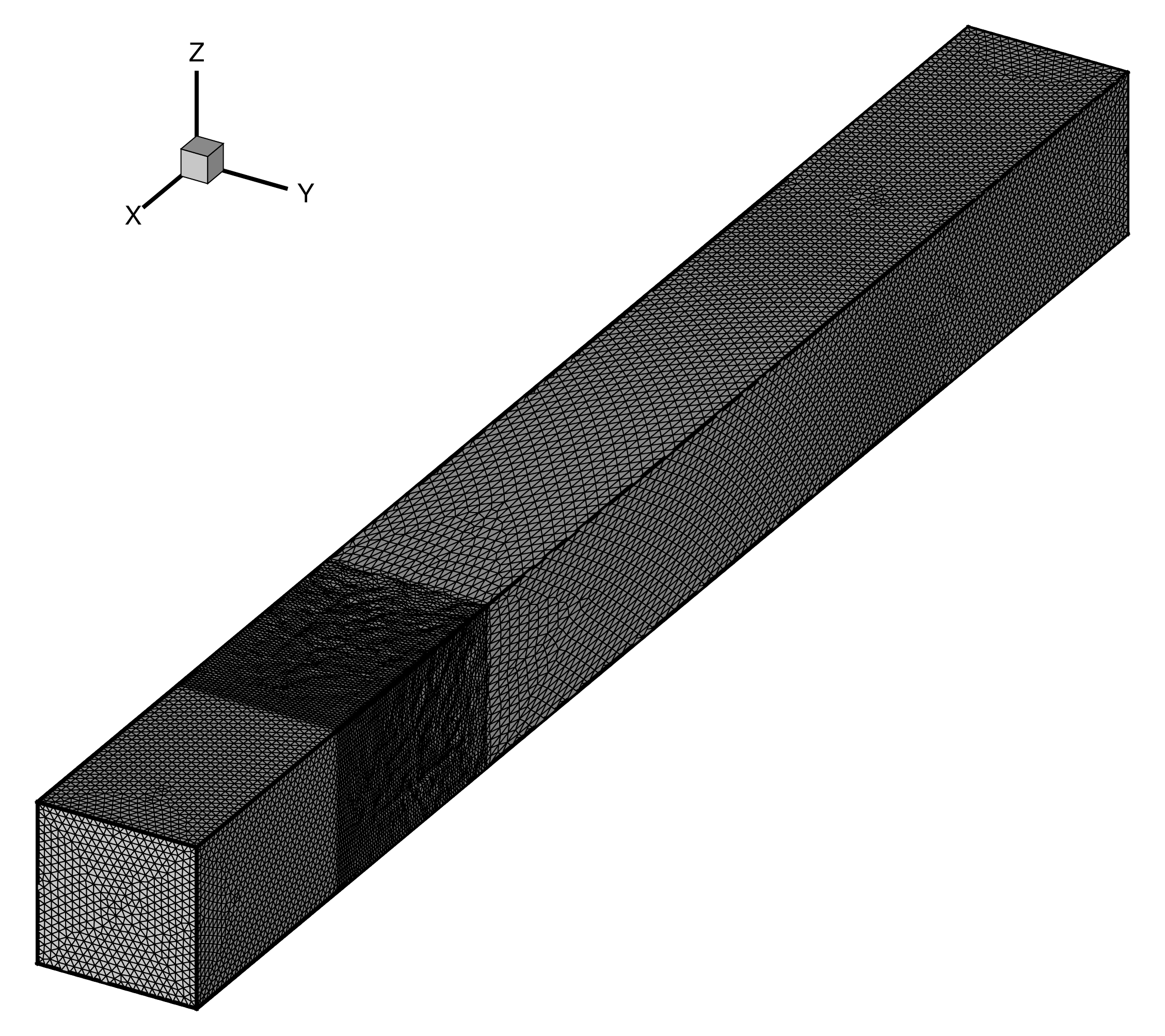}  &           
\includegraphics[width=0.47\textwidth]{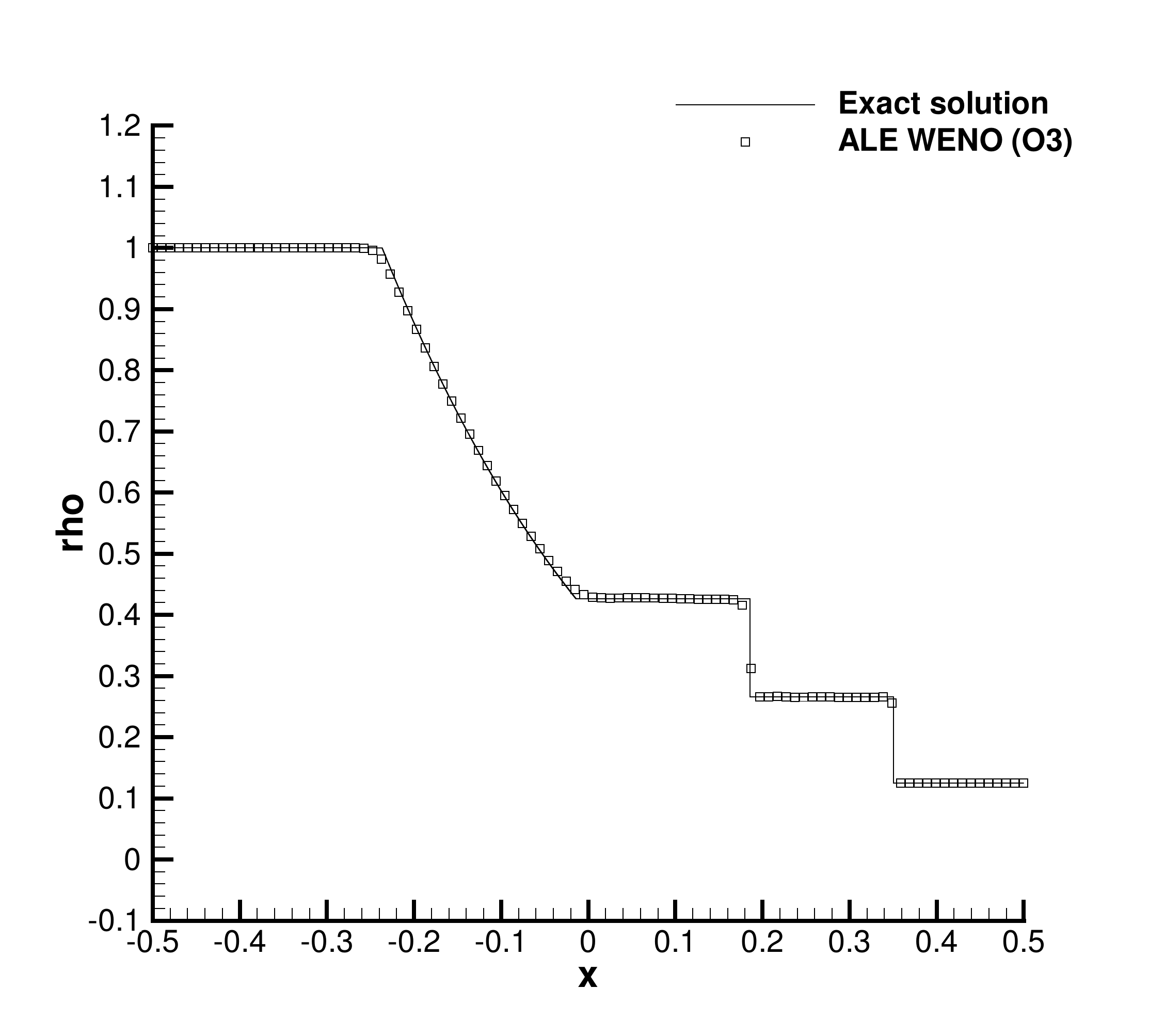} \\
\includegraphics[width=0.47\textwidth]{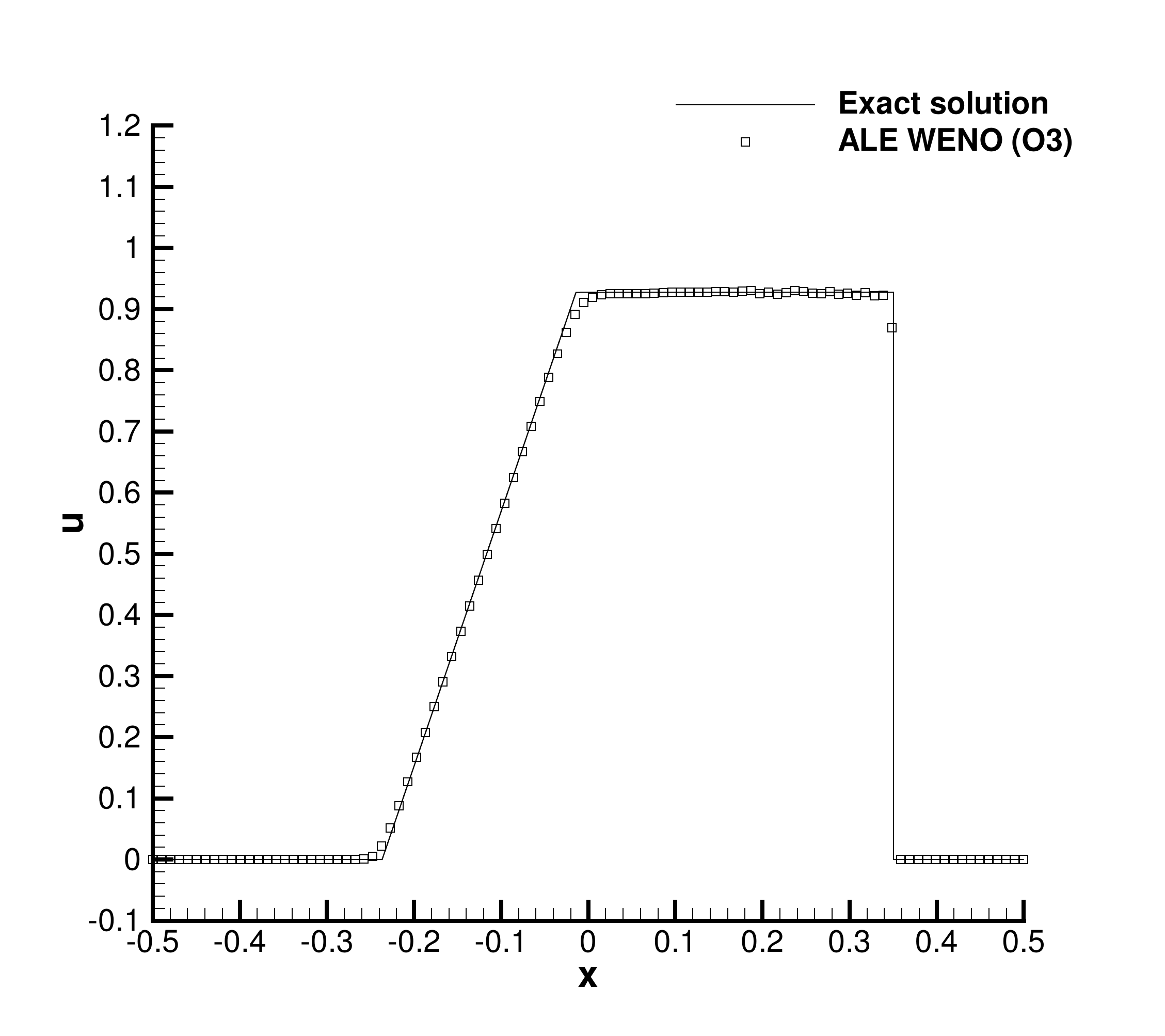}  &           
\includegraphics[width=0.47\textwidth]{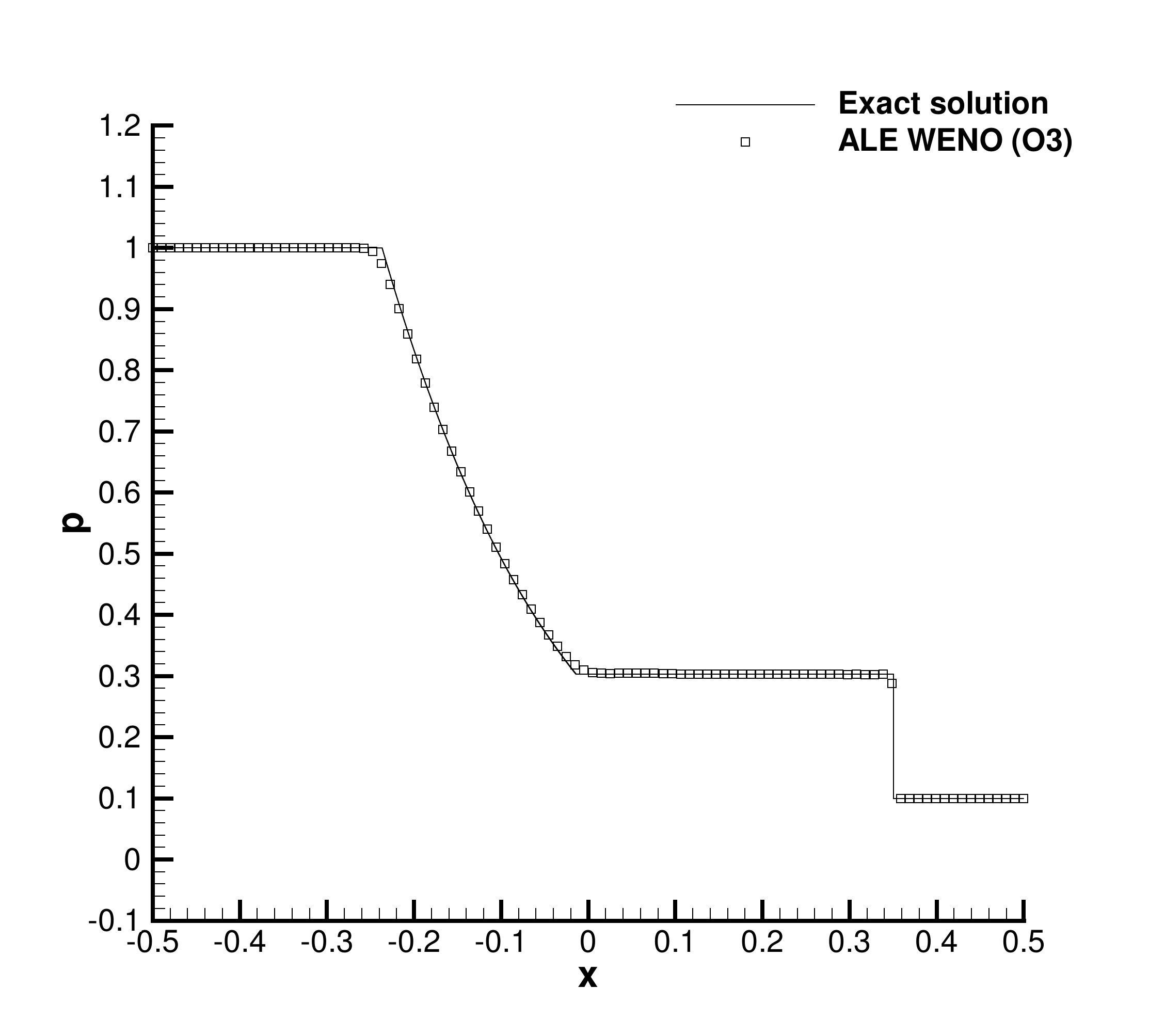} \\    
\end{tabular} 
\caption{Final 3D mesh configuration together with a 1D cut along the $x$-axis through the 
third order numerical results and comparison with exact solution for the three-dimensional Sod shock tube problem at time $t=0.2$.} 
\label{fig.RP3D}
\end{center}
\end{figure}

Here we solve in a three-dimensional setting the well-known Sod shock tube problem, which is a classical one-dimensional test problem that involves a rarefaction wave traveling towards the 
left boundary as well as a right-moving contact discontinuity and a shock wave traveling to the right. The initial computational domain is the box 
$\Omega(0)=[-0.5;0.5]\times[-0.05;0.05]\times[-0.05;0.05]$, 
which is discretized with a total number of $N_E=70453$ tetrahedral elements with a characteristic mesh size of $h=1/100$. We set periodic boundaries in the $y$ and $z$ directions, 
while transmissive boundaries are imposed along the $x$ direction. The final time of the simulation is chosen to be $t_f=0.2$. The initial condition consists in 
a discontinuity located at $x_0=0$ between two different states $\U_L$ and $\U_R$, where $\U=(\rho,u,v,w,p)$ denotes the vector of \textit{primitive} variables: 
\begin{equation}
\U(\x,0) = \left\{ \begin{array}{lll} \U_L = \left( 1.0,  0,0,0,1.0 \right), & \textnormal{ if } & x \leq x_0, \\ 
                                      \U_R = \left( 0.125,0,0,0,0.1 \right), & \textnormal{ if } & x > x_0. 
                      \end{array}  \right. 
\label{eq:Sod_IC}
\end{equation}
Although the Sod problem is a one-dimensional test case, it becomes multidimensional when applied to unstructured meshes, where in general the element faces are not aligned with the coordinate axis 
or the fluid motion. Hence, it is actually a non trivial test problem. Moreover, a contact wave is present in the solution, so that one can check how well it is resolved by the Lagrangian scheme. 
A third order scheme has been used together with the Osher-type flux \eqref{eqn.osher}. The computational results are shown in Figure \ref{fig.RP3D} and compared with the exact solution obtained
with the exact Riemann solver presented in \cite{ToroBook}. The contact wave has been resolved very well with only one intermediate point and overall a very good agreement with the exact solution 
is achieved for density, as well as for pressure and for the horizontal velocity component $u$. 

\subsubsection{Three-dimensional explosion problem} 
\label{sec.EP3D}

The explosion problem can be seen as a fully three-dimensional extension of the Sod problem presented in Section \ref{sec.RP3D} before. 
The initial domain is the sphere of radius $R_o=1$, i.e. $\Omega(0)=\left\{ \x : \left\| \x \right\| \leq R_o \right\}$, in which a sphere of radius 
$R=0.5$ separates two different constant states: 
\begin{equation}
  \U(\x,0) = \left\{ \begin{array}{ccc} \U_i=(1,0,0,0,1),       & \textnormal{ if } & \left\| \x \right\| \leq R, \\ 
                                        \U_o=(0.125,0,0,0,0.1), & \textnormal{ if } & \left\| \x \right\| > R.        
                      \end{array}  \right. 
\end{equation}
The \textit{inner state} $\U_i$ and the \textit{outer state} $\U_o$ correspond to the ones of the 1D Sod problem. 
%
%
For spherically symmetric problems, the multidimensional Euler system \eqref{eqn.pde.nc}-\eqref{eulerTerms} can be simplified to a one-dimensional system with geometric source terms, 
see \cite{ToroBook,Lagrange2D}. It reads 
\begin{equation}
\label{inhomEuler}
\Q_t + \F(\Q)_r = \S(\Q),
\end{equation}  
with
\begin{equation}
\label{matrixInhomEuler}
\Q = \left(\begin{array}{c} \rho \\ \rho u \\ \rho E \end{array}\right), \quad \F = \left(\begin{array}{c} \rho u \\ \rho u^2 + p \\  u(\rho E+p) \end{array}\right), \quad \S = -\frac{d-1}{r}\left(\begin{array}{c} \rho u \\ \rho u^2 \\  u(\rho E+p) \end{array}\right).
\end{equation}
The radial direction is denoted as usual by $r$, while $u$ represents the radial velocity and $d$ is the number of space dimensions. In order to compute a suitable reference solution we set $d = 3$ 
and a classical second order TVD scheme \cite{ToroBook} with Rusanov flux has been used to solve the inhomogeneous system of equations \eqref{inhomEuler} on a one-dimensional mesh of 15000 points 
in the radial interval $r \in [0;1]$. 

\begin{figure}[!htbp]
\begin{center}
\begin{tabular}{cc} 
\includegraphics[width=0.47\textwidth]{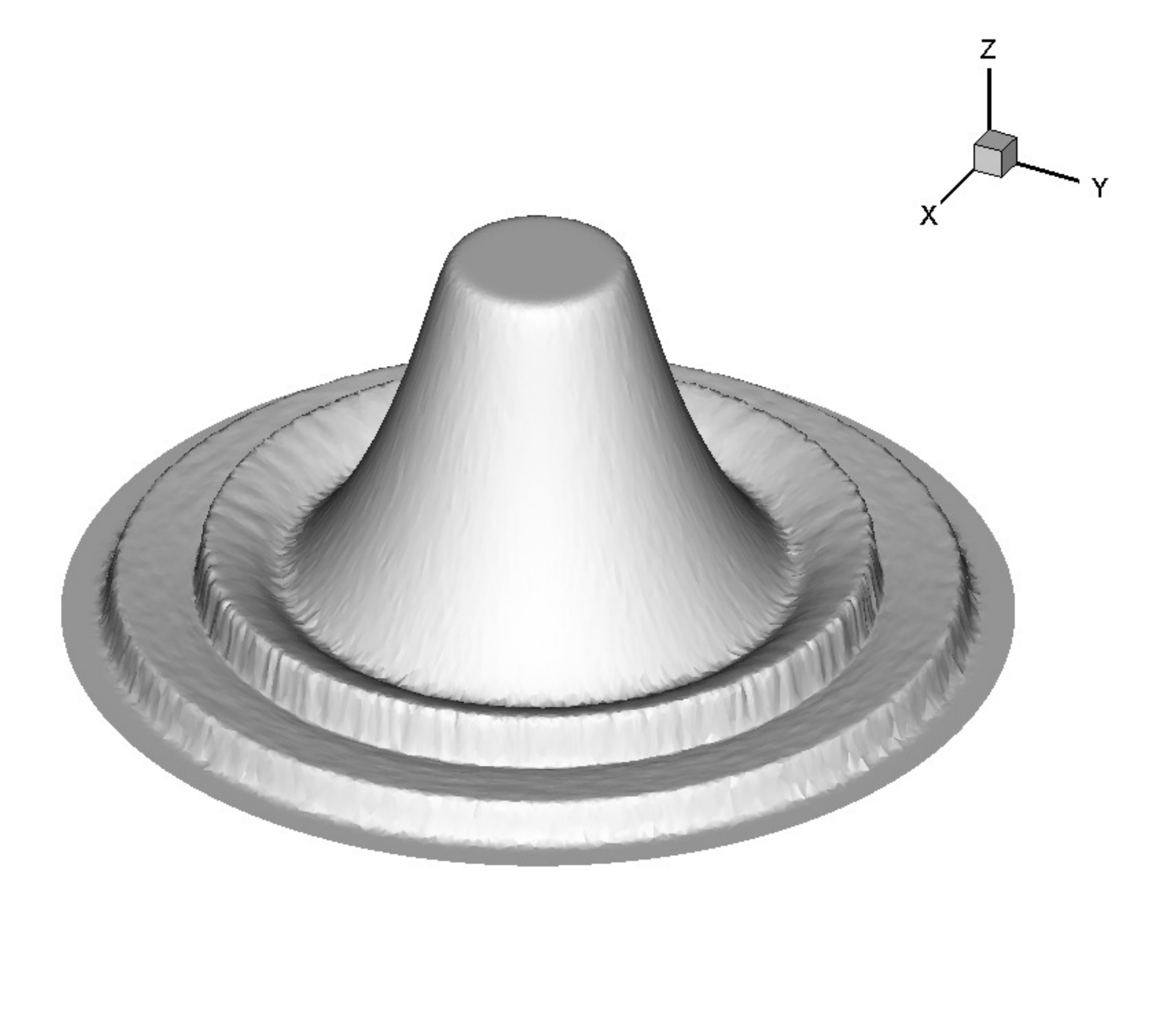}  &           
\includegraphics[width=0.47\textwidth]{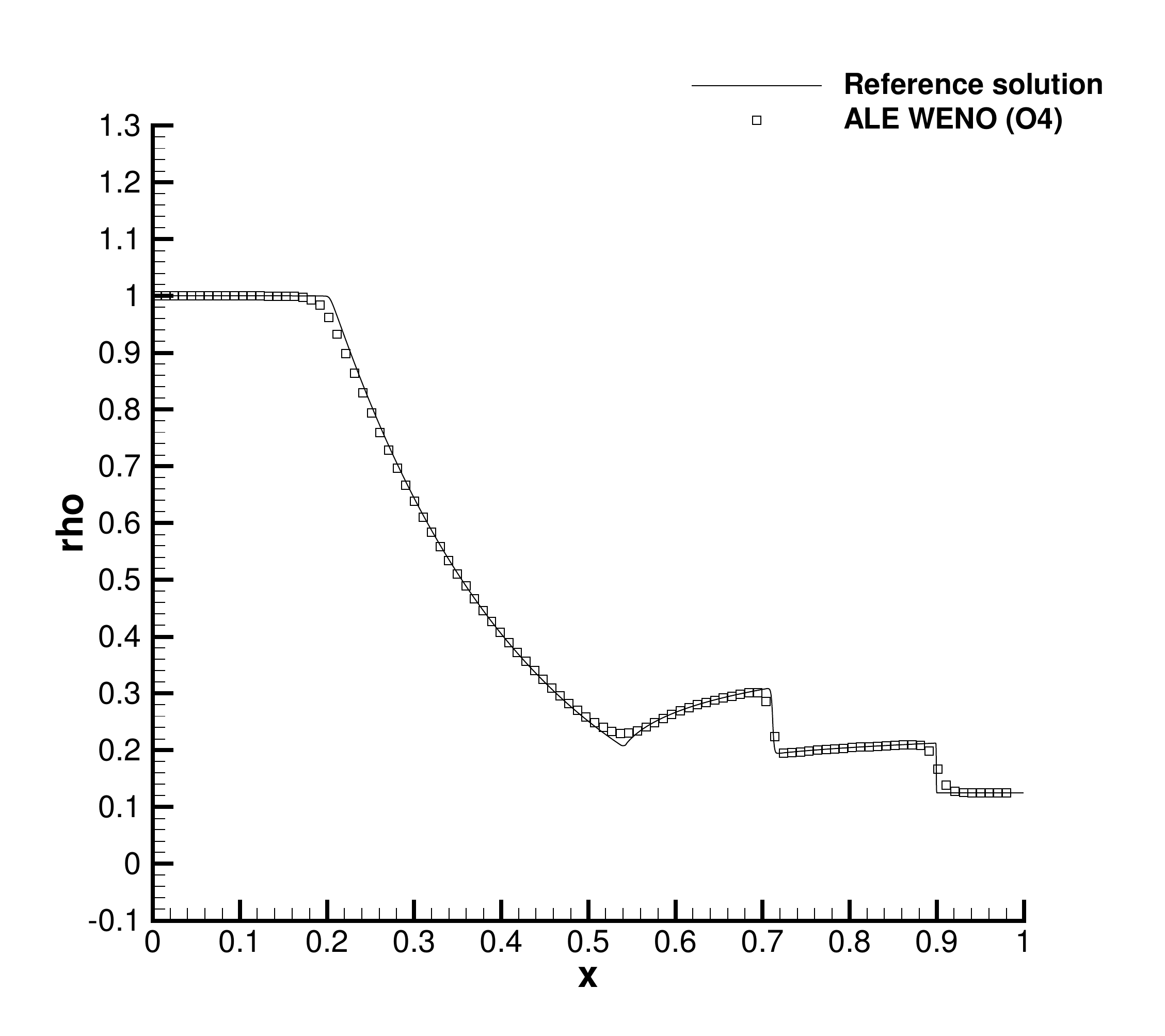} \\
\includegraphics[width=0.47\textwidth]{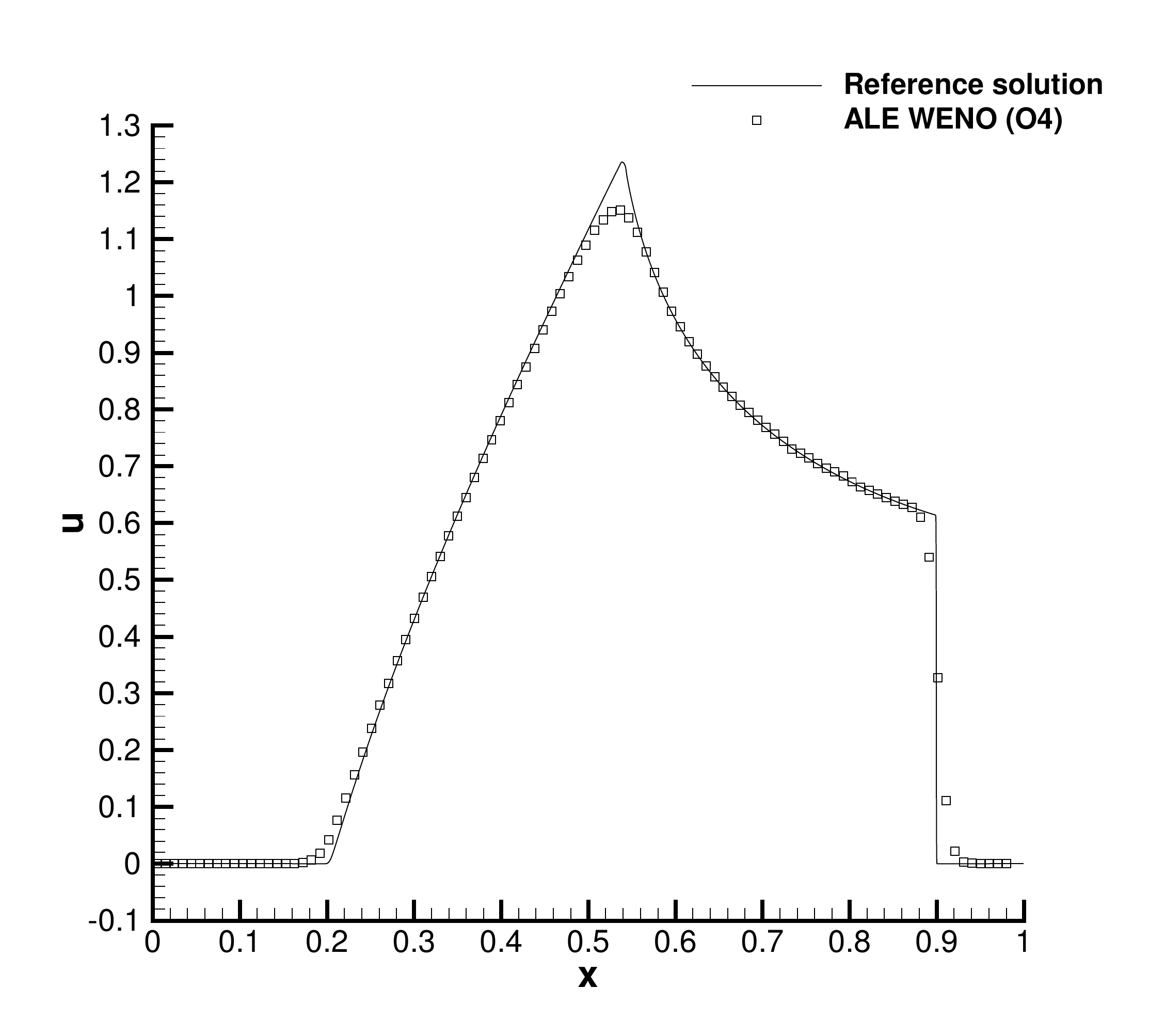}  &           
\includegraphics[width=0.47\textwidth]{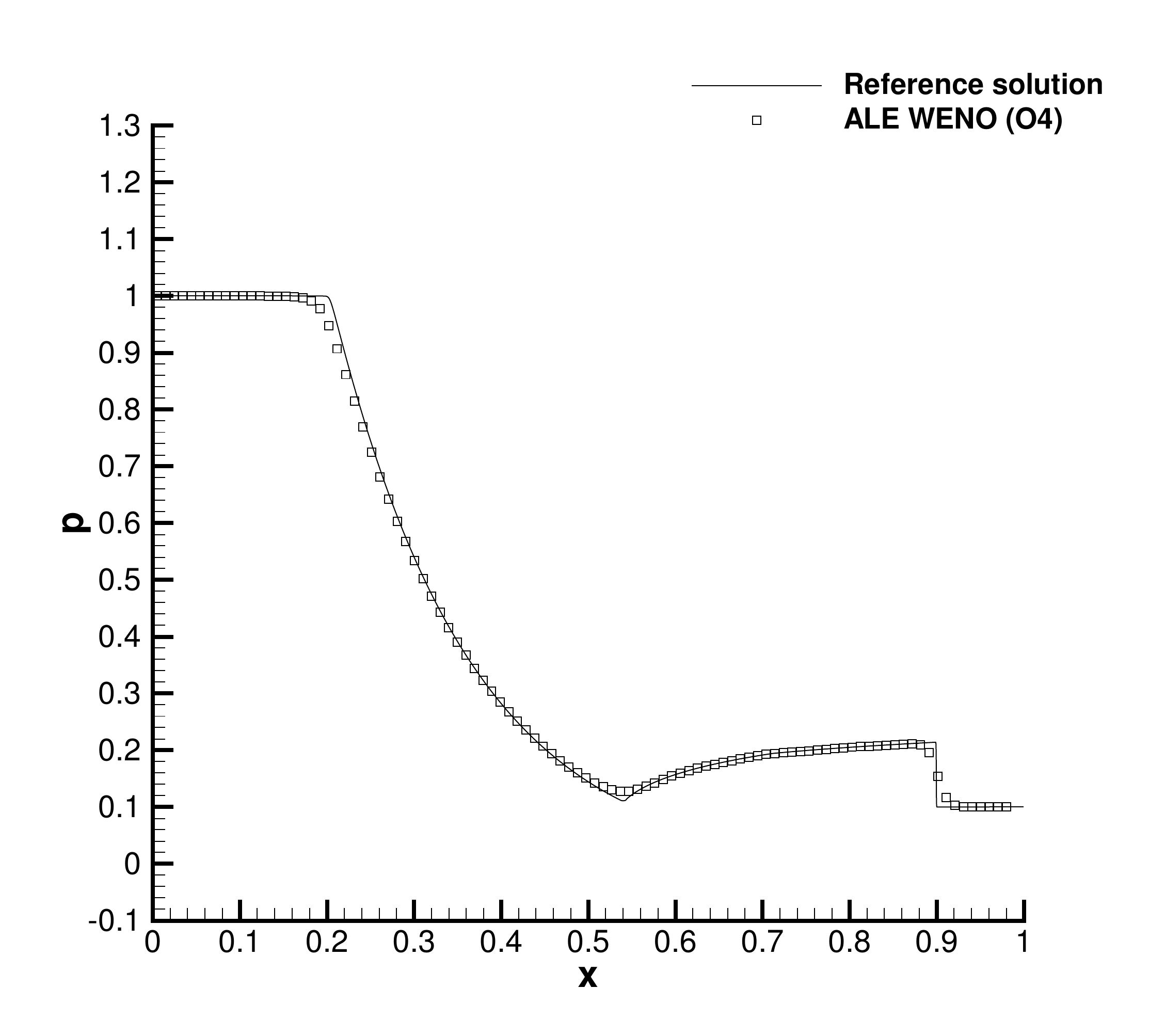} \\    
\end{tabular} 
\caption{Fourth order numerical results and comparison with the reference solution for the three-dimensional explosion problem at time $t=0.25$.} 
\label{fig.EP3D}
\end{center}
\end{figure}

\begin{figure}[!htbp]
\begin{center}
\begin{tabular}{cc} 
\includegraphics[width=0.47\textwidth]{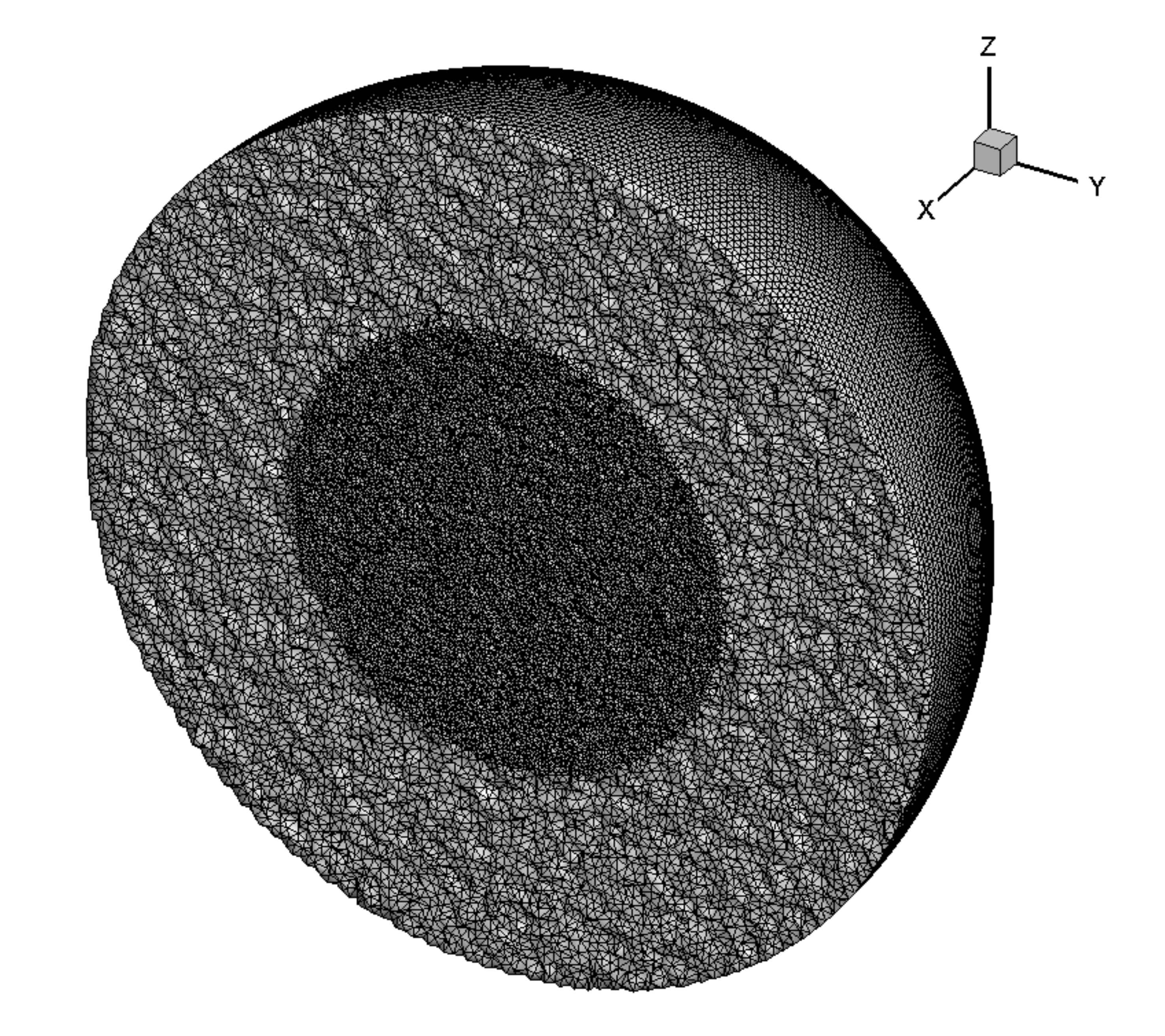}  &           
\includegraphics[width=0.47\textwidth]{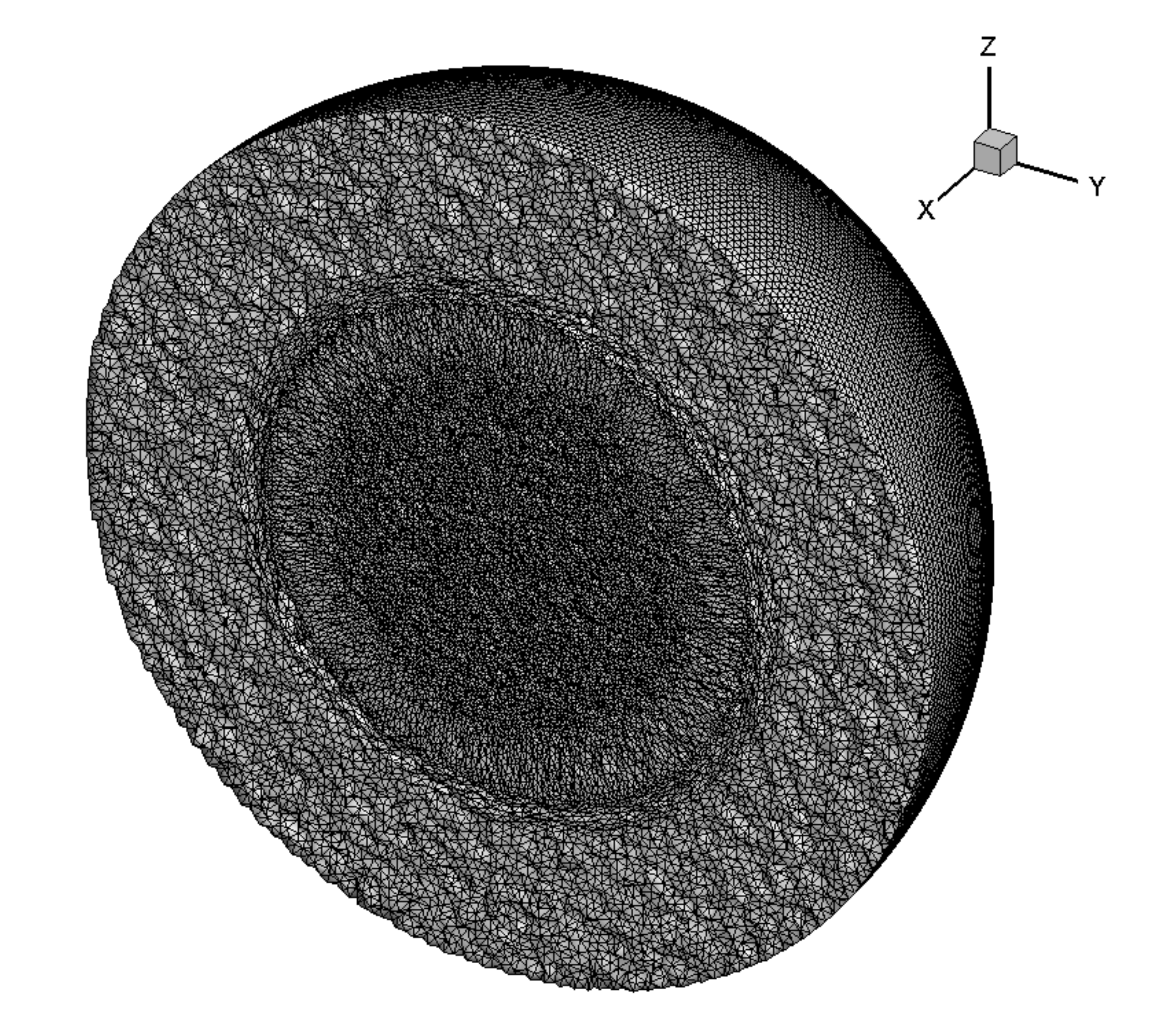} \\
\includegraphics[width=0.47\textwidth]{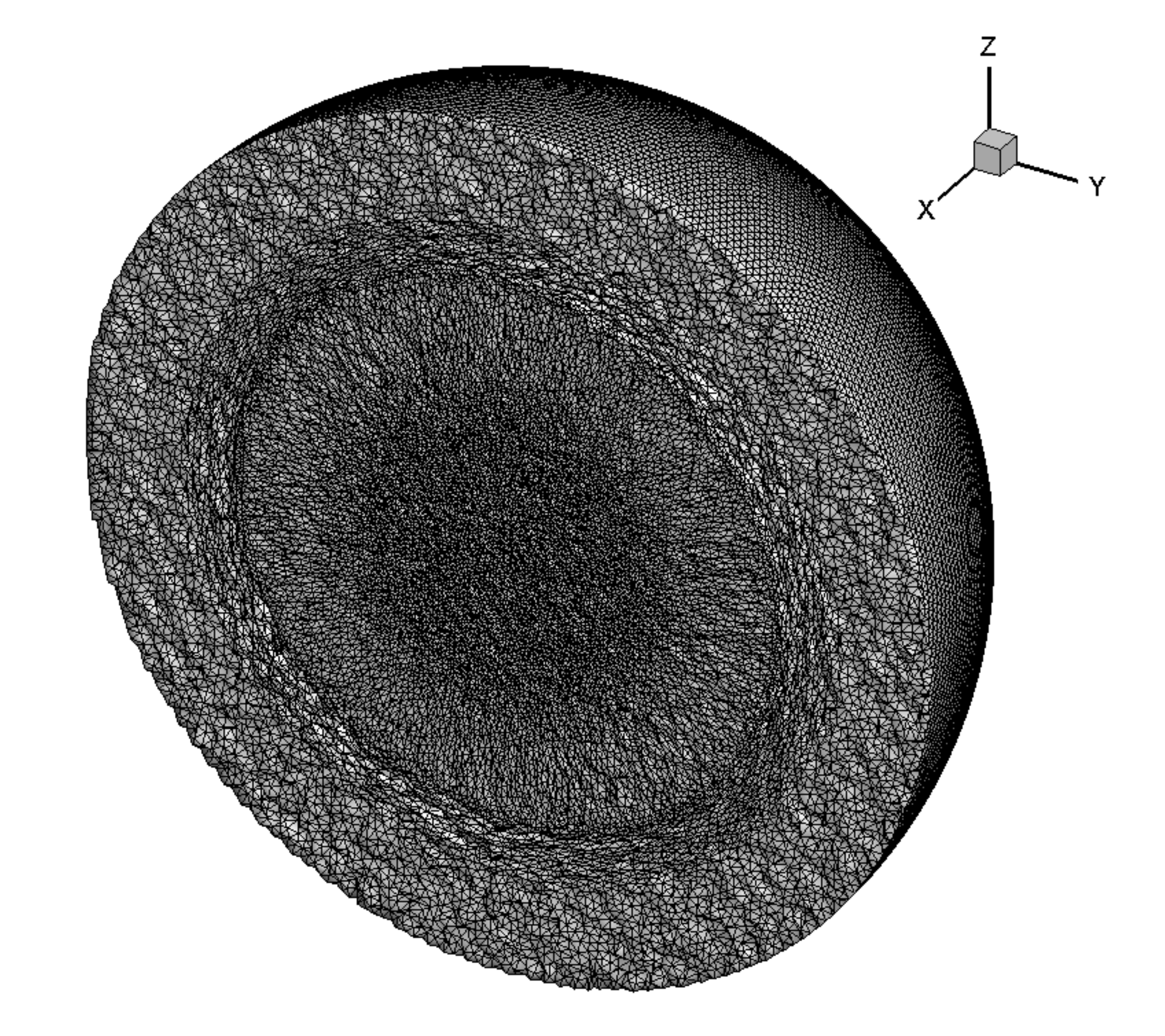}  &           
\includegraphics[width=0.47\textwidth]{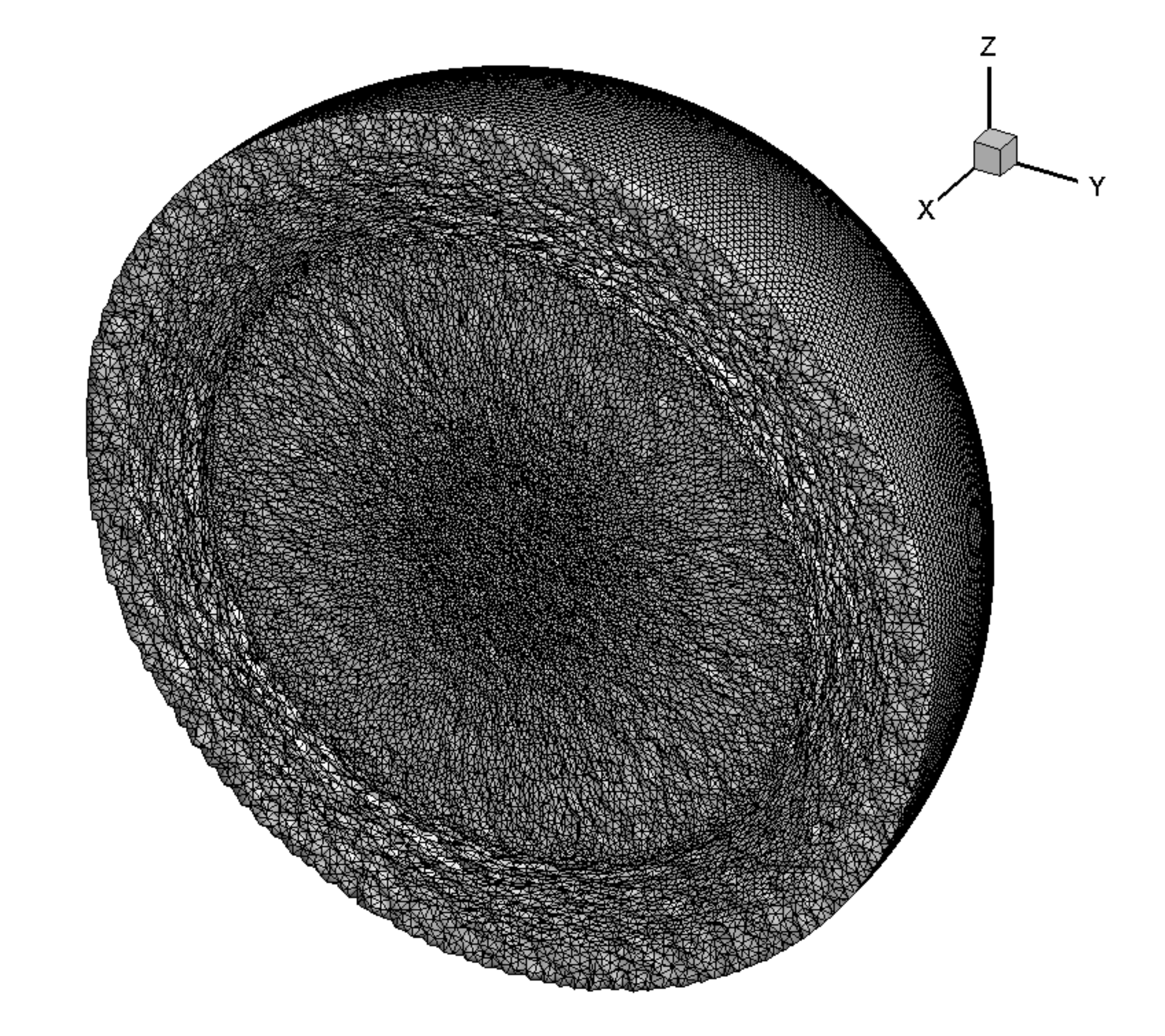} \\    
\end{tabular} 
\caption{Mesh configuration for the explosion problem at times $t=0.00$, $t=0.08$, $t=0.16$ and $t=0.25$.} 
\label{fig.EP3Dgrid}
\end{center}
\end{figure}

The ratio of specific heats is assumed to be $\gamma=1.4$ and the final time is $t_f=0.25$. The computational domain is discretized with a total number of elements of $N_E=7225720$ and transmissive boundary conditions have been imposed on the external boundary. Figure \ref{fig.EP3D} shows a comparison between the reference solution and the numerical results, computed using the fourth order version of our ALE WENO schemes together with the Osher-type numerical flux \eqref{eqn.osher}. The solution involves three different waves, namely one spherical shock wave traveling towards the external boundary of the domain, the rarefaction fan which is moving to the opposite direction and the contact wave in between, that is very well resolved due to the Lagrangian approach together with the use of the little  diffusive Osher-type numerical flux. A slice of the entire mesh configuration at four different output times is depicted in Figure \ref{fig.EP3Dgrid}, where the progressively compression of the tetrahedra located at the shock frontier can be clearly identified.

\subsubsection{The Kidder problem} 
\label{sec.Kidder}

In \cite{Kidder1976} Kidder proposed this test problem, which has become a classical benchmark for Lagrangian schemes \cite{Maire2009,Despres2009}. It consists in an isentropic compression of a portion of a shell filled with a prefect gas which is assigned with the following initial condition: 
\begin{equation}
\left( \begin{array}{c} \rho_0(r) \\ \mathbf{v}_0(r) \\ p_0(r) \end{array}  \right) = \left( \begin{array}{c}  \left(\frac{r_{e,0}^2-r^2}{r_{e,0}^2-r_{i,0}^2}\rho_{i,0}^{\gamma-1}+\frac{r^2-r_{i,0}^2}{r_{e,0}^2-r_{e,0}^2}\rho_{e,0}^{\gamma-1}\right)^{\frac{1}{\gamma-1}} \\ 0 \\ s_0\rho_0(r)^\gamma \end{array}  \right), 
\label{eq:KidderIC}
\end{equation}
where $r=\sqrt{x^2+y^2+z^2}$ represents the general radial coordinate, $\left(r_i(t),r_e(t)\right)$ are the time-dependent internal and external frontier that delimit the shell, $\rho_{i,0}=1$ and $\rho_{e,0}=2$ are the corresponding initial values of density and $\gamma=\frac{5}{3}$ is the ratio of specific heats. Furthermore $s_0$ denotes the initial entropy distribution, that is assumed to be uniform, i.e. $s_0= \frac{p_0}{\rho_0^\gamma} = 1$.

\begin{figure}[!htbp]
	\begin{center}
	\begin{tabular}{cc} 
	\includegraphics[width=0.47\textwidth]{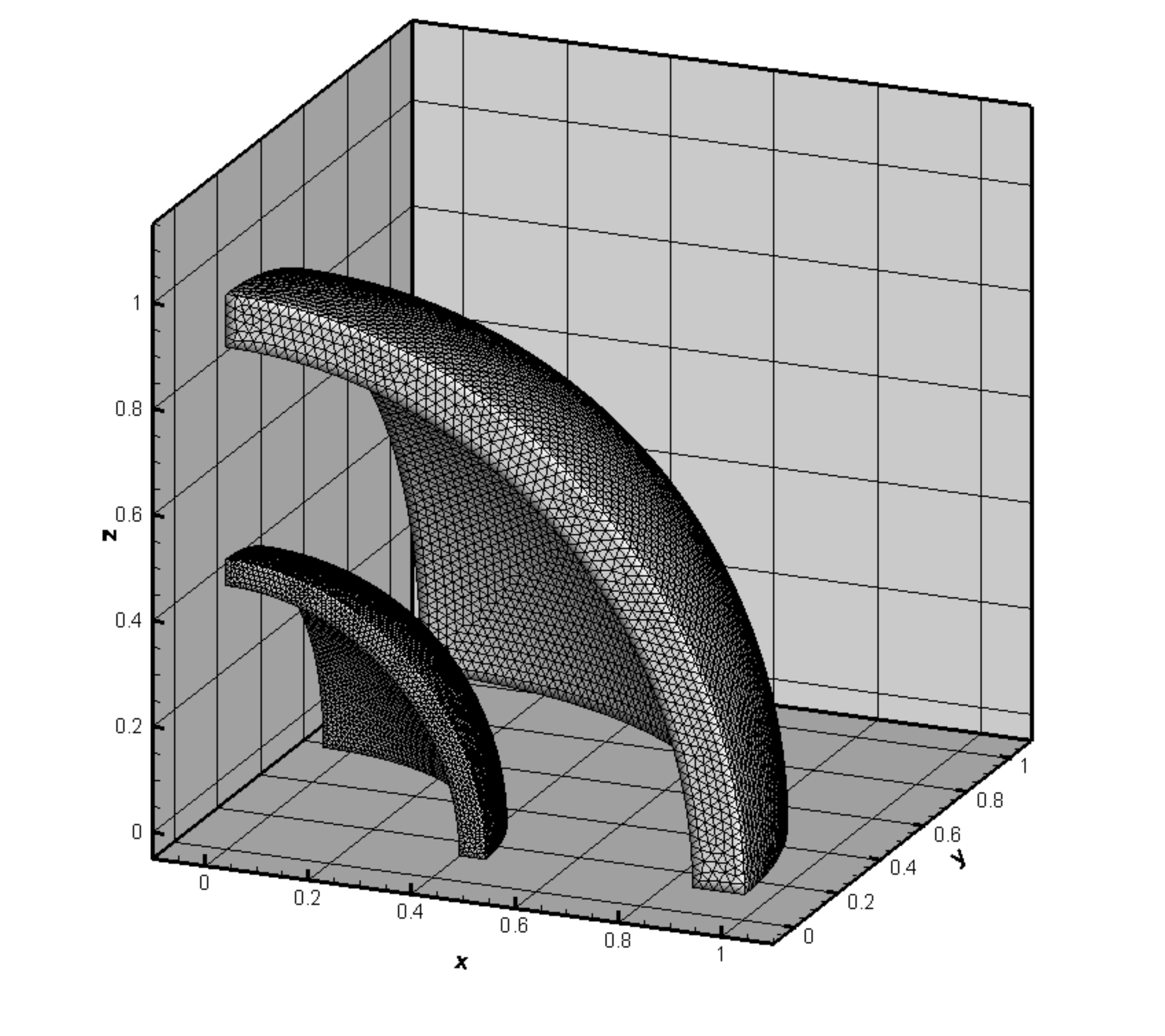}  &           
	\includegraphics[width=0.47\textwidth]{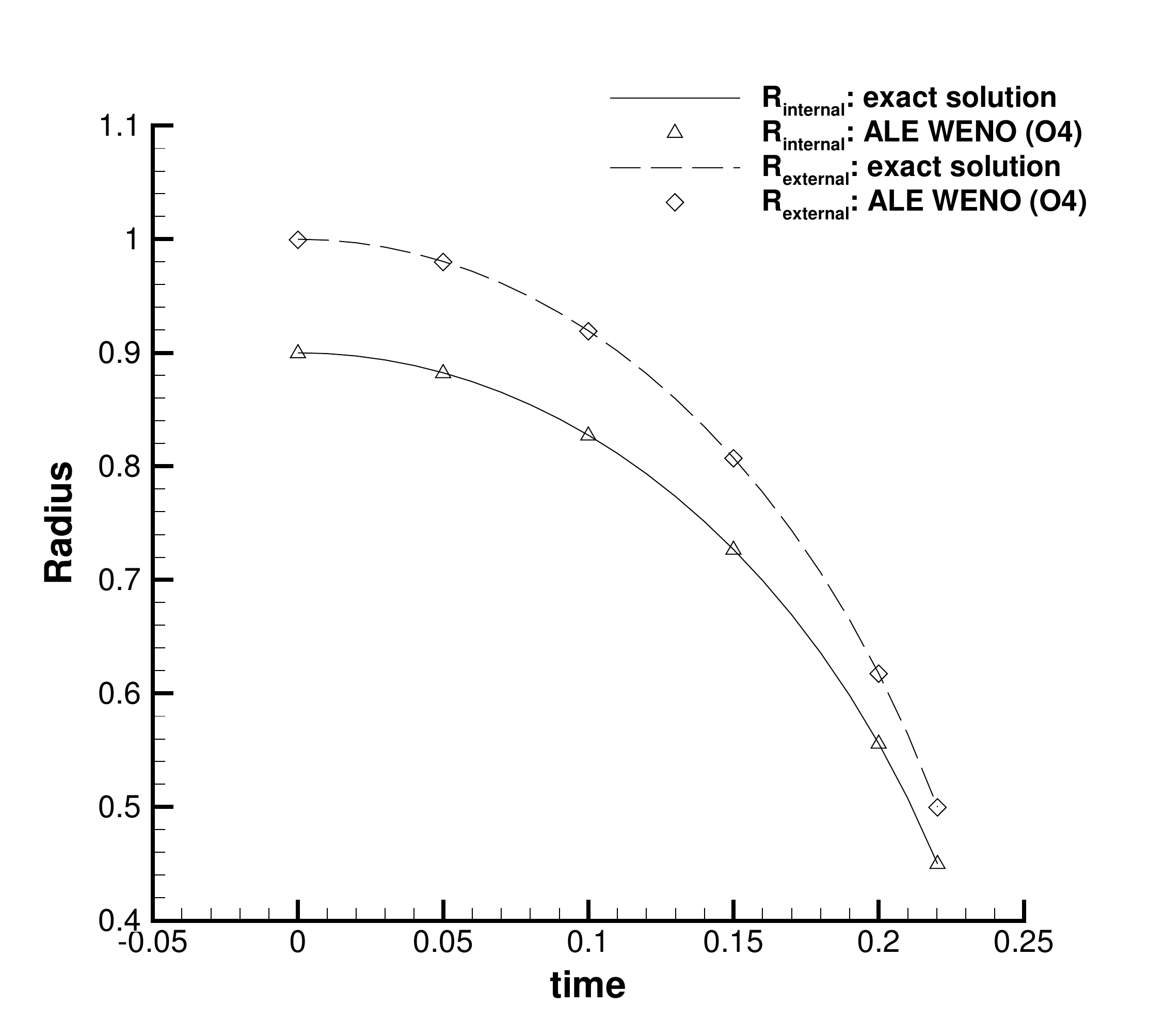} \\
	\end{tabular}
	\caption{Left: position and mesh configuration of the shell at times $t=0$ and at $t=t_f$. Right: Evolution of the internal and external radius of the shell and comparison between analytical and numerical solution.}
	\label{fig:Kidder}
	\end{center}
\end{figure}

The initial computational domain $\Omega(0)$ is one eighth of the entire shell and is depicted in Figure \ref{fig:Kidder} on the left. Sliding wall boundary conditions are imposed on the lateral faces and on the bottom, while on the internal and on the external frontier a space-time dependent state is assigned according to the exact analytical solution $R(r,t)$ \cite{Kidder1976}, which is defined at the general time $t$ for a fluid particle initially located at radius $r$ as a function of the radius and the homothety rate $h(t)$, i.e.
\begin{equation}
  R(r,t) = h(t)r, \qquad h(t) = \sqrt{1-\frac{t^2}{\tau^2}},
\label{eqKidderEx}
\end{equation}
where $\tau$ is the focalisation time 
\begin{equation}
\tau = \sqrt{\frac{\gamma-1}{2}\frac{(r_{e,0}^2-r_{i,0}^2)}{c_{e,0}^2-c_{i,0}^2}}
\end{equation}
with $c_{i,e}=\sqrt{\gamma\frac{p_{i,e}}{\rho_{i,e}}}$ representing the internal and external sound speeds. As done in \cite{Despres2009,Maire2009}, the final time of the simulation is chosen in such a way that the compression rate is $h(t_f)=0.5$, hence $t_f=\frac{\sqrt{3}}{2}\tau$ and the the exact location of the shell is bounded with $0.45 \leq R \leq 0.5$.

The computational domain is discretized with a total number of $N_E=111534$ elements and we use the fourth order version of our ALE ADER-WENO scheme together with the Osher-type flux \eqref{eqn.osher}. Figure \ref{fig:Kidder} shows the initial and the final density distribution of the shell as well as the evolution of the internal and external frontier location during the simulation. Furthermore Table \ref{tab:radiusKidder} reports the associated absolute error $|err|$, that has been evaluated as the difference between the analytical and the numerical location of the internal and external radius at the final time $t_f$. 
 
\begin{table}[!htbp]
	\begin{center}
		\begin{tabular}{|c|c|c|c|}
		\hline
		    		  			& $r_{ex}$ 		& $r_{num}$  & $|err|$     \\
		\hline
		\textit{Internal radius}	& 0.450000 	& 0.449765  & 2.35E-04 \\
		\hline
		\textit{External radius}	& 0.500000 	& 0.499727  & 2.73E-04 \\
		\hline
		\end{tabular}
	\end{center}
	\caption{Absolute error for the internal and external radius location between exact ($r_{ex}$) and numerical ($r_{num}$) solution.}
	\label{tab:radiusKidder}
\end{table}

\subsubsection{The Saltzman problem} 
\label{sec.Saltzman}

\begin{figure}[!htbp]
	\begin{center}
	\begin{tabular}{cc} 
	\includegraphics[width=0.47\textwidth]{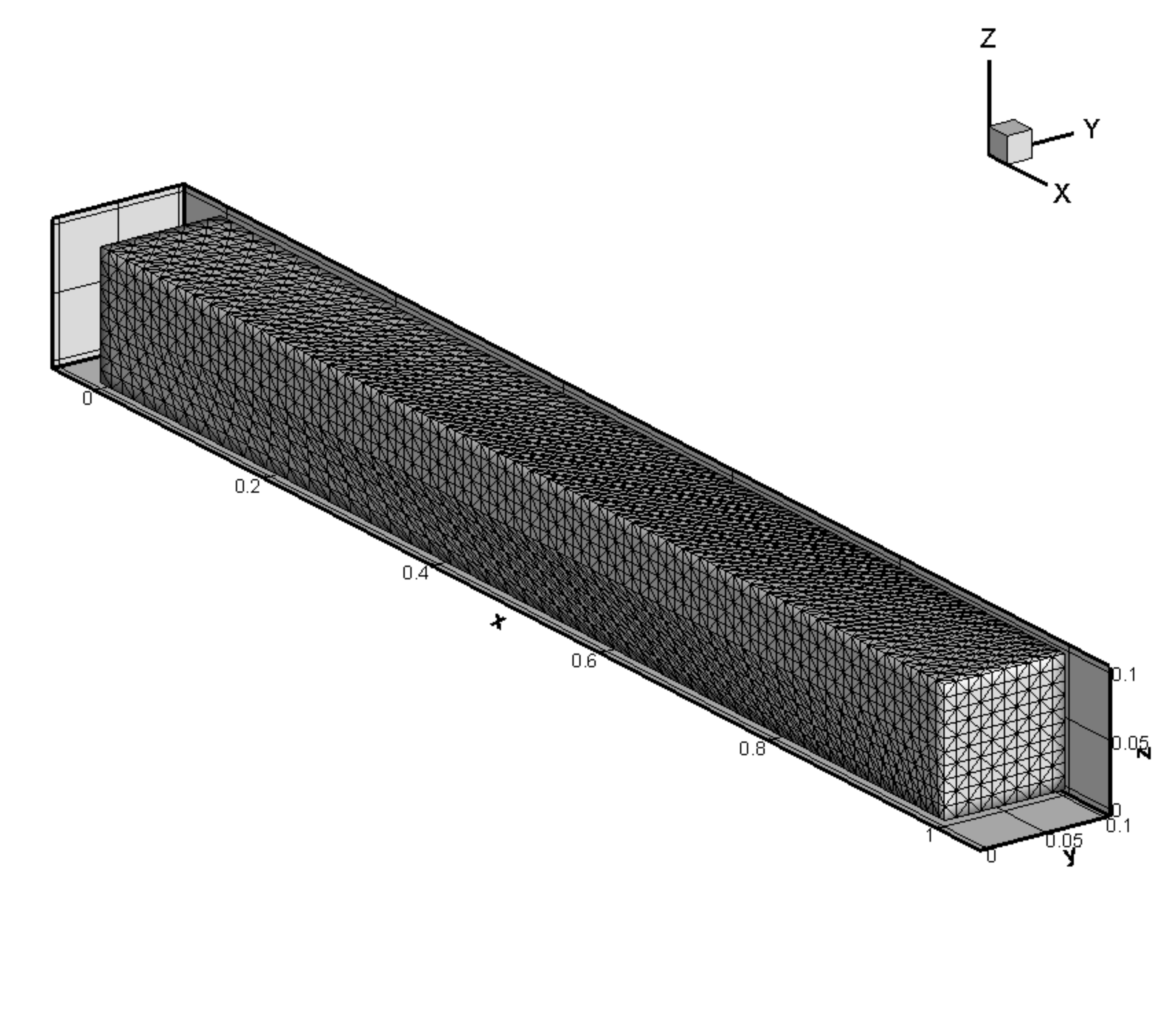}  &           
	\includegraphics[width=0.47\textwidth]{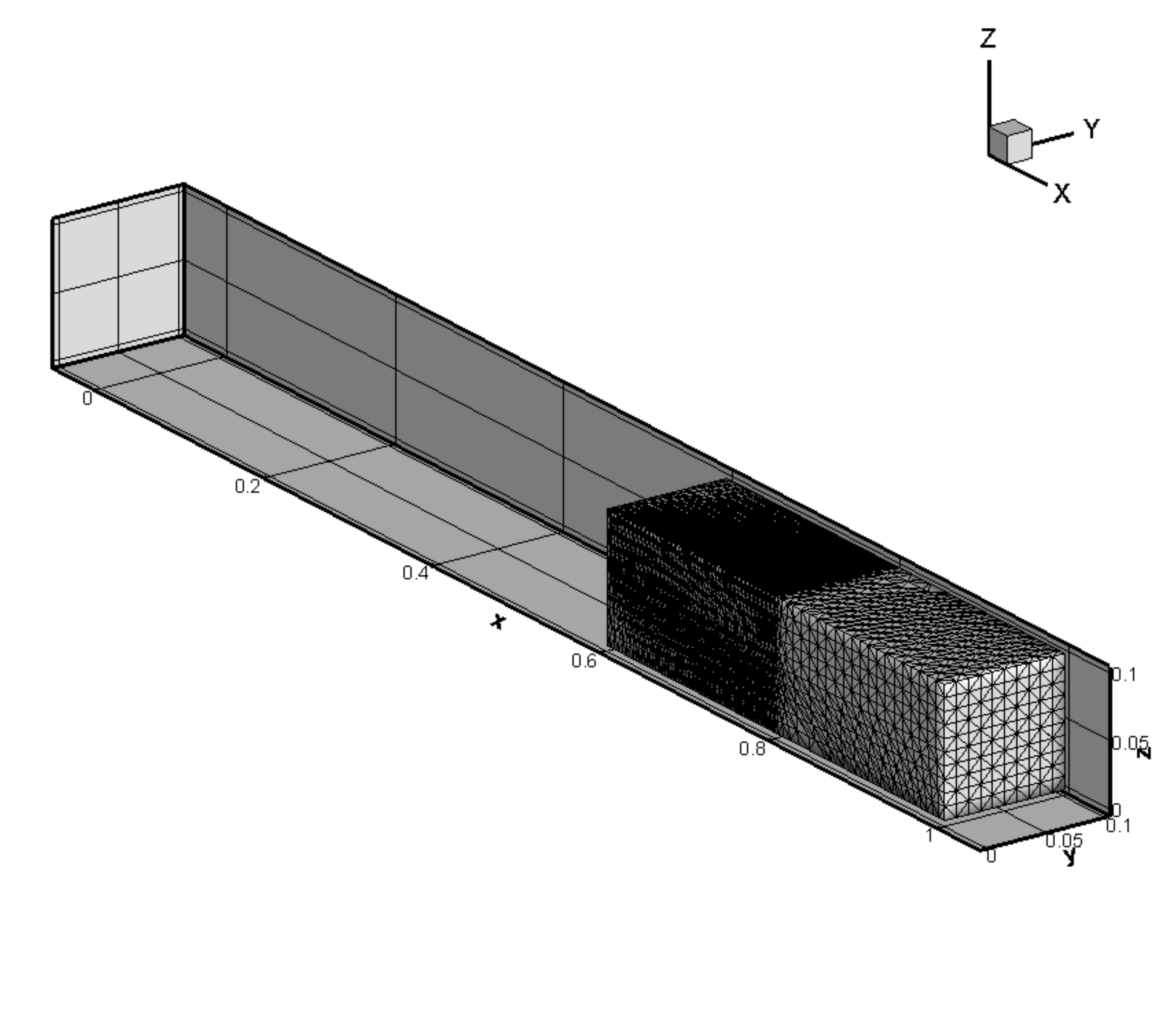} \\
	\end{tabular}
	\caption{Initial and final mesh configuration for the Saltzman problem.}
	\label{fig:SaltzGrid}
	\end{center}
\end{figure}

The Saltzman problem involves a strong shock wave that is caused by the motion of a piston traveling along the main direction of a rectangular box. This test case was first proposed in \cite{SaltzmanOrg} for a two-dimensional Cartesian grid that has been skewed and it represents a very challenging test problem that allows the robustness of any Lagrangian scheme to be validated, because the mesh is not aligned with the fluid motion. According to \cite{Maire2009b}, we consider the three-dimensional extension of the original problem \cite{SaltzmanOrg,SaltzmanOrg3D}, hence the initial computational domain is the box $\Omega(0)=[0;1]\times[0;0.1]\times[0;0.1]$ which is discretized with a total number of $N_E=50000$ tetrahedral elements. The computational mesh is obtained as follows: 
\begin{itemize} 
	\item the domain is initially meshed with a uniform Cartesian grid composed by $100 \times 10 \times 10$ cubic elements, as done in \cite{Maire2009b};
	\item each cube is then split into five tetrahedra;
	
	\item finally we use the mapping given in \cite{SaltzmanOrg3D,Maire2009b} to transform the uniform grid, defined by the coordinate vector $\mathbf{x}=(x,y,z)$, to the skewed configuration $\mathbf{x'}=(x',y',z')$:
	  \begin{eqnarray}
    	x' &=& x + \left( 0.1 - z \right) \left( 1 - 20y \right) \sin(\pi x) \quad \textnormal{for} \quad 0 \leq y \leq 0.05 \nonumber, \\
    	x' &=& x + z\left( 20y - 1 \right) \sin(\pi x) \quad \textnormal{for} \quad 0.05 < y \leq 0.1 \nonumber, \\
    	y' &=& y \nonumber, \\ 
    	z' &=& z.
	   \label{eqSaltzSkew}
	  \end{eqnarray}
	  The initial mesh configuration as well as the final mesh configuration are depicted in Figure \ref{fig:SaltzGrid}.
\end{itemize} 

According to \cite{chengshu2}, the computational domain is filled with a perfect gas with the initial state $\Q_{0}$ given by
\begin{equation}
 \Q_{0} = \left( 1, 0, 0, 0, \epsilon \right).
\label{eqSaltz_ini}
\end{equation}
The ratio of specific heats is taken to be $\gamma = \frac{5}{3}$, $\epsilon = 10^{-4}$ and the final time is set to $t_f=0.6$. The piston is traveling from the left to the right side of the domain with velocity $\mathbf{v}_p = (1,0,0)$ and it starts moving at the initial time while the gas is at rest. In the initial time steps the scheme must obey a geometric $\textnormal{CFL}$ condition, i.e. the piston must not move
more than one element per time step. Sliding wall boundary conditions have been set everywhere, except for the piston, which has been assigned with moving slip wall boundary condition.

The exact solution $\Q_{ex}$ for the Saltzman problem can be computed by solving a one-dimensional Riemann problem, see \cite{Lagrange2D,ToroBook} for details. It reads
\begin{equation}
  \Q_{ex}(\x,t_f) = \left\{ \begin{array}{ccc} \left( 4, 4, 0, 0, 4     \right) & \textnormal{ if } & x \leq x_f, \\
                                               \left( 1, 0, 0, 0, \epsilon \right) & \textnormal{ if } & x > x_f,        
                      \end{array}  \right. 
\end{equation}
where $x_f=0.8$ is the shock location at time $t_f = 0.6$.

\begin{figure}[!htbp]
	\begin{center}
	\begin{tabular}{cc} 
	\includegraphics[width=0.47\textwidth]{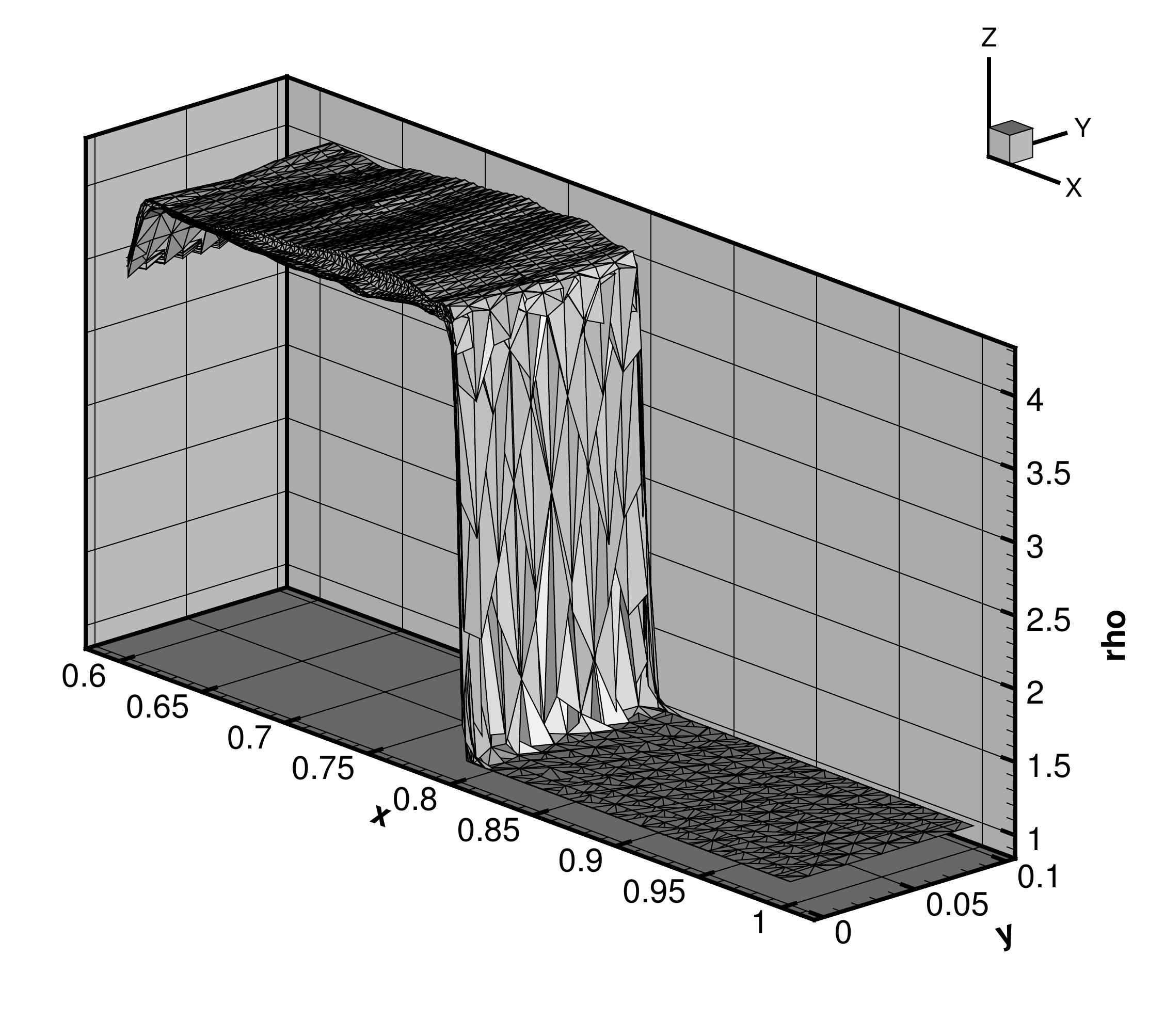}  &           
	\includegraphics[width=0.47\textwidth]{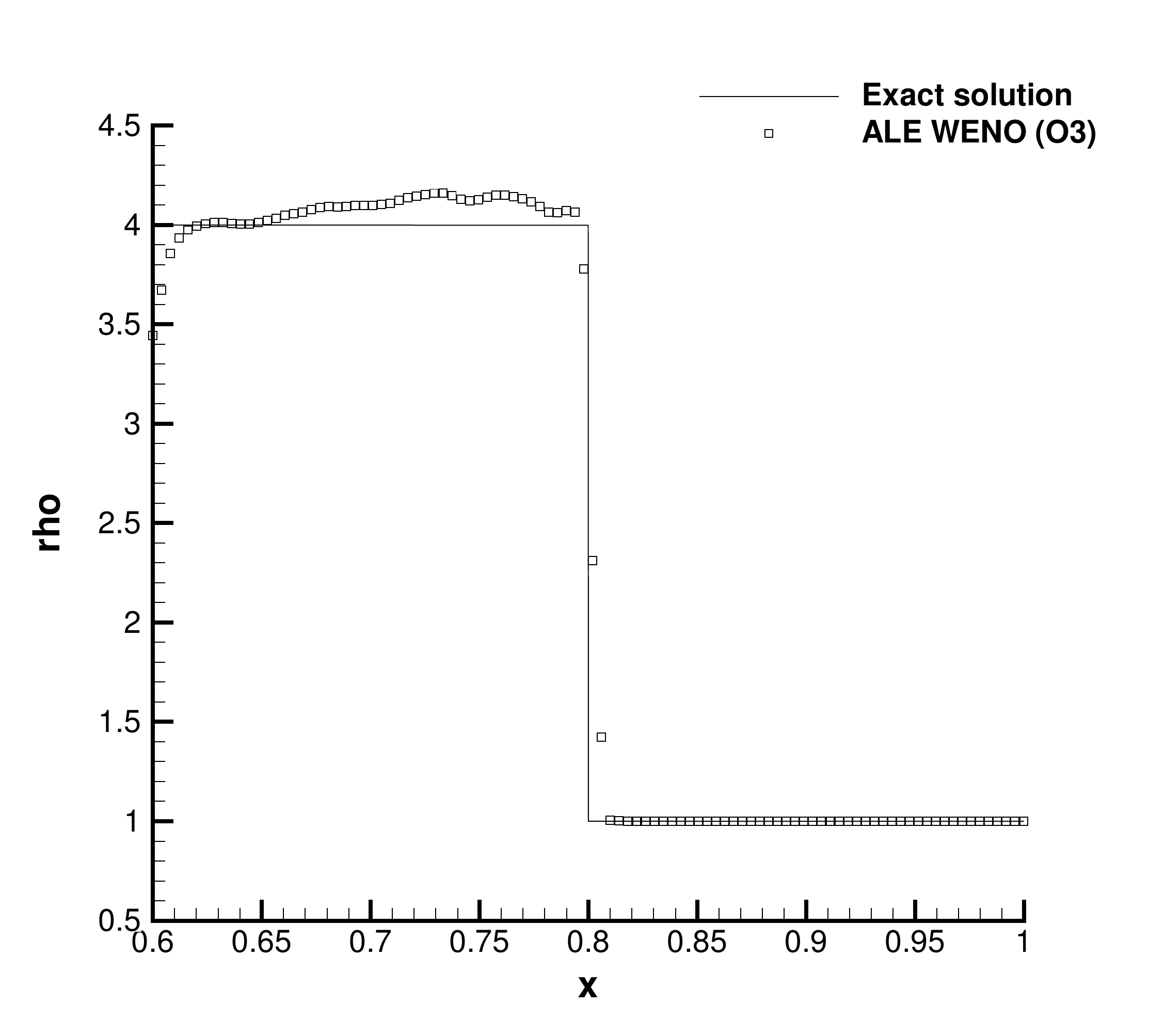} \\
	\includegraphics[width=0.47\textwidth]{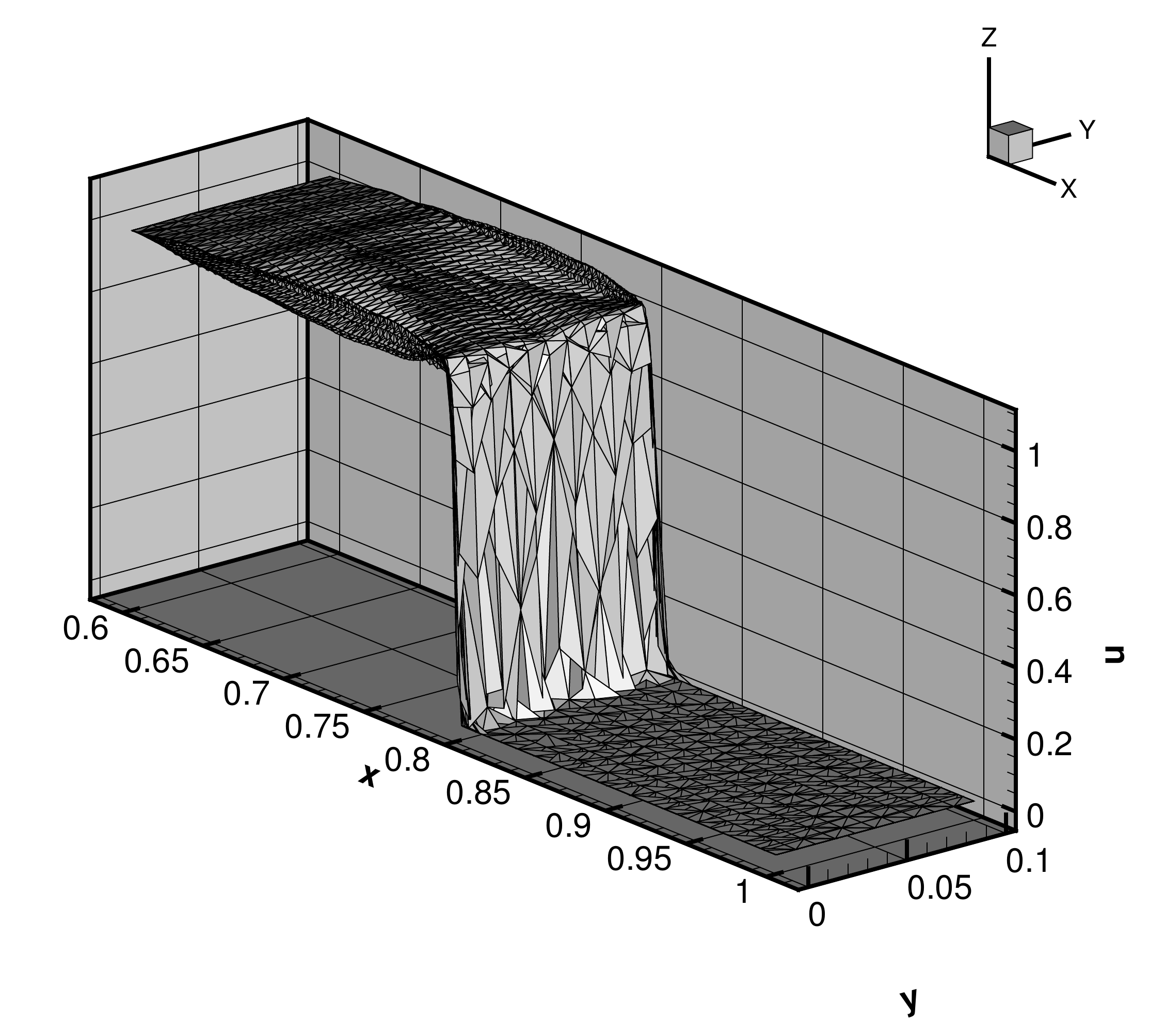}  &           
	\includegraphics[width=0.47\textwidth]{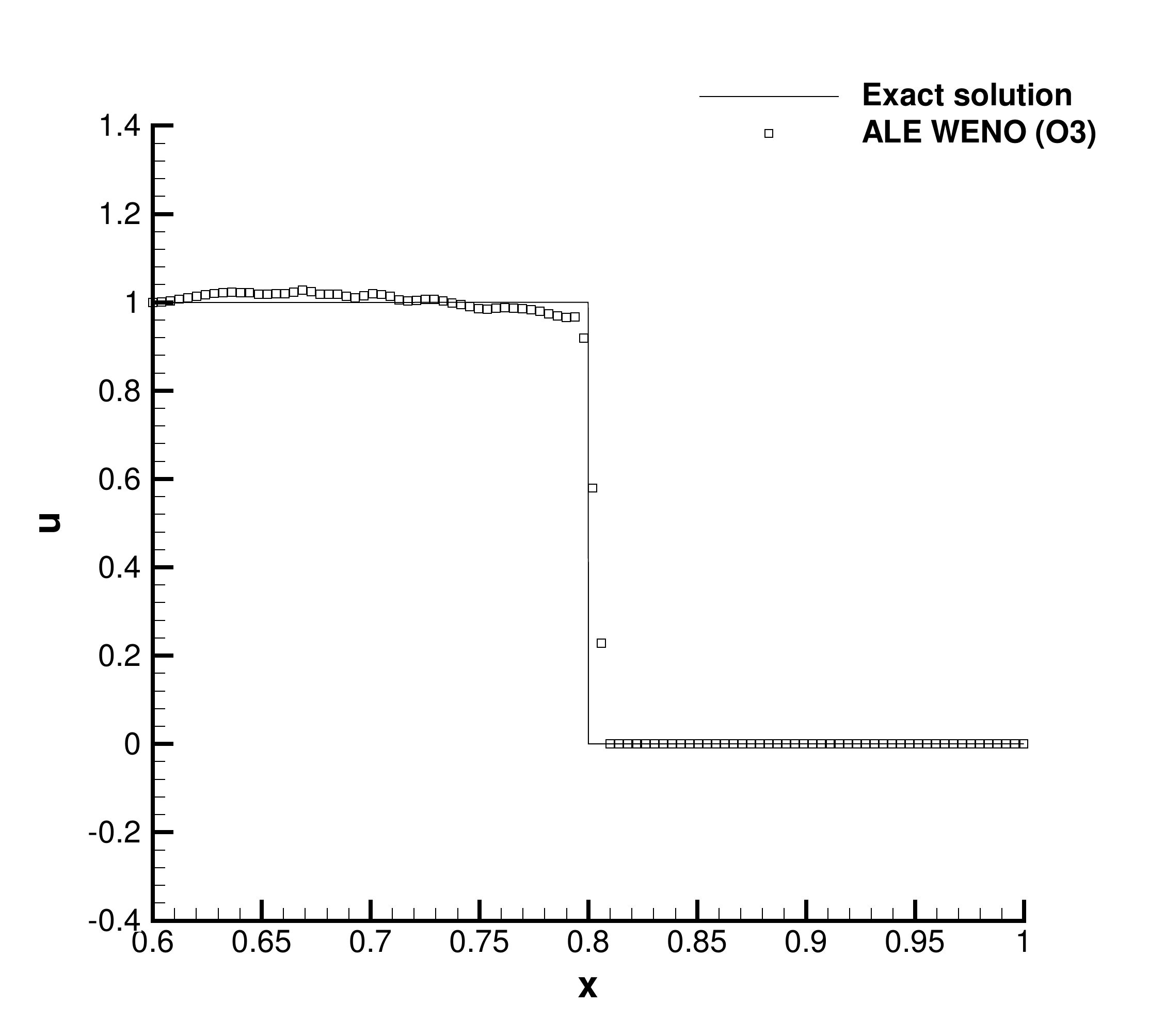} \\
	\end{tabular}
	\caption{Third order numerical results for the Saltzman problem: density (top) and velocity (bottom) distribution and comparison with analytical solution at time $t=0.6$.}
	\label{fig:Saltz3D}
	\end{center}
\end{figure}

The numerical results have been obtained with the third order ALE WENO scheme using a robust Rusanov-type numerical flux \eqref{eqn.rusanov} and they are depicted in Figure \ref{fig:Saltz3D}. A good agreement with the exact solution can be noticed regarding both density and velocity distribution at the final time $t_f=0.6$. The decrease of density near the piston is due to the well 
known \textit{wall-heating problem}, see \cite{toro.anomalies.2002}. The positivity preserving technique presented in Section \ref{sec.flattener} has been used to smear out some unphysical 
oscillations occurring at the shock.

\subsubsection{The Sedov problem} 
\label{sec.Sedov}

\begin{figure}[!htbp]
	\begin{center}
	\begin{tabular}{cc} 
	\includegraphics[width=0.47\textwidth]{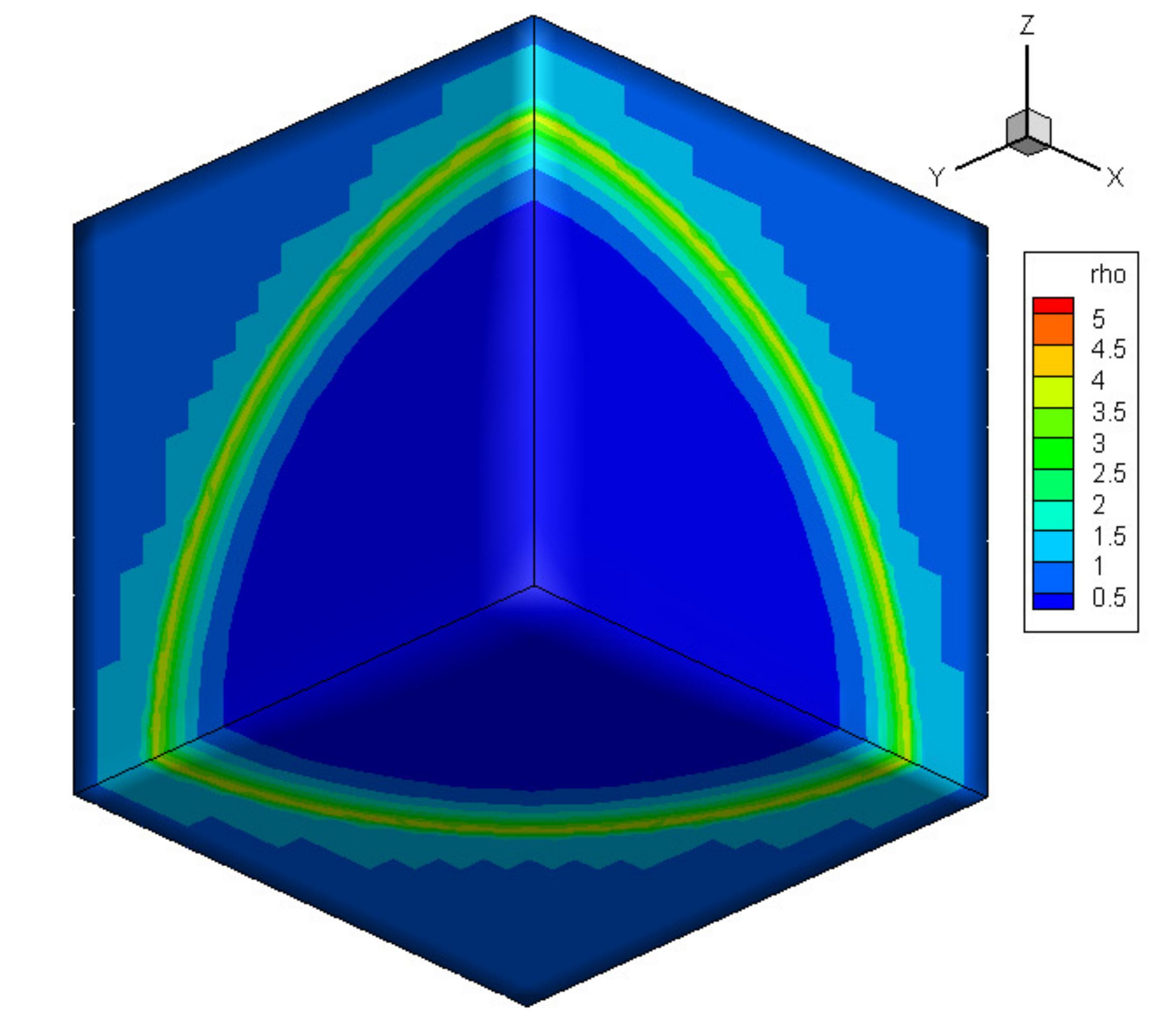}  &           
	\includegraphics[width=0.47\textwidth]{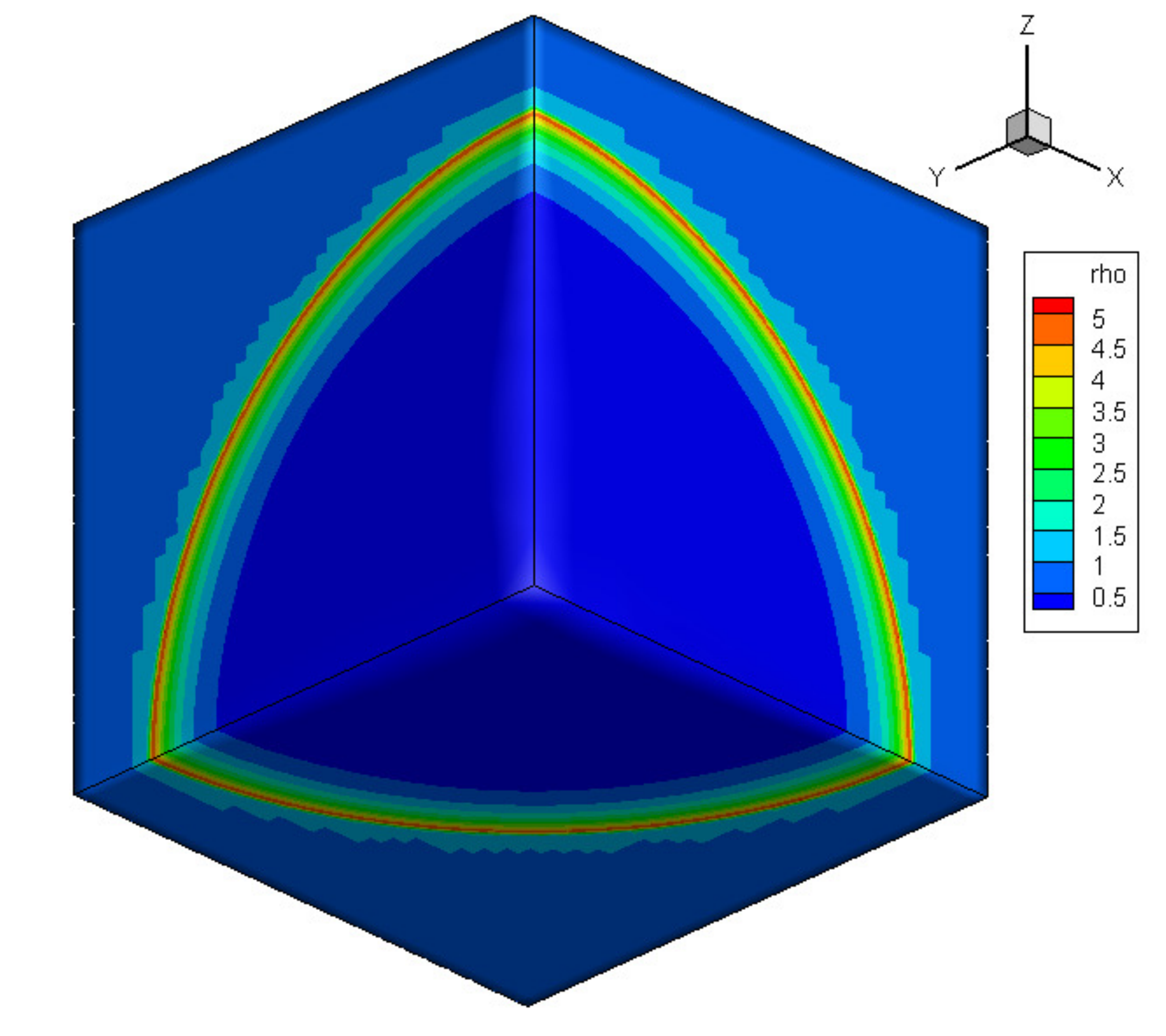} \\
	\includegraphics[width=0.47\textwidth]{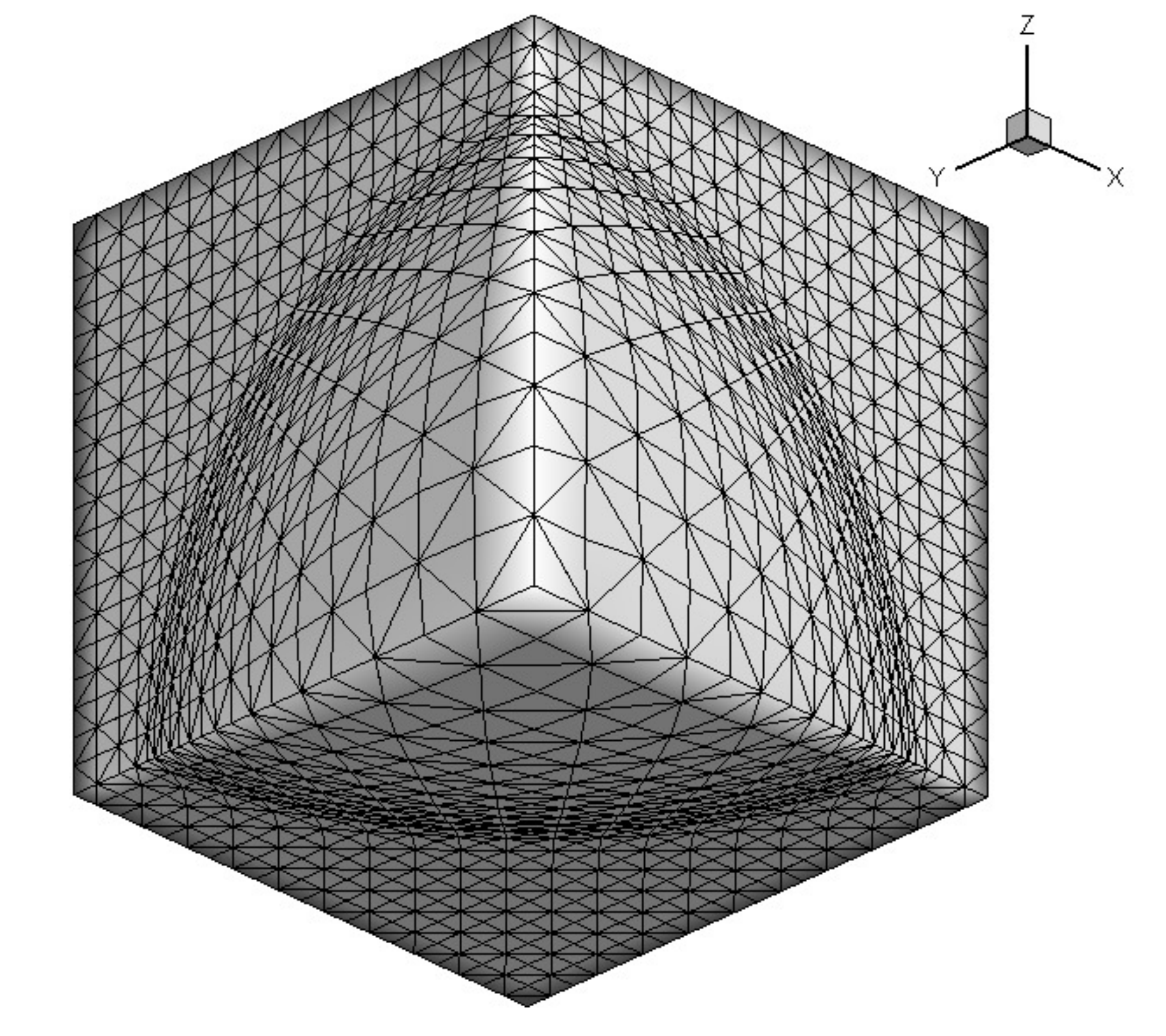}  &           
	\includegraphics[width=0.47\textwidth]{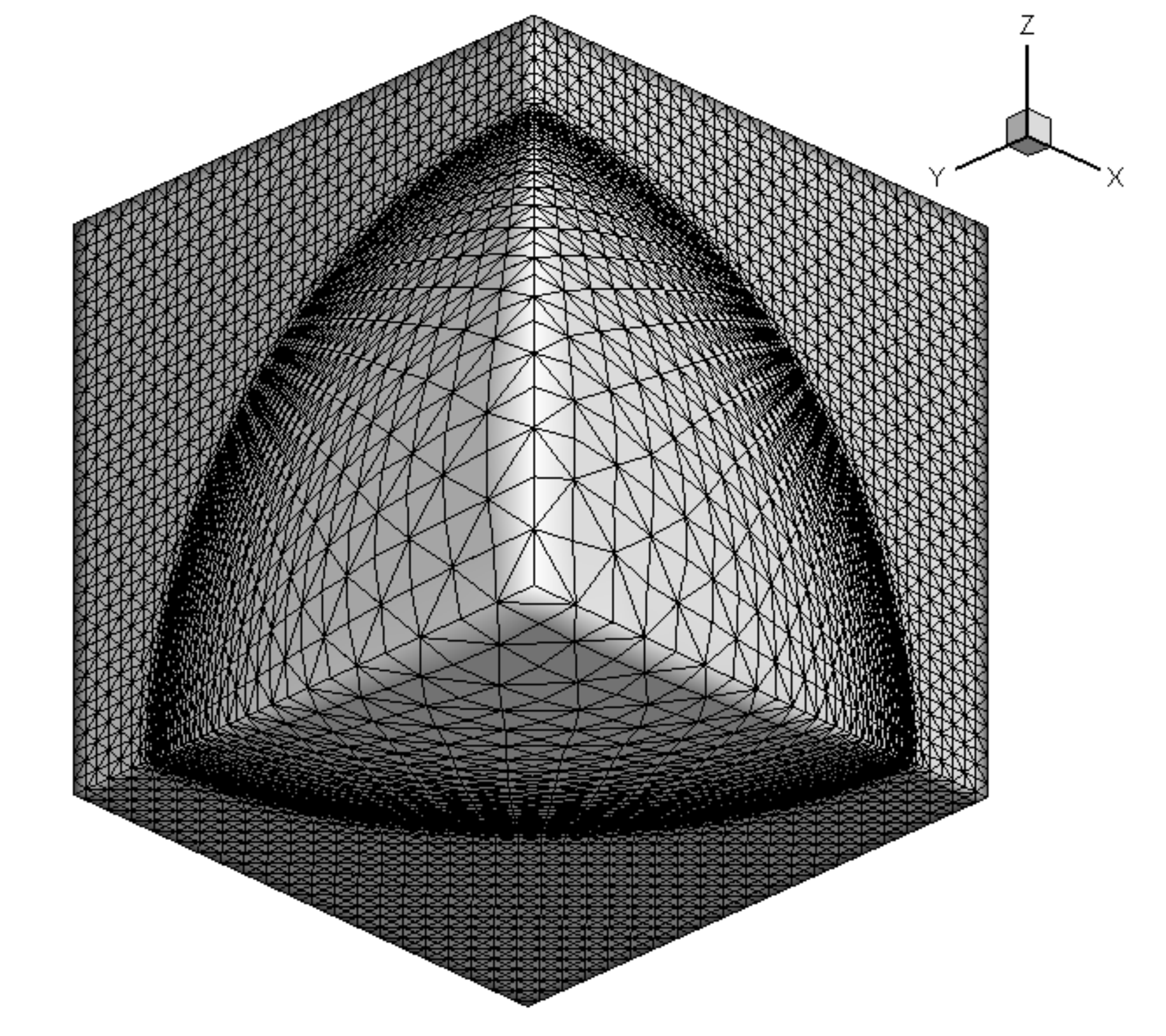} \\
	\includegraphics[width=0.47\textwidth]{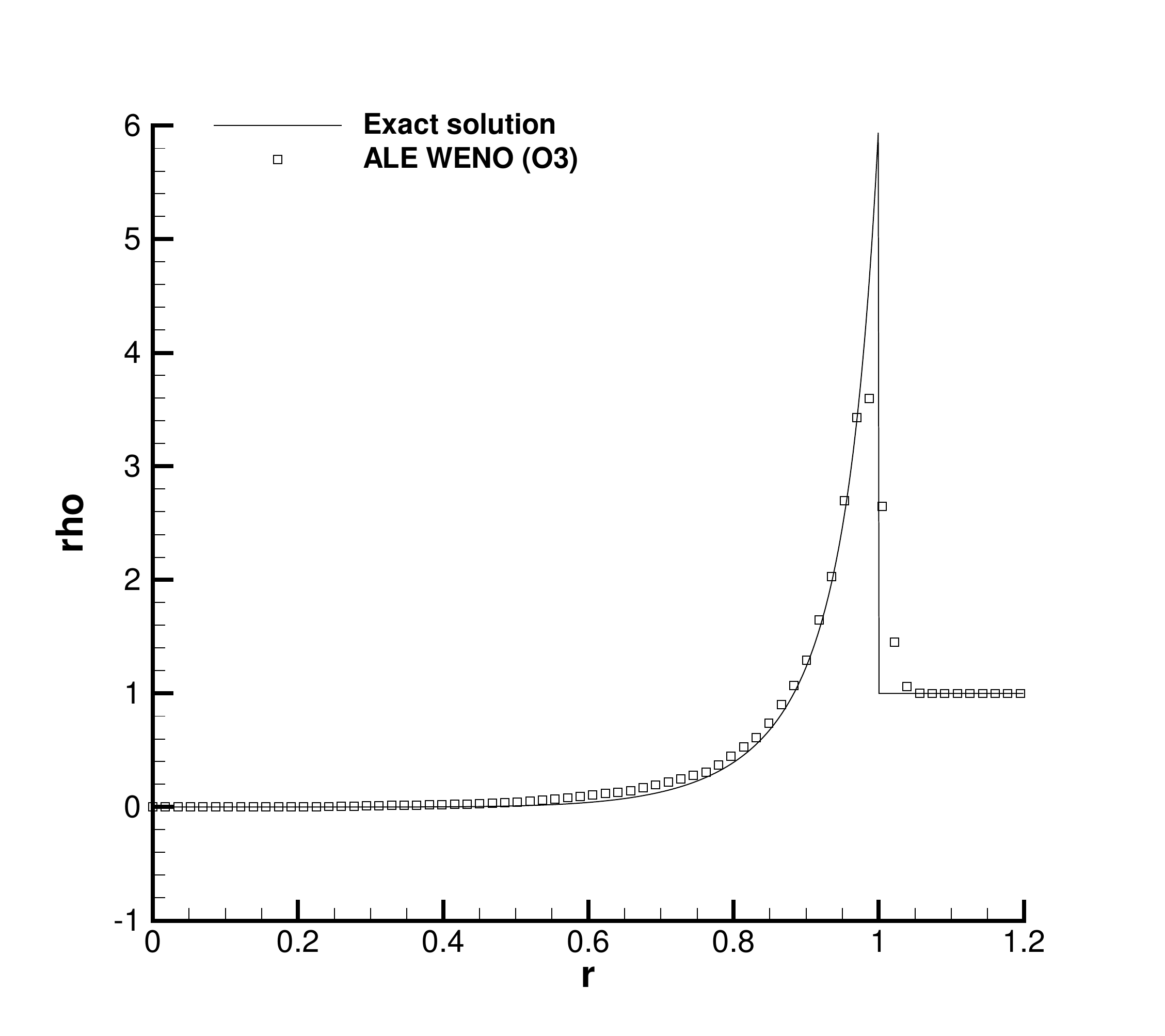}  &           
	\includegraphics[width=0.47\textwidth]{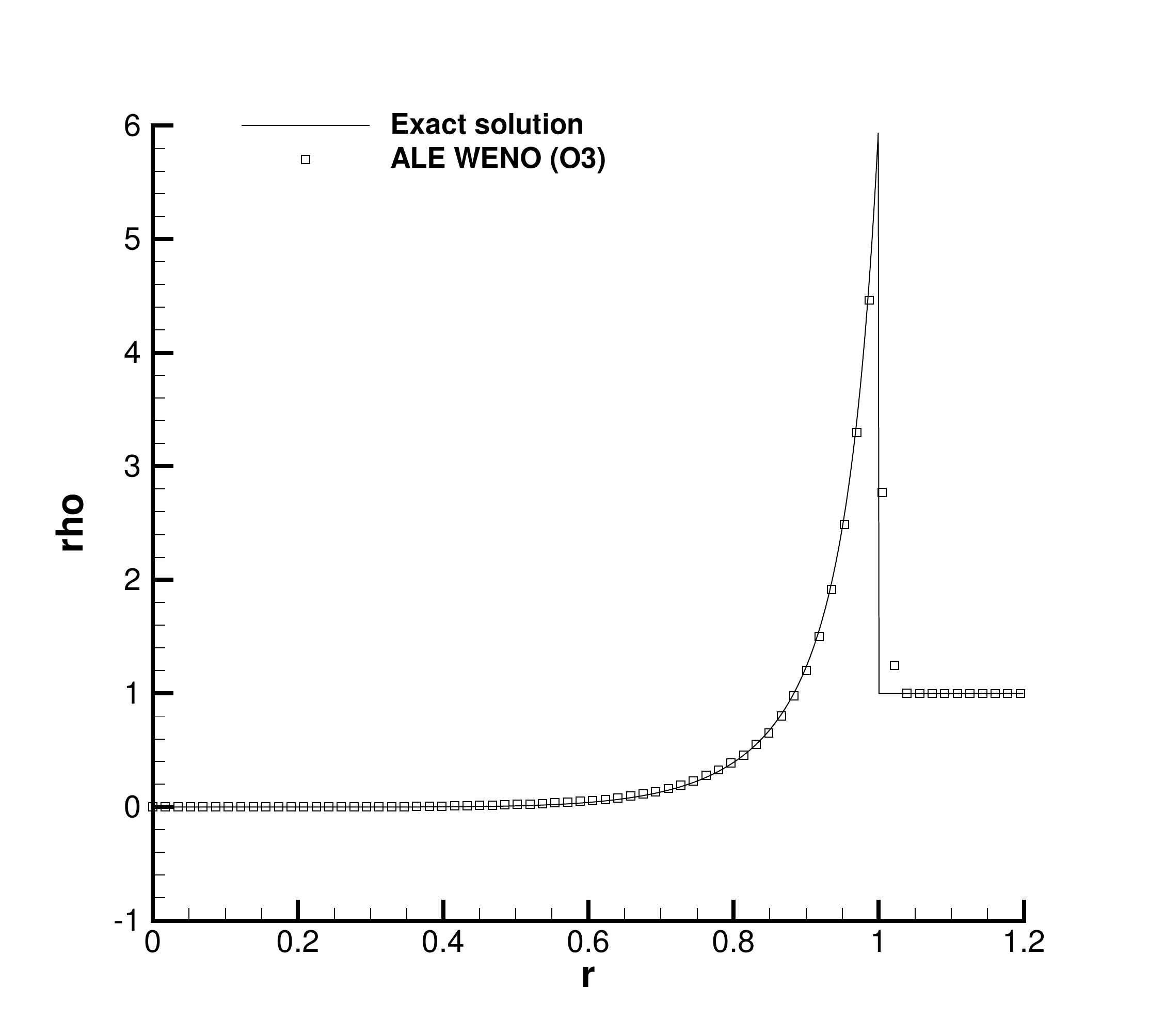} \\
	\end{tabular}
	\caption{Third order results for the Sedov problem on the coarse grid $m_1$ (left column) and on the fine grid $m_2$ (right column). From top to bottom: solution for density at the final time of the simulation (top row), mesh configuration at the final time $t_f=1.0$ (middle row) and comparison between analytical and numerical density distribution along the diagonal straight line that crosses the cubic computational domain (bottom row).}
	\label{fig:Sedov3D}
	\end{center}
\end{figure}

In this section we consider the spherical symmetric Sedov problem, which describes the evolution of a blast wave generated at the origin $\mathbf{O}=(x,y,z)=(0,0,0)$ of the initial cubic computational domain $\Omega(0)=[0;1.2]\times[0;1.2]\times[0;1.2]$. It is a well-known test case for Lagrangian schemes \cite{Maire2009,Maire2009b,LoubereSedov3D} that becomes very challenging in the three-dimensional case. An analytical solution which is based on self-similarity arguments is furthermore available from the work of Kamm et al. \cite{SedovExact}. As done in \cite{LoubereSedov3D} we consider two different meshes, the first one $m_1$ is composed by $20\times 20 \times 20$ cubes, while the second one $m_2$ involves $40\times 40 \times 40$ elements. Each cube is then split into five tetrahedra for a total number of elements of $N_{E,1}=40000$ and $N_{E,2}=320000$. The computational domain is filled with a prefect gas with $\gamma=1.4$, which is initially at rest and is assigned with a uniform density $\rho_0=1$. The total energy $E_{tot}$ is concentrated only in the cell $c_{or}$ containing the origin $\mathbf{O}$, therefore the initial pressure is given by
\begin{equation}
p_{or} = (\gamma-1)\rho_0 \frac{E_{tot}}{8 \cdot V_{or}},
\label{eqn.p0.sedov}
\end{equation}
where $V_{or}$ is the volume of the cell $c_{or}$,  which is composed by five tetrahedra, and the factor $\frac{1}{8}$ takes into account the spherical symmetry, since the computational domain $\Omega(0)$ is only the eighth part of the entire domain, which would have to be considered if we did not assume the spherical symmetry. According to \cite{LoubereSedov3D} we set $E_{tot}=0.851072$, while in the rest of the domain the initial pressure is $p_0=10^{-6}$. At the final time of the simulation $t_f=1.0$ the exact solution is a symmetric spherical shock wave located at radius $r=1$ with a density peak of $\rho=6$. We use the third order accurate version of the ALE WENO schemes presented in this paper together with the Rusanov-type numerical flux \eqref{eqn.rusanov} and the positivity preserving algorithm illustrated in Section \ref{sec.flattener}. The numerical solution for the Sedov problem has been computed on both meshes $m_1$ and $m_2$. Figure \ref{fig:Sedov3D} shows the solution for density at the final time of the simulation as well as the mesh configuration and a comparison between the numerical and the exact density distribution along the radial direction.   

\subsection{The magnetohydrodynamics (MHD) equations} 
The equations of ideal classical magnetohydrodynamics (MHD) constitute a more complicated hyperbolic conservation law compared to the Euler equations used so far, especially because this 
system introduces an additional constraint regarding the divergence of the magnetic field that must remain zero in time, i.e. 
\begin{equation}
\nabla \cdot \mathbf{B} = 0.
\label{eq:divB}
\end{equation}  
If the magnetic field $\mathbf{B}$ is initialized with data that are guaranteed to be divergence-free, then Eqn. \eqref{eq:divB} is always satisfied for the exact solution. 
The difficulty appears at the discrete level, where the numerical divergence-free constraint has to be carefully taken into account and properly treated. As done for the 
two-dimensional ALE WENO finite volume schemes presented in \cite{LagrangeMHD}, we adopt the  hyperbolic version of the generalized Lagrangian multiplier (GLM) divergence cleaning 
approach proposed by Dedner et al. \cite{Dedneretal}, hence adding to the MHD system one more variable $\Psi$ as well as one more linear scalar PDE that aims at transporting the 
divergence errors out of the computational domain with an artificial divergence cleaning speed $c_h$. The augmented MHD system can be cast into form \eqref{eqn.pde.nc} and reads 
\begin{equation}
\label{MHDTerms} 
  \frac{\partial}{\partial t} \left( \begin{array}{c} 
  \rho \\   \rho \v \\   \rho E \\ \B \\ \psi \end{array} \right)   
  + \nabla \cdot \left( \begin{array}{c} 
  \rho \v  \\ 
  \rho \v \v + p_t \mathbf{I} - \frac{1}{4 \pi} \B \B \\ 
  \v (\rho E + p_t ) - \frac{1}{4 \pi} \B ( \v \cdot \B ) \\ 
  \v \B - \B \v + \psi \mathbf{I} \\
  c_h^2 \B 
  \end{array} \right) = 0. 
\end{equation} 
The non-conservative part of the ideal MHD system is zero, the velocity vector is denoted by $\mathbf{v}=v_i=(u,v,w)$ and similarly the vector of the magnetic field is addressed with $\mathbf{B}=B_i=(B_x,B_y,B_z)$. The system is then closed by the equation of state
\begin{equation}
p = \left(\gamma - 1 \right) \left(\rho E - \frac{1}{2}\mathbf{v}^2 - \frac{\mathbf{B}^2}{8\pi}\right),
\label{MHDeos}
\end{equation}
with $\gamma$ representing the ratio of specific heats and the total pressure being defined as $p_t = p + \frac{1}{8 \pi} \B^2$.

We define also the fastest magnetosonic speed, needed for the node solver $\mathcal{NS}_m$. It reads
\begin{equation}
c = \sqrt{ \frac{1}{2} \left( \frac{\gamma p}{\rho} + (B_x+B_y+B_z) + \sqrt{ \left(\frac{\gamma p}{\rho} + (B_x+B_y+B_z)\right)^2 - 4 \frac{\gamma p}{\rho} \frac{B_x^2}{4\pi\rho}} \right) }.
\label{eq:cMHD}
\end{equation}

\subsubsection{The MHD rotor problem} 
\label{sec.MHDRotor}

\begin{figure}[!htbp]
\begin{center}
\begin{tabular}{cc} 
\includegraphics[width=0.47\textwidth]{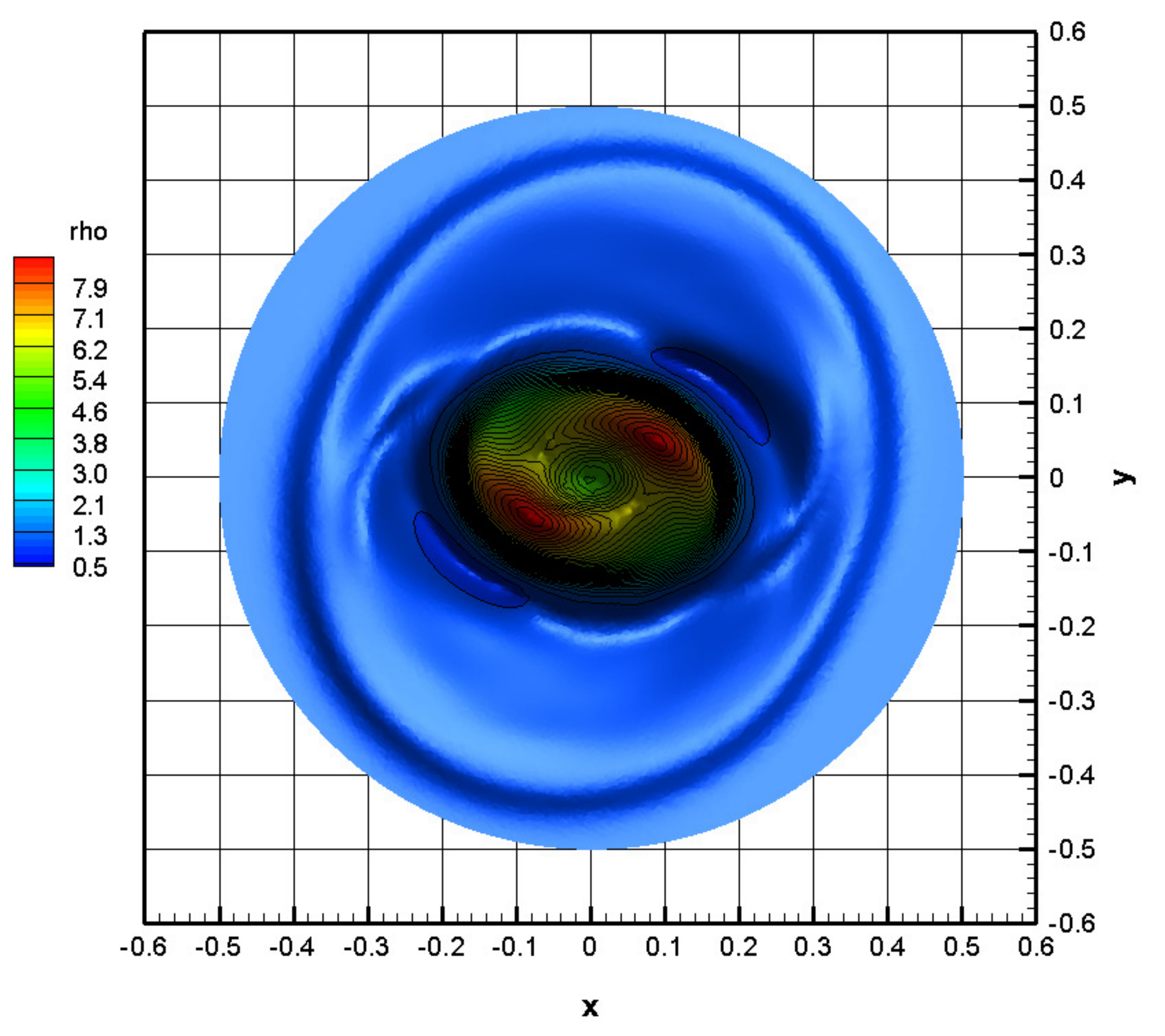}  &           
\includegraphics[width=0.47\textwidth]{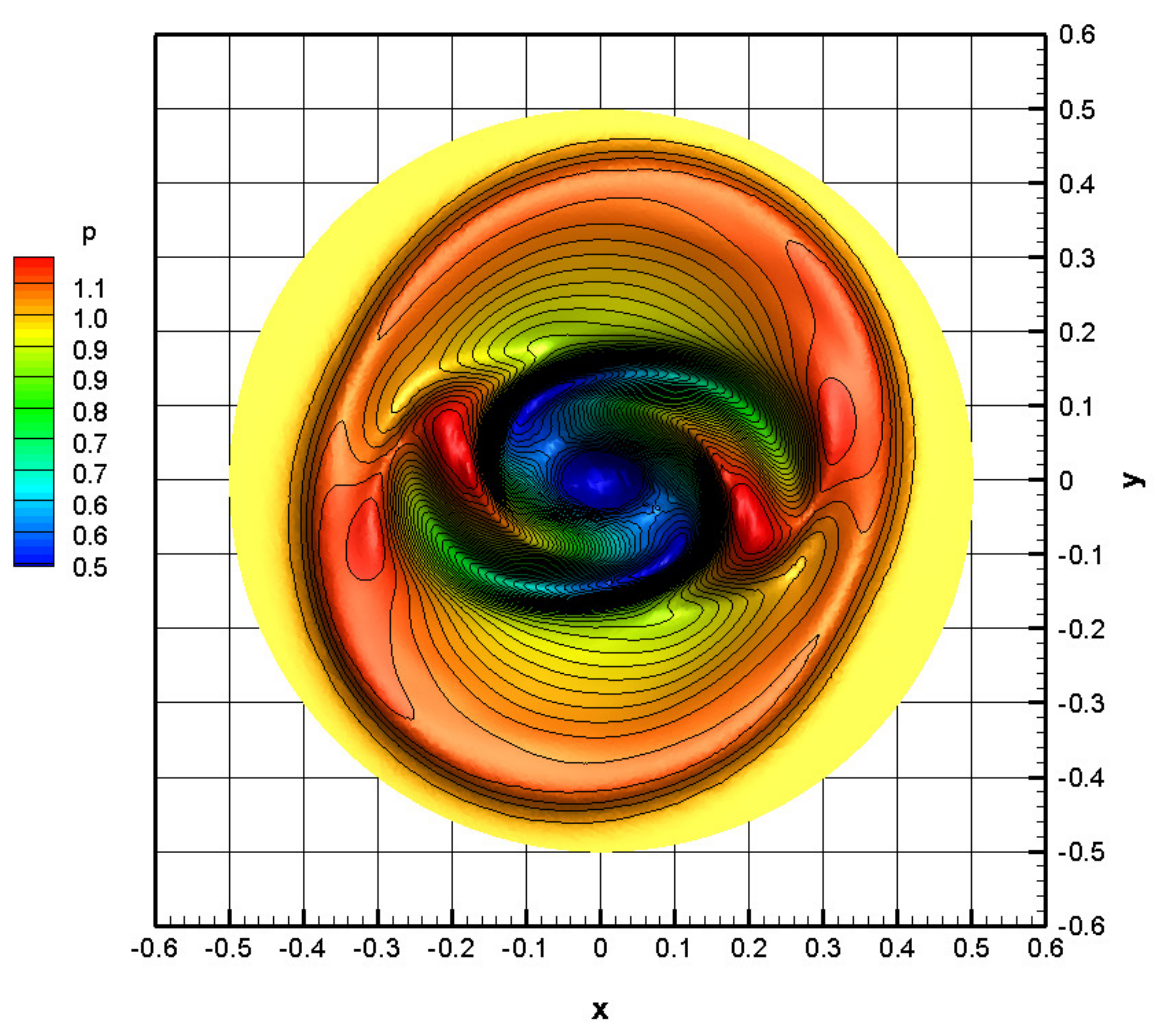} \\
\includegraphics[width=0.47\textwidth]{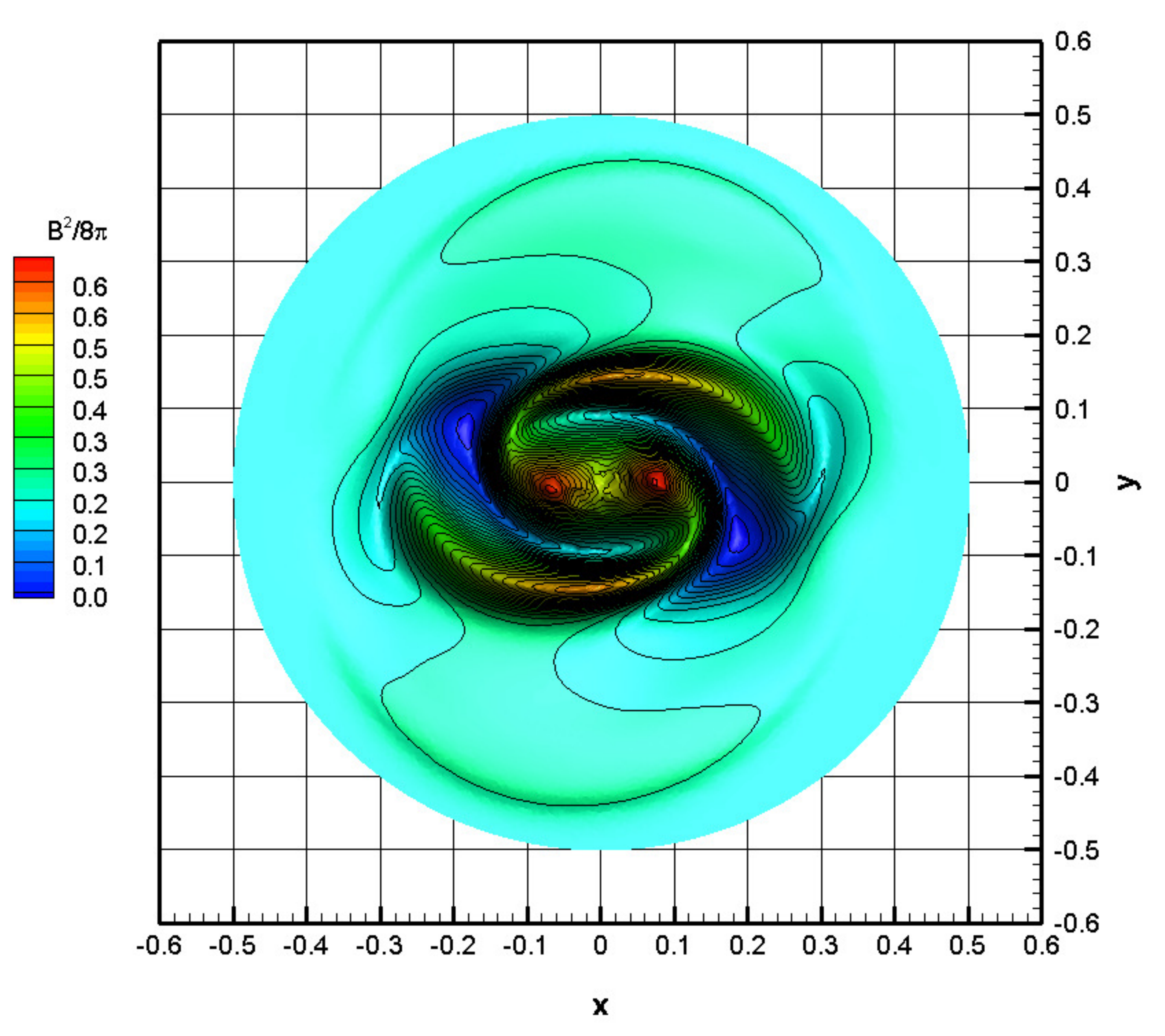} & 
\includegraphics[width=0.47\textwidth]{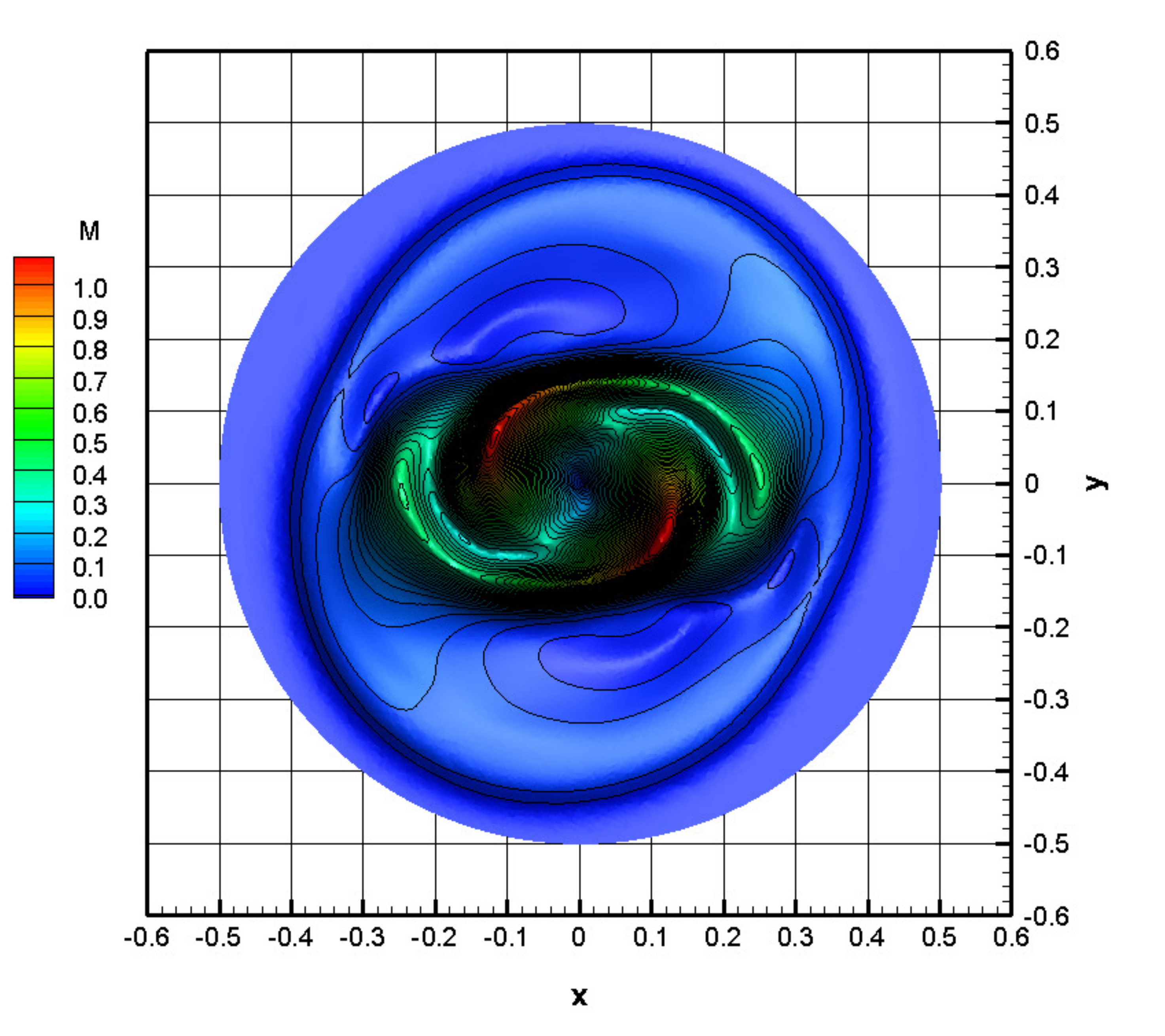}  \\  
\end{tabular} 
\caption{Third order numerical results for the ideal MHD rotor problem at time $t=0.25$. Top: density and pressure. Bottom: magnitude of the magnetic field and Mach number.} 
\label{fig.MHDRotor}
\end{center}
\end{figure}

\begin{figure}[!htbp]
\begin{center}
\begin{tabular}{cc} 
\includegraphics[width=0.47\textwidth]{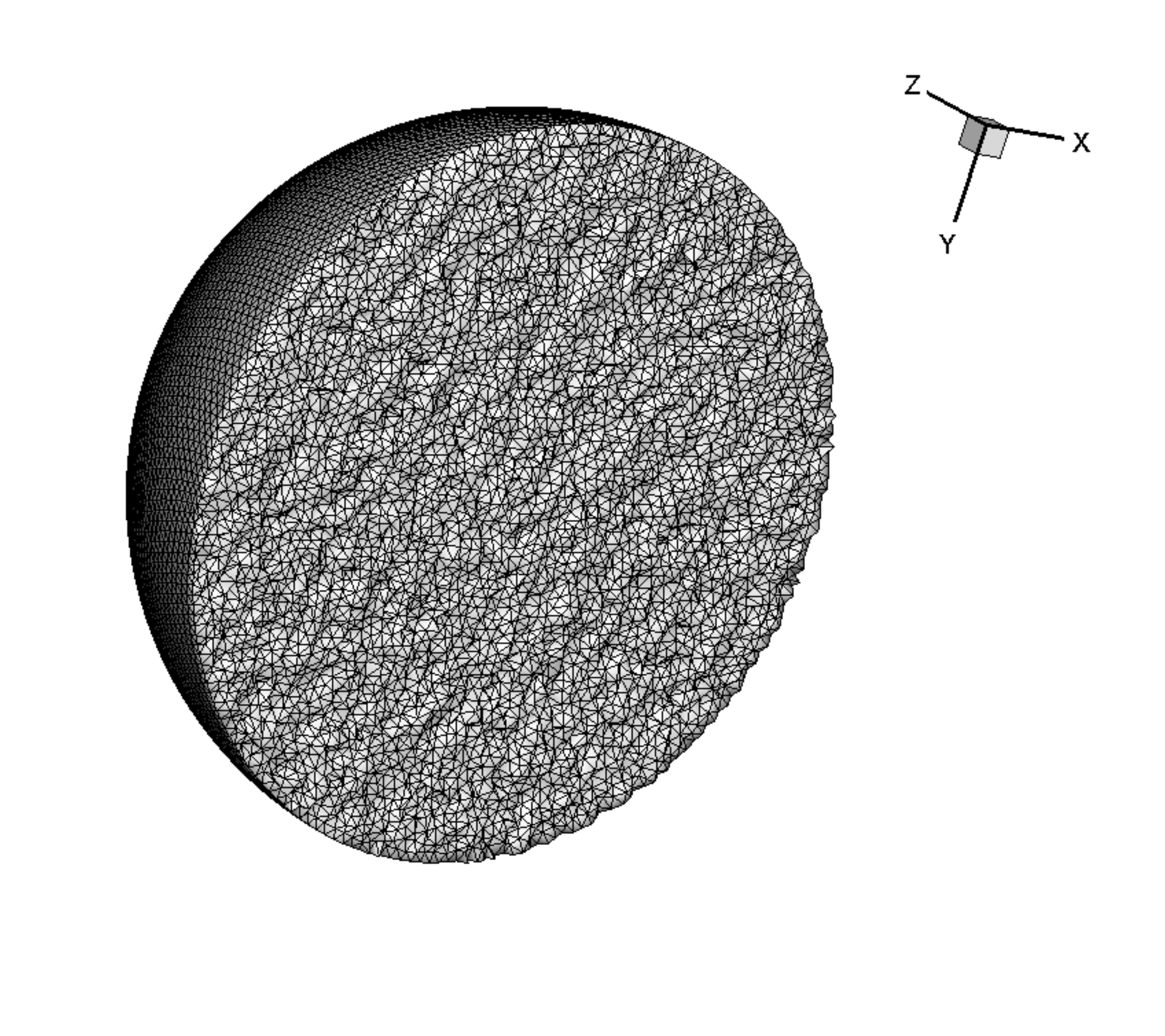}  &           
\includegraphics[width=0.47\textwidth]{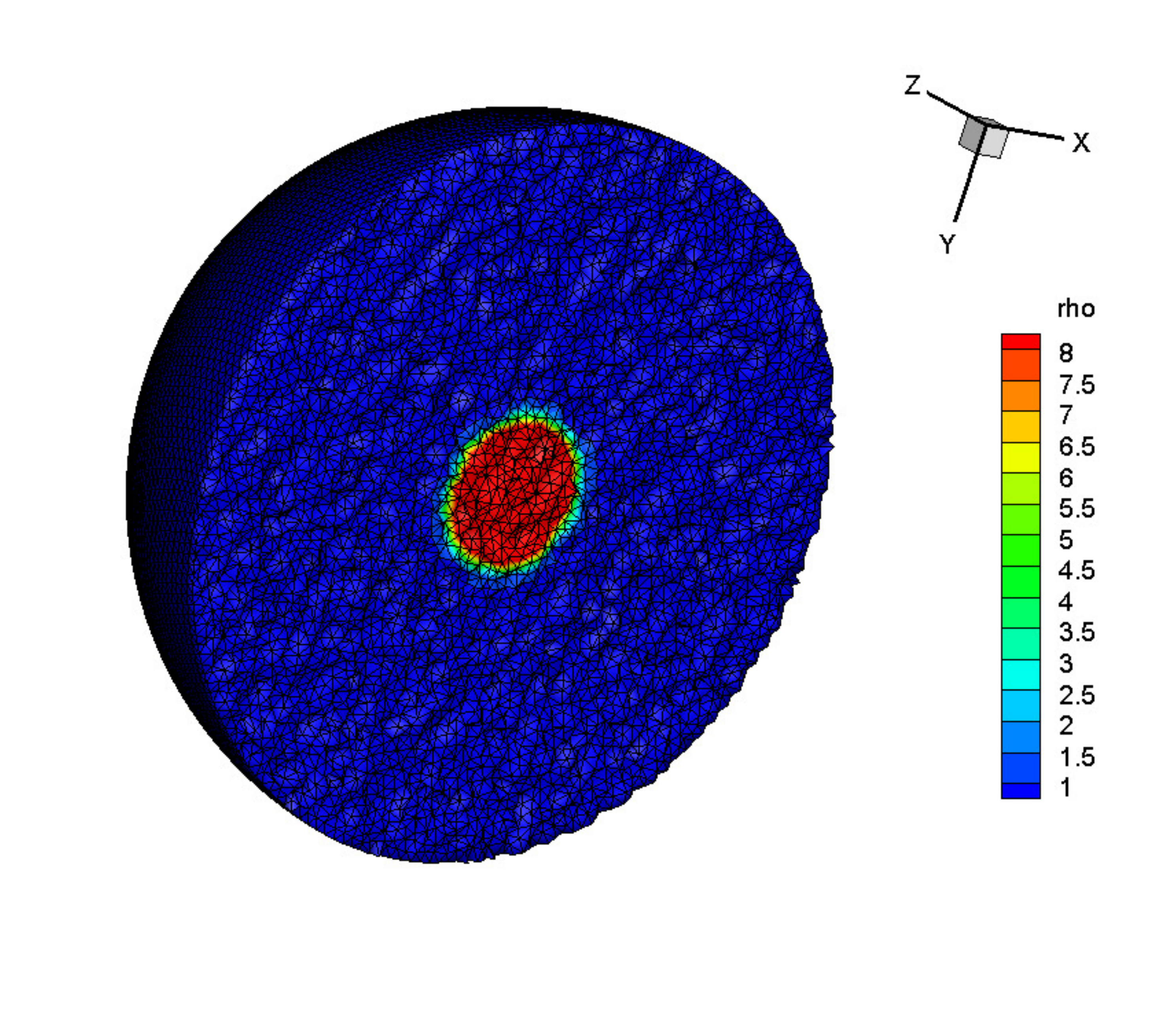} \\
\includegraphics[width=0.47\textwidth]{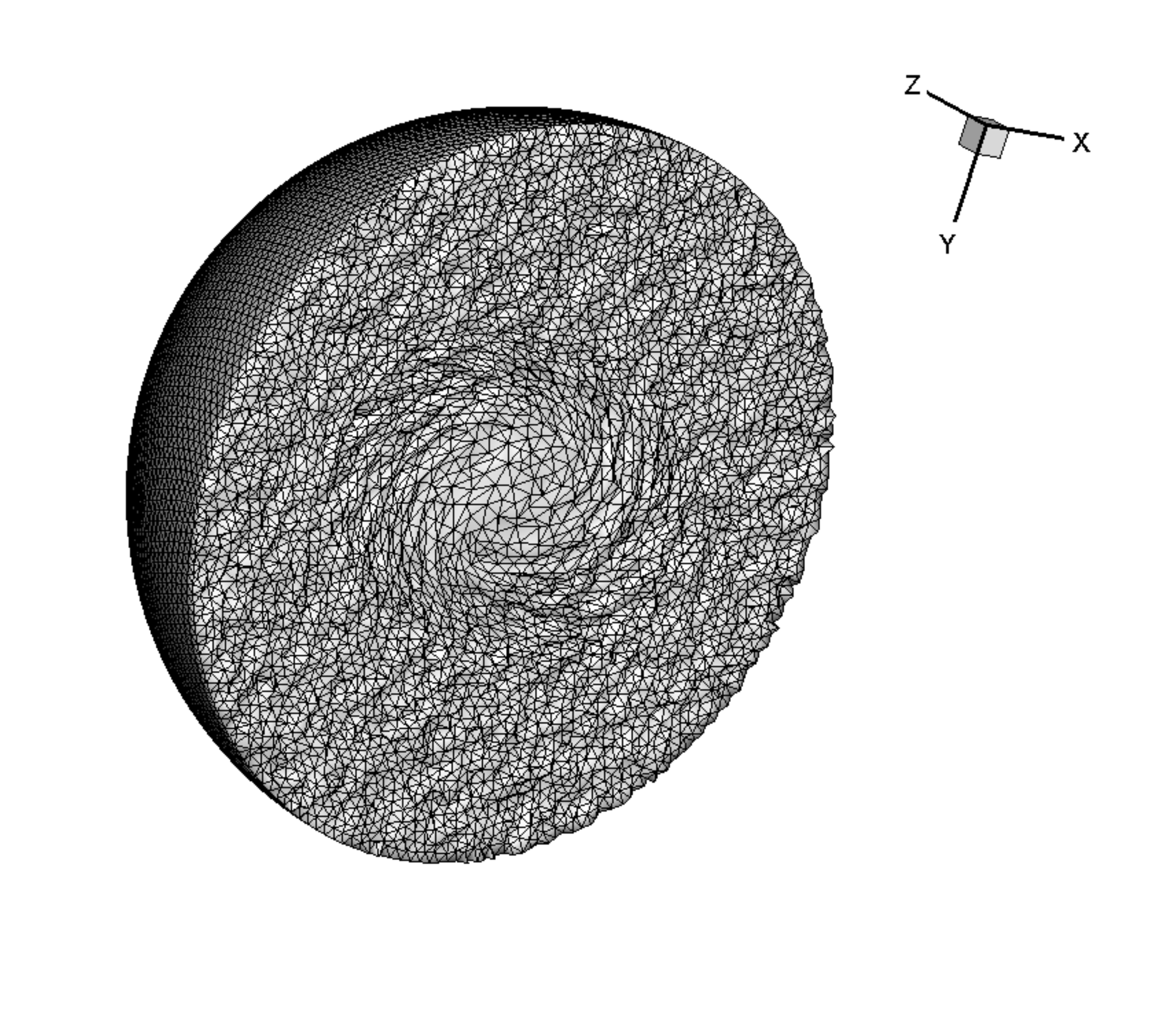} & 
\includegraphics[width=0.47\textwidth]{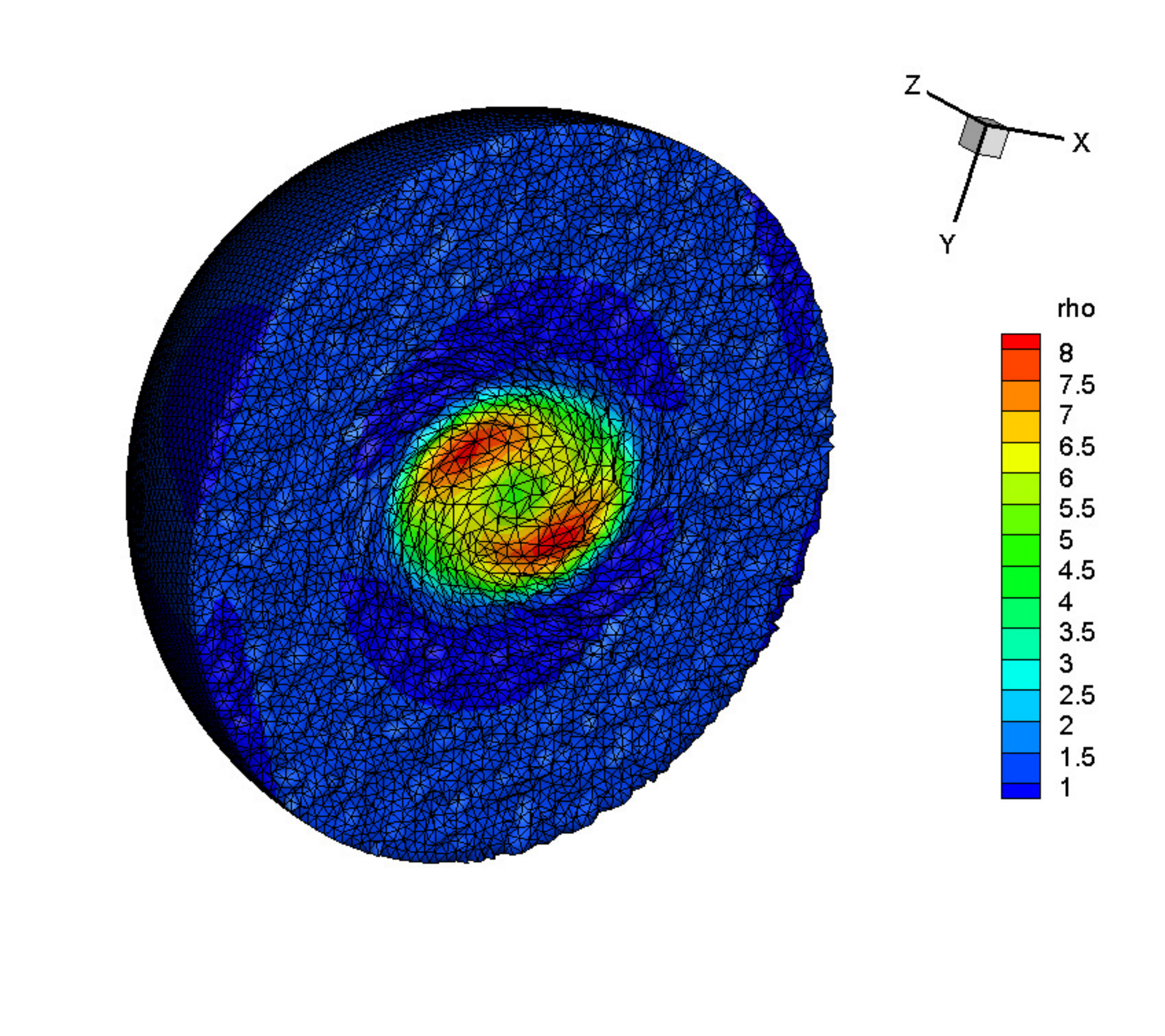}  \\  
\end{tabular} 
\caption{Mesh configuration and density distribution for the MHD rotor problem at the initial time $t=0.0$ (top) and at the final time $t=0.25$ (bottom).} 
\label{fig.MHDRotorGrid}
\end{center}
\end{figure}

The first test case for the ideal classical MHD equations is the MHD rotor problem proposed by Balsara et al. in \cite{BalsaraSpicer1999}. It consists in a fluid of high density that is rotating very quickly, surrounded by a fluid at rest with low density. The initial computational domain $\Omega(0)$ is a sphere of radius $R_0=0.5$ which is discretized with a total number of tetrahedra of $N_E=1089071$. The generic radial position is denoted by $r=\sqrt{x^2+y^2+z^2}$ and at radius $R=0.1$ the \textit{inner} region with the high density fluid is separated by the \textit{outer} region. Therefore the initial density distribution is $\rho=10$ for $0 \leq r \leq R$ and $\rho=1$ in the rest of the domain, while the angular velocity $\omega$ of the rotor is assumed to be constant and it is chosen 
in such a way that at $r=R$ the toroidal velocity is $v_t=\omega \cdot R =1$. The initial discontinuity for density and velocity occurring at the frontier $r=R$ is smeared out according to \cite{BalsaraSpicer1999}, where a linear taper bounded by $0.1 \leq r \leq 0.13$ is applied in such a way that the internal values for density and velocity match exactly those ones of the outer region. The pressure is $p=1$ in the whole computational domain and a constant magnetic field $\mathbf{B}=(2.5,0,0)^T$ is imposed everywhere. The divergence cleaning velocity is taken to be $c_h=2$, while the ratio of specific heats is set to $\gamma=1.4$ and the final time is $t_f=0.25$. Transmissive boundary conditions have been imposed at the external boundary. The numerical results for the MHD rotor problem have been obtained using the third order version of the ALE WENO schemes presented in this paper with the Rusanov-type flux \eqref{eqn.rusanov} and they are depicted in Figure \ref{fig.MHDRotor}. Although the mesh adopted for the simulation is coarser than the one used in \cite{BalsaraSpicer1999}, we can note a good qualitative agreement with the solution presented in \cite{BalsaraSpicer1999} (note that the
present simulation is carried out in 3D). 
The rezoning procedure described in Section \ref{sec.meshMot} allows the mesh to be reasonably well shaped, even with the strong deformations produced by the velocity field of the rotor. 
Figure \ref{fig.MHDRotorGrid} shows the initial and the final mesh configuration and the corresponding density distribution. 

\subsubsection{The MHD blast wave problem}
\label{sec.MHDblast}

\begin{figure}[!htbp]
\begin{center}
\begin{tabular}{cc} 
\includegraphics[width=0.47\textwidth]{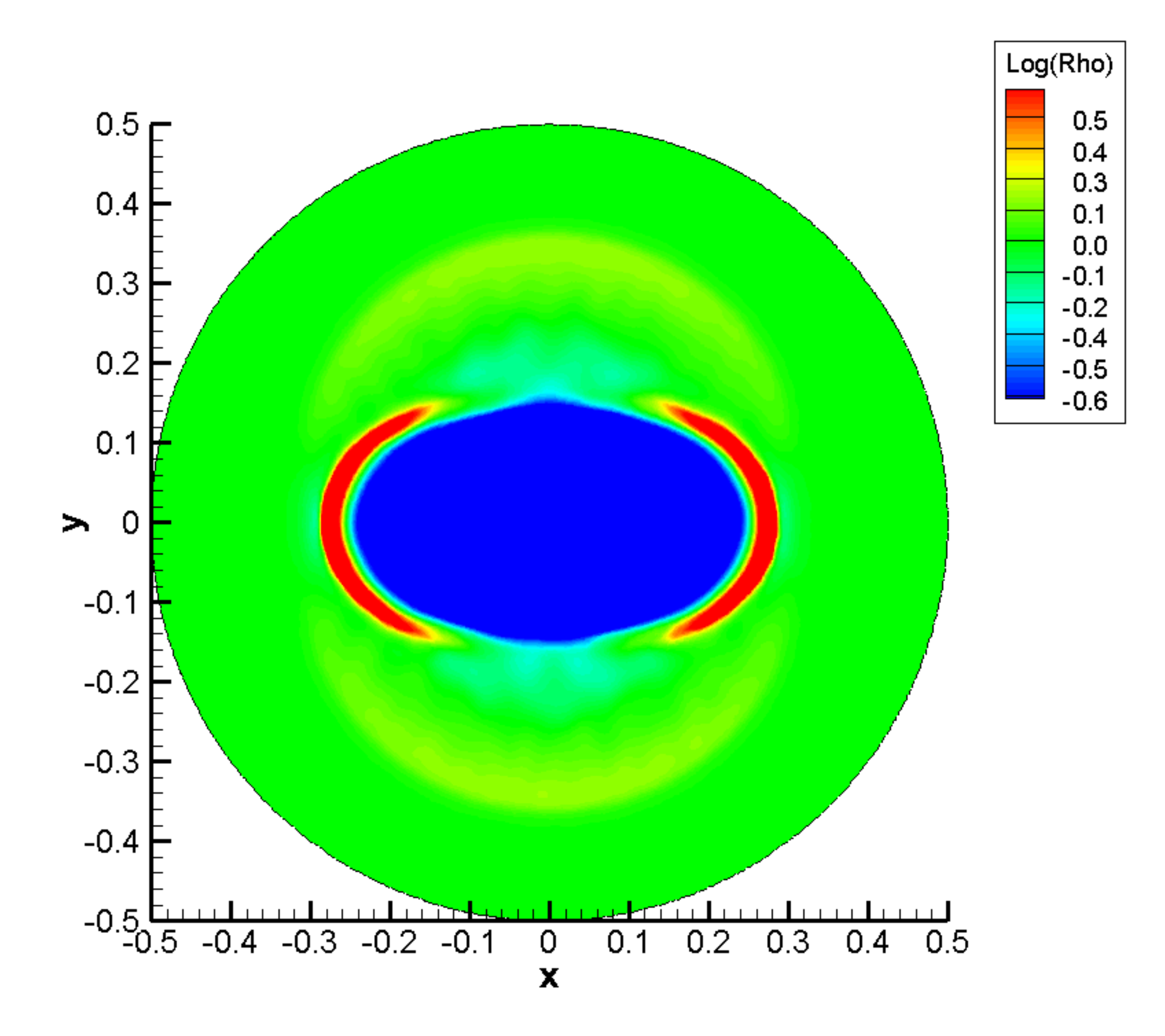}  &           
\includegraphics[width=0.47\textwidth]{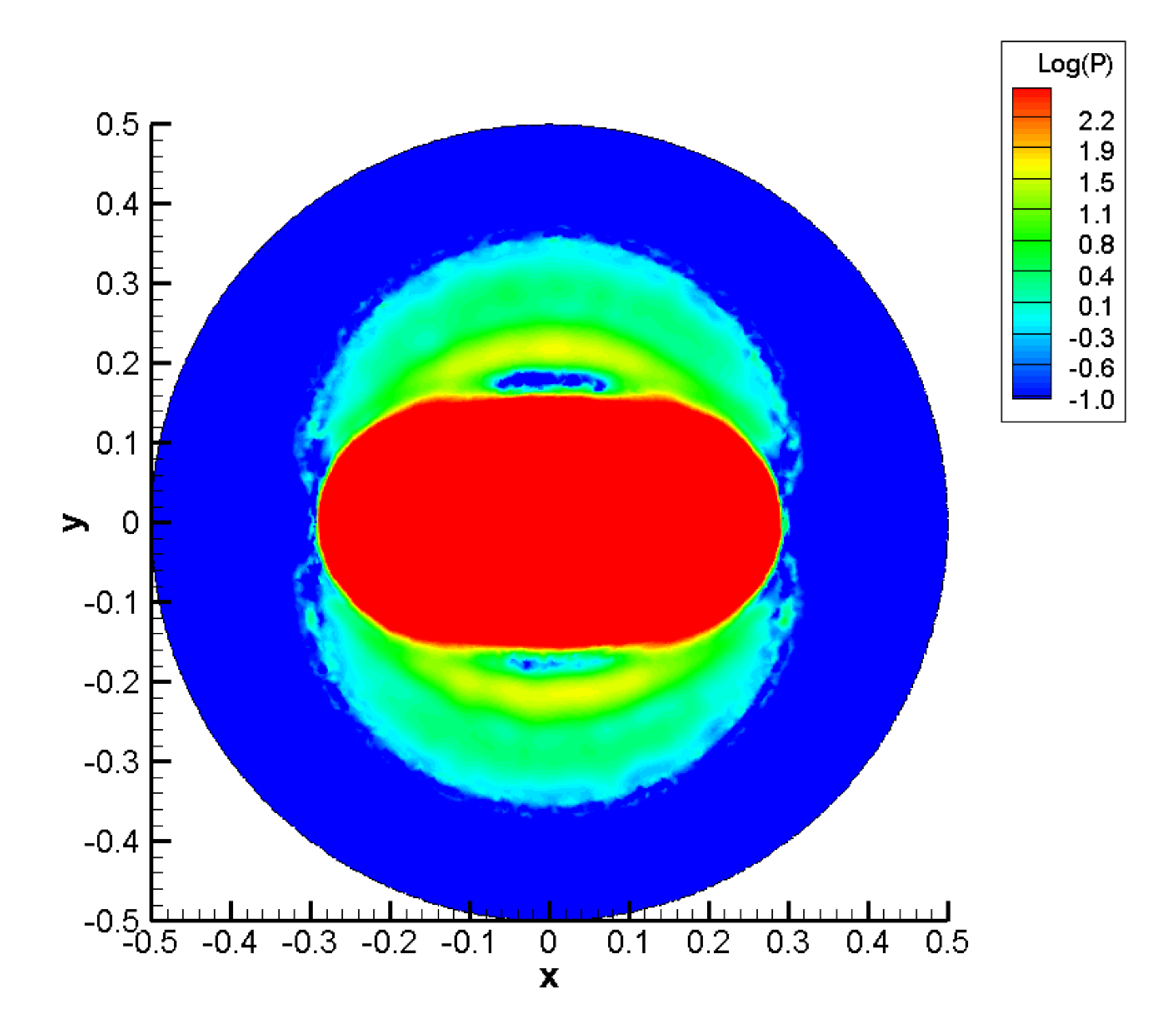} \\ 
\includegraphics[width=0.47\textwidth]{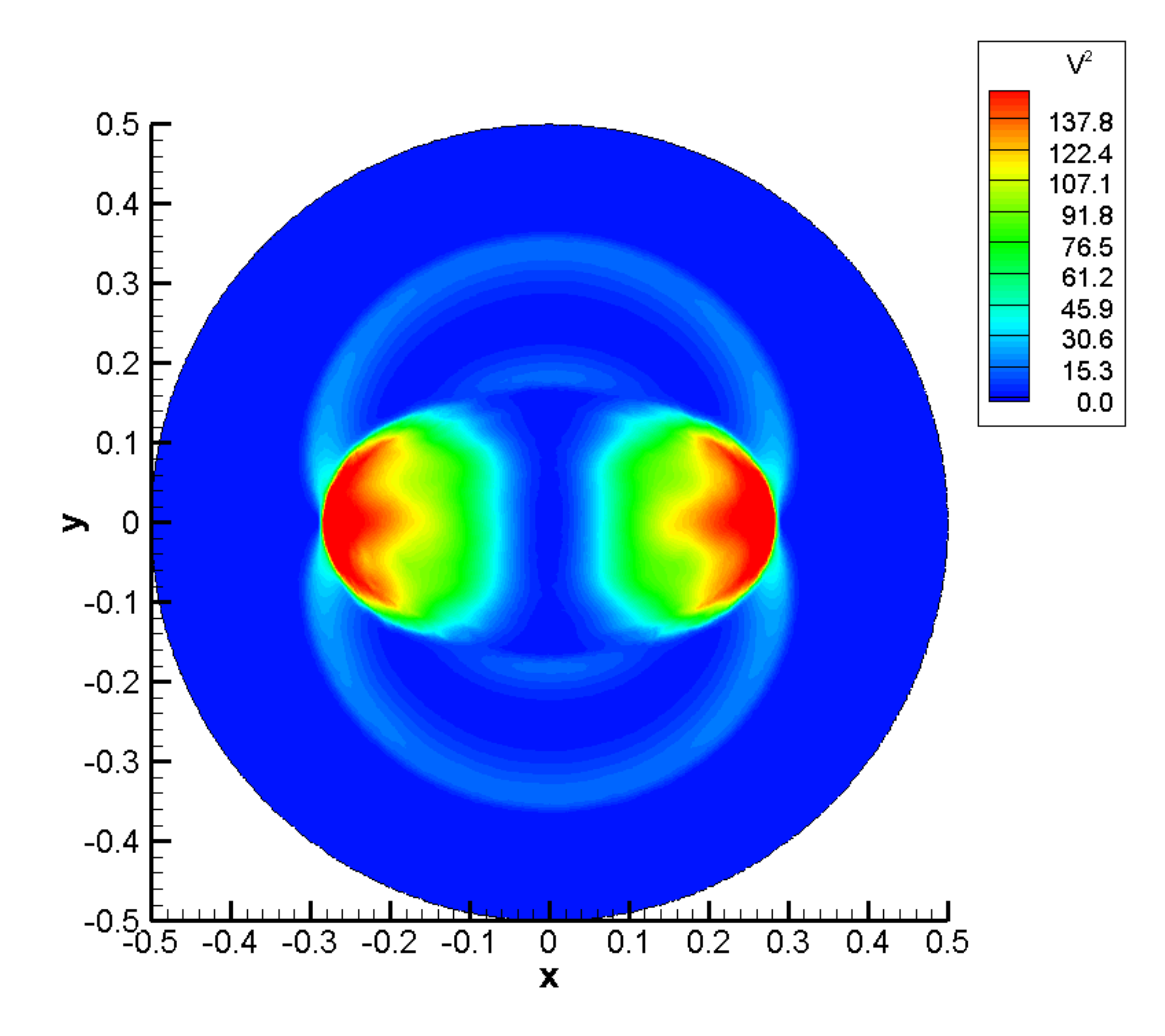}  &           
\includegraphics[width=0.47\textwidth]{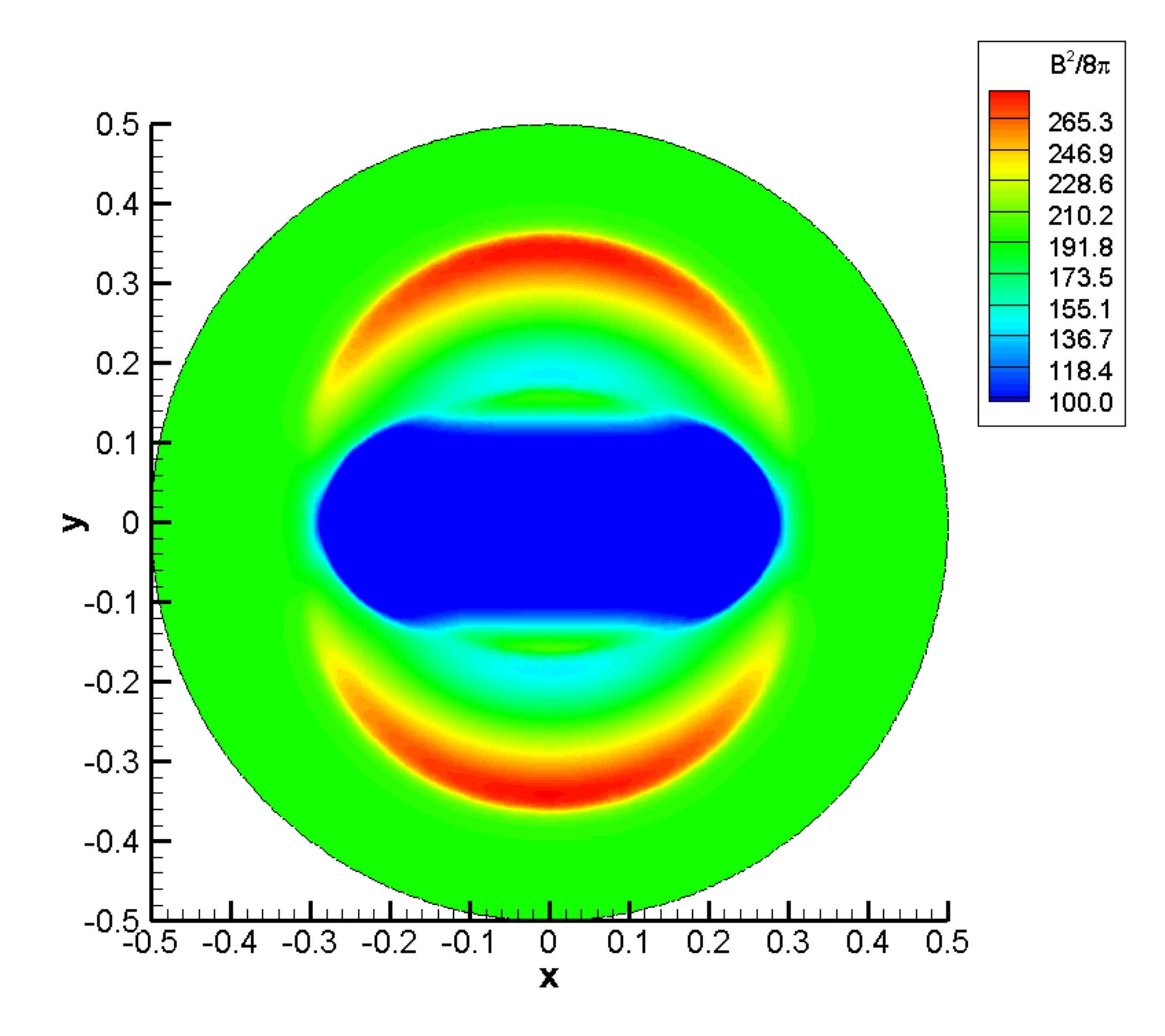} \\  
\end{tabular} 
\caption{Third order numerical results for the Blast problem at time $t=0.01$. Top: logarithm (base 10) of the density and logarithm (base 10) of the pressure. Bottom: magnitude of the velocity field and the magnetic field.} 
\label{fig.Blast}
\end{center}
\end{figure}

The blast wave problem constitutes a benchmark in magnetohydrodynamics. A strong circular fast magnetosonic shock wave is traveling from the center to the boundaries of the initial computational domain $\Omega(0)$, which is a sphere of radius $R_0=0.5$. The frontier delimited by radius $R=0.1$ splits the domain into two parts, hence defining an \textit{inner} state $\U_i$ and an \textit{outer} state $\U_o$, that are initially assigned in terms of primitive variables $\U=(\rho,u,v,w,p,B_x,B_y,B_z,\psi)$ as 
\begin{equation}
  \U(\x,0) = \left\{ \begin{array}{ccc} \U_i = \left( 1.0, 0.0, 0.0, 0.1, 70, 0.0, 0.0 \right) & \textnormal{ if } & r \leq R, \\ 
                                        \U_o = \left( 1.0, 0.0, 0.0, 1000, 70, 0.0, 0.0 \right) & \textnormal{ if } & r > R,        
                      \end{array}  \right. 
\end{equation}
where $r=\sqrt{x^2+y^2+z^2}$. We use the same mesh adopted for the MHD rotor problem described in the previous section and we set transmissive boundary conditions at the external boundary. The final time of the computation is $t_f=0.01$ and the ratio of specific heats is taken to be $\gamma=1.4$. The numerical results obtained with the third order version of the scheme using the Rusanov-type flux 
\eqref{eqn.rusanov} are depicted in Figure \ref{fig.Blast}, where the logarithm of density and pressure, as well as the magnitude of both the velocity and the magnetic field are reported. 
The solution is in qualitative agreement with the results shown in \cite{LagrangeMHD,LagrangeMDRS}, where the two-dimensional version of our Lagrangian WENO algorithm has been used to run this 
test case. The tetrahedral mesh at the final time $t=0.01$ is depicted in Figure \ref{fig.Blast.grid}. 

%
\subsection{The Baer-Nunziato model of compressible two-phase flows}
We consider here the Baer-Nunziato model for compressible two-phase flows, which has been developed by Baer and Nunziato in \cite{BaerNunziato1986} for describing detonation waves in solid-gas combustion processes. The first phase is normally addressed as the \textit{solid} phase, while the second one as the \textit{gas} phase and in this paper we will use the subscripts $1$ and $2$ to define them. We will use equivalently also the subscripts $s$ and $g$ to denote the solid and the gas phase. Let $k=1,2$ be the phase number and $\phi_k$ be the volume fraction of phase $k$ with the condition 
$\phi_1 + \phi_2 = 1$, while $\rho_k$ and $p_k$ represent the corresponding density and pressure, respectively. Let furthermore the velocity vector of each phase be addressed with 
$\textbf{u}_k=(u_k,v_k,w_k)$. The full seven-equation Baer-Nunziato model with relaxation source terms results in a \textit{non-conservative system} of nonlinear hyperbolic PDE that can be written as 
\begin{equation}\left.\begin{array}{l}
\label{ec.BN}


	\frac{\partial}{\partial t}\left(\ar\right)+\nabla\cdot\left(\ar\ub\right)=0,
	
	\\


\frac{\partial}{\partial t}\left(\ar\ub\right)
+\nabla\cdot\left(\phi_1\rho_1\textbf{u}_\textbf{1} \textbf{u}_\textbf{1}\right)+\nabla\phi_1p_1 
= p_I\nabla\phi_1 - \lambda\left(\textbf{u}_\textbf{1}-\textbf{u}_\textbf{2} \right),

\\


\frac{\partial}{\partial t}\left(\phi_1\rho_1E_1\right)
+\nabla\cdot\left(\left(\phi_1\rho_1E_1+\phi_1p_1\right)\textbf{u}_\textbf{1}\right) = 
-p_I\partial_t\phi_1 - \lambda \, \mathbf{u_I} \cdot \left(\mathbf{u_1}-\mathbf{u_2}\right),

\\


\frac{\partial}{\partial t}\left(\phi_2\rho_2\right)+\nabla\cdot\left(\phi_2\rho_2\textbf{u}_\textbf{2}\right)=0,

\\ 


\frac{\partial}{\partial t}\left(\phi_2\rho_2\textbf{u}_\textbf{2}\right)
+\nabla\cdot\left(\phi_2\rho_2\textbf{u}_\textbf{2} \textbf{u}_\textbf{2}\right)+\nabla\phi_2p_2=p_I\nabla\phi_2 
- \lambda \, \left(\textbf{u}_\textbf{2}-\textbf{u}_\textbf{1}\right),

\\


\frac{\partial}{\partial t}\left(\phi_2\rho_2E_2\right)
+\nabla\cdot\left(\left(\phi_2\rho_2E_2+\phi_2p_2\right)\textbf{u}_\textbf{2}\right)=p_I\partial_t\phi_1 
- \lambda \, \mathbf{u_I} \cdot \left(\mathbf{u_2}-\mathbf{u_1}\right),

\\


\frac{\partial}{\partial t}\phi_1+\textbf{u}_\textbf{I}\nabla\phi_1 = \mu (p_1 - p_2),
\end{array}\right\}
\end{equation} 
where only strongly simplified interphase drag and pressure relaxation source terms are considered. Further details on the choice and the formulation of such terms can be found in \cite{Kapila2001}. 
The so-called stiffened gas equation of state is then used for each of the two phases to close the system: 
\begin{equation}
\label{eqn.eosBN} 
   e_k = \frac{p_k + \gamma_k \pi_k}{\rho_k (\gamma_k -1 )}.
\end{equation}
The specific total energy of each phase is $E_k = e_k + \halb \mathbf{u_k}^2$ with $e_k$ denoting the corresponding internal energy, while in the system \eqref{ec.BN} $\mu$ is a parameter which characterizes pressure relaxation and $\lambda$ is related to the friction between the phases. According to \cite{BaerNunziato1986,Kapila2001} the velocity at the interface $I$ is taken to be the 
\textit{solid} velocity, while for the interface pressure we choose the \textit{gas} pressure, hence 
\begin{equation}
\mathbf{u_I} = \mathbf{u_1} \qquad p_I = p_2.
\label{up_interface}
\end{equation} 
Other choices are possible,  see \cite{SaurelAbgrall,SaurelGavrilyukRenaud} for a detailed discussion.

The resolution of material interfaces, which are given by jumps in the volume fraction $\phi_k$, is a challenging task for the numerical methods applied to the Baer-Nunziato model \eqref{ec.BN}. In the following we present numerical results for some well-known test cases in order to validate our algorithm. In all cases the mesh velocity is chosen to be the interface velocity in order to resolve 
the material contact properly, hence we set $\mathbf{V}=\mathbf{u_I}$. 

\subsubsection{Riemann problems} 
\label{sec.BN-RP}
In this section we apply the three-dimensional high order ALE ADER-WENO finite volume schemes presented in this article to the seven-equation Baer-Nunziato model \eqref{ec.BN} by solving a set of 
Riemann problems. Those test cases have been taken from \cite{DeledicquePapalexandris,USFORCE2} and the initial condition for each of the Riemann problem is listed in Table \ref{tab.rpbn.ic}. 
The exact solution of the Riemann problem for the BN model has been derived in \cite{AndrianovWarnecke,Schwendeman,DeledicquePapalexandris} and is used here for comparison. 

\begin{figure}[!htbp]
	\begin{center}
	\begin{tabular}{cc} 
	\includegraphics[width=0.47\textwidth]{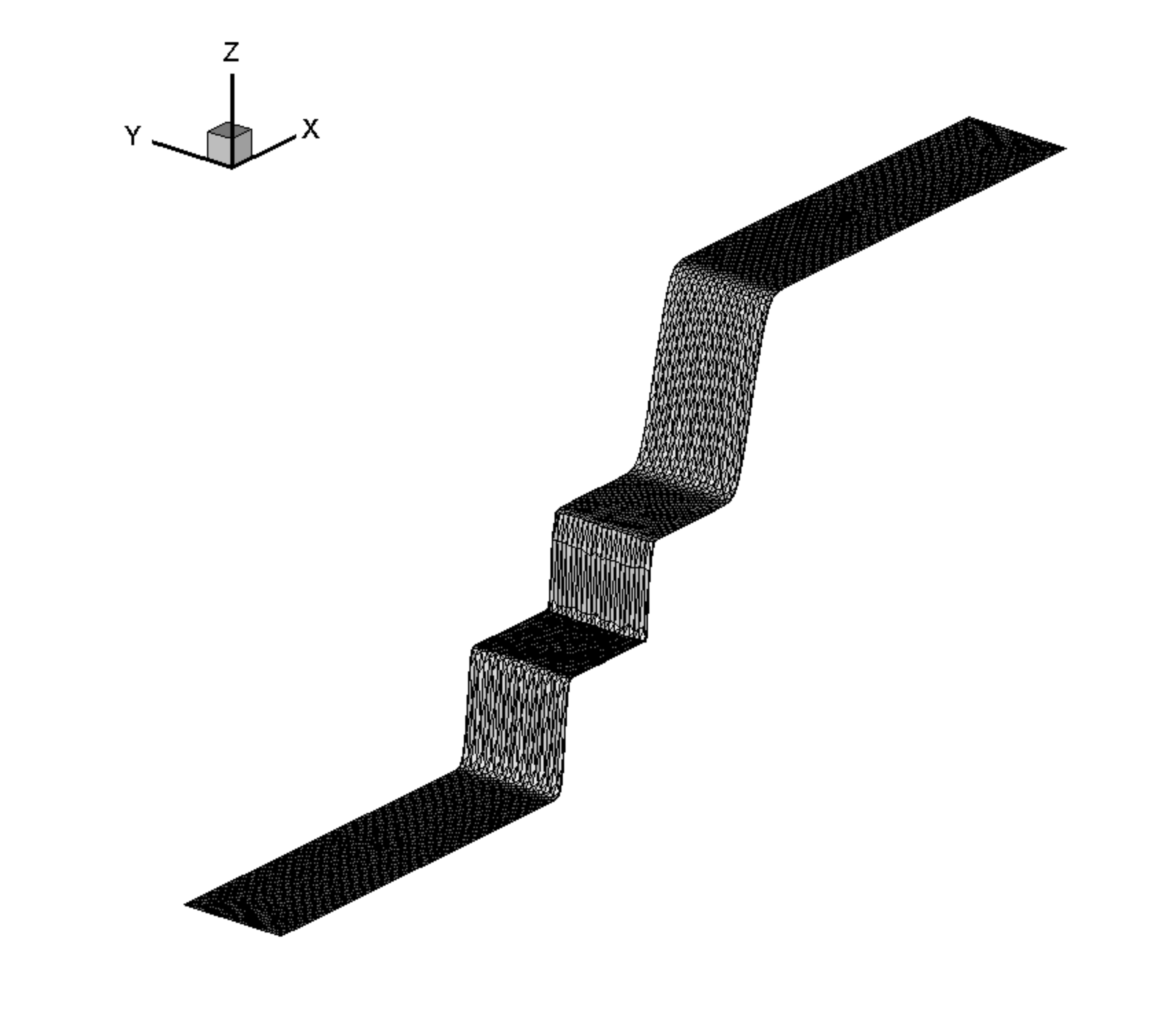}  &           
	\includegraphics[width=0.47\textwidth]{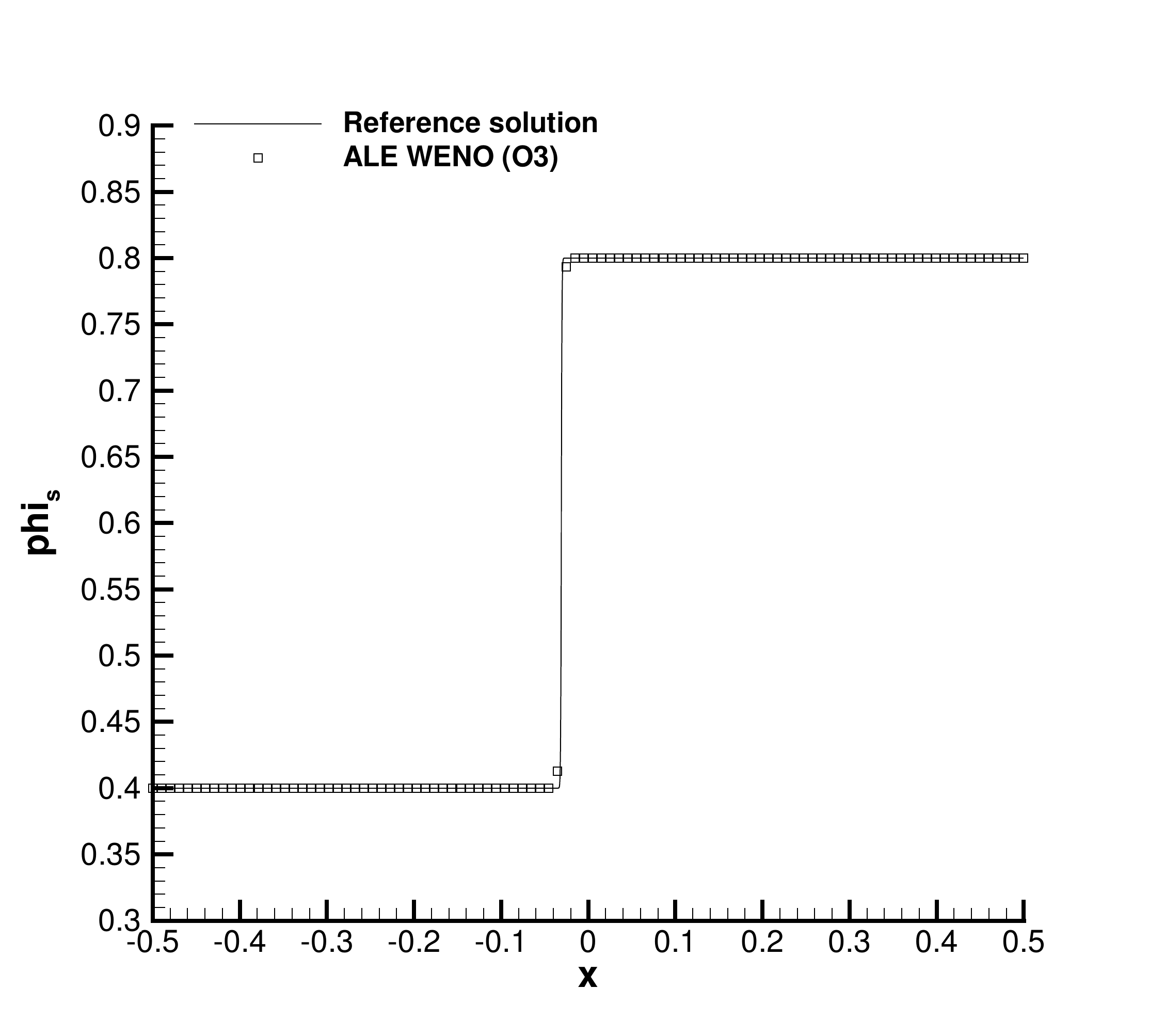} \\
	\includegraphics[width=0.47\textwidth]{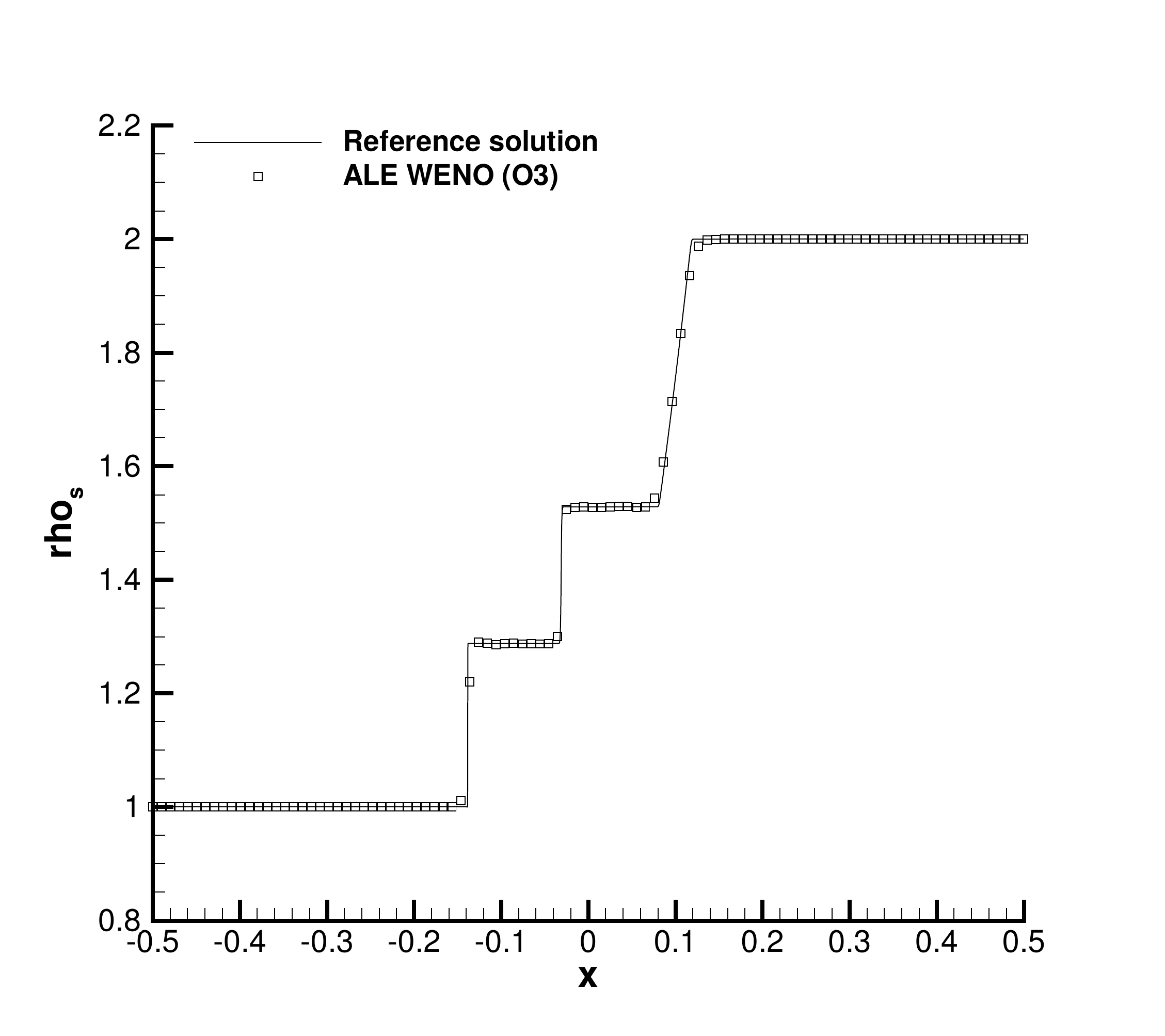}  &           
	\includegraphics[width=0.47\textwidth]{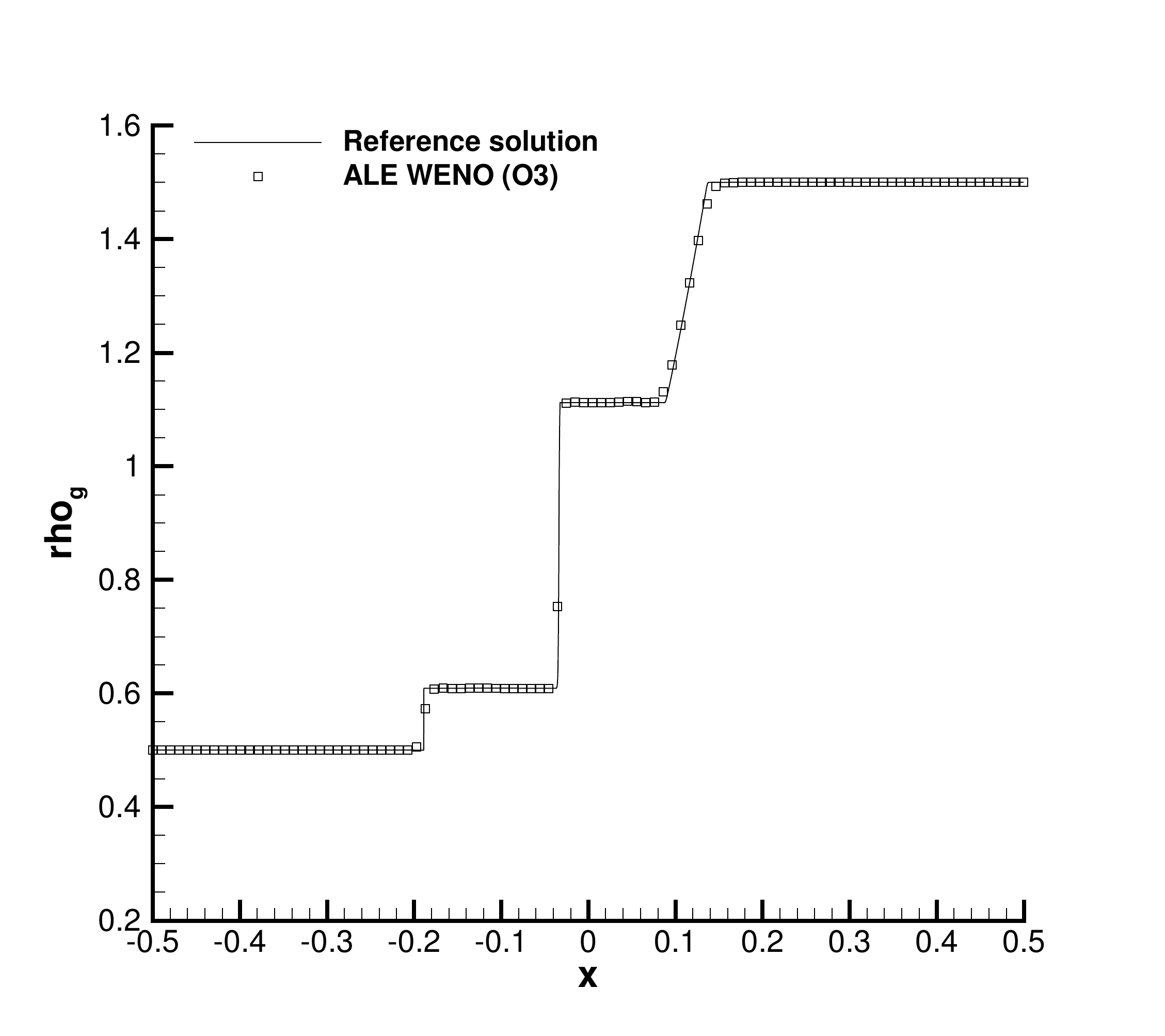} \\
	\includegraphics[width=0.47\textwidth]{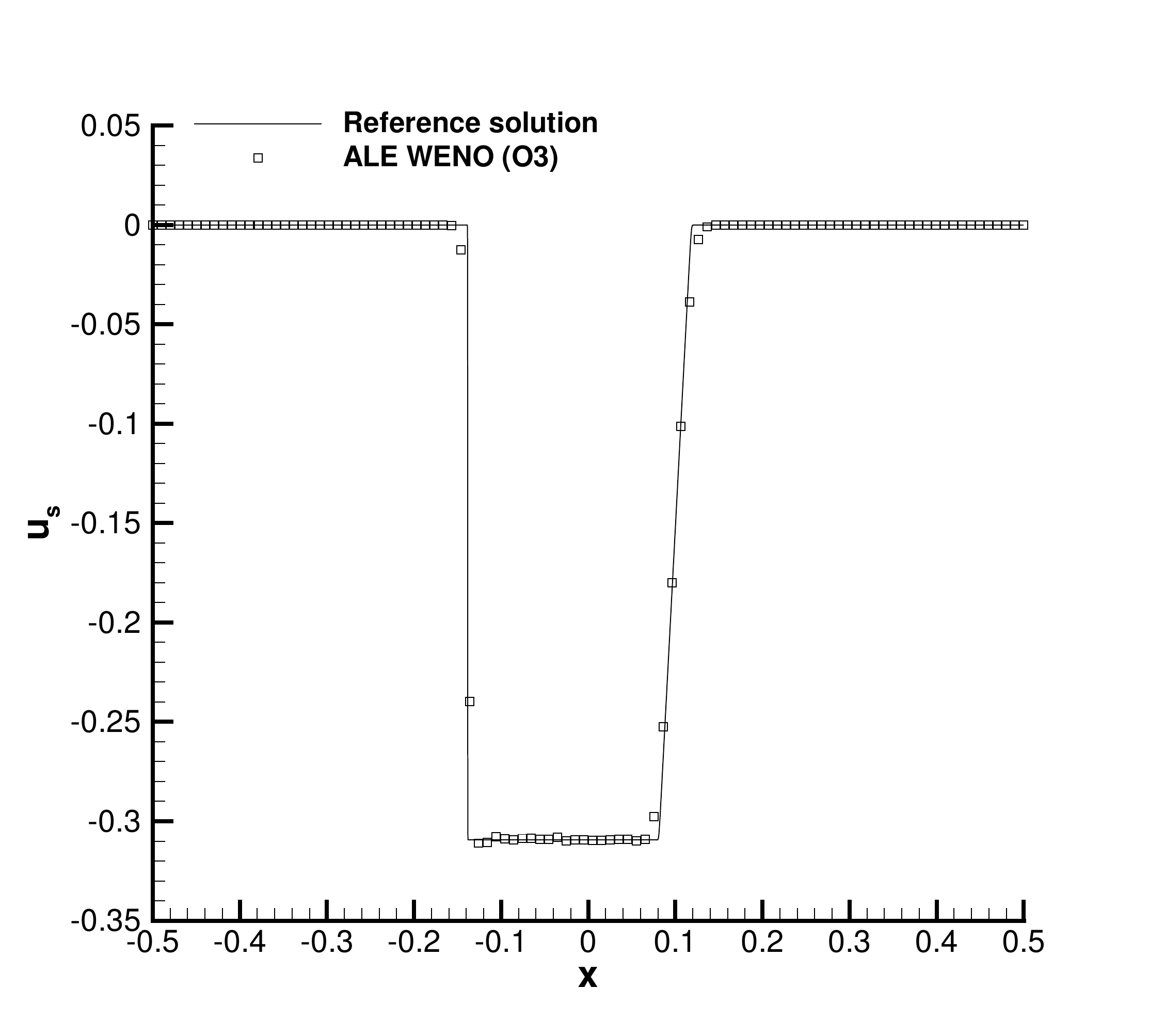}  &           
	\includegraphics[width=0.47\textwidth]{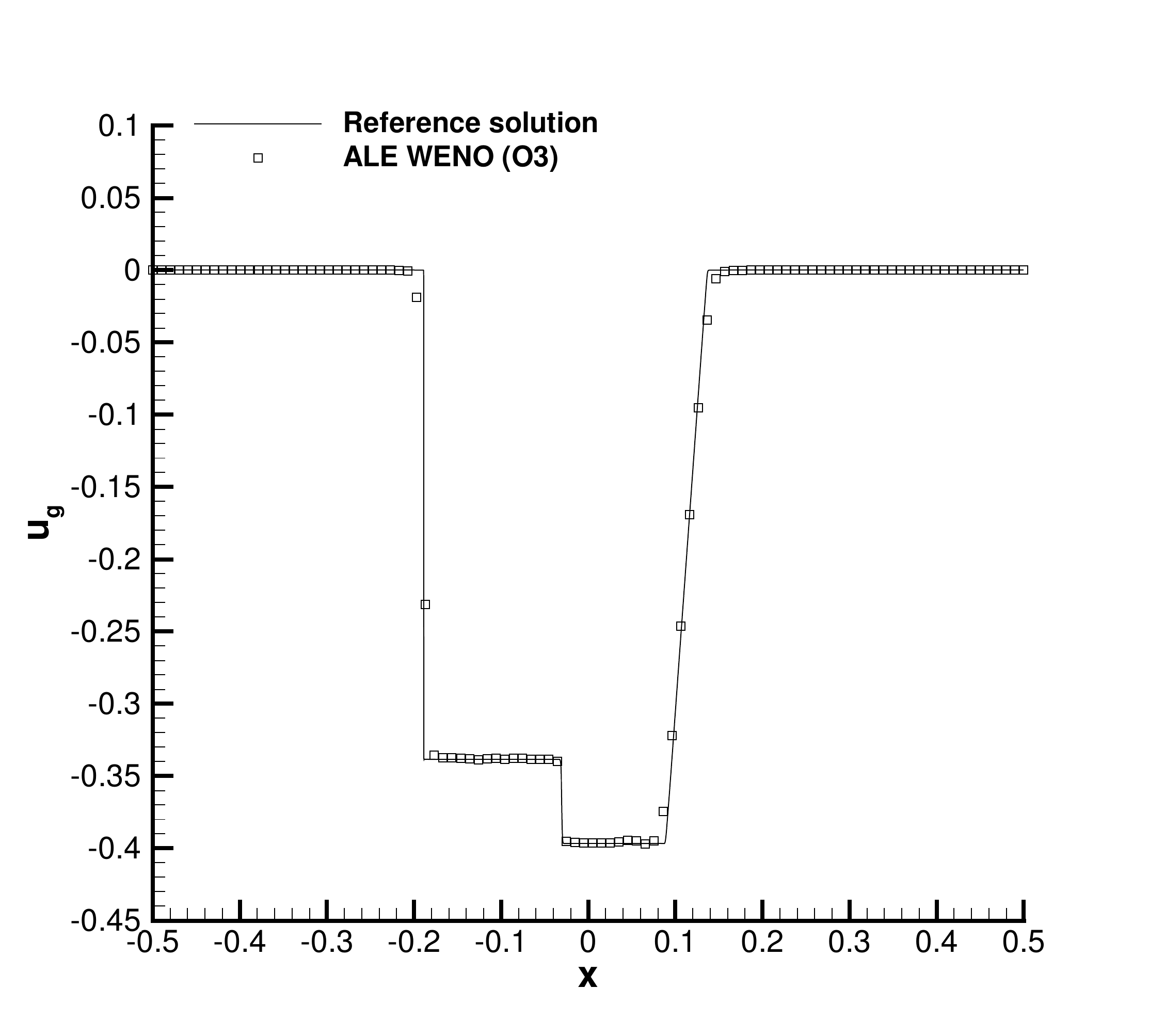} \\
	\includegraphics[width=0.47\textwidth]{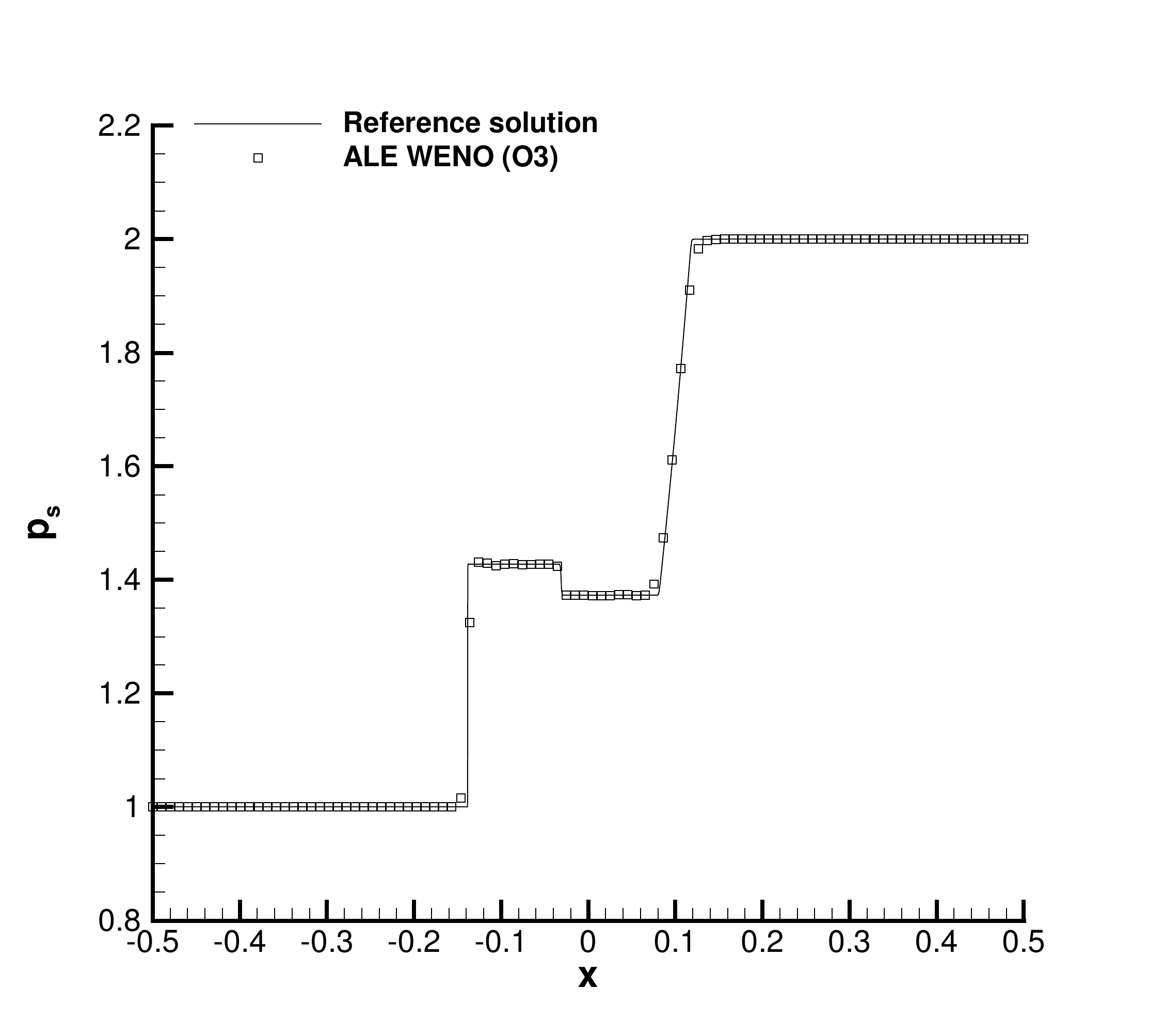}  &           
	\includegraphics[width=0.47\textwidth]{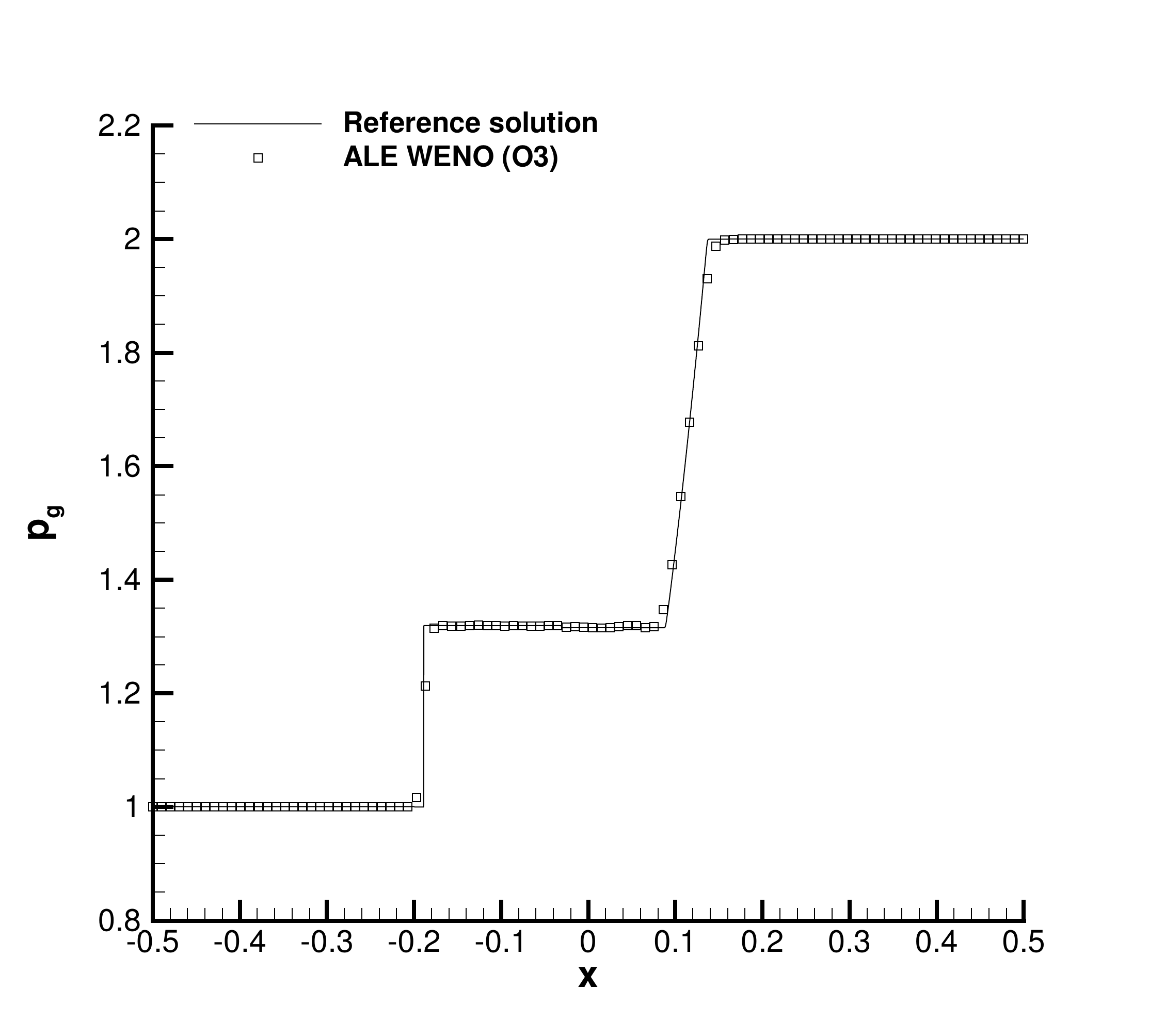} \\
	\end{tabular}
	\caption{Third order numerical results for Riemann problem RP1 of the seven-equation Baer-Nunziato model at time $t=0.1$ and comparison with the reference solution.}
	\label{fig.RP1}
	\end{center}
\end{figure}

\begin{figure}[!htbp]
	\begin{center}
	\begin{tabular}{cc} 
	\includegraphics[width=0.47\textwidth]{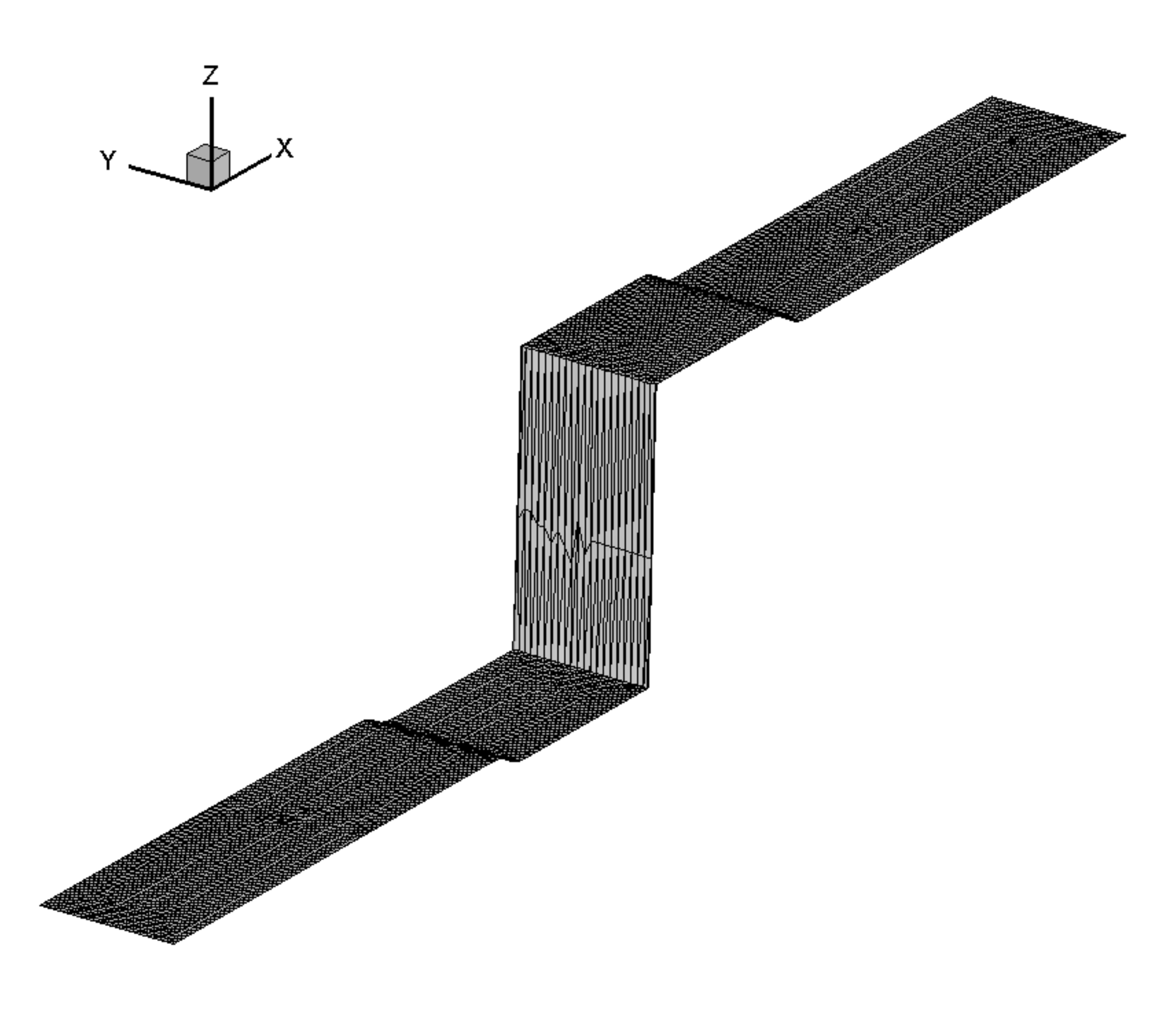}  &           
	\includegraphics[width=0.47\textwidth]{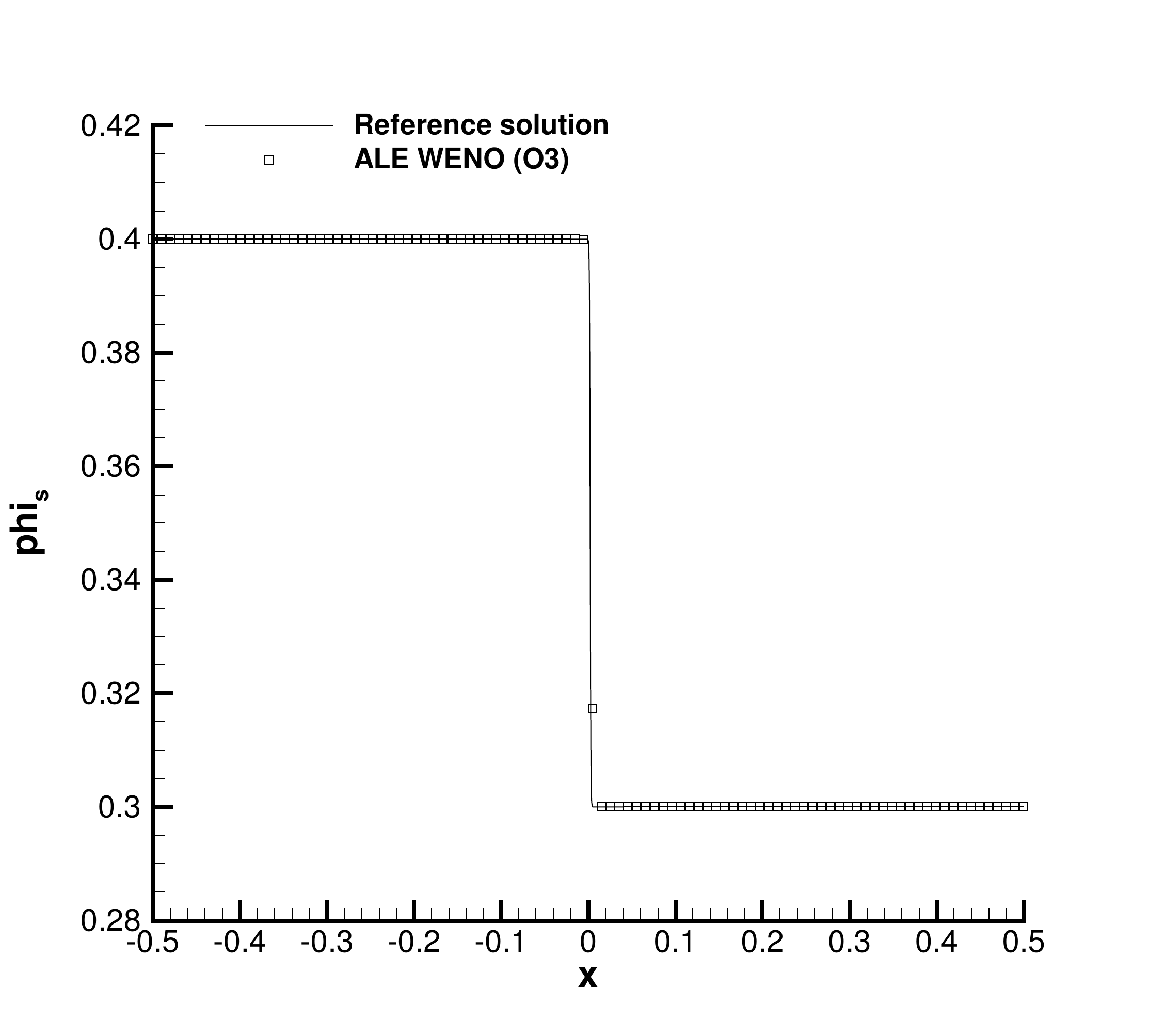} \\
	\includegraphics[width=0.47\textwidth]{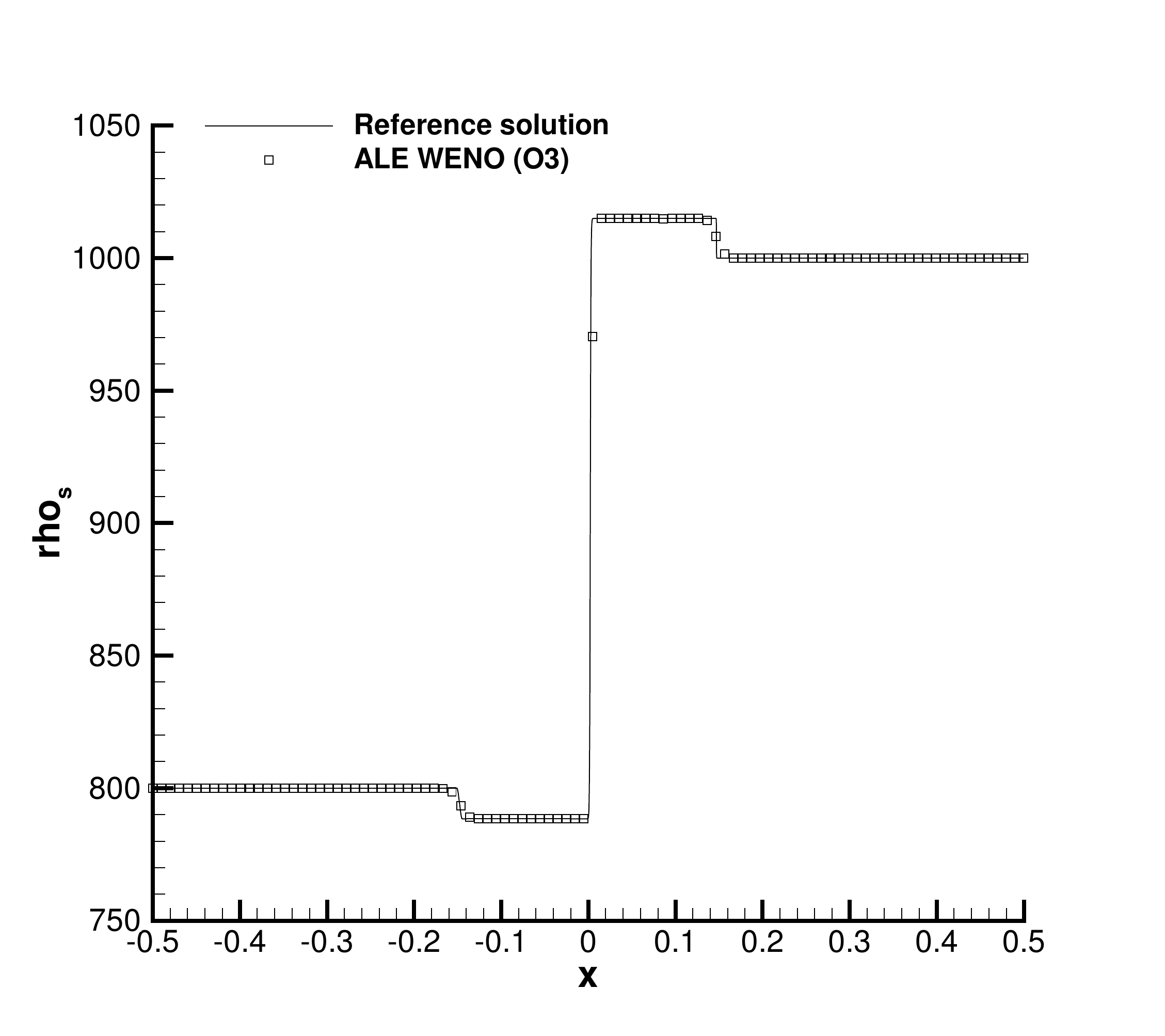}  &           
	\includegraphics[width=0.47\textwidth]{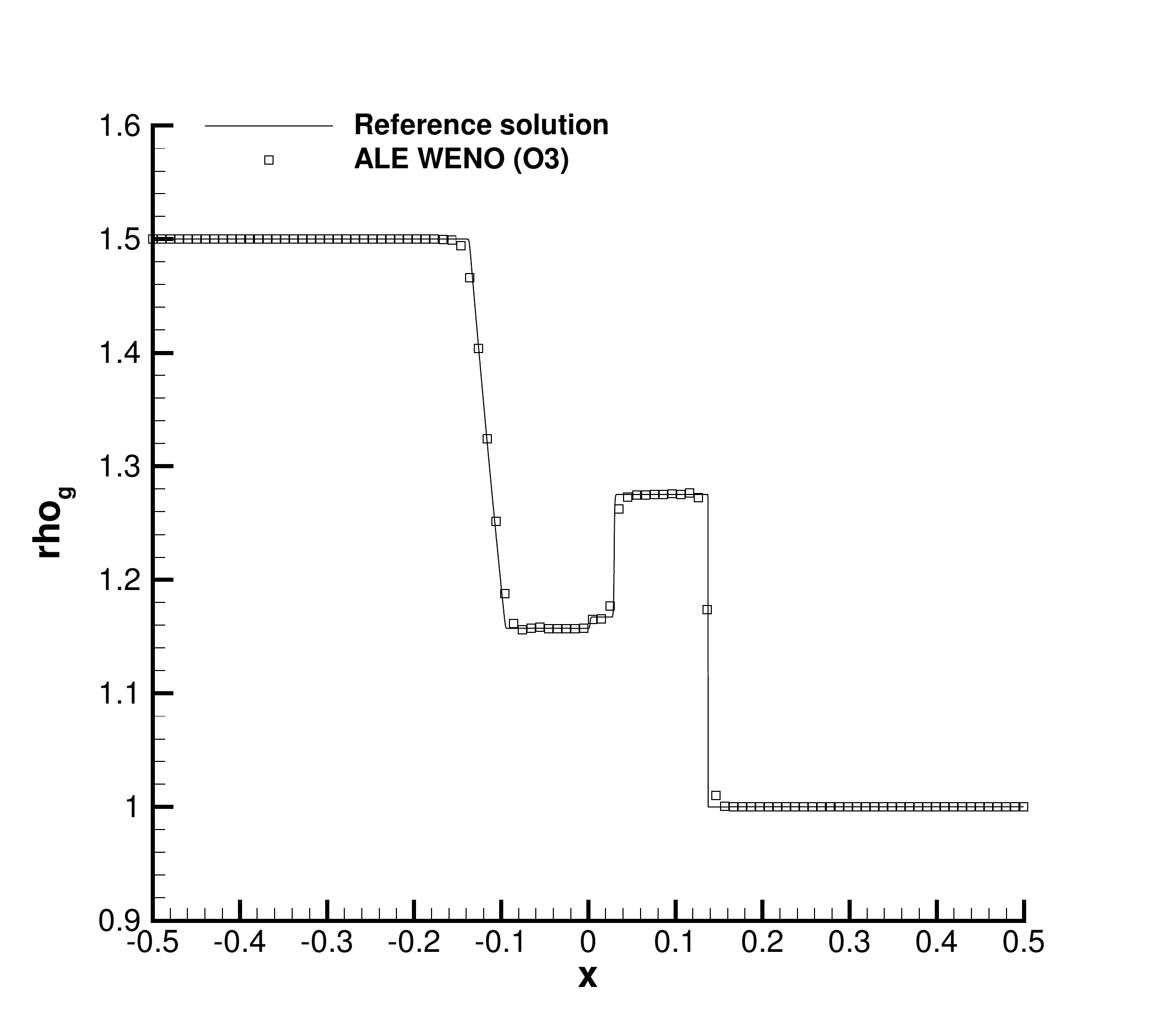} \\
	\includegraphics[width=0.47\textwidth]{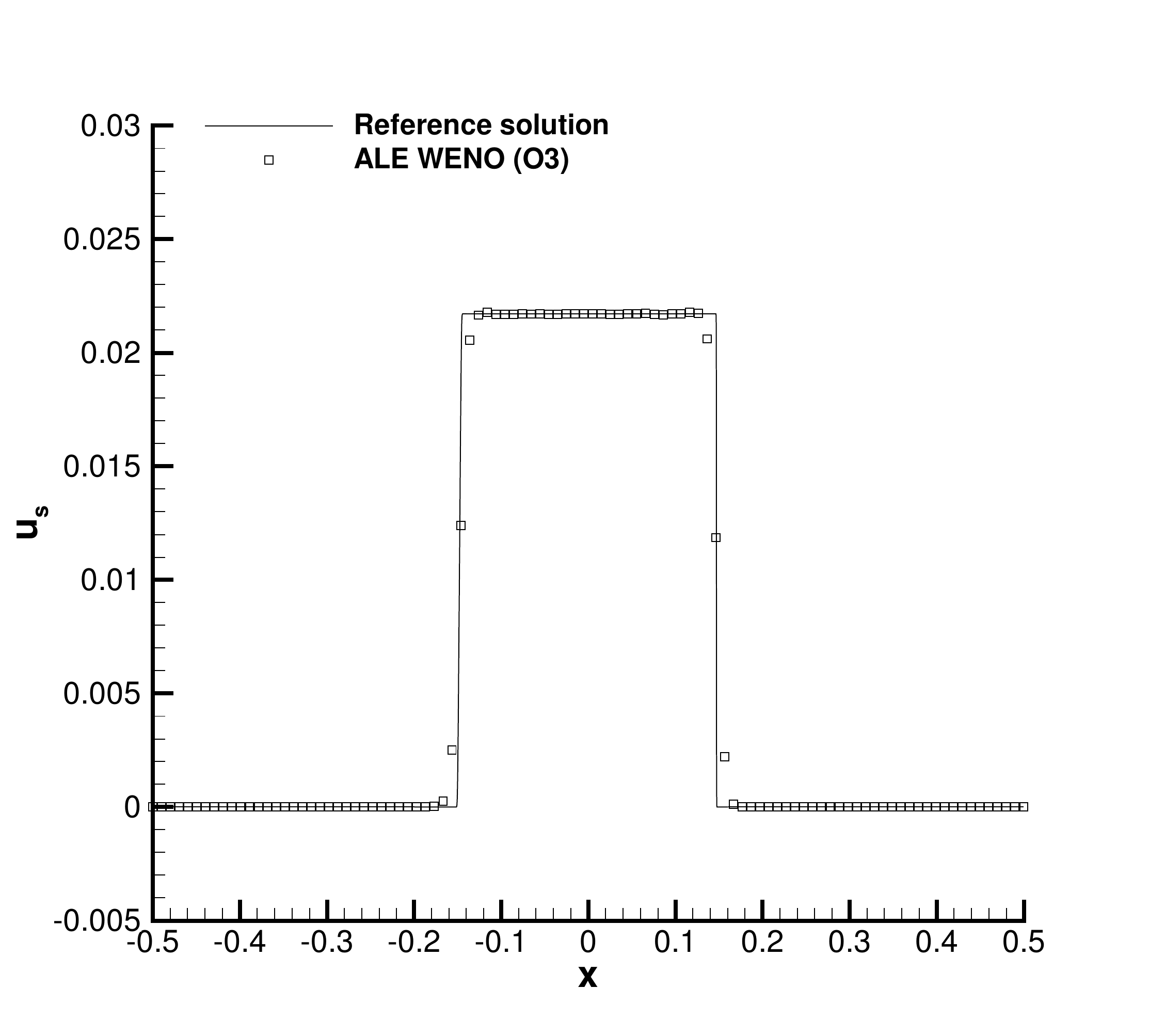}  &           
	\includegraphics[width=0.47\textwidth]{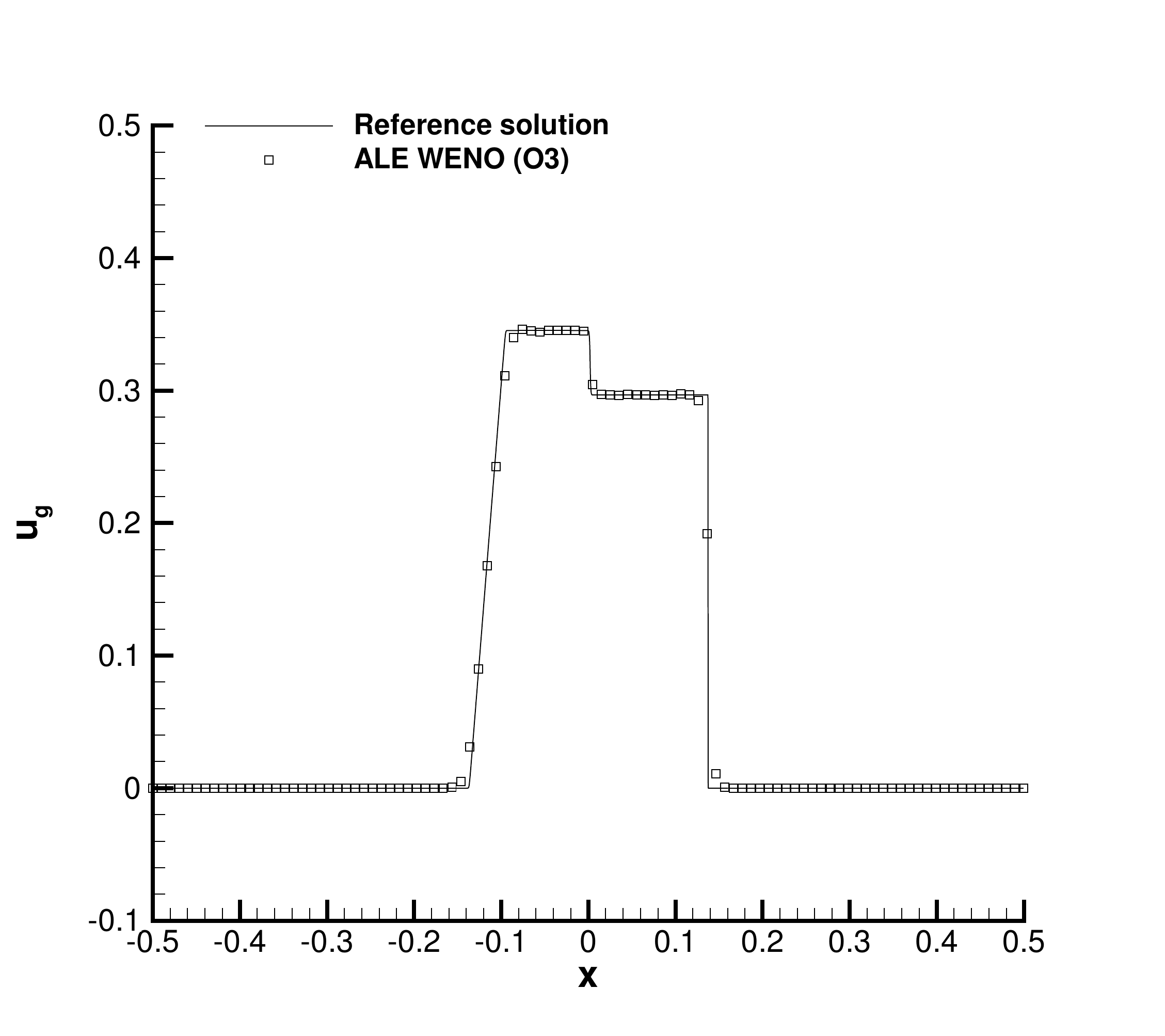} \\
	\includegraphics[width=0.47\textwidth]{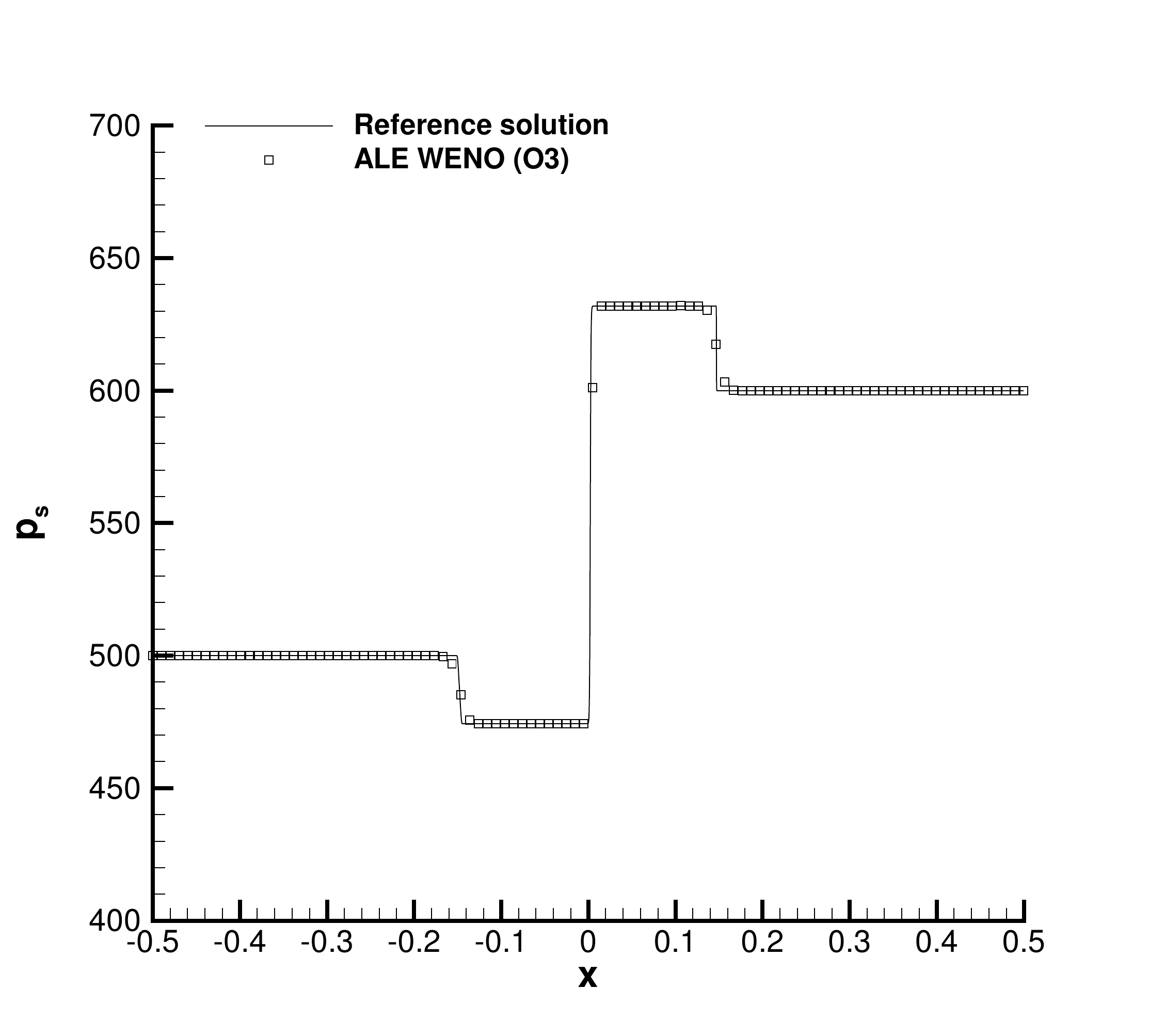}  &           
	\includegraphics[width=0.47\textwidth]{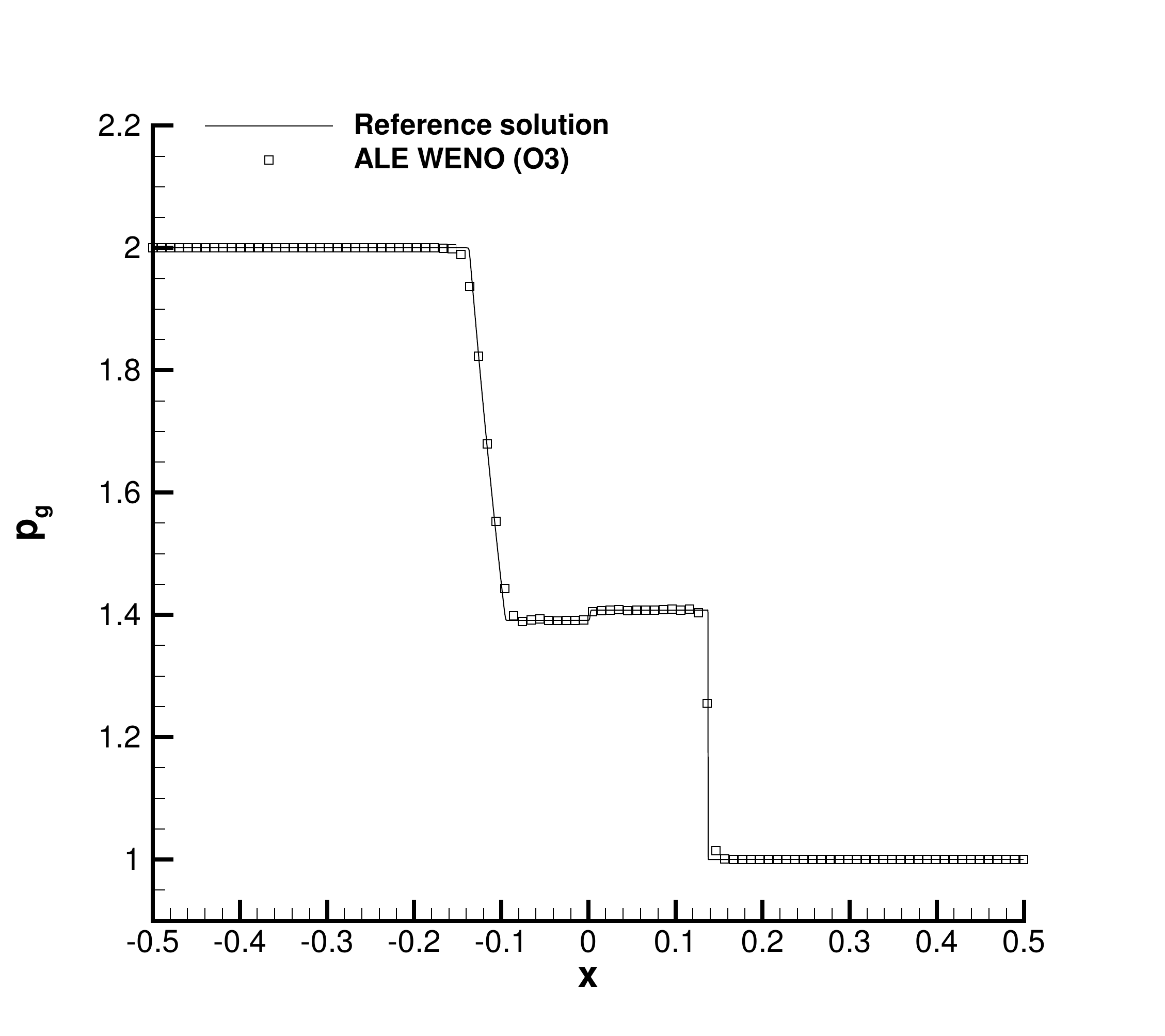} \\
	\end{tabular}
	\caption{Third order numerical results for Riemann problem RP2 of the seven-equation Baer-Nunziato model at time $t=0.1$ and comparison with the reference solution.}
	\label{fig.RP2}
	\end{center}
\end{figure}

\begin{figure}[!htbp]
	\begin{center}
	\begin{tabular}{cc} 
	\includegraphics[width=0.47\textwidth]{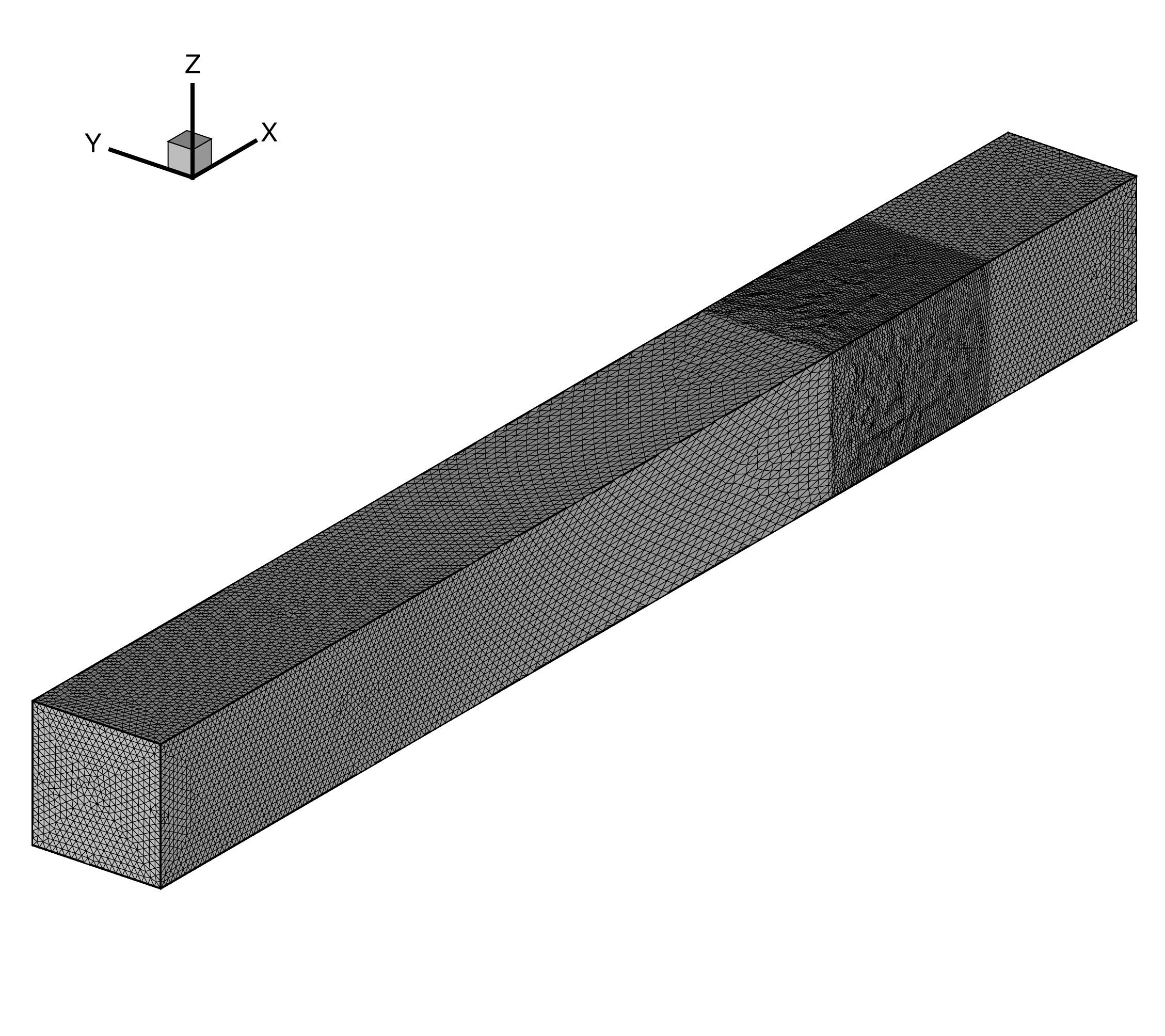}  &           
	\includegraphics[width=0.47\textwidth]{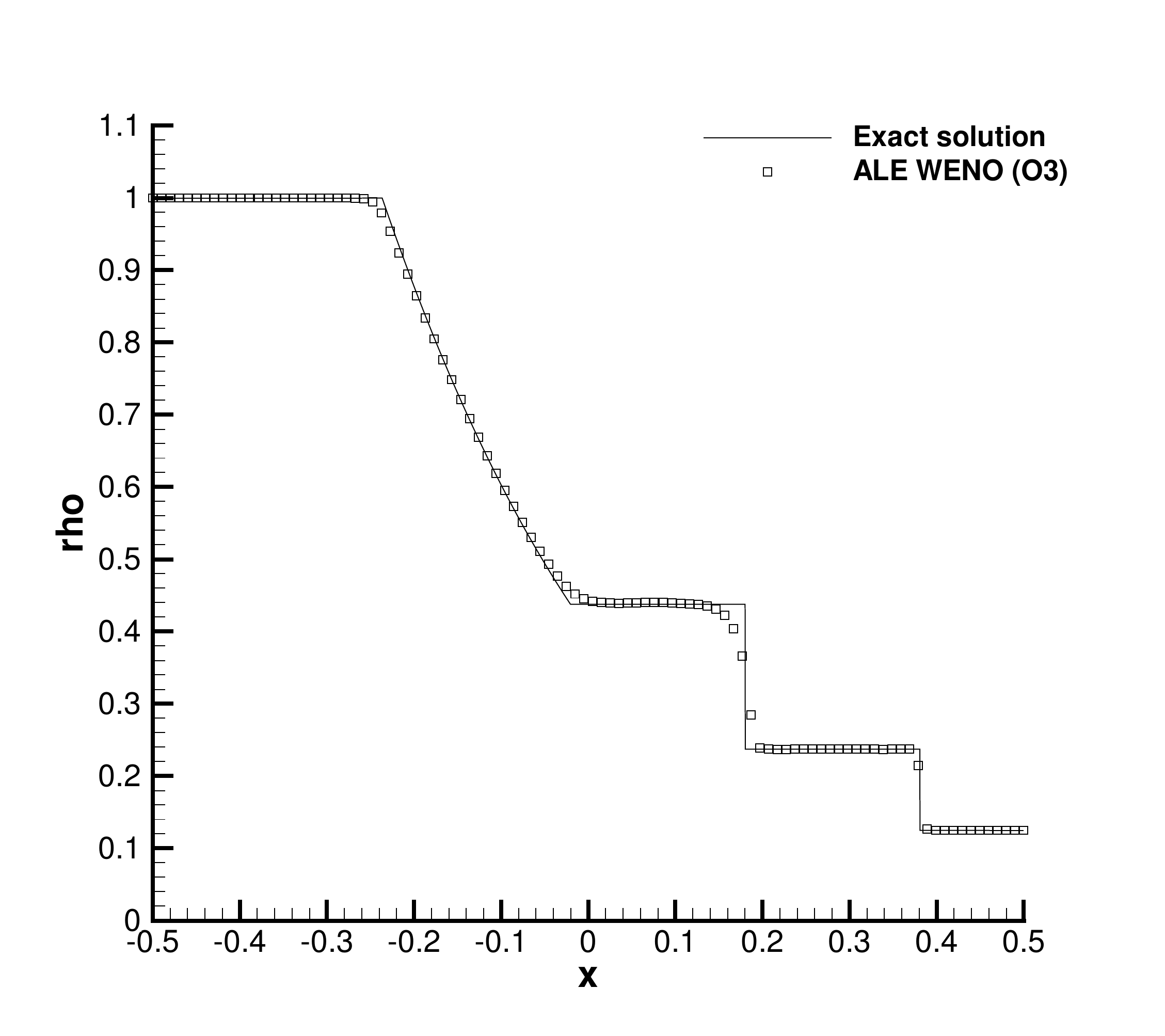} \\
	\includegraphics[width=0.47\textwidth]{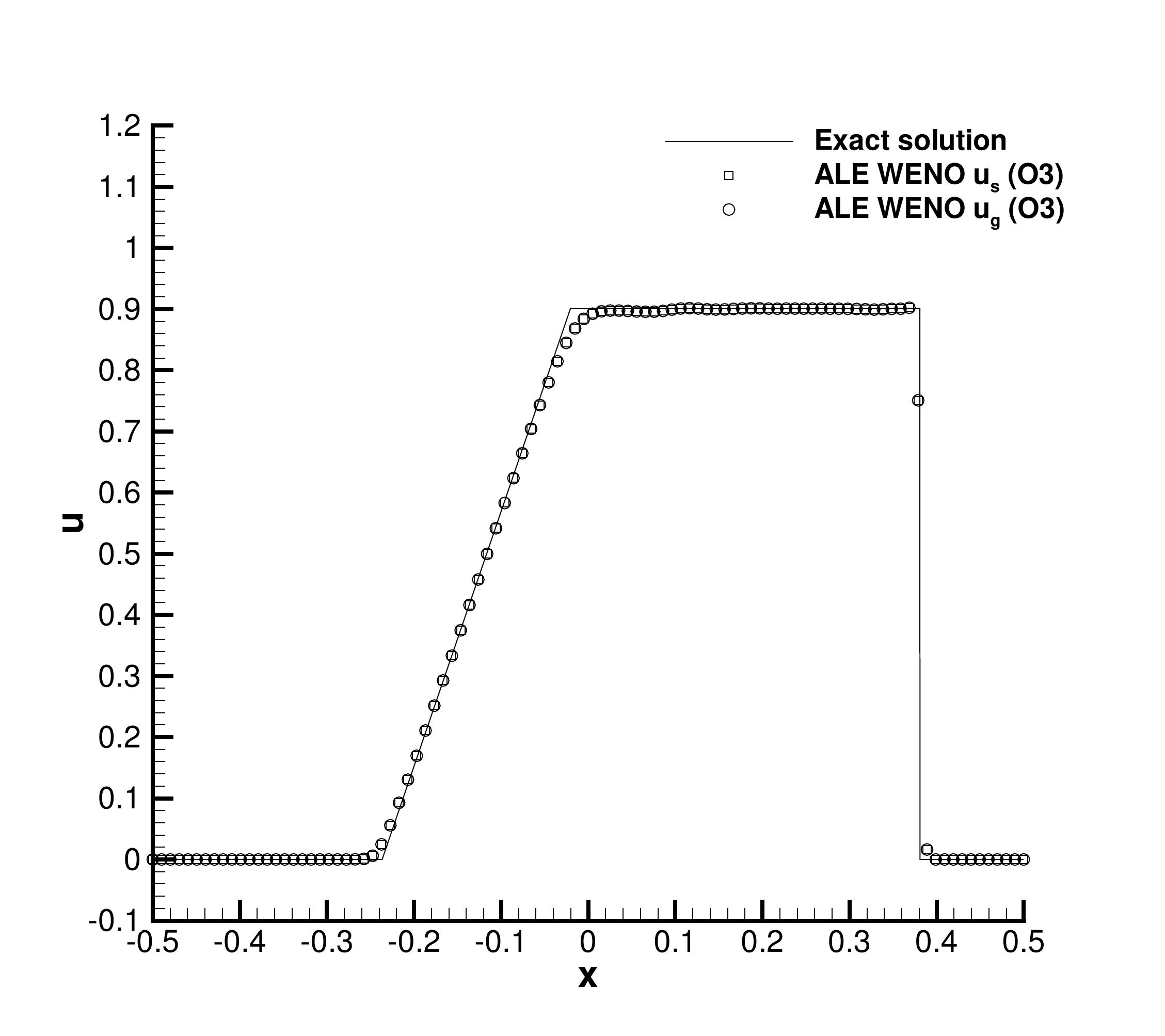}  &           
	\includegraphics[width=0.47\textwidth]{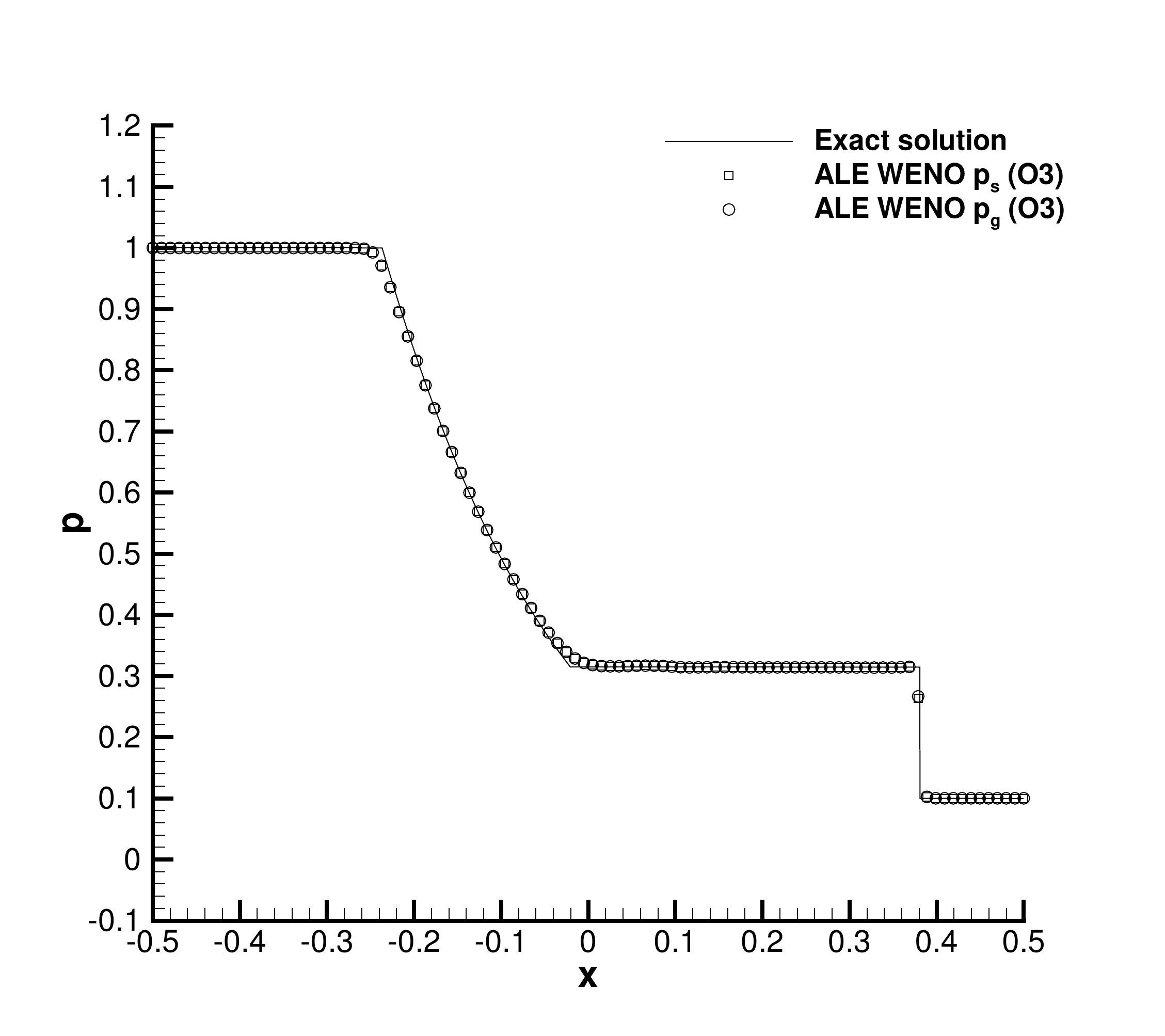} \\
	\end{tabular}
	\caption{Third order numerical results for Riemann problem RP3 of the seven-equation Baer-Nunziato model with drag and pressure relaxation ($\lambda=10^3, \mu=10^2$) 
	at time $t=0.2$ and comparison with the reference solution.} 
	\label{fig.RP3}
	\end{center}
\end{figure}

The initial computational domain is defined by $\Omega(0)=[-0.5;0.5]\times[-0.05;0.05]\times[-0.05;0.05]$ and the initial discontinuity between the left state $\Q_L$ and the right state $\Q_R$ is located at $x_0=0$. The domain is discretized using a characteristic mesh size of $h=1/200$, hence obtaining a total number of tetrahedra of $N_E=563357$. Periodic boundary conditions have been imposed in $y$ and $z$ directions, while we use transmissive boundaries along the $x$ direction.
Friction and pressure relaxation are neglected in the first two Riemann problems RP1 and RP2, while for RP3 we use a moderately stiff interphase drag $\lambda=10^3$ and pressure relaxation $\mu=10^2$. RP3 involves two almost pure ideal gases that differ in their value of $\gamma$. As done in \cite{USFORCE2,LagrangeNC} the exact solution for RP3 is computed using the exact Riemann solver for the Euler equations of compressible gas dynamics \cite{ToroBook} with two different values of $\gamma$ on the left and on the right of the contact discontinuity, respectively. The numerical results have been  obtained using the third order version of our ALE finite volume schemes together with the Osher-type method \eqref{eqn.osher} for RP1 and RP2, while the more robust Rusanov-type method 
\eqref{eqn.rusanov} has been adopted for RP3. Figures \ref{fig.RP1}-\ref{fig.RP3} show a comparison between the reference solution and a one-dimensional 
cut through the reconstructed numerical solution $\w_h$. For RP3 we show the \textit{mixture} density $\rho=\phi_s\rho_s + (1-\phi_s)\rho_g$ in Figure \ref{fig.RP3}. In all cases one can note a 
very good agreement between numerical solution and reference solution. The material contact is well resolved in all cases. 

\begin{table}
\caption{Initial condition for the left state (L) and the right state (R) for the Riemann problems solved with the 
Baer-Nunziato model. Values for $\gamma_k$, $\pi_k$ and the final time $t_f$ are also given.}
\renewcommand{\arraystretch}{1.0}
\begin{center}
\begin{tabular}{ccccccccc}
\hline
   & $\rho_s$ & $u_s$  & $p_s$ & $\rho_g$ & $u_g$ & $p_g$ & $\phi_s$ & $t_f$  \\
\hline 
\multicolumn{1}{l}{\textbf{RP1 \cite{DeledicquePapalexandris}:} } & 
\multicolumn{8}{c}{ $\gamma_s = 1.4, \quad \pi_s = 0, \quad \gamma_g = 1.4, \quad \pi_g = 0, \quad \lambda=\mu=0$}  \\
\hline 
L & 1.0    & 0.0   & 1.0  & 0.5 & 0.0   &  1.0 & 0.4 & 0.10 \\
R & 2.0    & 0.0   & 2.0  & 1.5 & 0.0   &  2.0 & 0.8 &      \\
\hline 
\multicolumn{1}{l}{\textbf{RP2 \cite{DeledicquePapalexandris}:}} & 
\multicolumn{8}{c}{ $\gamma_s = 3.0, \quad \pi_s = 100, \quad \gamma_g = 1.4, \quad \pi_g = 0, \quad \lambda=\mu=0$}  \\
\hline
L & 800.0   & 0.0   & 500.0  & 1.5 & 0.0   & 2.0 & 0.4 & 0.10  \\
R & 1000.0  & 0.0   & 600.0  & 1.0 & 0.0   & 1.0 & 0.3 &       \\
\hline 
\multicolumn{1}{l}{\textbf{RP3 \cite{USFORCE2}:}} & 
\multicolumn{8}{c}{ $\gamma_s = 1.4, \quad \pi_s = 0, \quad \gamma_g = 1.67, \quad \pi_g = 0, \quad \lambda=10^3, \,\, \mu=10^2$}  \\
\hline
L & 1.0    & 0.0   & 1.0 & 1.0   & 0.0  & 1.0 & 0.99 & 0.2   \\
R & 0.125  & 0.0   & 0.1 & 0.125 & 0.0  & 0.1 & 0.01 &       \\
\hline 
\end{tabular}
\end{center}
\label{tab.rpbn.ic}
\end{table}

\subsubsection{Spherical explosion problems} 
\label{sec.BN-EP}

\begin{figure}[!htbp]
	\begin{center}
	\begin{tabular}{cc} 
	\includegraphics[width=0.47\textwidth]{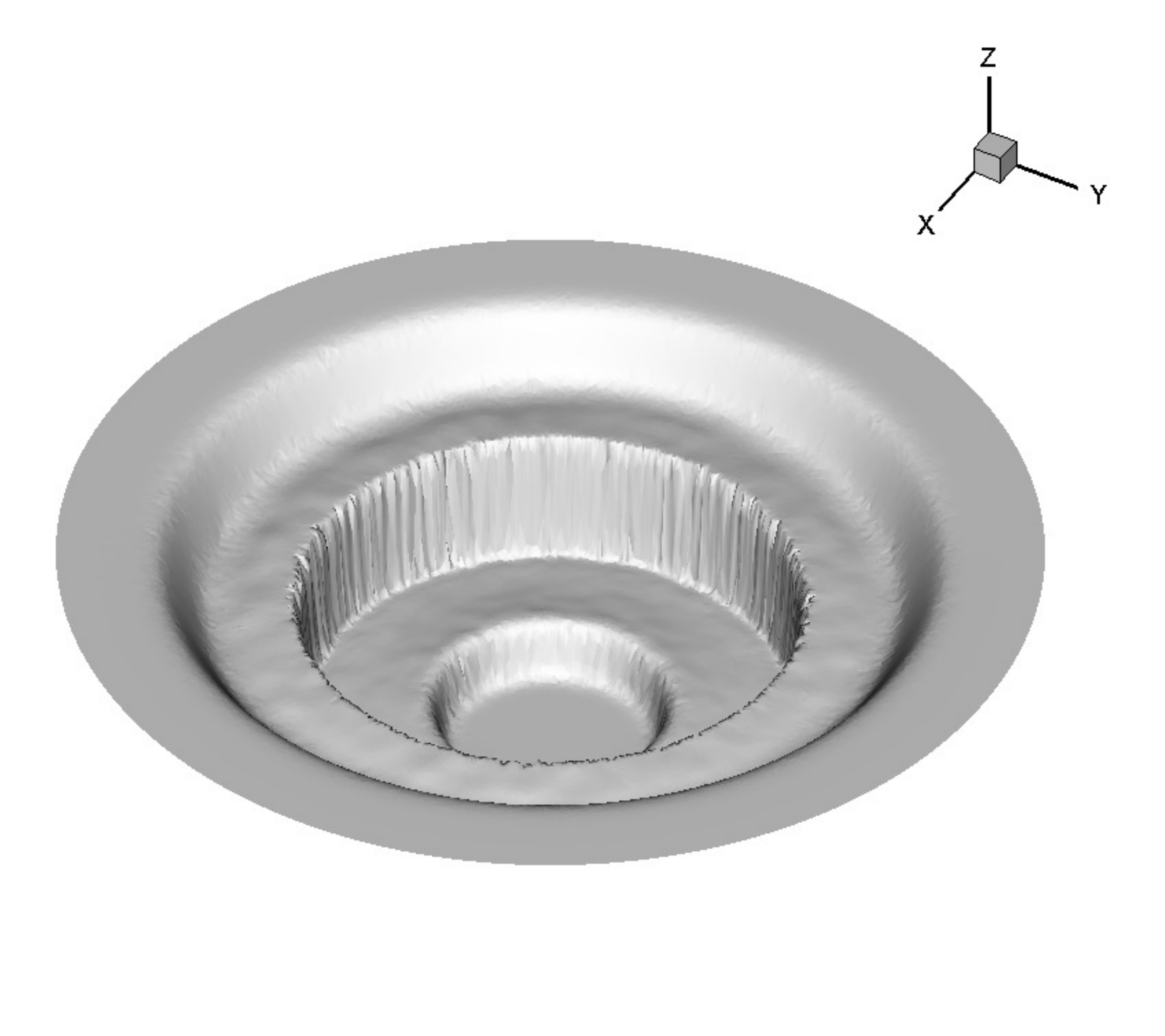}  &           
	\includegraphics[width=0.47\textwidth]{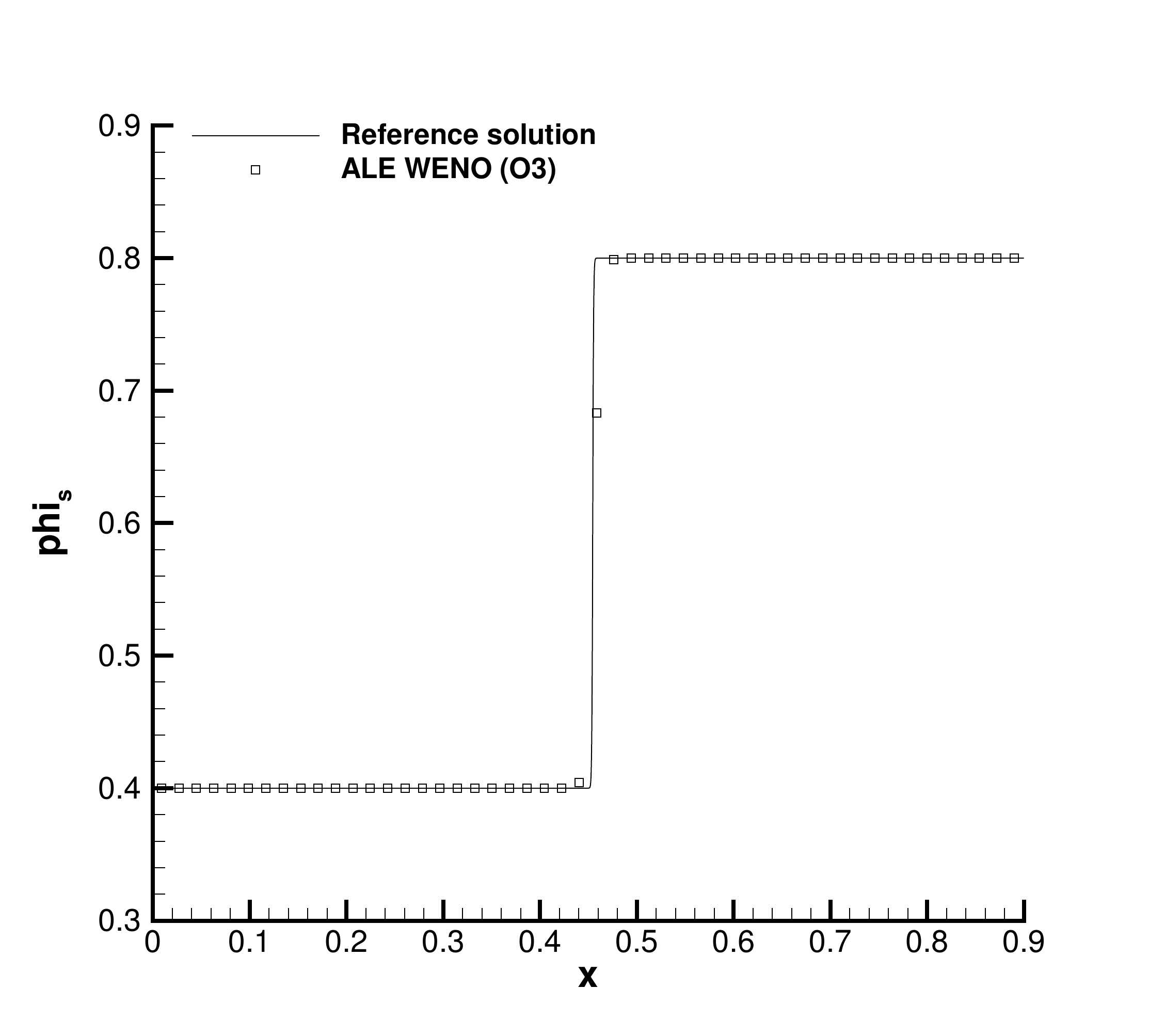} \\
	\includegraphics[width=0.47\textwidth]{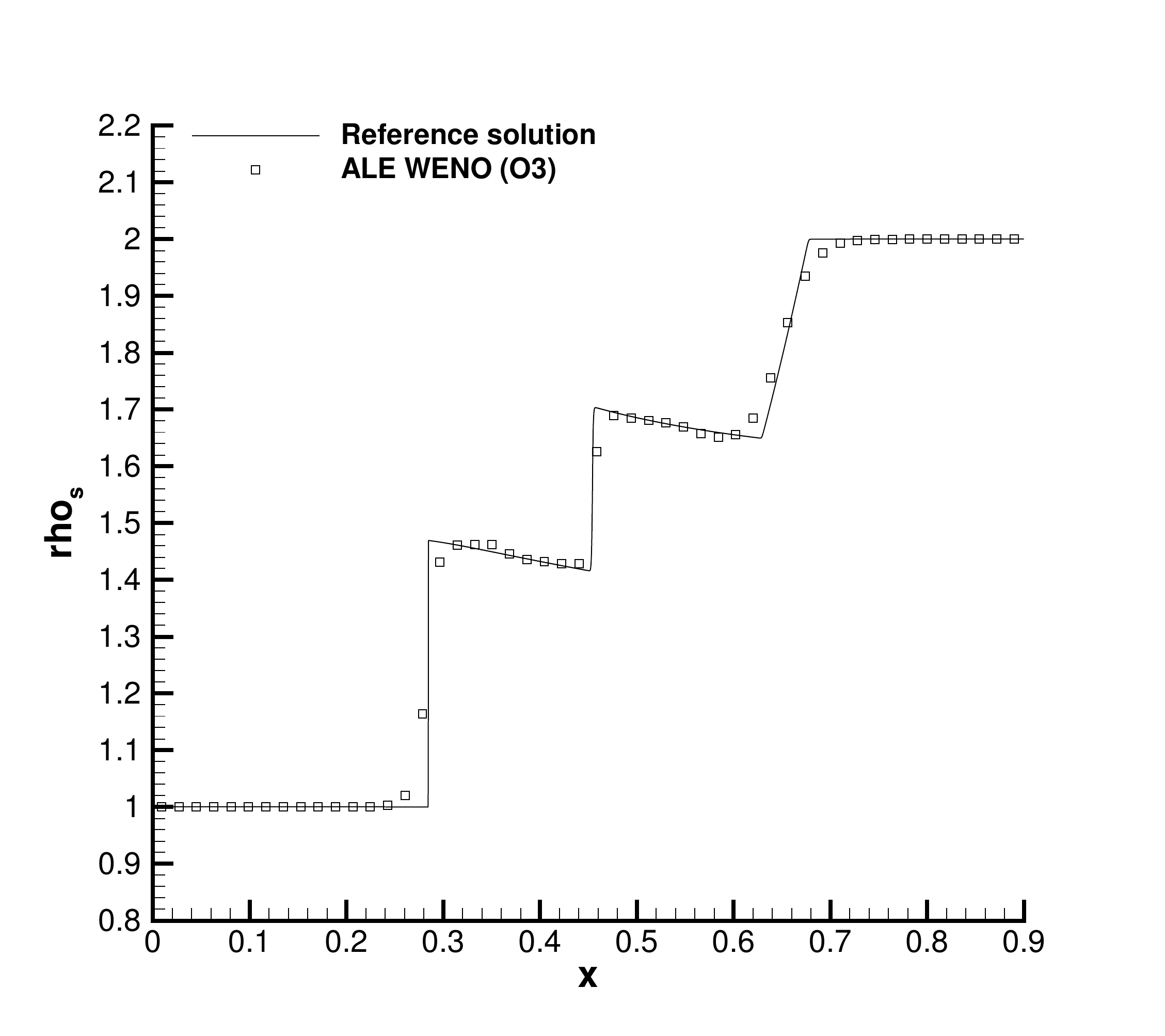}  &           
	\includegraphics[width=0.47\textwidth]{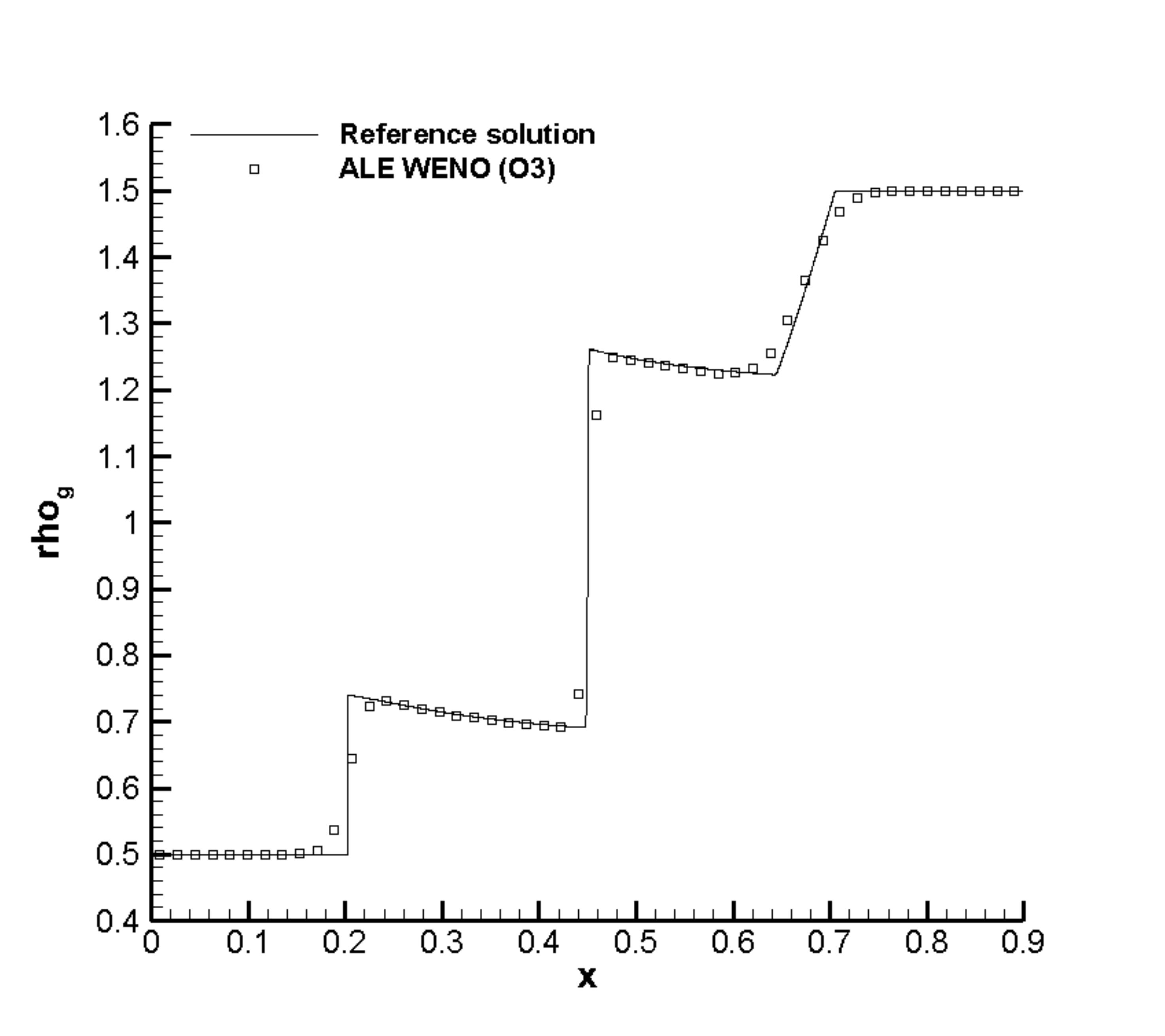} \\
	\includegraphics[width=0.47\textwidth]{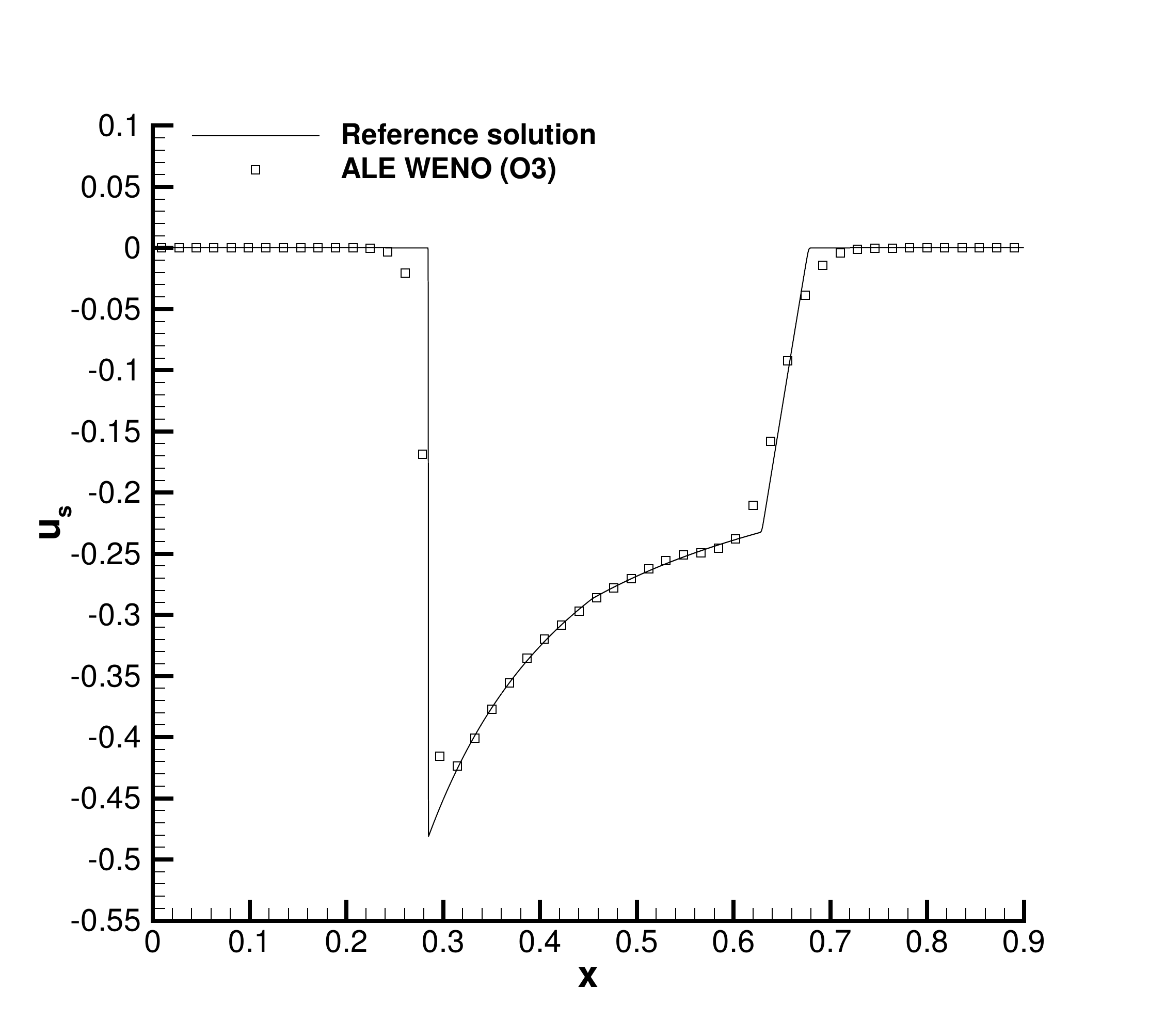}  &           
	\includegraphics[width=0.47\textwidth]{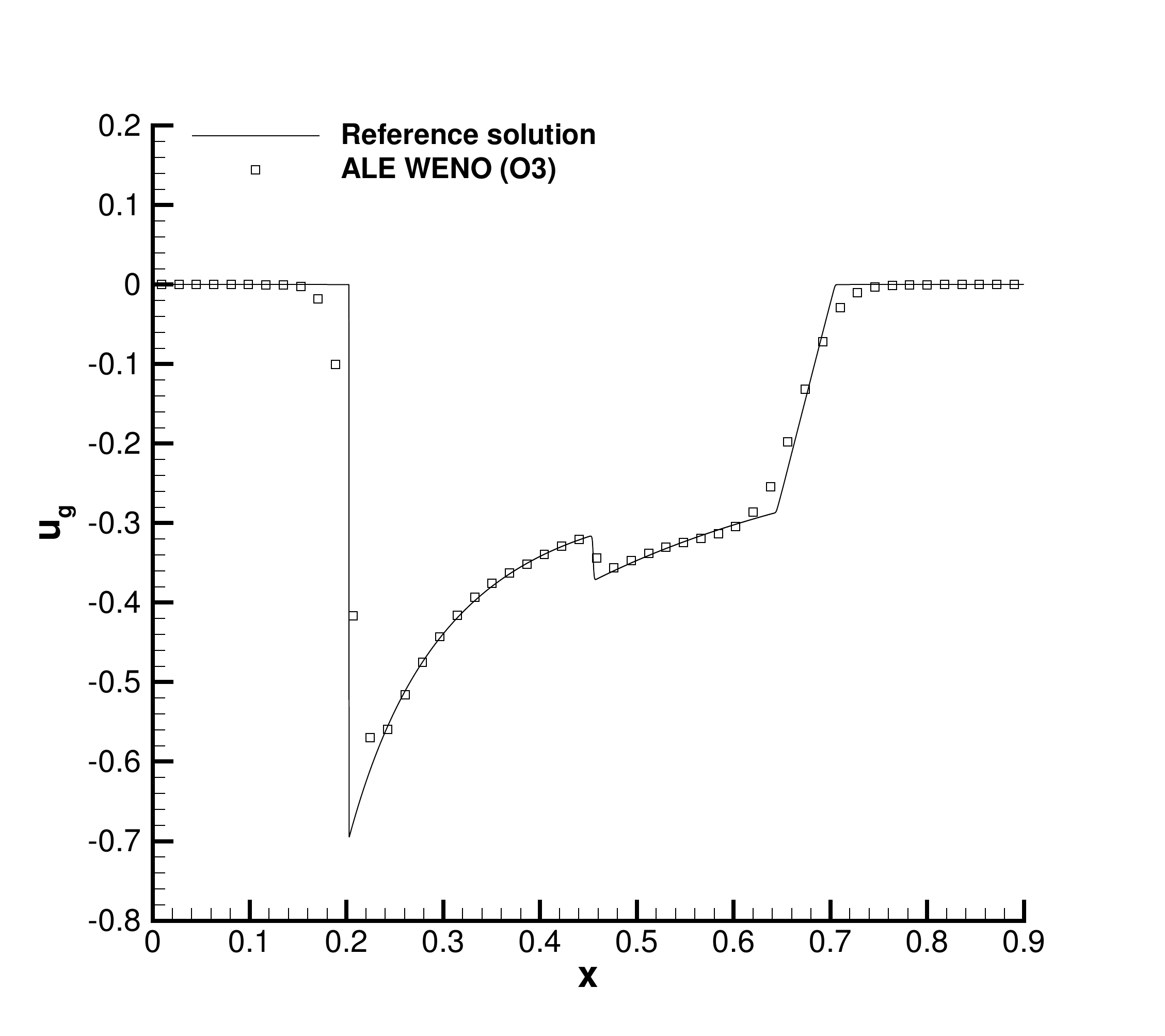} \\
	\includegraphics[width=0.47\textwidth]{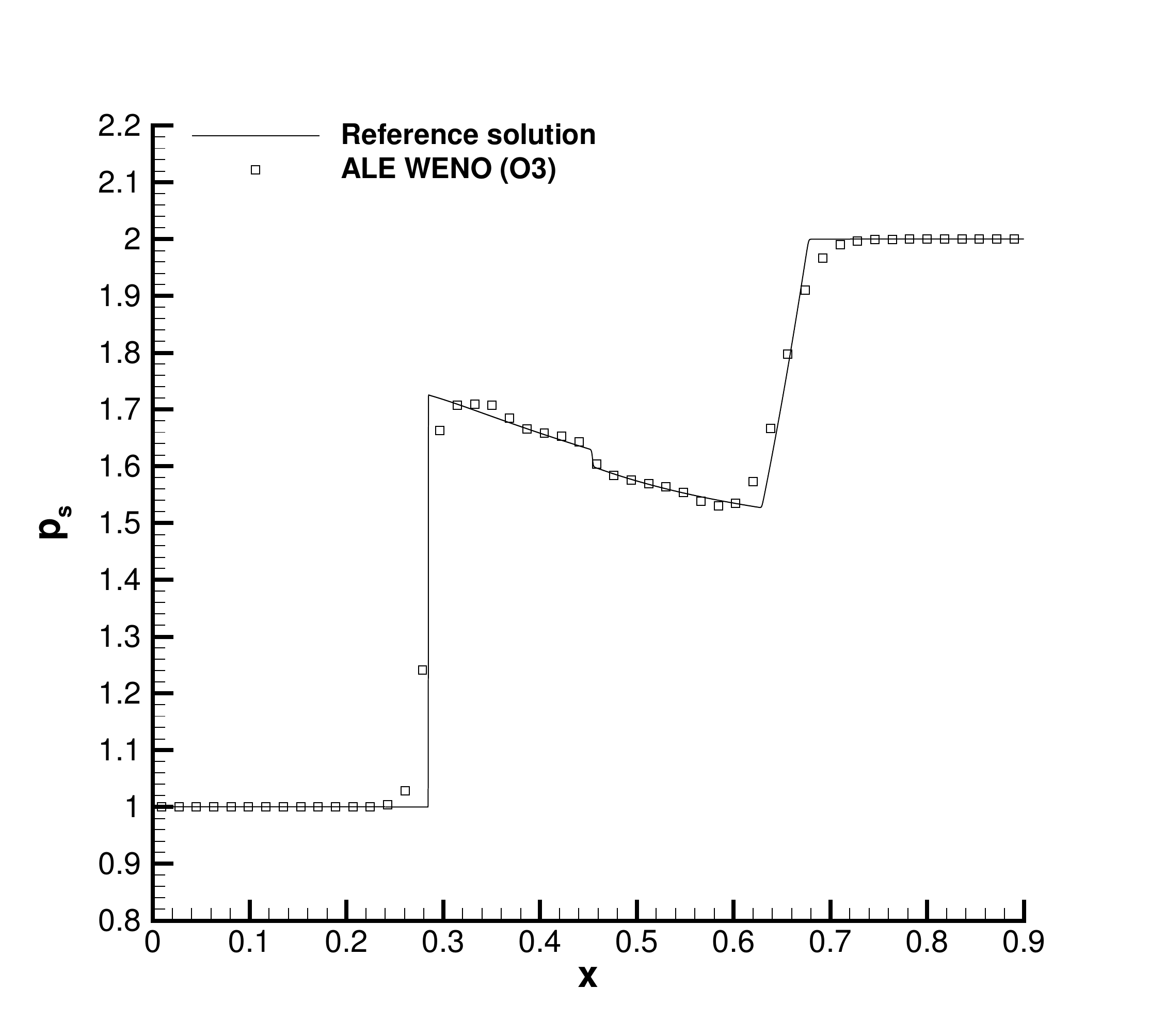}  &           
	\includegraphics[width=0.47\textwidth]{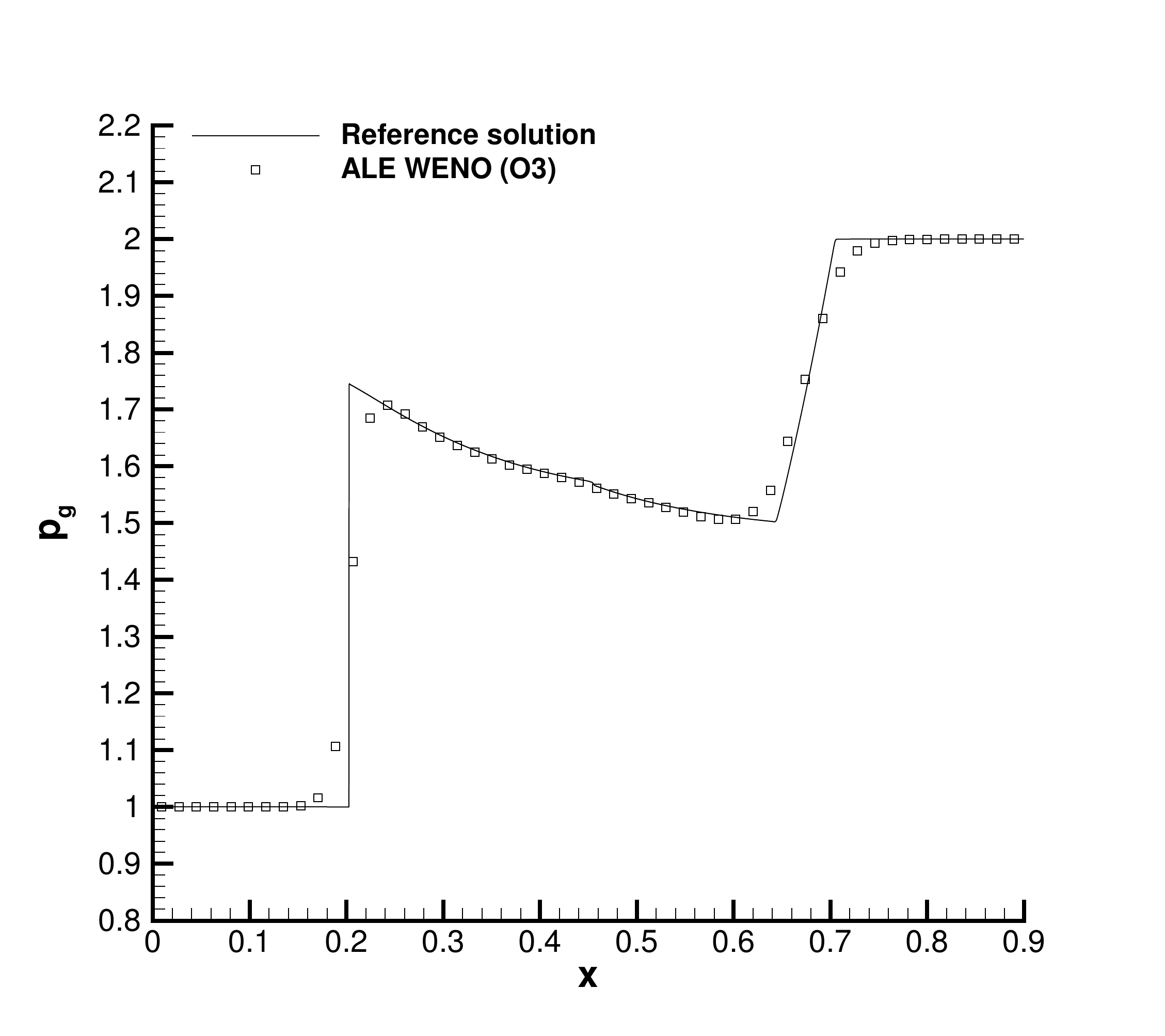} \\
	\end{tabular}
	\caption{Third order numerical results for the 3D explosion problem EP1 of the seven-equation Baer-Nunziato model at time $t=0.15$ and comparison with the reference solution.}
	\label{fig.EP1}
	\end{center}
\end{figure}

\begin{figure}[!htbp]
	\begin{center}
	\begin{tabular}{cc} 
	\includegraphics[width=0.47\textwidth]{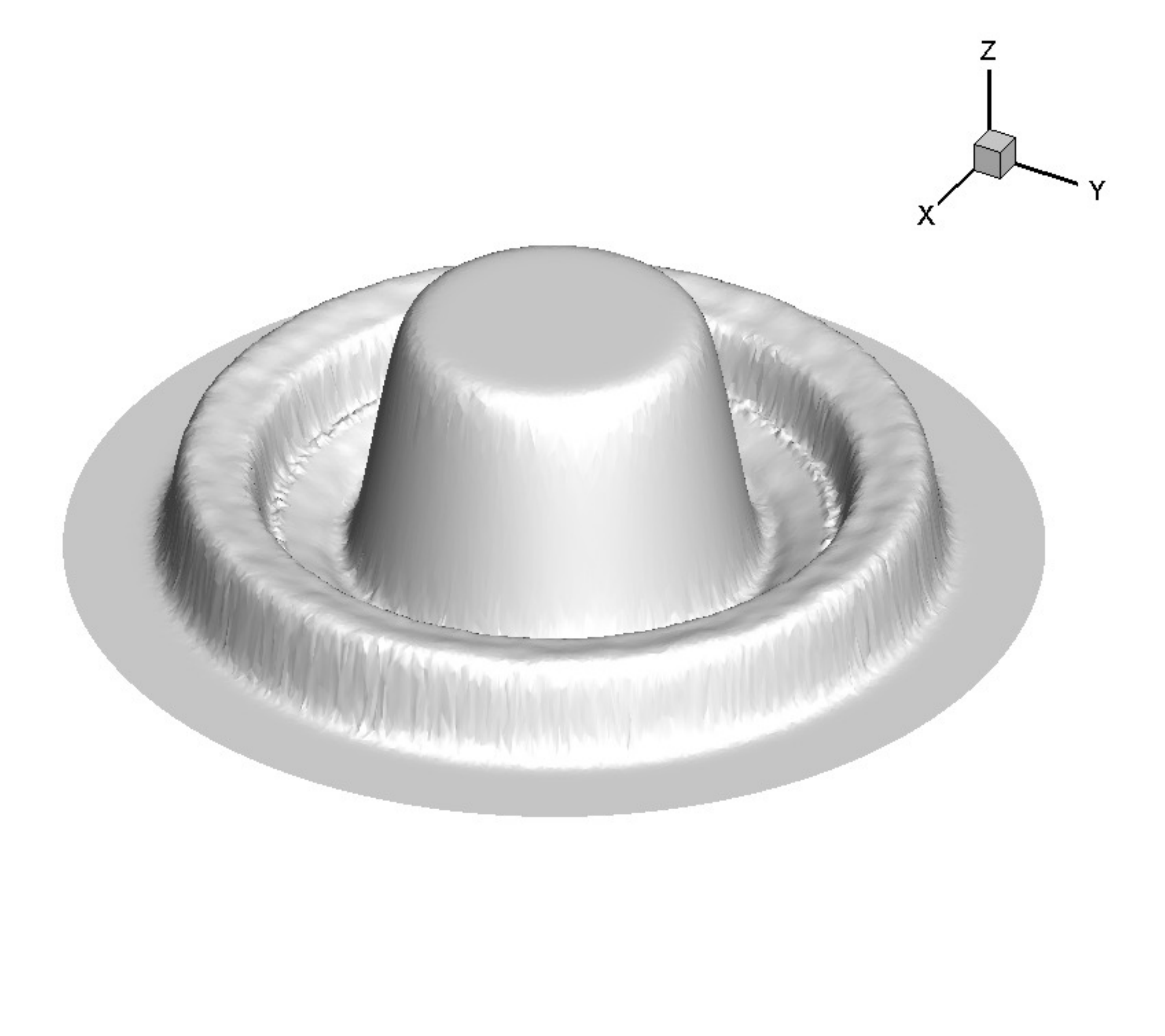}  &           
	\includegraphics[width=0.47\textwidth]{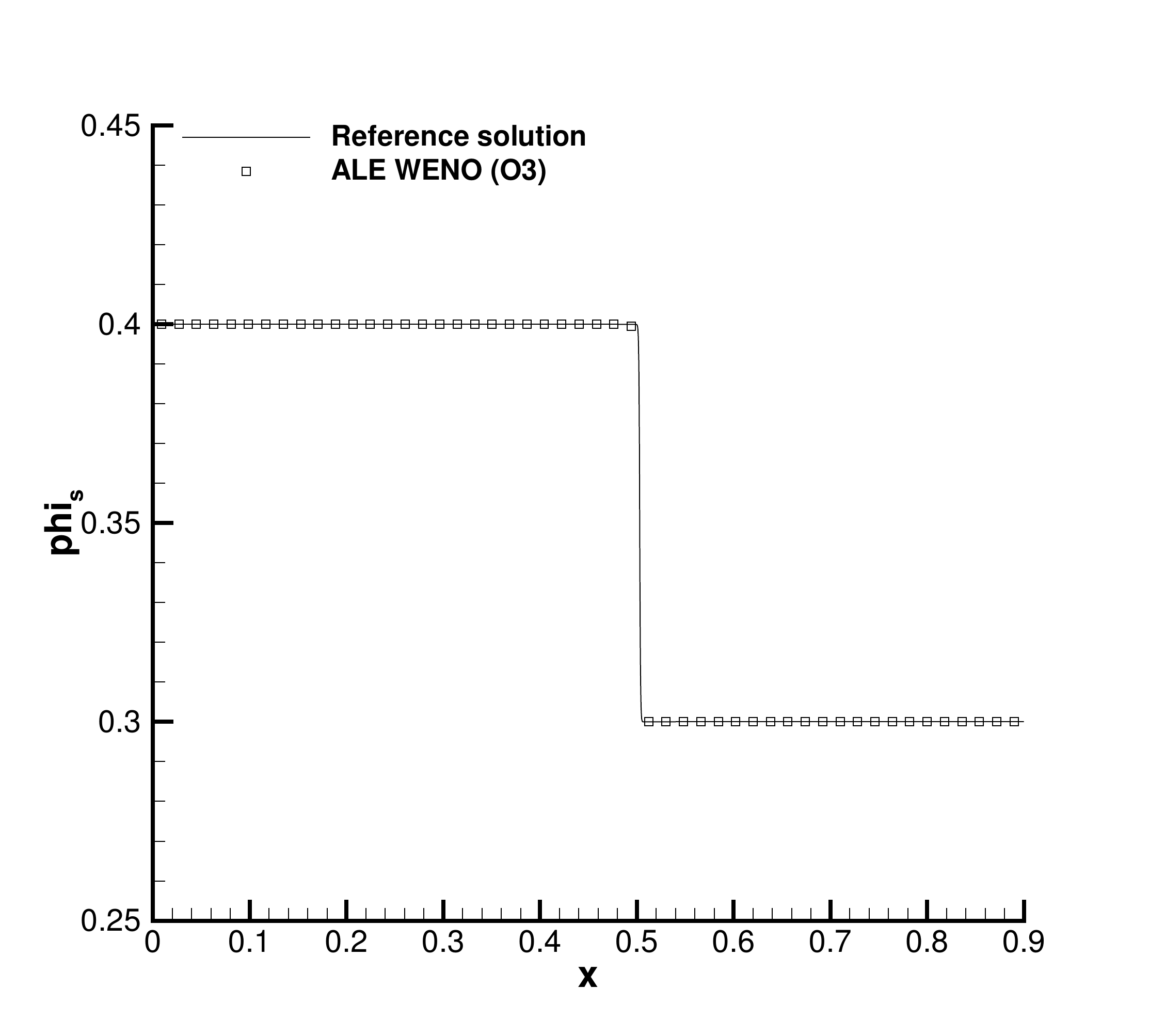} \\
	\includegraphics[width=0.47\textwidth]{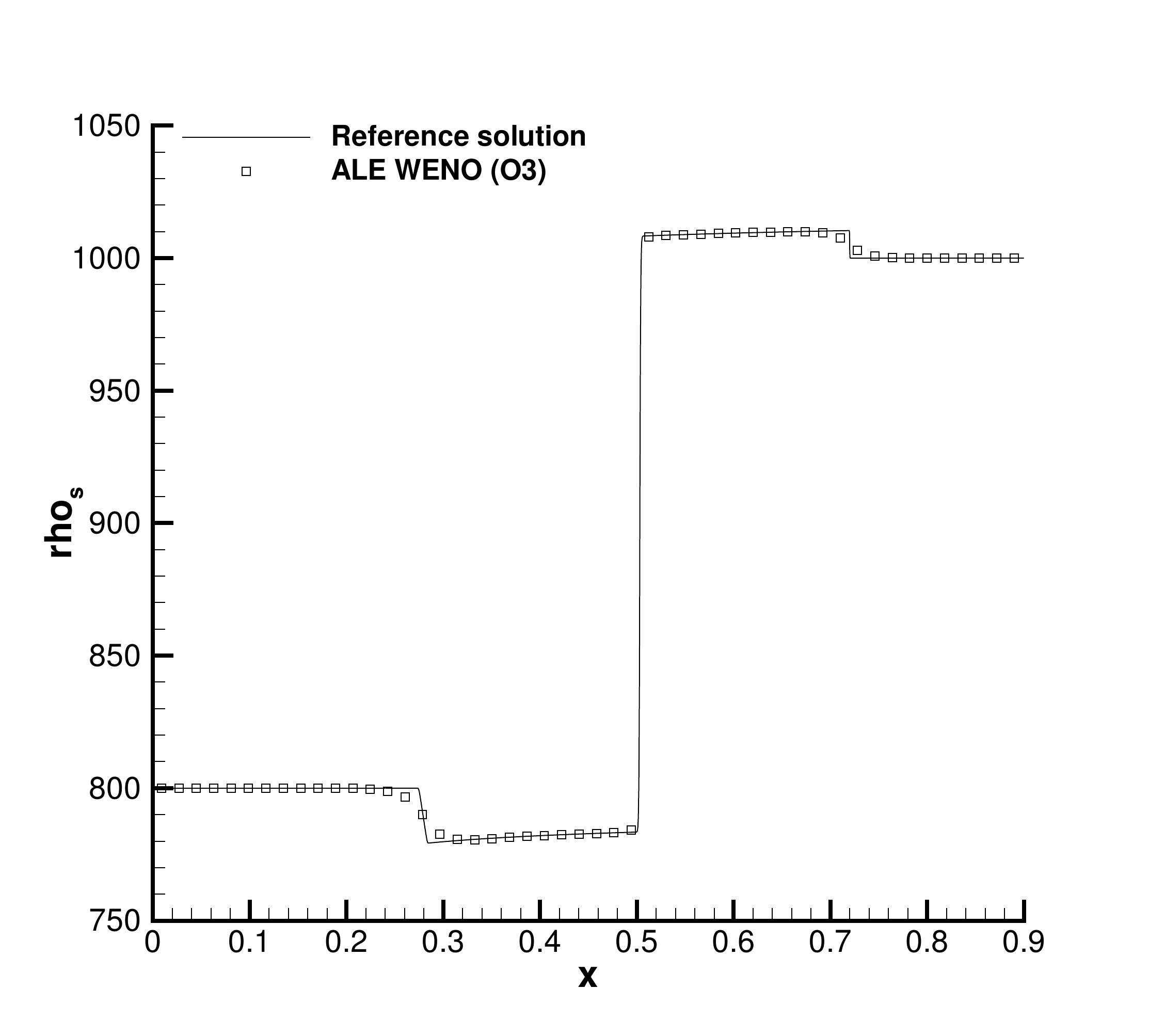}  &           
	\includegraphics[width=0.47\textwidth]{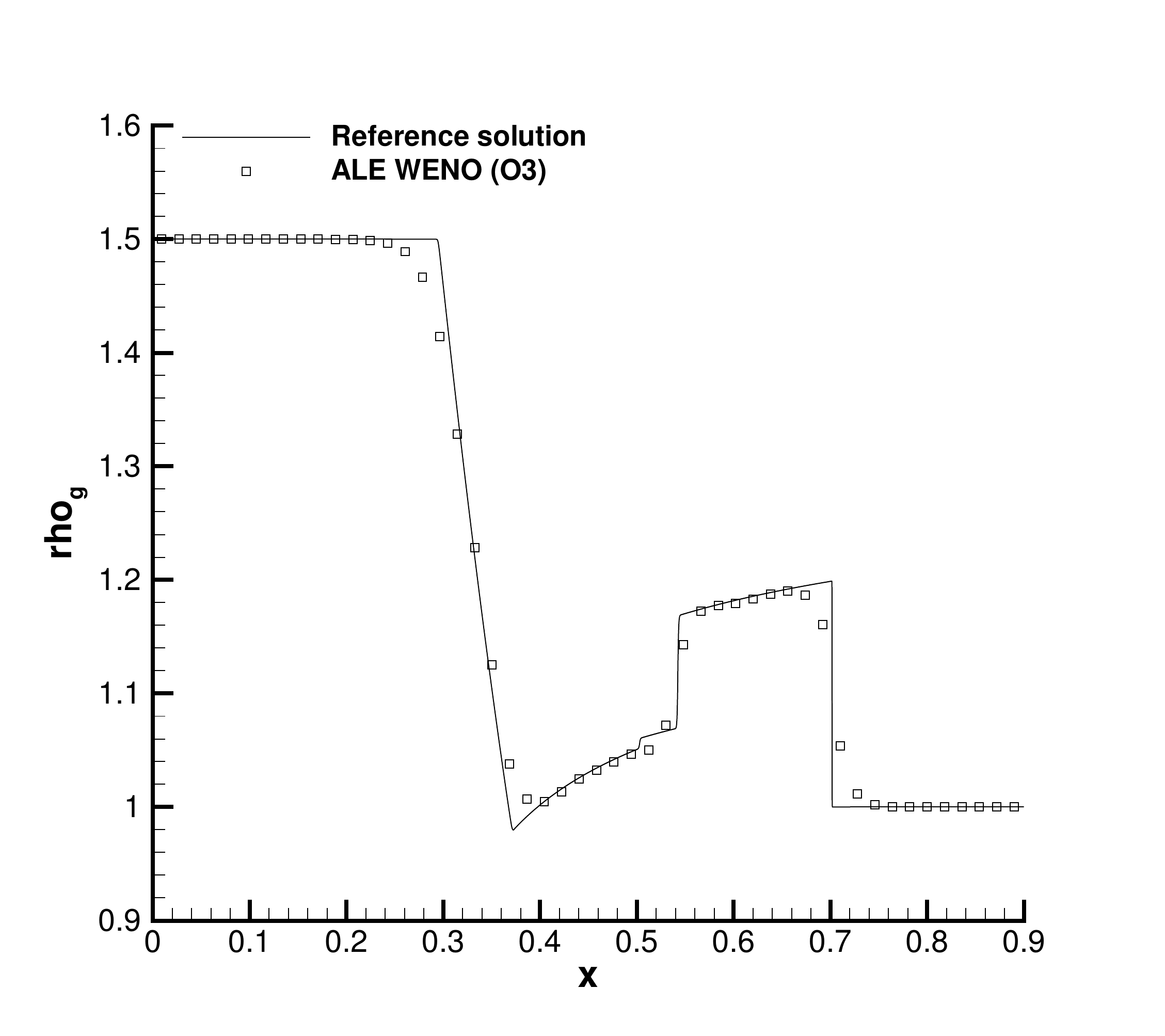} \\
	\includegraphics[width=0.47\textwidth]{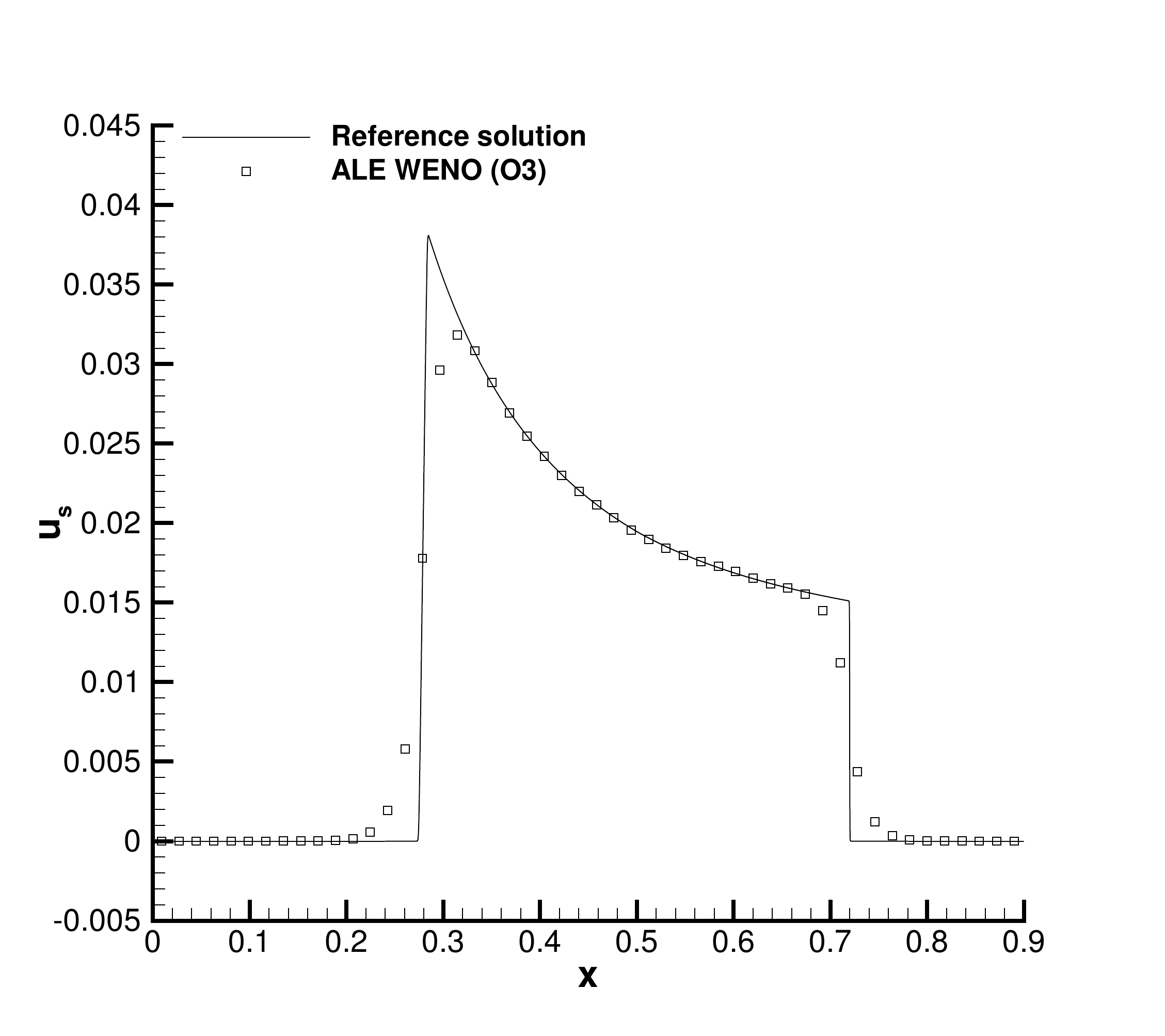}  &           
	\includegraphics[width=0.47\textwidth]{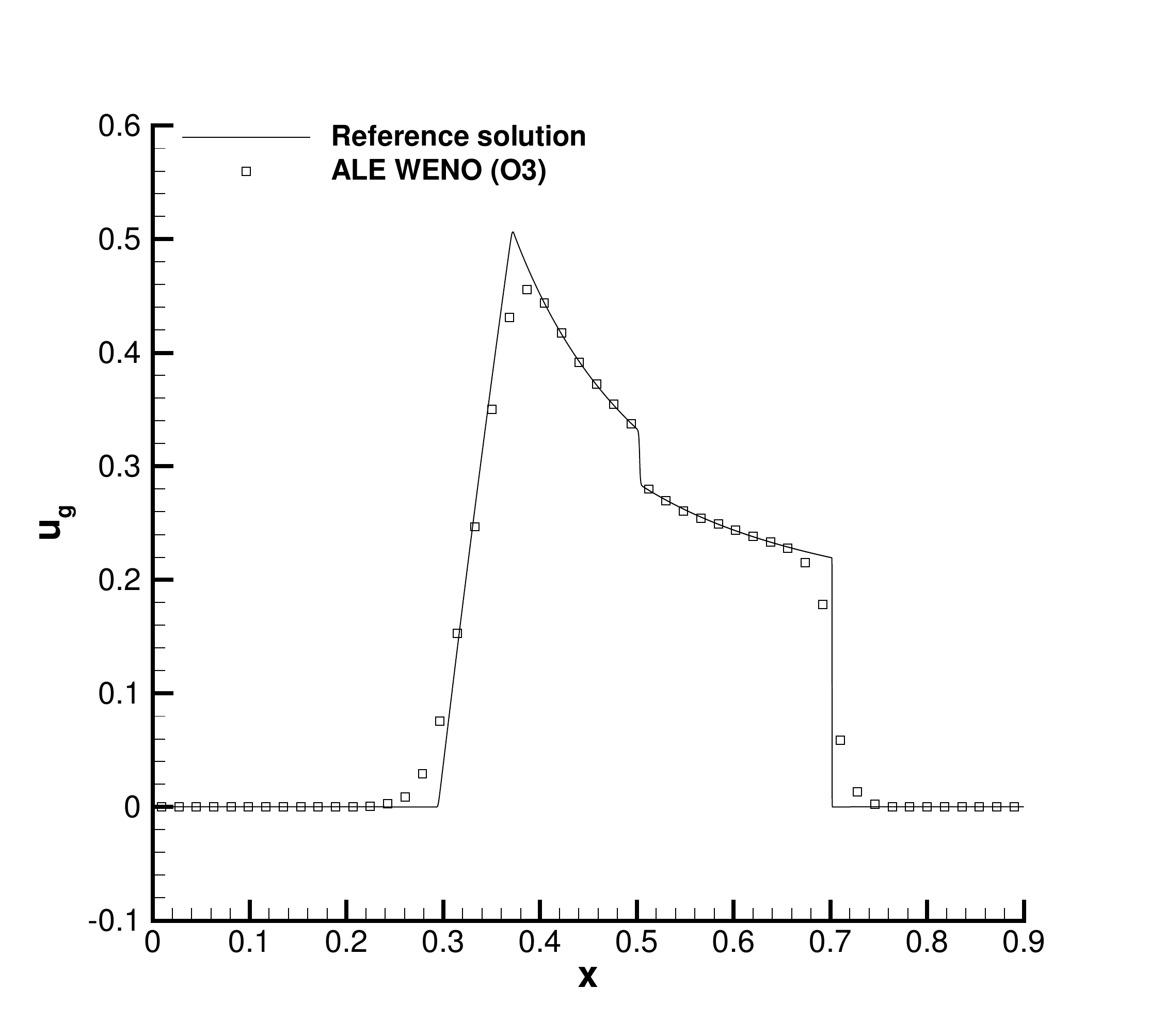} \\
	\includegraphics[width=0.47\textwidth]{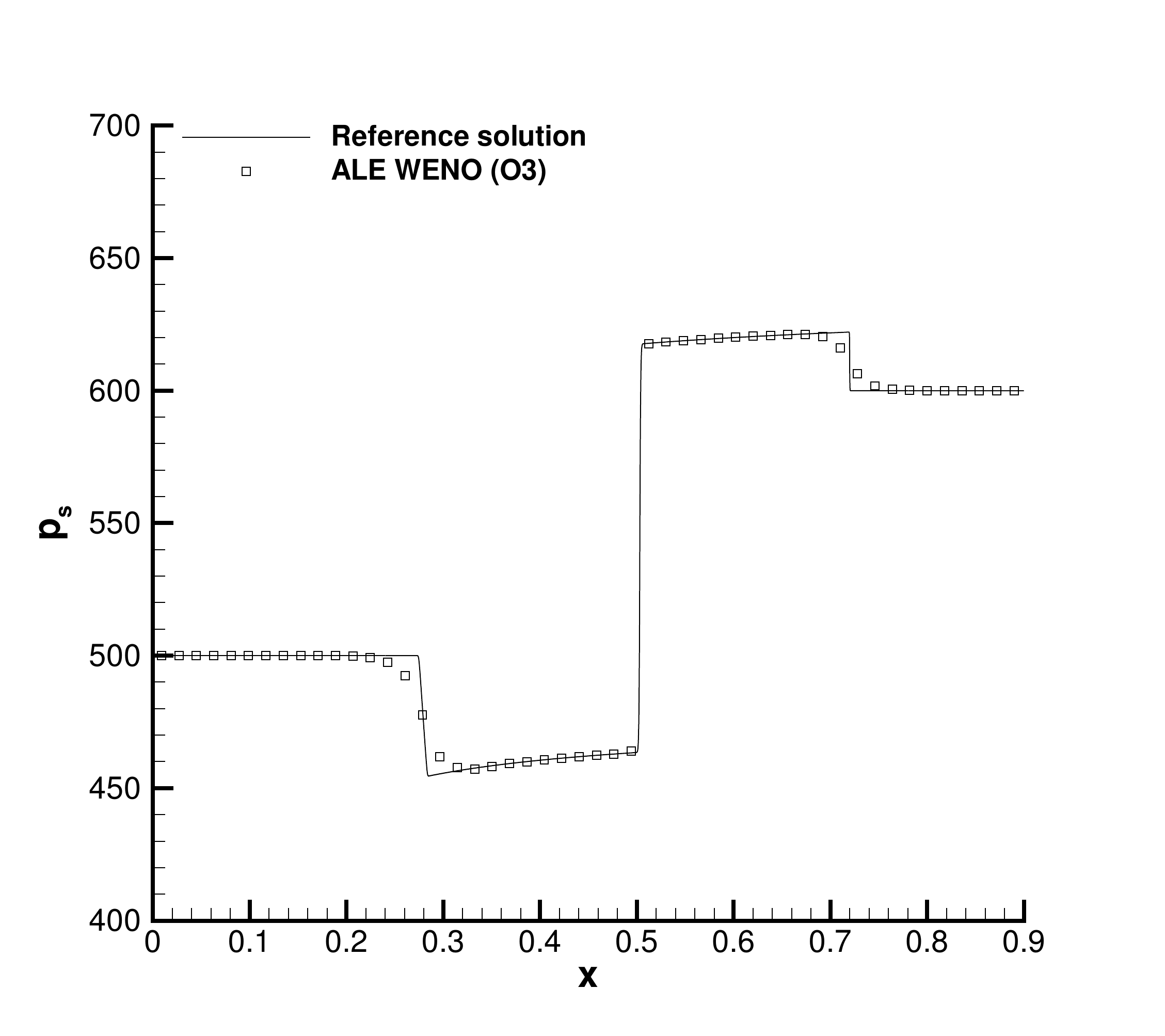}  &           
	\includegraphics[width=0.47\textwidth]{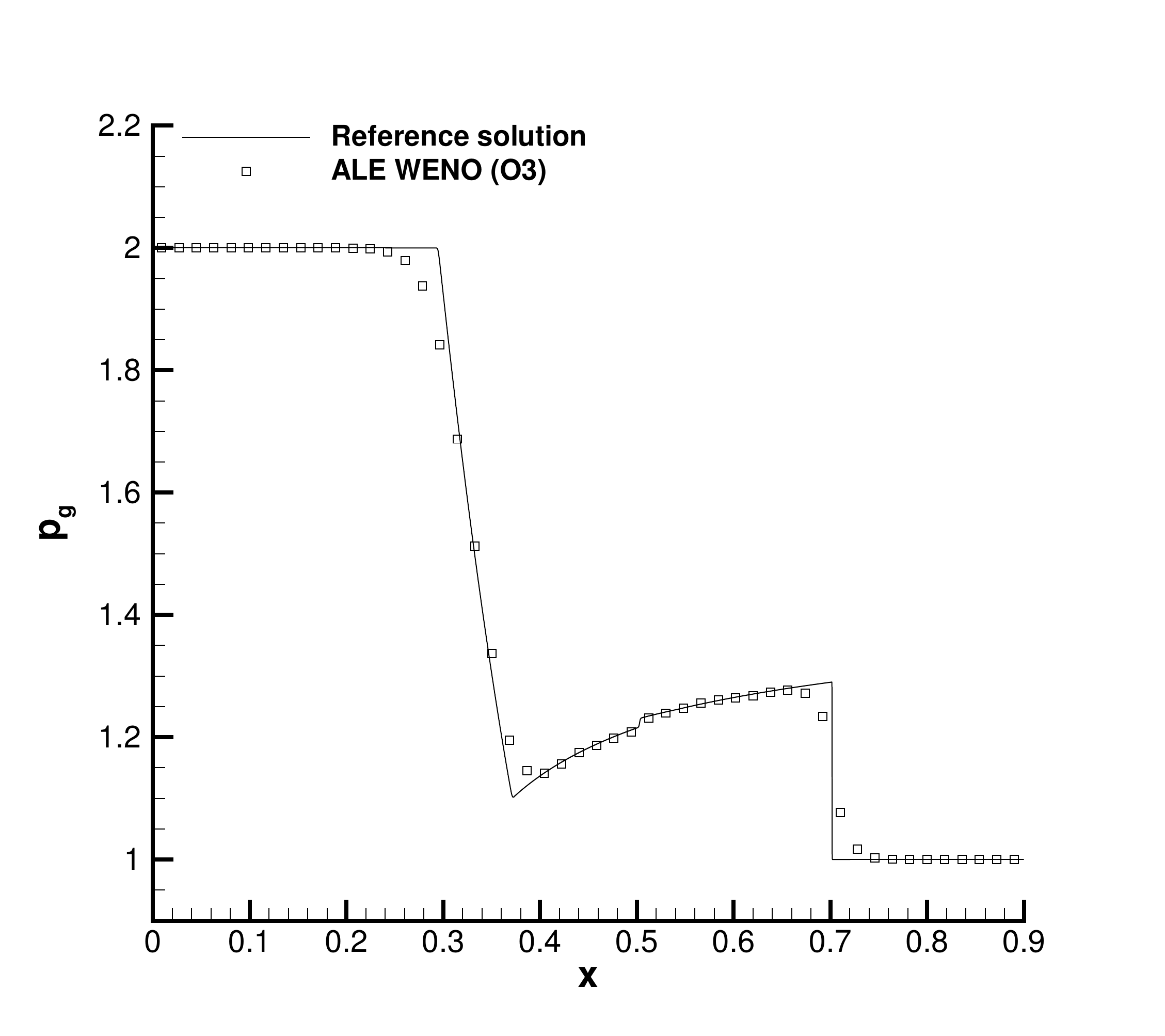} \\
	\end{tabular}
	\caption{Third order numerical results for the 3D explosion problem EP2 of the seven-equation Baer-Nunziato model at time $t=0.15$ and comparison with the reference solution.}
	\label{fig.EP2}
	\end{center}
\end{figure}

\begin{figure}[!htbp]
	\begin{center}
	\begin{tabular}{cc} 
	\includegraphics[width=0.47\textwidth]{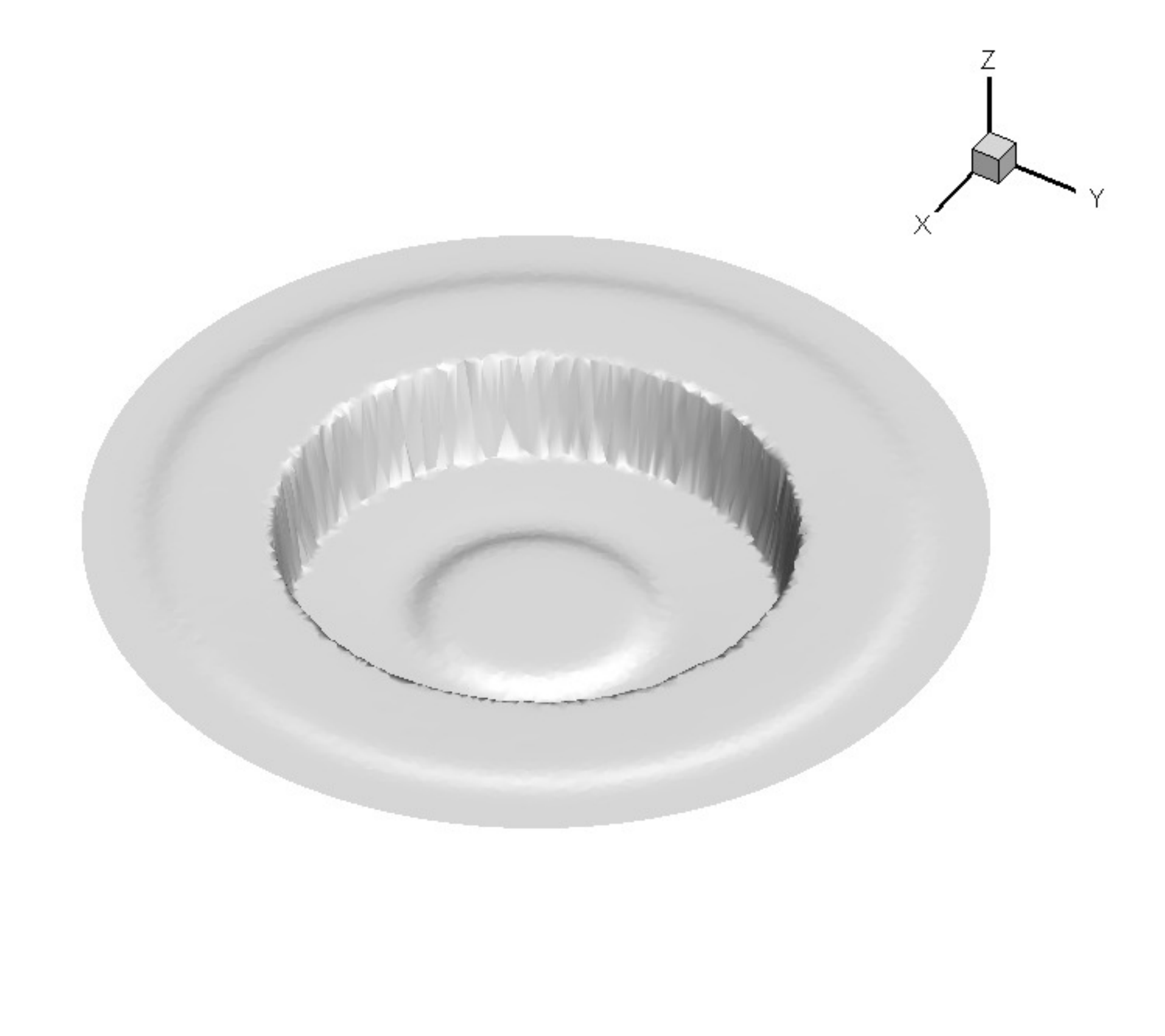}  &           
	\includegraphics[width=0.47\textwidth]{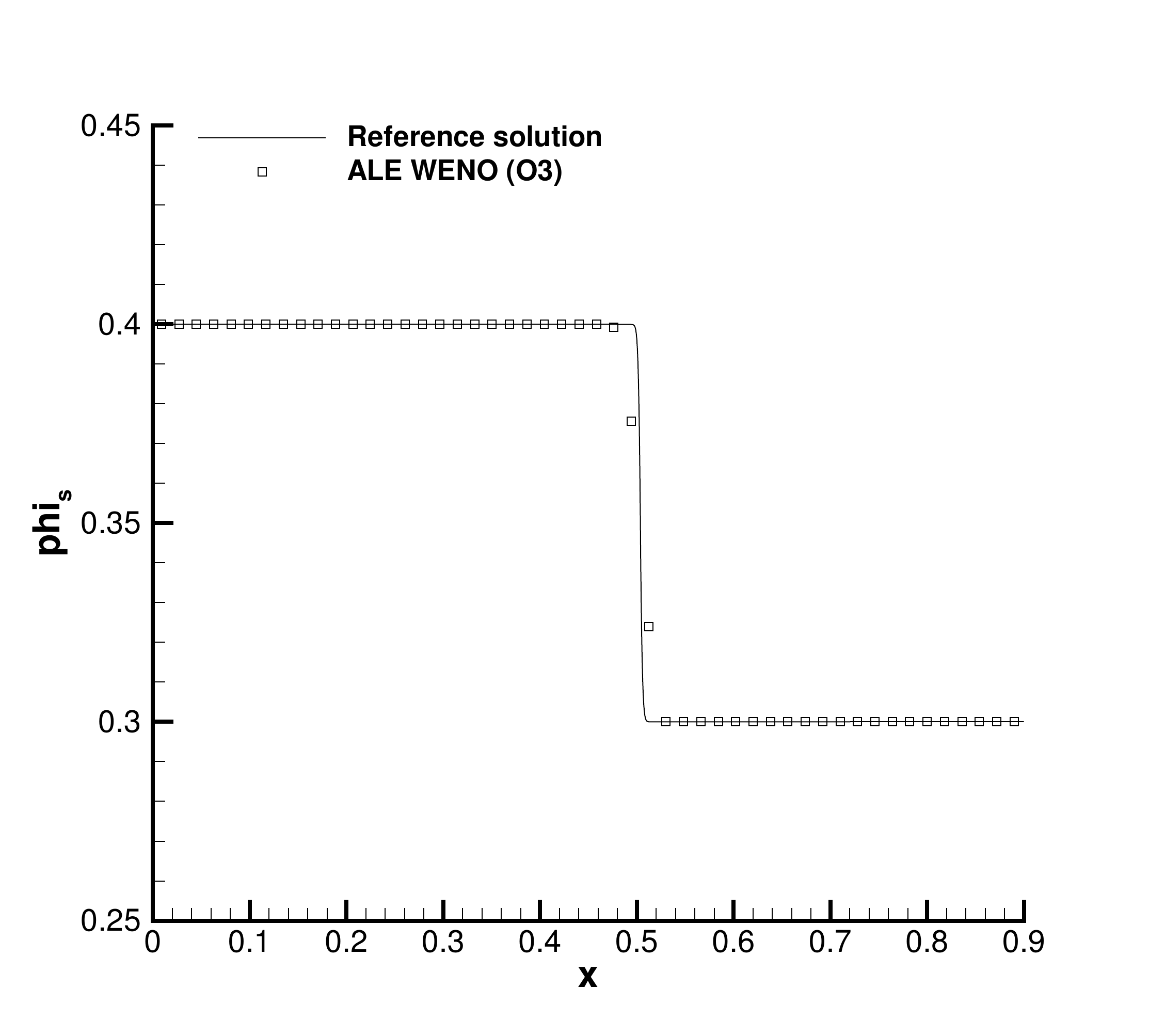} \\
	\includegraphics[width=0.47\textwidth]{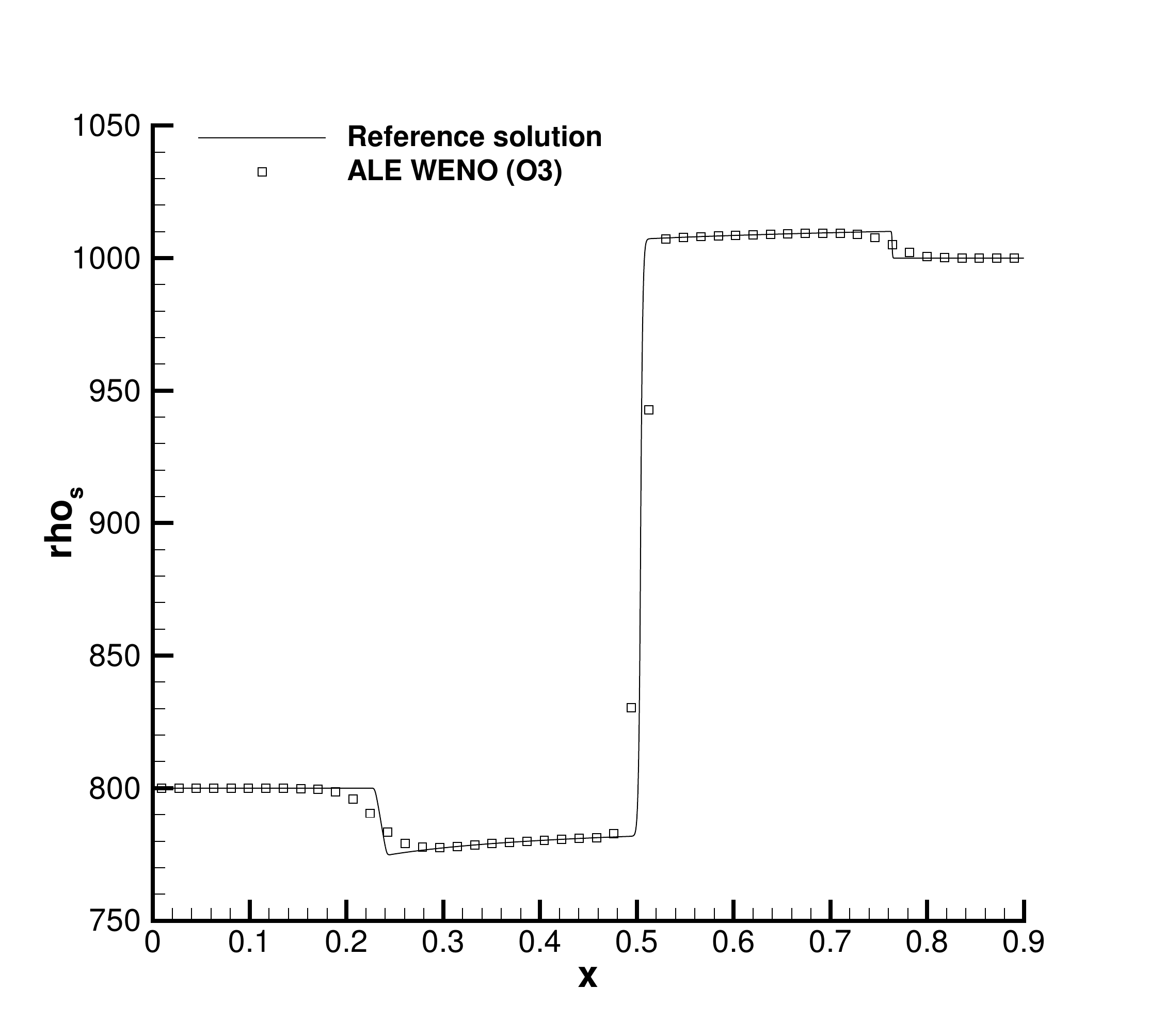}  &           
	\includegraphics[width=0.47\textwidth]{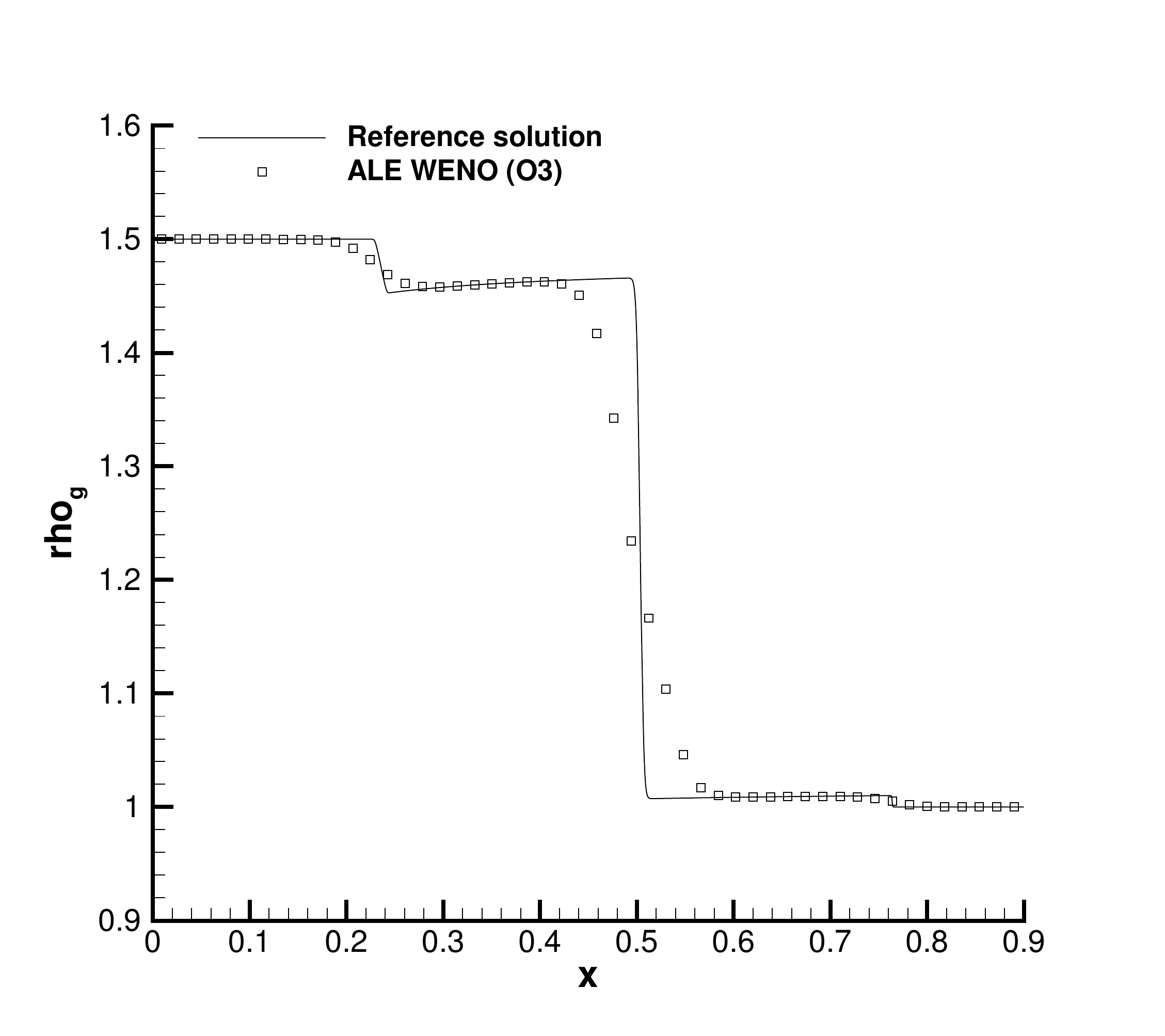} \\
	\includegraphics[width=0.47\textwidth]{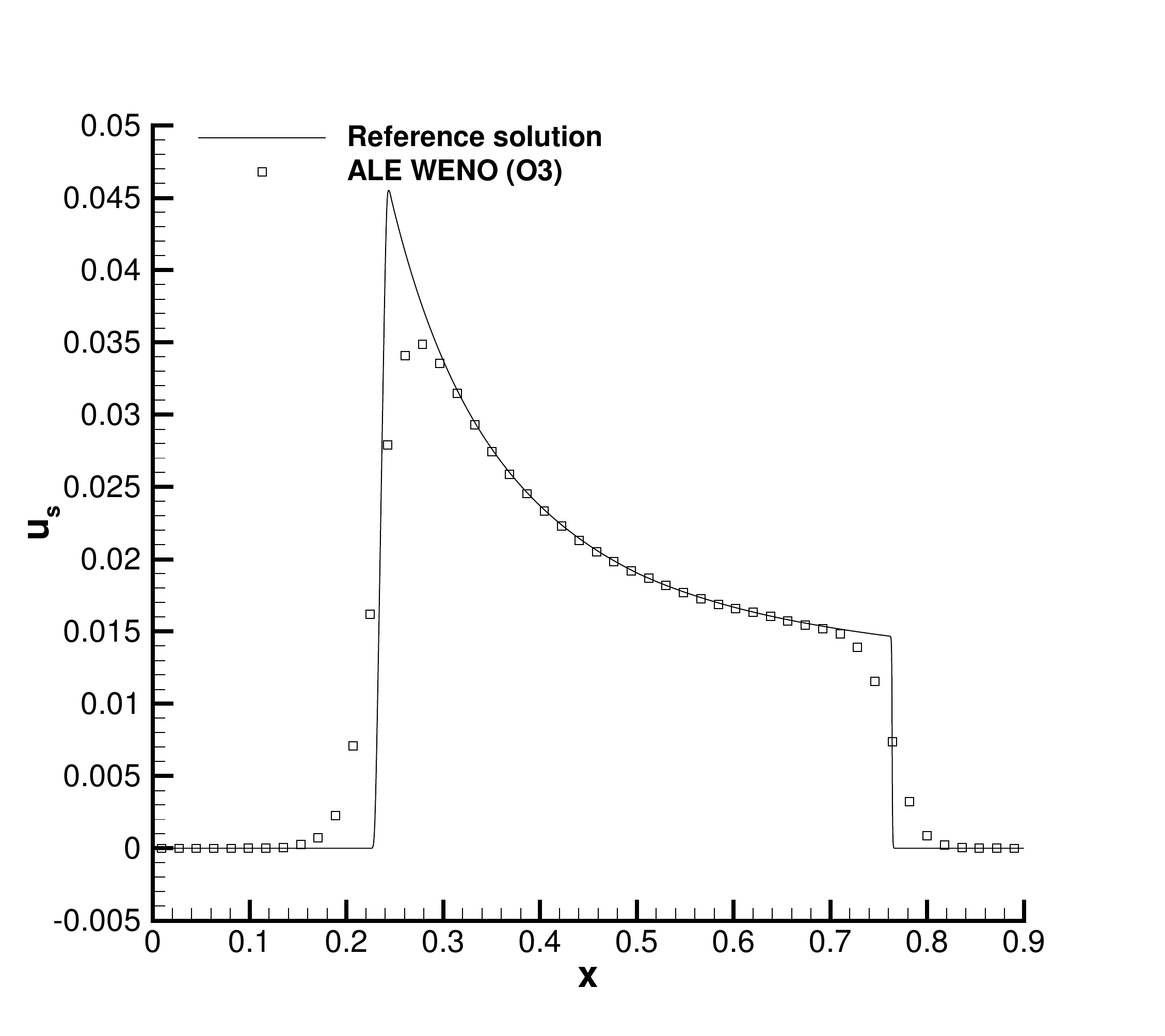}  &           
	\includegraphics[width=0.47\textwidth]{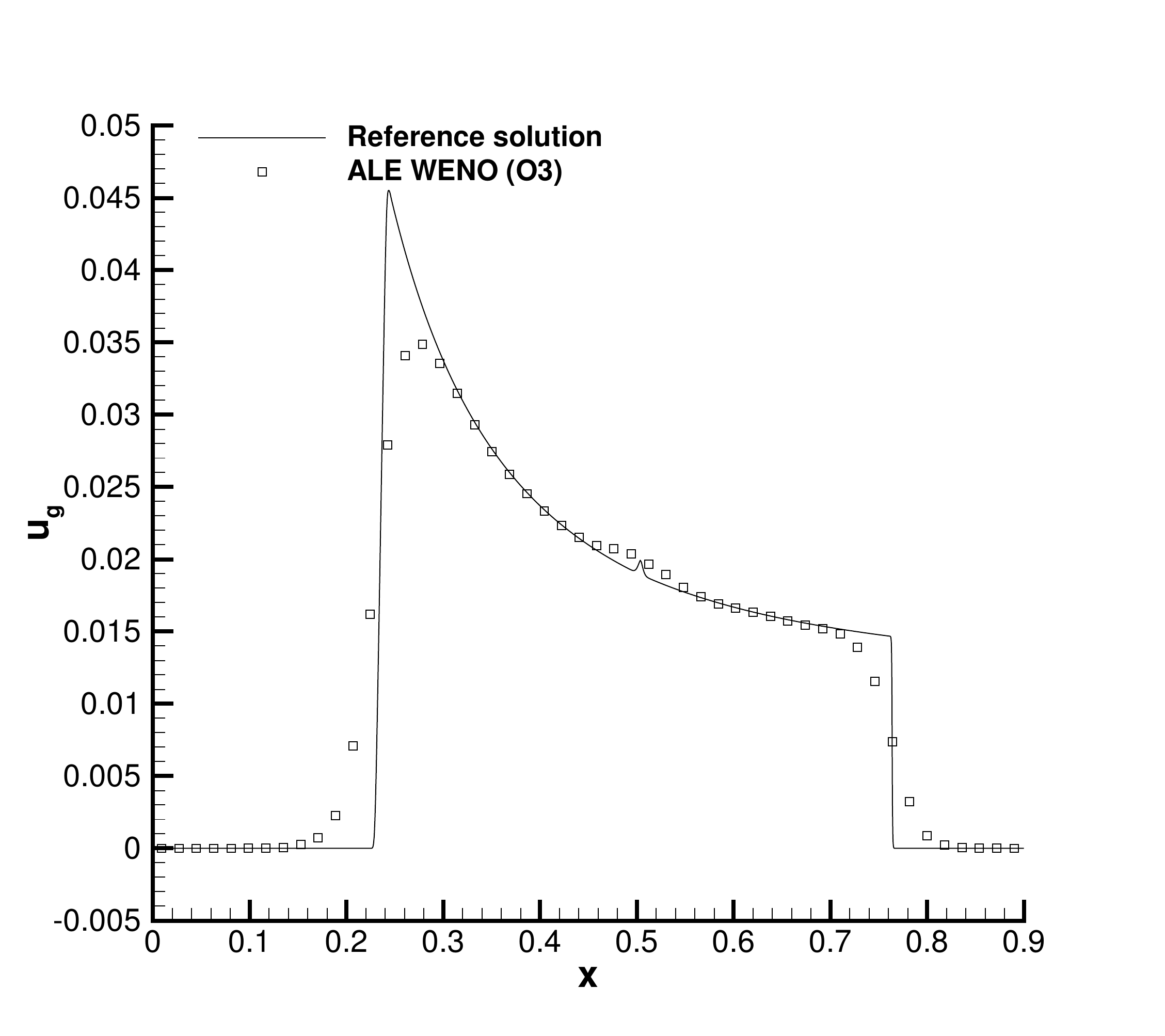} \\
	\includegraphics[width=0.47\textwidth]{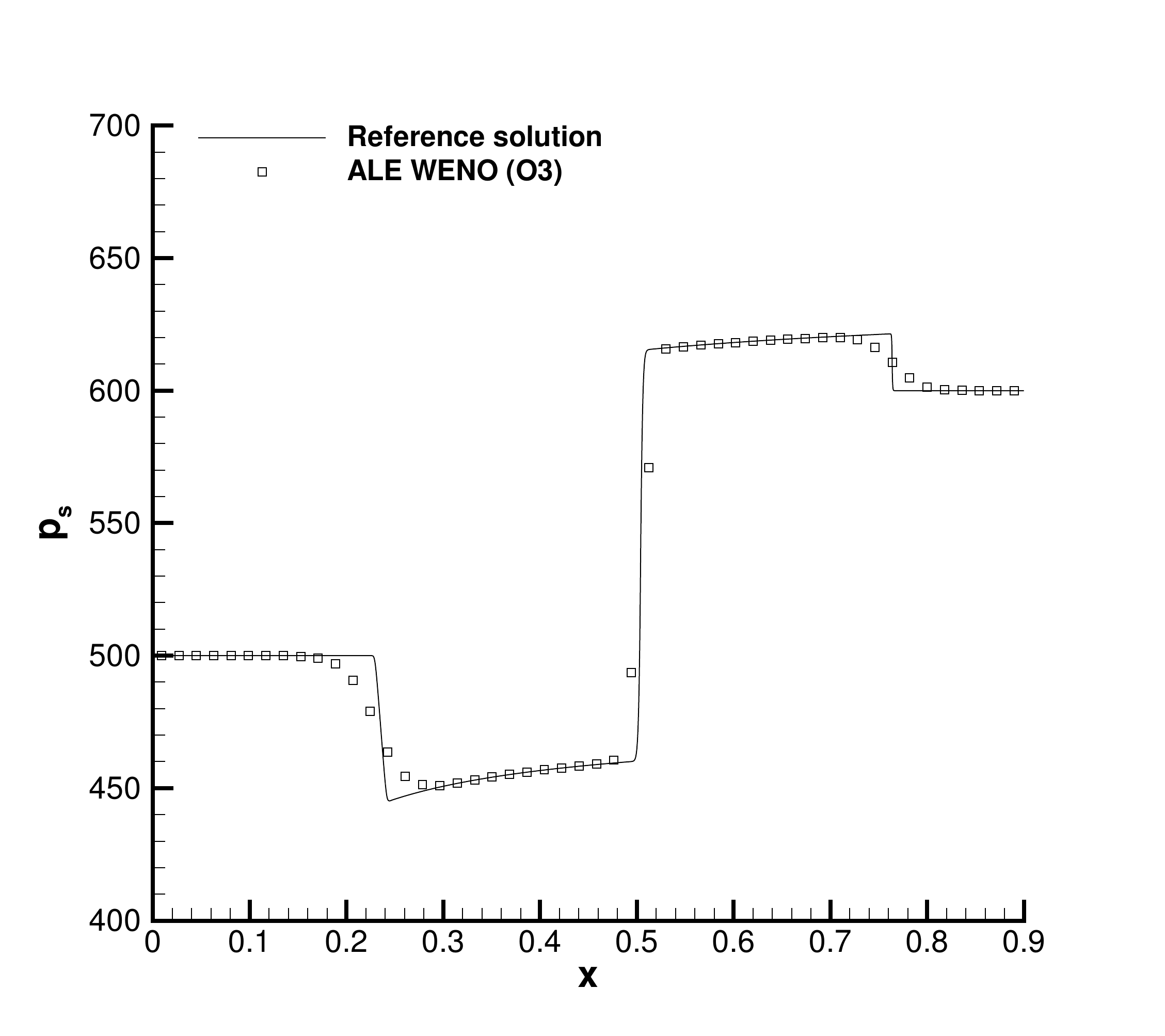}  &           
	\includegraphics[width=0.47\textwidth]{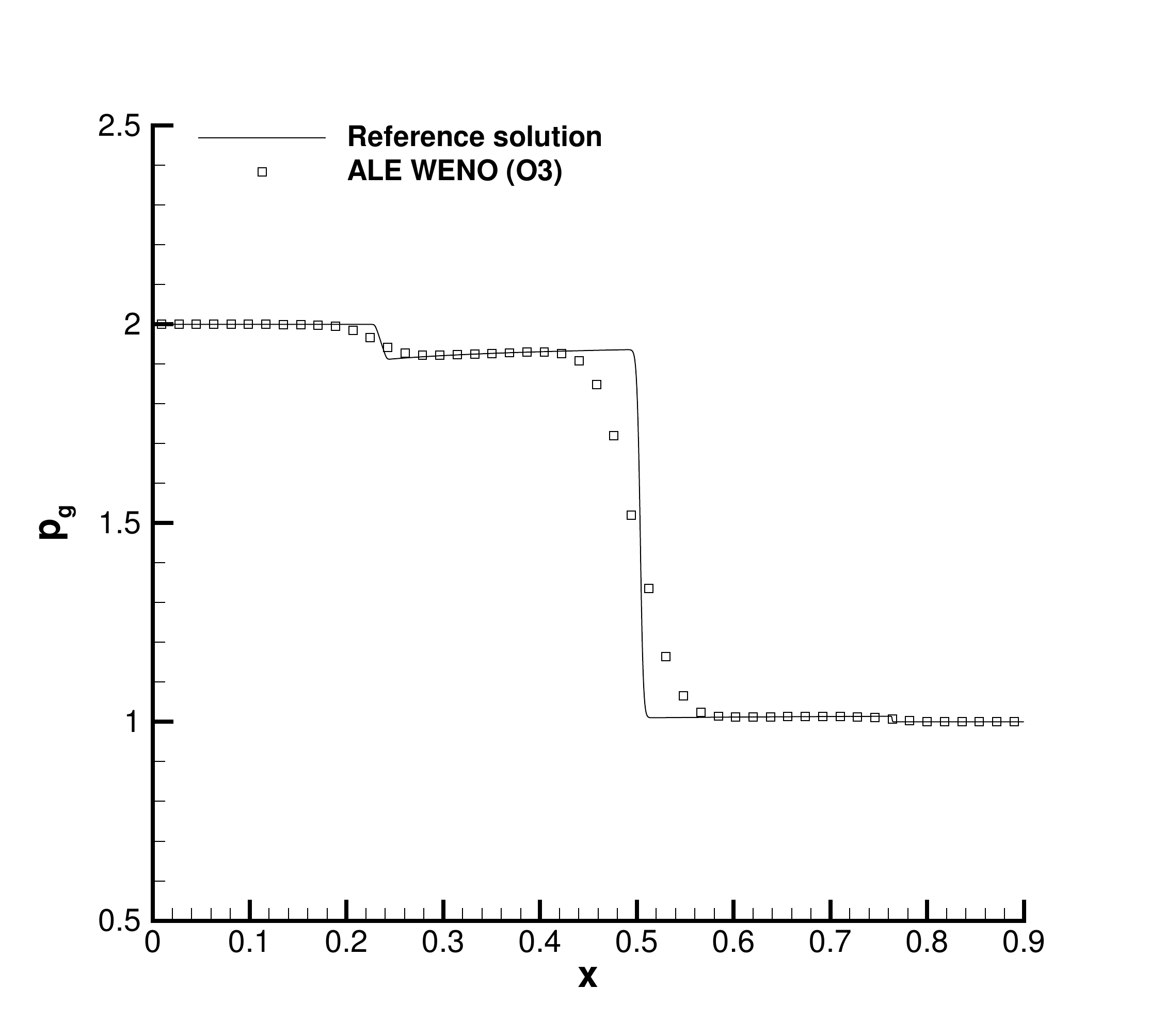} \\
	\end{tabular}
	\caption{Third order numerical results for the 3D explosion problem EP3 with $\lambda=10^5$ of the seven-equation Baer-Nunziato model at time $t=0.18$ and comparison with the reference solution.}
	\label{fig.EP3}
	\end{center}
\end{figure}

We use the same initial condition given for the Riemann problems in Table \ref{tab.rpbn.ic} to solve three different spherical explosion problems with the compressible Baer-Nunziato model \eqref{ec.BN}. The computational domain $\Omega(0)$ is initially the sphere of radius $R=0.9$, which has been discretized with a characteristic mesh size of $h=1/100$ for $r \leq r_c$ and $h=1/50$ for $r > r_c$ for a total number of elements of $N_E=2632305$. In all cases the initial state $\Q(\x,0)$ is assigned taking 
\begin{equation}
  \Q(\x,0) = \left\{ \begin{array}{ll} 
                            \Q_i, & \qquad \textnormal{ if } \quad | \x |  < r_c \\
                            \Q_o, & \qquad \textnormal{ else}          
                      \end{array} \right. ,  
\label{eqn.ep.ic}                       
\end{equation}
with $r_c=0.5$ representing the location of the initial discontinuity. The left state reported in Table \ref{tab.rpbn.ic} is assumed to be the inner state $\Q_i$, while the right state represents here the outer state $\Q_o$. In particular the initial condition of Riemann problem RP1 is used for the first explosion problem EP1 and the same applies for RP2 and the corresponding EP2. In the third explosion problem EP3 we use again the initial values of RP2 and we set $\lambda=10^5$ and $\mu=0$, hence adopting a stiff interphase drag. As done for the spherical explosion problem for the compressible Euler equations, the reference solution is obtained by solving an equivalent non-conservative one-dimensional PDE in radial direction with geometric reaction source terms using a path-conservative second 
order TVD scheme, see \cite{USFORCE2} for details. The final time is set to $t_f=0.15$ for EP1 and EP2, while $t_f=0.18$ is used for EP3. Figures \ref{fig.EP1} - \ref{fig.EP3} show a comparison between the numerical results obtained with a third order ADER-WENO ALE scheme and the one-dimensional reference solution. We use the the path-conservative Osher-type method \eqref{eqn.osher} since it is 
less dissipative than the Rusanov-type scheme \eqref{eqn.rusanov}, hence a better resolution of the material contact can be achieved. Since the mesh is moving with the interface velocity 
$\mathbf{u_I}$, i.e. $\mathbf{V}=\mathbf{u_I}=\mathbf{u}_1$, the contact discontinuity of the first phase $\phi_1$ is very well resolved in all cases. 

\section{Conclusions}
\label{sec.concl}
In this article we have presented a new high order Arbitrary-Lagrangian-Eulerian one-step ADER-WENO finite volume scheme on unstructured tetrahedral meshes in three space dimensions. Numerical convergence studies up to sixth order of accuracy in space and time have been shown and the algorithm is formulated in a very general manner so that it can be applied to both conservative and non-conservative hyperbolic systems, with and without stiff source terms. To the knowledge of the authors, this is the first better than second order accurate Lagrangian finite volume scheme ever presented on 
unstructured tetrahedral meshes. Several classical test problems have been run for the Euler equations of compressible gasdynamics, for the MHD equations and for the seven-equation Baer-Nunziato model of compressible multiphase flows. Where possible, the obtained numerical results have been carefully compared with exact or other numerical reference solutions. \\
Further work will regard the extension of the presented Lagrangian ADER-WENO finite volume schemes to the more general framework of the new $P_{N}P_{M}$ method proposed in \cite{Dumbser20088209}, which 
can deal with either pure finite volume or pure discontinuous Galerkin finite element methods, or with a hybridization of both. Future work will also concern an extension of the present method to the 
\textit{a posteriori} limiter paradigm MOOD \cite{MOODorg,MOODhighorder,ADERMOOD} and to the use of time-accurate local time stepping (LTS), see \cite{LTS}. Further research will also be necessary to  extend the multi-dimensional Riemann solvers used in \cite{BalsaraMultiDRS,LagrangeMHD,LagrangeMDRS} to the case of moving unstructured tetrahedral meshes. Last but not least, another important topic 
will be the application of the present scheme to more realistic real world simulations in engineering and physics.   

\section*{Acknowledgments}
The presented research has been financed by the European Research Council (ERC) under the European Union's Seventh Framework 
Programme (FP7/2007-2013) with the research project \textit{STiMulUs}, ERC Grant agreement no. 278267. The authors acknowledge 
PRACE for awarding us access to the SuperMUC supercomputer of the Leibniz Rechenzentrum (LRZ) in Munich, Germany. 

\bibliography{Lagrange3D}
\bibliographystyle{plain}

\end{document}